\documentclass[12pt]{amsart}

\usepackage{CJKutf8}

\usepackage{adjustbox}

\usepackage{ifthen}
\usepackage[utf8]{inputenc}
\usepackage{setspace}
\usepackage[english]{babel}
\usepackage{enumerate}
\usepackage{thmtools}
\usepackage[shortlabels]{enumitem}
\usepackage[T1]{fontenc}
\usepackage{tikz}
\usetikzlibrary{shapes.geometric,intersections,decorations.markings,snakes}
\usetikzlibrary{calc,intersections,through,backgrounds}
\usetikzlibrary{patterns}
\usepackage{stmaryrd}

\DeclareMathOperator{\wt}{wt}

\DeclareMathOperator{\zz}{\mathbb{Z}}

\makeatletter
\newcommand{\preceqdot}{\mathrel{\mathpalette\pr@ceqd@t\relax}}
\newcommand{\pr@ceqd@t}[2]{%
  \begingroup
  \sbox\z@{$#1\prec$}\sbox\tw@{$#1\preceq$}%
  \dimen@=\dimexpr\ht\tw@-\ht\z@\relax
  {\preceq}%
  \mkern-5mu
  \raisebox{\dimen@}{$\m@th#1\cdot$}%
  \endgroup
}
\makeatother

\usepackage{wasysym}
\usepackage{booktabs}

\makeatletter
\providecommand*\xrightarrowtriangle[2][]{%
  \ext@arrow 0055{\arrowfill@\relbar\relbar\rightarrowtriangle}{#1}{#2}}
  \providecommand*\wto[2][]{%
  \ext@arrow 0055{\arrowfill@\relbar\relbar\rightarrowtriangle}{#1}{#2}}
\makeatother

\usepackage{amsmath}
\usepackage{amsfonts}
\usepackage{xcolor}
\usepackage{tikz}
\usetikzlibrary{backgrounds}
\usetikzlibrary{arrows}
\usetikzlibrary{shapes,shapes.geometric,shapes.misc}
\tikzstyle{tikzfig}=[baseline=-0.25em,scale=0.5]
\pgfkeys{/tikz/tikzit fill/.initial=0}
\pgfkeys{/tikz/tikzit draw/.initial=0}
\pgfkeys{/tikz/tikzit shape/.initial=0}
\pgfkeys{/tikz/tikzit category/.initial=0}
\pgfdeclarelayer{edgelayer}
\pgfdeclarelayer{nodelayer}
\pgfsetlayers{background,edgelayer,nodelayer,main}

\tikzstyle{none}=[inner sep=0mm]

\tikzstyle{label}=[inner sep=0.15mm, font={\footnotesize}]
\tikzstyle{vert}=[fill=black, draw=white, shape=circle, line width=0.25mm, tikzit shape=circle, inner sep=0.35mm]
\tikzstyle{label-small}=[inner sep=0.1mm, font={\scriptsize}]
\tikzstyle{label-white}=[fill=white, draw=none, shape=circle, inner sep=0.07mm, font={\scriptsize}]
\tikzstyle{vert-white}=[fill=white, draw=black, shape=circle, inner sep=0.35 mm]
\tikzstyle{vert-red}=[fill=red, draw=white, shape=circle, inner sep=0.4 mm, font={\scriptsize}, text=black]

\tikzstyle{Vert}=[fill=black, draw=black, shape=circle, minimum size=0.4em, inner sep=0.4pt, scale=0.6]
\tikzstyle{label-red}=[fill=white, draw=black, shape=circle, scale=0.4, color=red]
\tikzstyle{label small}=[fill=none, draw=none, shape=circle, scale=0.7, inner sep = 0.1]
\tikzstyle{bullet}=[fill=black, draw=black, shape=circle, scale=0.5]
\tikzstyle{G-vert}=[fill=white, draw=black, shape=circle, inner sep=0.7pt, minimum size=0.4em, scale=0.6]
\tikzstyle{label}=[inner sep=0.15mm, font={\footnotesize}]
\tikzstyle{vert}=[fill=black, draw=white, shape=circle, line width=0.15mm, tikzit shape=circle, inner sep=1 mm, scale = 1.5]
\tikzstyle{nodes}=[fill=white, draw=none, shape=circle, inner sep=0.7pt, font={\scriptsize}, minimum size=11pt]
\tikzstyle{label-s}=[inner sep=0.1mm, font={\scriptsize}]

\tikzstyle{label-w}=[fill=white, draw=white, shape=circle, inner sep=0.07mm, font={\scriptsize}]
\tikzstyle{charge}=[fill=white, draw={rgb,255: red,16; green,77; blue,209}, shape=circle, inner sep=0.1 em, scale=0.8]

\tikzstyle{label}=[inner sep=0.15mm, font={\footnotesize}]
\tikzstyle{vert}=[fill=black, draw=white, shape=circle, line width=0.3mm, tikzit shape=circle, inner sep=0.35mm,minimum size=6pt]
\tikzstyle{label-s}=[inner sep=0.1mm, font={\scriptsize}]
\tikzstyle{label-w}=[fill=white, draw=none, shape=circle, inner sep=0.07mm, font={\scriptsize}]
\tikzstyle{charge}=[fill=white, draw={rgb,255: red,16; green,77; blue,209}, shape=circle, inner sep=0.2 em, scale=0.6]
\tikzstyle{particle}=[fill=black, draw=black, shape=circle, inner sep=0.76 mm]
\tikzstyle{vert-white}=[fill = white, draw = black, shape = circle, inner sep = 1.89 pt]

\tikzstyle{blue}=[-, draw={rgb,255: red,49; green,149; blue,255}, line width=0.8pt]
\tikzstyle{blue thick}=[-, tikzit draw=blue, draw=blue, line width=1.7pt]
\tikzstyle{thick}=[-, ultra thick]
\tikzstyle{dot}=[-, dotted]
\tikzstyle{fade}=[-, opacity=0.5]
\tikzstyle{yellow}=[-, draw={rgb,255: red,208; green,208; blue,42}, line width=1.5]
\tikzstyle{arrow}=[->, line width=2em]
\tikzstyle{faded}=[-, opacity=0.5]
\tikzstyle{matching}=[-, draw=blue, line width=2.5pt, opacity=0.5]
\tikzstyle{orange}=[-, draw={rgb,255: red,255; green,128; blue,0}]
\tikzstyle{red arrow}=[draw=red, ->]
\tikzstyle{double edge}=[-, double, double distance=1 pt, draw=blue, line width=1.3, opacity=1]
\tikzstyle{single edge}=[-, opacity=1, draw=blue, line width=1.3]
\tikzstyle{dot}=[-, dotted]
\tikzstyle{red arrow thick}=[->, line width=0.8pt, draw=red]
\tikzstyle{dimer}=[-, line width=1.75 pt, draw=red]

\tikzstyle{dashed}=[-, densely dotted]
\tikzstyle{virtual}=[-, double]
\tikzstyle{RED}=[-, draw=red]
\tikzstyle{CYAN}=[-, draw=cyan]
\tikzstyle{BLUE}=[-, draw=blue]
\tikzstyle{LGREEN}=[-, draw=green]
\tikzstyle{DGREEN}=[-, draw={rgb,255: red,0; green,128; blue,128}]
\tikzstyle{PURPLE}=[-, draw={rgb,255: red,128; green,0; blue,128}]
\tikzstyle{ORANGE}=[-, draw={rgb,255: red,255; green,128; blue,0}]
\tikzstyle{MAGENTA}=[-, draw=magenta]
\tikzstyle{BLUE(dashed)}=[-, draw=blue, densely dotted]
\tikzstyle{GREEN}=[-, draw={rgb,255: red,0; green,208; blue,6}]
\tikzstyle{MAGENTA(dashed)}=[-, draw=magenta, densely dotted]
\tikzstyle{ORANGE(dashed)}=[-, draw={rgb,255: red,255; green,128; blue,0}, densely dotted]
\tikzstyle{RED(dashed)}=[-, draw={rgb,255: red,191; green,0; blue,64}, densely dotted]
\tikzstyle{PURPLE(dashed)}=[-, draw={rgb,255: red,128; green,0; blue,128}, densely dotted]
\tikzstyle{GREEN(dashed)}=[-, draw={rgb,255: red,0; green,208; blue,6}, densely dotted]
\tikzstyle{shade}=[-, opacity=0.4, draw={rgb,255: red,128; green,128; blue,128}, line width=6.5, fill=none, line cap=round,rounded corners]
\tikzstyle{hyper}=[-, fill={rgb,255: red,48; green,255; blue,214}, fill opacity=0.6, draw={rgb,255: red,18; green,229; blue,85}, tikzit fill={rgb,255: red,48; green,255; blue,214}]
\tikzstyle{hyper1}=[-, fill={rgb,255: red,230; green,138; blue,9}, draw={rgb,255: red,255; green,128; blue,0}, fill opacity=0.5]
\tikzstyle{hyper2}=[-, fill={rgb,255: red,128; green,179; blue,255}, draw={rgb,255: red,46; green,87; blue,115}, fill opacity=0.55]
\tikzstyle{new edge style 0}=[-, fill={rgb,255: red,245; green,255; blue,39}, draw={rgb,255: red,168; green,170; blue,22}, fill opacity=0.6]
\tikzstyle{blue}=[-, draw={rgb,255: red,49; green,149; blue,255}, line width=0.7 pt]
\tikzstyle{red}=[-, draw=red, line width= 0.7pt]
\tikzstyle{blue thick}=[-, tikzit draw=blue, draw=blue, line width=1.7pt]
\tikzstyle{thick}=[-, ultra thick]
\tikzstyle{dot}=[, dotted]
\tikzstyle{fade}=[-, opacity=0.5]
\tikzstyle{shade}=[-, opacity=0.23, fill={rgb,255: red,191; green,191; blue,191}, draw=none]

\tikzset{%
every path/.append style={line width = 0.8 pt}
}
\usepackage{ytableau}
\usepackage{mathtools, amssymb, old-arrows}
\usepackage[margin = 1in]{geometry}
\include{theorems}
\usepackage{tikz-cd}

\renewcommand{\leq}{\leqslant}

\usepackage{mathbbol}
\usepackage{amssymb}            

\DeclareSymbolFontAlphabet{\amsmathbb}{AMSb}

\newlength\friezelen 

\usepackage{mathrsfs} 
\usepackage{mathtools}
\usepackage{verbatim}
\usepackage{comment}
\usetikzlibrary{arrows}
\tikzset {->-/.style={decoration={markings, mark=at position .5 with {\arrow{latex}}}, postaction={decorate}}}
\tikzset {-->-/.style={decoration={markings, mark=at position .5 with {\arrow[scale=2]{latex}}}, postaction={decorate}}}

\usepackage{comment}
\usepackage{graphicx}
\usepackage{color}
\usepackage{etoolbox}

\usepackage{amsmath,amsthm,amssymb,graphicx,tikz,tikz-cd}
\usepackage{hyperref}
\hypersetup{
    colorlinks = true,
    citecolor = orange,%
}
\usepackage[noabbrev]{cleveref}
\usepackage{subcaption}

\makeatletter
\patchcmd{\@settitle}{\uppercasenonmath\@title}{}{}{}
\patchcmd{\@setauthors}{\MakeUppercase}{}{}{}
\patchcmd{\section}{\scshape}{}{}{}
\makeatother
\makeatletter
\patchcmd{\@maketitle}
  {\ifx\@empty\@dedicatory}
  {\ifx\@empty\@date \else {\vskip2ex 
  \centering\footnotesize\@date\par\vskip1ex}\fi
   \ifx\@empty\@dedicatory}
  {}{}
\patchcmd{\@adminfootnotes}
  {\ifx\@empty\@date\else \@footnotetext{\@setdate}\fi}
  {}{}{}
\makeatother

\DeclareMathOperator{\fH}{\mathfrak{H}}

\theoremstyle{plain}
\newtheorem{theorem}{Theorem}[section]

\newtheorem{lemma}[theorem]{Lemma}
\newtheorem{prop}[theorem]{Proposition}

\newtheorem{corollary}[theorem]{Corollary}

\theoremstyle{definition}
\newtheorem{remark}[theorem]{Remark}
\newtheorem{example}[theorem]{Example}
\newtheorem{definition}[theorem]{Definition}

\DeclareMathOperator{\spin}{spin}
\DeclareMathOperator{\id }{id}
\newcommand{\fh}[0]{\mathfrak{h}}

\newcommand{\fff}[0]{\mathbb{F}}
\newcommand{\bbb}[0]{\mathbb{B}}
\renewcommand{\geq}{\geqslant}

\newcommand{\bS}[0]{\overleftarrow{\mathfrak{S}}}

\newcommand{\sz}[0]{S_{\mathbb{Z}}}

\DeclareMathOperator{\bb}{{\bf B}}
\DeclareMathOperator{\ff}{{\bf F}}
\DeclareMathOperator{\qq}{\mathbb{Q}}
\DeclareMathOperator{\cc}{\mathbb{C}}
\DeclareMathOperator{\xx}{{\bf{x}}}
\DeclareMathOperator{\ms}{\mathfrak{S}}
\DeclareMathOperator{\mf}{\mathfrak{F}}
\DeclareMathOperator{\bR}{\overleftarrow{R}}
\renewcommand{\emptyset}{\varnothing}
\DeclareMathOperator{\len}{len}
\newcommand{\rtile}[0]{
\begin{tikzpicture}[scale=0.6, tikzfig, baseline = -0.5 em]
	\begin{pgfonlayer}{nodelayer}
		\node [style=none] (204) at (15.25, -1) {};
		\node [style=none] (205) at (15.25, 0.25) {};
		\node [style=none] (206) at (16.5, 0.25) {};
	\end{pgfonlayer}
	\begin{pgfonlayer}{edgelayer}
		\draw [style=red, rounded corners=0.2cm] (204.center) to (205.center)  to (206.center);
	\end{pgfonlayer}
\end{tikzpicture}
}

\newcommand{\jtile}[0]{
\begin{tikzpicture}[scale=0.6, tikzfig, baseline = 0.5 em]
	\begin{pgfonlayer}{nodelayer}
		\node [style=none] (201) at (14.25, 0.25) {};
		\node [style=none] (202) at (14.25, 1.5) {};
		\node [style=none] (203) at (13, 0.25) {};
	\end{pgfonlayer}
	\begin{pgfonlayer}{edgelayer}
		\draw [style=red, rounded corners=0.2cm] (203.center) to (201.center) to (202.center);
	\end{pgfonlayer}
\end{tikzpicture}
}

\title[Schubert Calculus and the Heisenberg Algebra]{\large Schubert Calculus and the Heisenberg Algebra}

\author[Sylvester W. Zhang]{Sylvester W. Zhang  \\  \begin{CJK}{UTF8}{gkai}
张文泽
\end{CJK}}

\thanks{\url{swzhang@umn.edu} University of Minnesota, School of Mathematics.}
\date{September 23, 2024}

\begin{document}

\begin{abstract}

We show that the Hilbert space with basis indexed by infinite permutations and the cohomology ring of the infinite flag variety can be seen as representations of the Heisenberg algebra, which are isomorphic using the back-stable Schubert polynomials. We give a model for infinite permutations as certain two dimensional fermions, generalizing the Maya diagram construction for partitions. Under this framework, the pipedream model for Schubert polynomials can be viewed as the Hamiltonian time evolution of the 2D fermions.  \end{abstract}

\maketitle
\setcounter{tocdepth}{1}
\tableofcontents
\setlength{\parindent}{0em}
\setlength{\parskip}{0.618em}

\tableofcontents
\section*{Introduction}
Originally appeared in string theory, the boson-fermion correspondence has found linkage to symmetric function theory, particularly through its application by the Kyoto school for deriving soliton solutions of the KP equations \cite{jimbo1983solitons}. In this framework, the Hilbert space spanned by the set of all Young diagrams $\mathbb{Y}$ is conceived as the fermionic Fock space $\ff$, while the ring of symmetric functions $\Lambda$ serves as the bosonic Fock space $\bb$. Then what Kac \cite{kac2013bombay} called the second part of the boson-fermion correspondence asserts that the map sending a partition $\lambda$ to its Schur function $s_\lambda$ forms an isomorphism as $\fh$-modules, where $\fh$ is the Heisenberg algebra.

Schur polynomials play an important rolein the context of Schubert calculus, representing the cohomology classes for Schubert varieties in a Grassmanian. They are special cases of \emph{Schubert polynomials}, which corresponds to Schubert cycles in the full flag variety $G/B$. We attempt to extend the story of boson-fermion correspondence from Grassmanians to full flag varieties.
It turns out that the most natural candidate for our special polynomials is the back-stable Schubert polynomials $\bS_w$ (introduced in  \cite{lam2021back}\footnote{see \cite{billey1995schubert,lee2019combinatorial,anderson2021schubert,anderson2021infinite} for similar constructions in other Lie types}), which correspond to Schubert classes in an infinite flag variety. The $\bS_w$'s are certain power series in infinitely many variables, which generate a ring $\bR$ called \emph{back-symmetric functions}. We define the analogue of fermionic Fock space to be the space spanned by $\sz$, the set of infinite permutations with finite support, denoted $\fff$, and define the analogue of bosonic Fock space to be $\bbb=\bR$.

We consider a generalized Heisenberg algebra $\fH$, containing the usual Heisenberg algebra $\fh$ as a subalgebra, is considered. The algebra $\fH$ is the analogue of ``bosons''. In fact, we will define bosonic operators on $\fff$ and $\bbb$ via actions of $\fH$. Our main result is the following, generalizing the work of \cite{nenashev2020differential} and \cite{lam2006combinatorial}.
\begin{theorem}
	Both $\fff$ and $\bbb$ carry a faithful $\fH$-action and the map $\Psi:w\mapsto \bS_w$ is an isomorphism of $\fH$-representations.
\end{theorem}

Let $H$ be the Hamiltonian operator, a standard result in the theory of vertex algebras is that the map sending $\ff\to \bb\otimes \ff$ given by $|\lambda\rangle \mapsto \langle \text{vac}|\exp(H)|\lambda\rangle$ preserves $\fh$-actions (see \cite{jimbo1983solitons}).
We adapt this approach, and realize the Stanley symmetric function as $\langle \id |\exp(H)|w\rangle$ where $H$ is a map $\fff\to \bb\otimes \fff^*$. This establishes an intermediate step between the classical BF correspondence and our generalization. This connection is depicted in the following commutative diagram.
\begin{equation}\label{eq:comm-diagram}
	\begin{tikzcd}[scale = 1.2]
$\ff=\mathbb{Y}$ \arrow[rrrrr,leftrightarrow, "\Phi:\lambda\leftrightarrow s_\lambda"] \arrow["\mathfrak{h}"', loop, distance=2em, in=125, out=55]                                                                                                                                                          &  &  &   &  &$\bb = \Lambda$ \arrow["\mathfrak{h}"', loop, distance=2em, in=125, out=55]                                                                    \\                                                                                                                                                                                                                                                                                          &  &  &   &  &                                                                                                                                               \\                                                                                                                                                                                                                                                                                        &  &  &   &  &                                                                                                                                              \\
$\fff=\sz$ \arrow[uuu, "\text{Edelman-Greene}"'] \arrow["\mathfrak{h}\subset { \mathfrak{H}}"', loop, distance=2em, in=305, out=235] \arrow[rrrrruuu, "\tilde\Phi:w\mapsto F_w"', two heads] \arrow[rrrrr, leftrightarrow,"\Psi:w\leftrightarrow \overleftarrow{\mathfrak{S}}_w"', shift right=2] &  &  &  &  &   $\bbb=\overleftarrow{R}$ \arrow["\mathfrak{h}\subset {\mathfrak{H}}"', loop, distance=2em, in=305, out=235] \arrow[uuu, "\eta_0"]
\end{tikzcd}
\end{equation}

The space $\ff$ contains 1-dimensional free fermions, and their evolutions give rise to fermionic lattices which can be viewed as an integrable vertex model for Schur polynomials \cite{zinn2009six}. Elements of $\fff$ are naturally modeled by special families of two dimensional fermions. These 2D fermions in $\fff$ evolve according to almost the exact same rules as in $\ff$, and they give rise to certain 3-dimensional lattice paths. We will show that these 3D lattices can be simplified to the pipedream model of Schubert polynomials (\cite{bergeron1993rc,fomin1996yang}).

The paper is structured as follows. In \Cref{sec:pre}, we review the basics of symmetric functions and Schubert calculus. In \Cref{sec:schur_fermion}, we review the classical BF correspondence for Schur polynomials (the arrow $\Phi$ in \Cref{eq:comm-diagram}). In \Cref{sec:H_acts_on_sz}, we define bosonic operators on $\sz$, and establish the map $\tilde \Phi:\fff\to\bb$ via the Hamiltonian operator, with proof deferred to \Cref{sec:proof}. In \Cref{def:2d-fermion}, we introduce 2D Maya diagrams for $\fff$, and show that the time evolution of these 2D fermions recovers the pipe-dream formula. In \Cref{sec:back-stable}, we introduce $\bR$ and its $\fH$-action, and finally complete the commutative diagram in \Cref{eq:comm-diagram}.
		
\section{Preliminaries}\label{sec:pre}

Let $S_n$ denote the symmetric group of permutations of $[n]$. Let $\sz$ denote the infinite symmetric group of permutations of $\mathbb{Z}$ with finite support, i.e. for any $w\in \sz$, there exists $M<N$ such that $w(i)=i$ for all $i\leq M$ and $i\geq N$. The group $\sz$ is generated by \emph{simple reflections} $\{s_i:i\in\mathbb{Z}\}$ with relations
\begin{align*}
	s_i^2 &= 0&&\forall i\\
	s_is_j &= s_js_i&&\forall |i-j|>1\\
	s_is_{i+1}s_i &= s_{i+1}s_i s_{i+1}&&\forall i
\end{align*}

We say $a_1a_2\cdots a_n$ is a \emph{word} of $w\in\sz$ if $w=s_{a_1}\cdots s_{a_n}$, and a word is said to be \emph{reduced} if it uses minimal number of generators. Define $R(w)$ the set of all reduced words of $w$. If $w=s_{a_1}\cdots s_{a_n}$ is reduced then we define the length of $w$ to be $n$, denoted $\ell(w) = n$. 

In terms of one-line notations of permutations, right multiplication by simple reflections  $w s_i$ corresponds to swapping the $i$-th and $(i+1)$-th positions of $w$. Analogously, we define \emph{reflections} $t_{ij}$ to be the permutation that swaps the $i$-th and $j$-th position.

We use the notation $w\doteq u_1\cdots u_r$ for $w=u_1\cdots u_r$ and $\ell(w)=\ell(u_1)+\cdots+\ell(u_r)$.
Define two partial orders on $\sz$. We say that $v$ is covered by $w$ in the \emph{weak Bruhat order}, denoted $v\lessdot w$ if $w \doteq s_i v$. We write $v<w$ if $v$ is less that $w$ in the weak order. We say that $v$ is covered by $w$ in the \emph{strong Bruhat order}, denoted $v\lessdot_{\text{strong}}w$ if $w = t_{ij}v$ and $\ell(w)=\ell(v)+1$. We further define the $k$-Bruhat order to be the sub-order of the strong order as follows. We say that $v$ is covered by $w$ in the $k$-Bruhat order, denoted $v\lessdot_k w$, if $w = t_{ij} v $ with $\ell(w)=\ell(v)+1$, and $i\leq k <j$.

\subsection{Schubert Calculus}\label{sec:polynomials}

Let $G=GL_n(\mathbb{C})$ and $B\subset G$ the Borel subgroup consisting of upper triangular matrices of non-zero determinant. The flag variety $G/B=\{V_1\subset V_2\subset\cdots \subset V_n = \cc^n:\dim(V_i)=i\}$ has finitely many $B$-orbits, which are parameterized by elements of the Weyl group $S_n$. For $w\in S_n$, define the Schubert variety $X_w$ to be the closure of the $B$-orbits corresponding to $w$.

The cohomology ring of the flag variety $H^*_\mathbb{Z}(G/B)$ has basis given by cohomology classes of the Schubert varieties $[X_w],w\in S_n$. The structure of $H^*_\mathbb{Z} (G/B)$ was first revealed by Berstein-Gelfand-Gelfand \cite{bernstein1973schubert} and Demazure \cite{demazure1974desingularisation} geometrically using divided difference operators $\partial_i$. The ring $H^*_\mathbb{Z}(G/B)$ is isomorphic to the quotient of the polynomial ring $R_n:=\mathbb{Z}[x_1,\cdots,x_n]/I$ where $I=\langle e_1,\cdots,e_n\rangle$ is the ideal generated by elementary symmetric functions (defined in the next section). Lascoux and Sch\"utzenberger \cite{alain1982polynomes} identified polynomial representatives of Schubert classes called \emph{Schubert polynomials}, which can be studied purely combinatorially. 

Let $\partial_i:\zz[x_1,x_2,\cdots]\to \zz[x_1,x_2,\cdots]$ be the divided difference operators defined by
\[\partial_i(f) = {f-s_i(f)\over x_i - x_{i+1}}\]
where $s_if$ is the polynomial obtained from $f$ by swapping $x_i\leftrightarrow x_{i+1}$. Deferring their definition to the next section, Schubert polynomials are the only family of polynomials `compatible' with these operators, i.e. they satisfy the following conditions
\begin{align*}
\mathfrak{S}_{\id} &= 1\\	
\partial_i\mathfrak{S}_{w} & =  \begin{cases}
	\mathfrak{S}_{w s_i}&\text{if }w s_i<w\\
	0&\text{otherwise}
\end{cases}\\
\mathfrak{S}_w& \text{ is homogeneous of degree }\ell({w}).
\end{align*}

In case of the Grassmanian $Gr_k(\cc^n) = \{V\subset \cc^n:\dim(V)=k\}$, its Schubert cells are indexed by $k$-\emph{Grassmanian permutations}, which are permutations with at most one descent at position $k$. Let $\lambda(w)$ denote the bijection between Grassmanian permutations and partitions. The Schubert polynomials of a $k$-Grassmanian permutation are certain Schur polynomials in the first $k$ variables, i.e. $\mathfrak{S}_w(x_1,\cdots,x_{n-1}) = s_{\lambda(w)}(x_1,\cdots,x_k)$.

\subsection{Tableaux and Symmetric Functions}
We now define the relevant polynomials, beginning with some basic combinatorial objects.
\ytableausetup{smalltableaux}

We define a partition of $n$ to be a sequence of weakly decreasing numbers $\lambda = \lambda_1\geq \lambda_2\geq \cdots$ summing up to $n$, denoted $\lambda\vdash n$. A partition is represented by a Young diagram, in which the $i$-th row contains $\lambda_i$ boxes. For example, the Young diagram of $\lambda = (3,3,1)$ is 
$$\ydiagram{3,3,1}$$
We denote $\mathbb{Y}$ the vector space spanned by all partitions.

If $\mu$ is a sub-shape of $\lambda$, define the skew Young diagram $\lambda/\mu$ to be the set-theoretic difference $\lambda -\mu$. For example, the following is the skew shape $(3,3,1)/(1,1)$.
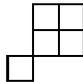
\[\begin{ytableau}
	\none& & \\
	\none& & \\
	\\
\end{ytableau}\]
A skew shape that is connected and does not contain a $2\times 2$ box is called a \emph{ribbon}, such as
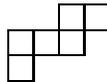
\[\begin{ytableau}
	\none& \none & & \\
 & & \\
	\\
\end{ytableau}\]
A skew shape that does not contain two boxes in the same column is called a \emph{horizontal strip}, e.g.
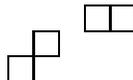
\[\begin{ytableau}
	\none& \none &\none & &\\
 \none & & \none\\
	\\
\end{ytableau}\]

A polynomial in $\qq[x_1,x_2,\cdots,x_n]$ is symmetric if it is invariant under swapping any pair of variables $x_i\leftrightarrow x_j$. We denote $\Lambda(\xx_n)$ to be the ring of symmetric polynomials in $n$ variables, and define $\Lambda(\xx)$ the ring of symmetric functions to be the limit of $\Lambda(\xx_n)$ as $n$ goes to infinity. The ring $\Lambda(\xx)$ has several important basis. 
\begin{enumerate}
	\item The elementary symmetric functions
	\[e_k(\xx)=\sum_{i_1<i_2<\cdots<i_k}x_{i_1}\cdots x_{i_k}\quad \quad e_\lambda = e_{\lambda_1}e_{\lambda_2}\cdots e_{\lambda_n}\]
	\item The complete homogeneous symmetric functions
	\[h_k(\xx)=\sum_{i_1\leq i_2\leq \cdots \leq i_n}x_{i_1}\cdots x_{i_k}\quad\quad h_\lambda = h_{\lambda_1}h_{\lambda_2}\cdots h_{\lambda_n}\]
	\item The power-sum symmetric functions\[p_k(\xx)=x_1^k+x_2^k+x_3^k\cdots\quad\quad\ \ \ \ p_\lambda = p_{\lambda_1}p_{\lambda_2}\cdots p_{\lambda_n}\]
\end{enumerate}
The homogeneous symmetric functions expand in the power-sum basis in the following way.
\[h_m=\sum_{\alpha\vdash m} {p_\alpha\over z_\alpha}\]
where the sum range over partitions of $m$, and $z_\alpha=1^{m_1(\alpha)}m_1(\alpha)!2^{m_2(\alpha)} m_2(\alpha)!\cdots$ where $m_i(\alpha)=\#\{i:\lambda_j=i\}$.

Last but not the least, Schur functions are generating function of \emph{semistandard Young tableaux} ($SSYT$), defined as follows. A SSYT of shape $\lambda$ is a filling of the Young diagram of $\lambda$ with numbers such that every row is weakly increasing and every column is strictly increasing. Denote $\text{SSYT}(\lambda)$ the set of all such tableau. For $T\in\text{SSYT}(\lambda)$, define $\xx^{\wt(T)}:=\prod_i x_i^{\text{number of }i\text{ in }T}$. For example, the following is a SSYT with $\xx^{\wt(T)}=x_1^2x_2^2x_3^3$.
\[
\begin{ytableau}
1& 1&2\\
2& 3&3\\
3
\end{ytableau}
\]
Then the Schur functions are defined to as follows
\[s_{\lambda}(\xx)=\sum_{T\in\text{SSYT}(\lambda)}\xx^{\wt(T)}\]

\subsection{(Stable) Schubert Polynomials}
We next define Schubert polynomials and their stable limit. They were originally defined using divided difference operators by Lascoux-Sch\"utzenberger. Here we will present an alternative definition using increasing factorizations.

We say that a permutation $v$ is increasing if it has a reduced word $v=s_{a_1}s_{a_2}\cdots s_{a_n}$ such that $a_1>a_2>\cdots >a_n$. We say that $v_1v_2\cdots v_n$ is a increasing factorization of $w$ if $\ell(w)=\ell(v_1)+\cdots+\ell(v_n)$ and each of $v_i$ is increasing. Denote $I(w)$ the set of all increasing factorization of $w$.Then we define $\mathfrak{S}_w$ the Schubert polynomial of $w\in S_n$ as follows.
\[\mathfrak{S}_w=\sum_{\substack{ v_1v_2\cdots v_{n-1} \in I(w) \\
v_i(j)=j\text{ if } j<i} }x_1^{\ell(v_1)}x_2^{\ell(v_2)}\cdots x_{n-1}^{\ell(v_{n-1})}\]

Define the shifting operator on permutations $\tau: S_n\to S_{n+1}$ by
$\tau(w)_1=1$ and $\tau(w)_i=w_i+1$ for $i\geq 2$. We can extend the definition of $\tau$ to be an automorphism on $S_\mathbb{Z}$. The Stanley symmetric functions, or stable Schubert polynomials $\mathfrak{F}_w$ are the stable limit of Schubert polynomials:
\[\mf_w = \lim_{n\to\infty} \ms_{\tau^n(w)}\]

So far (stable) Schubert polynomials are defined for $S_{\mathbb{Z}_>0}=\bigoplus_{n=1}^{\infty} S_n$. Following Lam-Lee-Shimozono's notion of back stable Schubert polynomials \cite{lam2021back}, we can extend the definition of $\mf_w$ to $S_{\zz}$. In particular, we define
\[\mf_w = \sum_{v_1 v_2\cdots v_r \in I(w)} x_1^{\ell(v_1)}x_2^{\ell(v_2)}\cdots x_{r}^{\ell(v_{r})}\]
for $w\in S_{\zz}$.

\begin{remark}Note that Schubert polynomials are, in general, not symmetric polynomials. 
Unlike Schubert polynomials, in the definition of Stanley symmetric functions, the length of the increasing factorizations are unrestricted, therefore it is a symmetric polynomial in infinitely many variables. In the language of Lam-Lee-Shimozono's back stable Schubert polynomials, our definition of Stanley symmetric function $\mf_w$ is the same as taking the back stable Schubert polynomial $\bS_w$, forgetting the non-symmetric part, then reversing the negative variables to positive variables, i.e. $\mf_w=\eta_0(\bS_w)|_{x_i\to x_{1-i}}$. Although later on we will switch back to negative variables in \Cref{sec:back-stable} concerning back symmetric functions, at this point we do not supply definitions regarding back-stable functions to avoid confusion.
\end{remark}

\begin{remark}
	Here the definitions of Schubert polynomials and Stanley functions are inverse to what is in the literature. Our definition of $\mathfrak{S}_w$ in the standard convention is the Schubert polynomial of $w^{-1}$. For Stanley symmetric functions, our convention is inverse to the convention of stable limit definition but is the same as Stanley's original definition. The benefit of using inverse permutation will become apparent in \Cref{sec:2dfermion,sec:proof} regarding Maya diagrams of 2D fermions.\end{remark}

\section{Classical Boson-Fermion Correspondence}\label{sec:schur_fermion}
In this section, we briefly survey the boson-fermion correspondence, and we refer to \cite{miwa2000solitons,zinn2009six} for details. 
To start, we define abstractly bosons and fermions as follows.
\begin{definition}[fermions]The \emph{Clifford algebra} $\mathfrak{C}$ is generated by $\{\psi_i:i\in\zz\}\cup \{\psi^*_i:i\in\zz\}$ satisfying anti-commutation relations
\[\{\psi_i,\psi_j\}=0\quad \{\psi^*_i,\psi^*_j\}=0\quad \{\psi_i,\psi_j^*\} = \delta_{i+j,0}. \]
The elements $\varphi_i,\varphi_i^*$ are called \emph{fermions}.
\end{definition}
\begin{definition}[bosons]
Let $\mathfrak{W}$ be the algebra generated by $\{\alpha_i,i\in\zz\}$ and central element $c$ with relations $[\alpha_i,\alpha_j]=c\delta_{i+j,0}.$
	The \emph{Heisenberg algebra} $\fh$ is defined to be $\mathfrak{W}/(c-1)$. The $\fh$ elements $\{\alpha_i\}$ are \emph{bosons}.	\end{definition}

 Let ${\bf B} =\mathbb{C}[t_1,t_2,\cdots]$ be the bosonic Fock space and define the fermionic Fock space to be $\ff=\Lambda = \{|\lambda\rangle : \lambda\vdash n, n \in \mathbb{N}\}$ the space of all partitions. Note that these are both the charge-0 part of what are traditionally know as the bosonic and fermionic spaces, and we restrict to this case for simplicity.
 We also define a dual space $\ff^*=\{\langle \lambda|:\lambda \vdash n, n\in\mathbb{N}\}$ and define the pairing $\langle \lambda |\mu\rangle = \delta_{\lambda,\mu}$. 
 
 Both $\bb$ and $\ff$ are $\mathcal H$-modules, where the actions are given as follows. 
 Let $\alpha_i^{\bb}$ and $\alpha_i^{\ff}$ denote the action of $\fh$ on $\bb$ and $\ff$ respectively, then we have
 \[\alpha_i^{\bb}(f) = \begin{cases}
 	{\partial\over \partial t_m} f&m>0\\
 	-mt_{-m}f&m<0\\
 	f&m=0
 \end{cases}\]
 \[ \alpha_i^{\ff} |\lambda\rangle =\begin{cases}\displaystyle
 	\sum_{\lambda/\mu\text{ is a }k\text{-ribbon}}(-1)^{\text{ht}(\lambda/\mu)-1}|\mu\rangle&k>0\\
 	\displaystyle\sum_{\mu/\lambda\text{ is a }k\text{-ribbon}}(-1)^{\text{ht}(\mu/\lambda)-1}|\mu\rangle&k<0\\
 	|\lambda\rangle &k=0
 \end{cases}\]
 where $\text{ht}(\lambda/\mu)$ is the height of the ribbon $\lambda/\mu$.
 
 Define the Hamiltonian on $\ff$ to be $H[t] = \sum_{i>0} t_i\alpha_i^{\ff}$ and consider its exponential $e^{H[t]}$, which can be viewed as a map from $\ff$ to $\ff^*\otimes\bb$. We can extend the pairing $\langle\mu|\lambda  \rangle $ to be a map from $\ff\times \ff^*$ to $\bb$, so that $\langle \mu|e^{H[t]}|\lambda\rangle \in \bb$.

The following theorem what Kac called the second part of the boson-fermion correspondence.

\begin{theorem}[Jimbo-Miwa \cite{jimbo1983solitons}]\label{thm:bf-schur}
Let $e^{H[t]}$ be defined as above, we have
\begin{enumerate}
	\item[\textup{(1)}] The map $\Phi: \ff\to \bb$ defined by $\Phi(\lambda)=\langle 0|e^{H[t]}|\lambda\rangle $ is an $\fh$-module isomorphism.
	\item[\textup{(2)}] They are skew Schur functions up to the Miwa transformation $t_i\to p_i/i$, i.e.
	\[\langle \mu|e^{H[t]}|\lambda\rangle\big|_{t_i\mapsto p_i/i} =s_{\lambda/\mu}(x_1,x_2,\cdots)\]
\end{enumerate}
\end{theorem}
We will write $\tilde H(x)$ the Hamiltonian after Miwa transform, and define $\phi(x)=\sum_{k>0} {x^k\over k}\alpha^{\ff}_k$ so that $\tilde H(x) = \sum_{i >1} \phi(x_i)$. Thus $e^{\tilde H(x)} = \prod_{i\geq 1} e^{\phi(x_i)}$ and
\[s_\lambda(\xx) =\langle 0|\cdots e^{\phi(x_2)}e^{\phi(x_1)}|\lambda \rangle.\]
From this perspective, the symmetry of Schur functions is a consequence of the commutativity $e^{\phi(x_i)}e^{\phi(x_{j})}=e^{\phi(x_j)}e^{\phi(x_i)}$. 

The boson-fermion correspondence for Schur functions is equivalent to the Murnagan-Nakayama rule. 
\begin{prop}[Murnagan-Nakayama rule]
	\[p_ks_\lambda = \sum_{\mu/\lambda\text{ is a }k\text{-ribbon}}(-1)^{\textup{ht}(\mu/\lambda)-1}s_{\mu}\]
\end{prop} 
\begin{proof}
	Recall that for $k>0$, the action of $\alpha_{-k}^{\ff}$ is multiplication by $kt_k=p_k$. By \Cref{thm:bf-schur}, we have
	\[\alpha_{-k}^{\ff}s_\lambda=\langle \emptyset|e^{H[t]}|\alpha_{-k}^{\ff} \lambda\rangle\]
	which is exactly the right-hand-side.
\end{proof}

For a sequence of numbers $\sigma=(\sigma_1,\cdots,\sigma_n)$, denote $\alpha^{\ff}_\sigma=\alpha^{\ff}_{\sigma_1}\cdots \alpha^{\ff}_{\sigma_n}$ and $\alpha^{\ff}_{-\sigma}=\alpha^{\ff}_{-\sigma_1}\cdots \alpha^{\ff}_{-\sigma_n}$.

\section{Bosonic Operators on $S_\mathbb{Z}$.}
\label{sec:H_acts_on_sz}
Let us begin by declaring the Hilbert space with basis indexed by $\sz$ to be the \emph{2-dimensional fermionic Fock space}, and denote by $\mathbb{F}$. The fact that they are 2D particles will be illustrated later. Slightly out of the blue, we define a `larger' Heisenberg algebra $\mathfrak{H}$, with generators $\{\alpha_i:i\in\mathbb{Z}_{\geq 0}\}\cup\{\alpha_{-i,k}:i\in\zz_{\geq 0},k\in\zz\}$ subject to the relations
\[[\alpha_i,\alpha_j]=0=[\alpha_{-i,k},\alpha_{-j,l}]\quad [\alpha_i,\alpha_{-j,k}]=i\delta_{i,j}\]
for all $i,j\in\mathbb{Z}_{\geq 0}$ and $k,l\in\zz$.

\begin{definition}\label{def:ribbon}
	We say a reduced word is $s_{a_1}\cdots s_{a_l}$ a \emph{weak-ribbon} of size $l$, if	\begin{enumerate}
		\item $A=\{a_1,\cdots,a_l\}$ is a consecutive set, i.e. there is a surjection from $A$ to some interval.
		\item $a_1\cdots a_l$ avoids the patterns $2132$ and $2312$. 
	\end{enumerate}
Define the \emph{spin} of a ribbon to be $\spin(a_1\cdots a_l)=\#\{k: \exists (i>j) \text{ s.t. } a_i=k=a_j-1\}$.
\end{definition}

\begin{example}
	The word $s_1s_2s_1s_4s_3$ is a weak-ribbon of spin $2$. The word $s_3s_2s_5s_4s_1$ is a weak-ribbon of spin $1$. The word $s_1s_2s_4s_1$ is not a ribbon because the index set $\{1,2,4,1\}$ is not consecutive. The word $s_3s_1s_2s_1s_4s_3$ is not a weak-ribbon because it contains the pattern $2132$ ($\underline{s_3}s_1\underline{s_2}s_1\underline{s_4s_3}$).
\end{example}

The next lemma, whose proof is omitted, allows the extension of the definition of ribbon to permutations.
\begin{lemma}
\label{lem:ribbon_equivalent}
If a weak-ribbon $s_{a_1}\cdots s_{a_l}$ is a reduced word for $v$, then all reduced words of $v$ are weak-ribbons of equal spin. 
\end{lemma}
We next define a dual version of ribbon called strong ribbon.

\begin{definition}\label{def:co-ribbon}
Fix a $k\in\mathbb{Z}$. For $u,w\in\sz$,
	if $w^{-1}u$ is a $(r+1)$-cycle and $u = w t_{a_1,b_1}t_{a_2,b_2}\cdots t_{a_r,b_r}$ such that $\ell (wt_{a_1,b_1}\cdots t_{a_i,b_i})=\ell(w)+i$ with $a_i\leq k<b_i$ for all $1\leq i\leq r$, then we say that $u/ w $ is a $k$-strong-ribbon of size $r$. Define the \emph{spin} of $u/ w$ to be $\spin(u/ w) = \#\{i\leq k:w^{-1}u(i)\neq i\}-1$.
\end{definition}

\begin{remark}
	Fomin and Greene \cite{fomin1998noncommutative}
	defined \emph{hook words} to be certain $w$ which has a reduced word of the form $s_{a_1}\cdots s_{a_r}s_{b_1}\cdots s_{b_l}$ such that $a_1>a_2>\cdots a_r<b_1\leq b_2\leq \cdots \leq b_l$ and $\{a_1,\cdots,a_r,b_1,\cdots,b_l\}$ is consecutive. This is equivalent to our definition of weak-ribbons, which we left as an exercise to the curious reader. On the other hand, the definition of $k$-strong-ribbon first appeared in the work of Morrison-Sottilie \cite{morrison2014murgnahan}.
\end{remark}

\begin{definition}\label{def:operators}
	Fix $k\in\mathbb{Z}$. Define operators $\{\alpha_n^{\fff},\alpha_{-n,k}^{\fff}:n\in\mathbb{Z}_{>0},k\in\zz\}$ acting on $\fff=\text{span}({\sz})$ as follows. \begin{align}
\alpha_n^{\fff} |w\rangle &= \sum_{\substack{w = u\cdot v\\\ell(w)-\ell(u)=\ell(v)=n\\ v\text{ is a weak ribbon}}  } (-1)^{\spin(v)}|u\rangle\\
\alpha_{-n,k}^{\fff} |w\rangle &= \sum_{\substack{\ell(u) =\ell(w)+n \\ w/u \text{ is a $k$-strong-ribbon}}} (-1)^{\spin(w/u)}|u\rangle
\end{align}
And $\alpha_0^{\fff}|w\rangle =|w\rangle$.

\end{definition}

\begin{theorem}\label{thm:heisenbergS}
The operators $\alpha_{n}^\fff,\alpha_{-n,k}^\fff$ satisfy the relations of $\fH$. Therefore $\fff$ is a faithful $\fH$-representation.
\end{theorem}

\begin{remark}
	By choosing a subalgebra $\fh\subset \fH$, one can also view $\fff$ as an $\fh$-representation. Note also that we are abusing the notation $\alpha_i$ when $i\geq 0$ to denote both elements in $\fh$ and $\fH$, since the actions of $\alpha_i\in\fh$ and $\alpha_i\in\fH$ on $\fff$ or $\ff$ are the same. \end{remark}

	Note that $\bb=\Lambda$ as a $\fH$-representation is not faithful, since the actions of $\alpha_{-i,k}$ are identical for all $k$. It turns out that the natural choice for a faithful representation of $\fH$ would be the ring of ``back-stable symmetric functions''. We will expand on this perspective in \Cref{sec:back-stable}.

\begin{theorem}\label{thm:stanley_hamiltonian}
	Define the Hamiltonian $H[t] = \sum_{i=1}^{\infty} t_i\alpha_i^\fff$, then the polynomials
	\[\langle \id |e^{H[t]} | w\rangle \in \mathbb{Q}[t_1,t_2,\cdots]\]
	are equal to the Stanley symmetric functions after the Miwa transformation $t_i \mapsto {p_i / i}$. In other words, if $\phi(x) = \sum_{i=1}^{\infty} {x^i\over i}H_i$, so that $\tilde H(x) = \sum_{k=1}^{\infty} \phi(x_k)$. Then
	\[\langle \id |\cdots e^{\phi(x_2)}e^{\phi(x_1)}|w\rangle=\mf_w(\xx)\]equals to the Stanley symmetric function defined in \Cref{sec:polynomials}. Furthermore, the map $\Phi': \fff\to\bb,|w\rangle \mapsto \mf_w(\xx)$ is a surjective $\fh$-module homomorphism.
\end{theorem}

\subsection{Applications}
We first look at some immediate applications of the bosonic operators on $\mf_w$, and defer the proofs to \Cref{sec:proof}.

\subsubsection{Murnagan-Nakayama Rule}
Define
\[\chi_w^\alpha=\sum_{\substack{v_1v_2\cdots v_k = w\\\text{each }v_i\text{ is a weak ribbon}\\\ell(v_i)=\alpha_i}}(-1)^{\spin(v_1)+\cdots+\spin(v_k)}\]
where $w\in\sz$ and $\alpha$ a partition of $\ell(w)$.
By unravelling the exponential, we have the following Murnagan-Nakayama rule, first proved by Fomin and Greene \cite{fomin1998noncommutative}.
\begin{prop}
	\[\mf_w =\sum_{\alpha \vdash \ell(w)}{\chi_w^\alpha\over z_\alpha}p_\alpha\]
\end{prop}
\begin{proof}
	By direct calculation of $\langle 0|e^{\tilde H(x)}|w\rangle$.
\end{proof}

There is also a dual version of the Murnagan-Nakayama rule, which is the stable limit of Morrison's \cite{morrison2014murgnahan} MN rule for Schubert polynomials.
\begin{prop}
	\[p_m \mf_w=\sum_{\substack{u = w\cdot \sigma \\ \ell(u)-\ell(w) = m\\ \sigma\text{ is a }k\text{-strong-ribbon}}}(-1)^{\spin(\sigma)}\mf_u\]
\end{prop}
\begin{proof}
	Take any $k\in\zz$, we have $p_m \mf_w=\langle 0|e^{\tilde H(x)}\alpha_{-m,k}^{}|w\rangle$ which equals to the r.h.s by definition.
\end{proof}
\subsubsection{Edelman-Greene coefficients}
Edleman and Greene developed an insertion algorithm which proved that $\mf_w$ expands positively in the Schur basis. That is, there are $j_{\lambda,w}\in\mathbb{Z}_{\geq 0}$ such that
\[\mf_{w}=\sum_{\lambda}j_{\lambda,w}s_\lambda\]

We define the Stanley operator to be
\[F_w^{(k)}:=\sum_{\sigma\vdash\ell(w)}{\chi_w^\sigma\over z_\sigma}\alpha_{-\sigma,k}^{\fff},\]
where $\alpha_{-\sigma,k}^\fff=\alpha_{-\sigma_1,k}^\fff\cdots \alpha_{-\sigma_n,k}^\fff$.
Interplaying with the boson-fermion correspondence, we have
\[\mf_w\cdot 1=\langle \id|e^{\tilde H(x)}F_w^{(k)}|\id\rangle\]
Since applying $\alpha_{-i,k}$'s on the identity permutation only produces chains in strong $k$-Bruhat order, we know that applying the Stanley operator on the identity permutation $F_w^{(k)}|\id\rangle$ will be a linear combination of $k$-Grassmanian permutations. This is indeed expanding Stanley in Schur basis.
Let $S_{\zz}^{(k)}$ denote the set of $k$-Grassmanian permutations, then
	$$F_w^{(k)}|\id\rangle = \sum_{\sigma\vdash \ell(w)}{\chi_w^\sigma\over z_\sigma}\alpha_{-\sigma,k}^\fff|\id\rangle=\sum_{u\in S_{\mathbb{Z}}^{(k)}}j_{\lambda(u),w}|u\rangle$$

More specifically, let $j_{\lambda,w}^\mu$ be the structure constant of the product of a Schur function and a Stanley function expanded in the Schur basis, i.e.
\[s_\lambda \mf_w=\sum_{\mu}j_{\lambda,w}^\mu s_\mu,\]
Then we have the following proposition.
\begin{prop}
	For $w\in \sz$ and $u\in\sz^{(k)}$, we have $$F_w ^{(k)}|u\rangle  = \sum_{v\in\sz^{(k)}} j_{\lambda(u),w}^{\lambda(v)} |v\rangle $$
\end{prop}

One can also define the adjoint of the Stanley operators,
\[F_w^\perp=\sum_{\sigma\vdash\ell(w)}{\chi_w^{\sigma}\over z_\sigma}\alpha^\fff_\sigma\]
where $\alpha_\sigma^\fff=\alpha_{\sigma_1}^\fff\cdots\alpha_{\sigma_l}^\fff$. In case when $w$ is Grassmanian permutation, $F^\perp_w$ can be understood as taking derivative in the Schur basis of $\bb$ with respect to $s_{\lambda(w)}$. Therefore applying $F_w^\perp$ on $u$ will extract the Edelman-Greene coefficients.

\section{Two Dimensional Fermions}
\label{sec:2dfermion}
\subsection{Schur polynomials and the 5 vertex models} We will now review the classical results about vertex models for Schur polynomials. References for material in the section can be found in \cite{zinn2009six,brubaker2011schur}.

\begin{figure}[h]
\begin{tikzpicture}[scale=0.65]
	\begin{pgfonlayer}{nodelayer}
		\node [style=none] (0) at (-7, 2) {};
		\node [style=none] (1) at (-7, -3) {};
		\node [style=none] (2) at (-2, 2) {};
		\node [style=none] (3) at (-6, 1) {};
		\node [style=none] (4) at (9.5, 0) {};
		\node [style=none] (5) at (-6, 0) {};
		\node [style=none] (6) at (-4, 1) {};
		\node [style=none] (7) at (-4, 2) {};
		\node [style=none] (8) at (-7, 0) {};
		\node [style=vert] (9) at (-7, -2.5) {};
		\node [style=vert] (10) at (-7, -1.5) {};
		\node [style=vert] (11) at (-7, -0.5) {};
		\node [style=vert] (12) at (-6, 0.5) {};
		\node [style=vert] (13) at (-4, 1.5) {};
		\node [style=vert-white] (14) at (-6.5, 0) {};
		\node [style=vert-white] (15) at (-5.5, 1) {};
		\node [style=vert-white] (16) at (-4.5, 1) {};
		\node [style=vert-white] (17) at (-3.5, 2) {};
		\node [style=vert-white] (18) at (-2.5, 2) {};
		\node [style=none] (19) at (-0.5, 0) {};
		\node [style=none] (20) at (-7.5, 2.5) {};
		\node [style=none] (21) at (-3.5, -1.5) {};
		\node [style=none] (22) at (4.5, 0.75) {};
		\node [style=none] (23) at (4.5, -0.75) {};
		\node [style=vert] (24) at (4, 0) {};
		\node [style=vert] (25) at (2, 0) {};
		\node [style=vert] (26) at (1, 0) {};
		\node [style=vert] (27) at (0, 0) {};
		\node [style=vert-white] (28) at (3, 0) {};
		\node [style=vert-white] (29) at (5, 0) {};
		\node [style=vert-white] (30) at (6, 0) {};
		\node [style=vert-white] (31) at (8, 0) {};
		\node [style=vert-white] (32) at (9, 0) {};
		\node [style=vert] (33) at (7, 0) {};
	\end{pgfonlayer}
	\begin{pgfonlayer}{edgelayer}
		\draw (1.center) to (0.center);
		\draw (0.center) to (2.center);
		\draw (8.center) to (5.center);
		\draw (5.center) to (3.center);
		\draw (3.center) to (6.center);
		\draw (6.center) to (7.center);
		\draw (19.center) to (4.center);
		\draw [style=dot] (20.center) to (21.center);
		\draw [style=dot] (22.center) to (23.center);
	\end{pgfonlayer}
\end{tikzpicture}	
\caption{Example of a partition and its Maya diagram. The dotted line indicated the main diagonal of the Young diagram, and the center of the Maya diagram.}
\label{fig:partition-maya}
\end{figure}
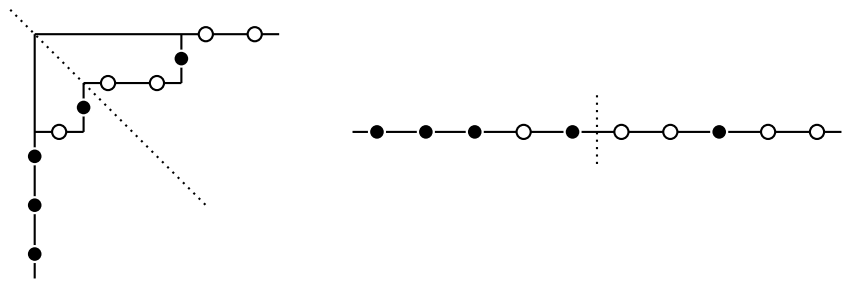

A \emph{Maya Diagram} is an infinite string of black and white nodes with a distinguished center, such that there are infinitely many black notes to the left, and infinitely many white nodes to the right. There is a bijection between Maya diagrams and partitions as follows. Take $\lambda$ Young diagram of a partition, along its boundary extended to infinity, we place a black node at each vertical edge and a white node at each horizontal edge and then a Maya diagram can be read off by straightening the boundary. We declare that the center of the Maya diagram is the main diagonal of the Young diagram. See \Cref{fig:partition-maya} for example.

More specifically, Maya diagrams can be viewed as maps from $\zz$ to $\{0,1\}$ where $1$ corresponds to black nodes and $0$ white. For a Maya diagram $f:\zz\to \{0,1\}$, there is a unique $k\in\zz$ such that $\#\{f(i)=0:i\leq k\}=\#\{f(j)=1:j>k\}$. Then $k$ is the center of $f$.
\begin{remark}
Recall that partitions are in bijection with $k$-Grassmanian permutations. From Maya diagrams, one can read off the corresponding inverse $k$-Grassmanian permutation easily. Fix a $k$, then label the black nodes of the Maya diagram from left to right by $\cdots, k-2,k-1,k$ and label the white nodes of the Maya diagram from left to right by $k+1,k+2,k+3,\cdots$. The resulting string of numbers will be the inverse of the $k$-Grassmanian permutation. For example the inverse $3$-Grassmanian permutation corresponding to $\lambda=(3,1)$ is:
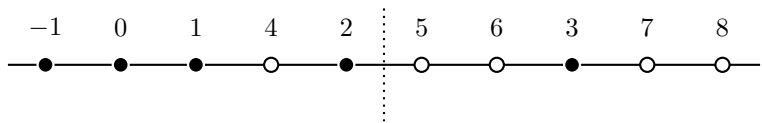
\[\begin{tikzpicture}
	\begin{pgfonlayer}{nodelayer}
		\node [style=none] (4) at (5, 0) {};
		\node [style=none] (19) at (-5, 0) {};
		\node [style=none] (22) at (0, 0.75) {};
		\node [style=none] (23) at (0, -0.75) {};
		\node [style=vert] (24) at (-0.5, 0) {};
		\node [style=vert] (25) at (-2.5, 0) {};
		\node [style=vert] (26) at (-3.5, 0) {};
		\node [style=vert] (27) at (-4.5, 0) {};
		\node [style=vert-white] (28) at (-1.5, 0) {};
		\node [style=vert-white] (29) at (0.5, 0) {};
		\node [style=vert-white] (30) at (1.5, 0) {};
		\node [style=vert-white] (31) at (3.5, 0) {};
		\node [style=vert-white] (32) at (4.5, 0) {};
		\node [style=vert] (33) at (2.5, 0) {};
		\node [style=label] (34) at (-4.5, 0.5) {$-1$};
		\node [style=label] (35) at (-3.5, 0.5) {$0$};
		\node [style=label] (36) at (-2.5, 0.5) {$1$};
		\node [style=label] (37) at (-1.5, 0.5) {$4$};
		\node [style=label] (38) at (-0.5, 0.5) {$2$};
		\node [style=label] (39) at (0.5, 0.5) {$5$};
		\node [style=label] (40) at (1.5, 0.5) {$6$};
		\node [style=label] (41) at (2.5, 0.5) {$3$};
		\node [style=label] (42) at (3.5, 0.5) {$7$};
		\node [style=label] (43) at (4.5, 0.5) {$8$};
	\end{pgfonlayer}
	\begin{pgfonlayer}{edgelayer}
		\draw (19.center) to (4.center);
		\draw [style=dot] (22.center) to (23.center);
	\end{pgfonlayer}
\end{tikzpicture}
\]
\end{remark}

Maya diagrams are models of 1D free fermions, where black nodes represent particles, and the white nodes represent empty spots (or ``holes''). We will next consider certain discrete time evolutions of fermions, and relate to the Hamiltonian interpretation of Schur functions.

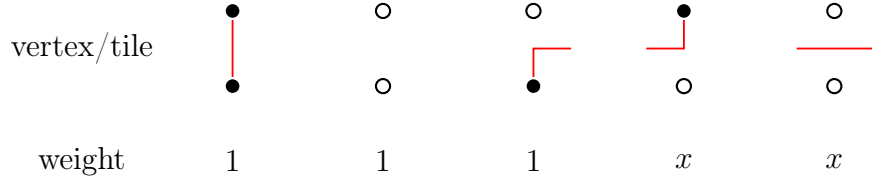
\begin{figure}[h]\centering
	\begin{tikzpicture}[baseline=-8.5em]
	\begin{pgfonlayer}{nodelayer}
		\node [style=vert] (52) at (-3, -2) {};
		\node [style=vert] (53) at (-3, -3) {};
		\node [style=vert-white] (54) at (-1, -2) {};
		\node [style=vert-white] (55) at (-1, -3) {};
		\node [style=vert-white] (56) at (1, -2) {};
		\node [style=vert-white] (57) at (3, -3) {};
		\node [style=vert] (58) at (1, -3) {};
		\node [style=vert] (59) at (3, -2) {};
		\node [style=vert-white] (60) at (5, -2) {};
		\node [style=vert-white] (61) at (5, -3) {};
		\node [style=none] (62) at (1, -2.5) {};
		\node [style=none] (63) at (1.5, -2.5) {};
		\node [style=none] (64) at (2.5, -2.5) {};
		\node [style=none] (65) at (3, -2.5) {};
		\node [style=none] (66) at (4.5, -2.5) {};
		\node [style=none] (67) at (5.5, -2.5) {};
		\node [style=none] (68) at (-3.5, -2.5) {};
		\node [style=none] (69) at (-2.5, -2.5) {};
		\node [style=none] (70) at (-1.5, -2.5) {};
		\node [style=none] (71) at (-0.5, -2.5) {};
		\node [style=none] (72) at (0.5, -2.5) {};
		\node [style=none] (73) at (3.5, -2.5) {};
		\node [style=none] (74) at (-5, -2.5) {\text{vertex}/\text{tile}};
		\node [style=none] (75) at (-5, -4) {\text{weight}};
		\node [style=none] (76) at (-3, -4) {$1$};
		\node [style=none] (77) at (-1, -4) {$1$};
		\node [style=none] (78) at (1, -4) {$1$};
		\node [style=none] (79) at (3, -4) {$x$};
		\node [style=none] (80) at (5, -4) {$x$};
	\end{pgfonlayer}
	\begin{pgfonlayer}{edgelayer}
		\draw [style=red] (52) to (53);
		\draw [style=red] (58)
			 to (62.center)
			 to (63.center);
		\draw [style=red] (64.center)
			 to (65.center)
			 to (59);
		\draw [style=red] (66.center) to (67.center);
	\end{pgfonlayer}
\end{tikzpicture}
\caption{Tiles for evolutions of free fermions with their weights shown below. Here $x$ is a spectral parameter.}
\label{fig:5-vertex-wt}
\end{figure}

We shall think of every Maya diagram, or equivalently $|\lambda\rangle \in\ff$, as a state of a discrete-time dynamics. To obtain the future state, each particle can move to the left by a number of steps, as long as they don't meet another particle. This can be summarized as certain `tilings' where the allowed tiles are summarized in \Cref{fig:5-vertex-wt}. Each tile has a weight associated to it, with $x$ being the spectral parameter.

\begin{figure}[h]\centering
\begin{tikzpicture}[]
	\begin{pgfonlayer}{nodelayer}
		\node [style=none] (4) at (5, 0) {};
		\node [style=none] (19) at (-5, 0) {};
		\node [style=vert] (24) at (-1.5, 0) {};
		\node [style=vert] (25) at (-2.5, 0) {};
		\node [style=vert] (26) at (-3.5, 0) {};
		\node [style=vert] (27) at (-4.5, 0) {};
		\node [style=vert-white] (28) at (-0.5, 0) {};
		\node [style=vert-white] (29) at (2.5, 0) {};
		\node [style=vert-white] (30) at (1.5, 0) {};
		\node [style=vert-white] (31) at (3.5, 0) {};
		\node [style=vert-white] (32) at (4.5, 0) {};
		\node [style=vert] (33) at (0.5, 0) {};
		\node [style=none] (34) at (5, 1) {};
		\node [style=none] (35) at (-5, 1) {};
		\node [style=vert] (38) at (-0.5, 1) {};
		\node [style=vert] (39) at (-2.5, 1) {};
		\node [style=vert] (40) at (-3.5, 1) {};
		\node [style=vert] (41) at (-4.5, 1) {};
		\node [style=vert-white] (42) at (-1.5, 1) {};
		\node [style=vert-white] (43) at (0.5, 1) {};
		\node [style=vert-white] (44) at (1.5, 1) {};
		\node [style=vert-white] (45) at (3.5, 1) {};
		\node [style=vert-white] (46) at (4.5, 1) {};
		\node [style=vert] (47) at (2.5, 1) {};
		\node [style=none] (48) at (2.5, 0.5) {};
		\node [style=none] (49) at (0.5, 0.5) {};
		\node [style=none] (50) at (-0.5, 0.5) {};
		\node [style=none] (51) at (-1.5, 0.5) {};
	\end{pgfonlayer}
	\begin{pgfonlayer}{edgelayer}
		\draw (19.center) to (4.center);
		\draw (35.center) to (34.center);
		\draw [style=red] (47)
			 to (48.center)
			 to (49.center)
			 to (33);
		\draw [style=red] (38)
			 to (50.center)
			 to (51.center)
			 to (24);
		\draw [style=red] (39) to (25);
		\draw [style=red] (40) to (26);
		\draw [style=red] (41) to (27);
	\end{pgfonlayer}
\end{tikzpicture}
\caption{Example of a one-step evolution of fermions, where the red paths represent movement of particles. The bottom row is the future state.}
\label{fig:vertex-one-step}
\end{figure}
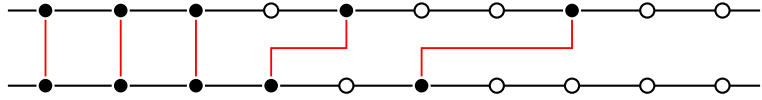

For example, $|(1)\rangle$ is a possible future state evolved from $|(3,1)\rangle$. The time evolution is represented by the path configuration in \Cref{fig:vertex-one-step}.

We now introduce linear operators $T(x)$ on $\ff$ called \emph{transfer matrices}. It can be simply described to be the weighted sum over all possible $|\mu\rangle$ that can be obtained by fermionic evolution from a given $|\lambda\rangle$. Equivalently, $\langle \mu|T(x)|\lambda \rangle$ is the weighted sum over all such configurations as in \Cref{fig:vertex-one-step}, where the top row is the Maya diagram of $\lambda$ and the bottom row is the Maya diagram of $\mu$. Now consider the evolution over multiple times, then the result will be a linear combination of configurations which are exactly the Gessel-Viennot's non-intersecting lattice paths. Therefore we have that $\langle \mu|T(x_n)\dots T(x_2)T(x_1)|\lambda\rangle=s_{\lambda/\mu}(x_1,\cdots,x_n)$ the skew Schur polynomials.

\begin{figure}[h]
\centering
\begin{tikzpicture}
	\begin{pgfonlayer}{nodelayer}
		\node [style=none] (4) at (3, 0) {};
		\node [style=none] (19) at (-2, 0) {};
		\node [style=vert] (24) at (-1.5, 0) {};
		\node [style=vert-white] (28) at (-0.5, 0) {};
		\node [style=vert-white] (29) at (2.5, 0) {};
		\node [style=vert-white] (30) at (1.5, 0) {};
		\node [style=vert] (33) at (0.5, 0) {};
		\node [style=none] (34) at (3, 1) {};
		\node [style=none] (35) at (-2, 1) {};
		\node [style=vert] (38) at (-0.5, 1) {};
		\node [style=vert-white] (42) at (-1.5, 1) {};
		\node [style=vert-white] (43) at (0.5, 1) {};
		\node [style=vert-white] (44) at (1.5, 1) {};
		\node [style=vert] (47) at (2.5, 1) {};
		\node [style=none] (48) at (2.5, 0.5) {};
		\node [style=none] (49) at (0.5, 0.5) {};
		\node [style=none] (50) at (-0.5, 0.5) {};
		\node [style=none] (51) at (-1.5, 0.5) {};
		\node [style=none] (52) at (3, -1) {};
		\node [style=none] (53) at (-2, -1) {};
		\node [style=vert] (54) at (-1.5, -1) {};
		\node [style=vert-white] (58) at (0.5, -1) {};
		\node [style=vert-white] (59) at (2.5, -1) {};
		\node [style=vert-white] (60) at (1.5, -1) {};
		\node [style=vert] (63) at (-0.5, -1) {};
		\node [style=none] (64) at (0.5, -0.5) {};
		\node [style=none] (65) at (-0.5, -0.5) {};
	\end{pgfonlayer}
	\begin{pgfonlayer}{edgelayer}
		\draw (19.center) to (4.center);
		\draw (35.center) to (34.center);
		\draw [style=red] (47)
			 to (48.center)
			 to (49.center)
			 to (33);
		\draw [style=red] (38)
			 to (50.center)
			 to (51.center)
			 to (24);
		\draw (53.center) to (52.center);
		\draw [style=red] (33)
			 to (64.center)
			 to (65.center)
			 to (63);
		\draw [style=red] (54) to (24);
	\end{pgfonlayer}
\end{tikzpicture}
\caption{Example of a vertex model configuration. The weight of this configuration is $x_1^3x_2$.}
\label{fig:vertex_model_eg}
\end{figure}
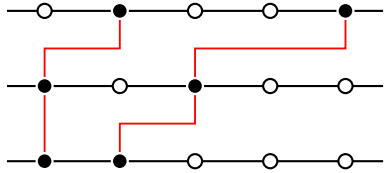

Since we always have infinitely many particles (resp. holes) to the left (resp. right), we may ignore the infinite redundancies and restrict ourselves in a finite window. If we allow multiple-step configurations like \Cref{fig:vertex-one-step}, we will have a family of configurations with fixed boundary conditions (such as \Cref{fig:vertex_model_eg}). We assign weights to these configurations according to \Cref{fig:5-vertex-wt} with the spectral parameter set to be $x_i$ for the $i$-th row. This is called a \emph{five vertex model}, and in this context, the weighted sum of all such configurations with fixed boundary is called the \emph{partition function}. As we have seen, the partition function of the $n$-row vertex model will be skew Schur polynomials in $n$ variables.

Recall that in the Hamiltonian interpretation of Schur functions, we have
\(s_\lambda = \langle 0|\prod e^{\phi(x_i)}|\lambda\rangle.\)
It turns out the $e^{\phi(x)}$'s are exactly the transfer matrices, i.e. $e^{\phi(x)}=T(x)$ as linear operators on $\ff$. From this equality we know that the transfer matrices commute, i.e. $T(x)T(y)=T(y)T(x)$. In the language of vertex models, the commutativity can be realized as the existence of solutions to the \emph{Yang-Baxter equations}. We will however not go into details about this perspective.

\subsection{Two Dimensional fermions}
We will now generalize this construction to $\tilde \ff = \sz$, the physical model here will be what we call \emph{2D fermions}.
 
Roughly speaking, a 2D fermionic configuration is an alignment of particles and holes on $\zz^2$, such that every row is an usual 1D fermionic configuration and every column is a vacuum with different centers. 

\begin{definition}\label{def:2d-fermion}
We define a \emph{2D free fermionic configuration} to be a map $f:\zz^2\to \{0,1\}$, satisfying the following conditions.
\begin{enumerate}[noitemsep]
	\item For every $j\in\mathbb{Z}$, the map $f(\underline{\ \ },j):\zz\to \{0,1\}$ is a Maya diagram centered at $j$. 
	\item For every $j\in\zz$, we have $\{i:f(i,j+1)=1\} = \{i:f(i,j)=1\} \sqcup\{k\}$ for some $k\in\zz$.
\end{enumerate}
\end{definition} 
Note that (2) of \Cref{def:2d-fermion} means that the $(j+1)$-th row of $f\in\tilde \ff$ has exactly one more particle than the $j$-th row. This guarantees that every column is a vacuum with different centers. Therefore we can associate each 2D Maya diagram $f$ a permutation $w\in\sz$ such that $w(i)$ is the center (i.e.  index of the top particle) of the $i$-th column. See \Cref{fig:2d-maya-eg} for example.
\begin{remark}
	For sake of readability, we will use bar-ed numbers in figures to represent negative numbers, e.g. $\bar 3 = -3$.
\end{remark}
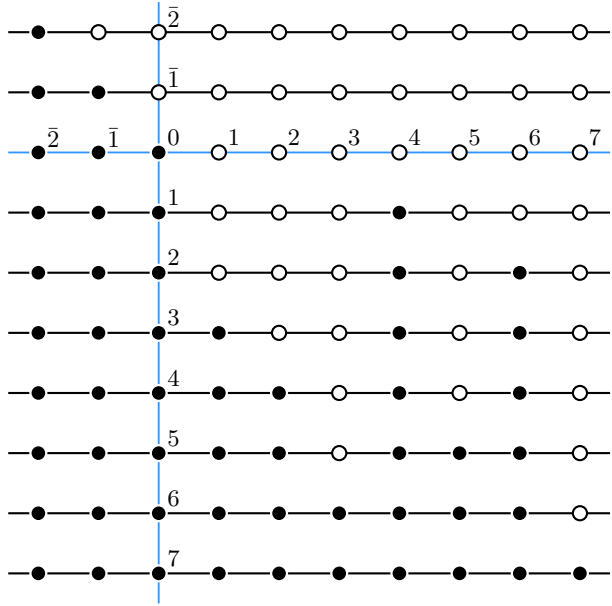
\begin{figure}[h]
\centering	
\begin{tikzpicture}[scale=0.8]
	\begin{pgfonlayer}{nodelayer}
		\node [style=label-small] (0) at (-0.75, 2.25) {$0$};
		\node [style=label-small] (1) at (0.25, 2.25) {$1$};
		\node [style=label-small] (2) at (1.25, 2.25) {$2$};
		\node [style=label-small] (3) at (2.25, 2.25) {$3$};
		\node [style=label-small] (4) at (3.25, 2.25) {$4$};
		\node [style=label-small] (5) at (4.25, 2.25) {$5$};
		\node [style=label-small] (8) at (-1.75, 2.25) {$\bar 1$};
		\node [style=label-small] (9) at (-2.75, 2.25) {$\bar 2$};
		\node [style=label-small] (14) at (-0.75, 3.25) {$\bar 1$};
		\node [style=label-small] (15) at (-0.75, 4.25) {$\bar 2$};
		\node [style=label-small] (24) at (-0.75, 1.25) {$1$};
		\node [style=label-small] (25) at (-0.75, 0.25) {$2$};
		\node [style=label-small] (26) at (-0.75, -0.75) {$3$};
		\node [style=label-small] (27) at (-0.75, -1.75) {$4$};
		\node [style=label-small] (28) at (-0.75, -2.75) {$5$};
		\node [style=none] (33) at (-3.5, 2) {};
		\node [style=none] (34) at (6.5, 2) {};
		\node [style=vert] (35) at (3, 0) {};
		\node [style=vert] (36) at (3, -1) {};
		\node [style=vert] (37) at (1, -2) {};
		\node [style=vert] (38) at (4, -3) {};
		\node [style=vert] (42) at (5, -3) {};
		\node [style=vert] (45) at (3, -3) {};
		\node [style=vert] (46) at (5, -2) {};
		\node [style=vert] (47) at (5, -1) {};
		\node [style=vert] (48) at (3, -2) {};
		\node [style=vert] (49) at (1, -3) {};
		\node [style=vert] (51) at (3, 1) {};
		\node [style=vert] (52) at (5, 0) {};
		\node [style=vert] (53) at (0, -1) {};
		\node [style=vert] (54) at (0, -2) {};
		\node [style=vert] (55) at (0, -3) {};
		\node [style=vert] (56) at (-1, 2) {};
		\node [style=vert] (57) at (-1, 1) {};
		\node [style=vert] (58) at (-1, 0) {};
		\node [style=vert] (59) at (-1, -1) {};
		\node [style=vert] (60) at (-1, -2) {};
		\node [style=vert] (61) at (-1, -3) {};
		\node [style=vert] (62) at (-2, 2) {};
		\node [style=vert] (63) at (-2, 1) {};
		\node [style=vert] (64) at (-2, 0) {};
		\node [style=vert] (65) at (-2, -1) {};
		\node [style=vert] (66) at (-2, -2) {};
		\node [style=vert] (67) at (-2, -3) {};
		\node [style=vert] (68) at (3, -4) {};
		\node [style=vert] (70) at (5, -4) {};
		\node [style=vert] (71) at (4, -4) {};
		\node [style=vert] (72) at (2, -4) {};
		\node [style=vert] (73) at (1, -4) {};
		\node [style=vert] (74) at (0, -4) {};
		\node [style=vert] (75) at (-1, -4) {};
		\node [style=vert] (76) at (4, -5) {};
		\node [style=vert] (78) at (6, -5) {};
		\node [style=vert] (79) at (5, -5) {};
		\node [style=vert] (80) at (3, -5) {};
		\node [style=vert] (81) at (2, -5) {};
		\node [style=vert] (82) at (1, -5) {};
		\node [style=vert] (83) at (0, -5) {};
		\node [style=vert] (84) at (-2, -4) {};
		\node [style=vert] (85) at (-1, -5) {};
		\node [style=vert] (86) at (-2, -5) {};
		\node [style=vert] (88) at (-2, 3) {};
		\node [style=vert] (90) at (-3, 3) {};
		\node [style=vert] (91) at (-3, 2) {};
		\node [style=vert] (92) at (-3, 1) {};
		\node [style=vert] (93) at (-3, 0) {};
		\node [style=vert] (94) at (-3, -1) {};
		\node [style=vert] (95) at (-3, -2) {};
		\node [style=vert] (96) at (-3, -3) {};
		\node [style=vert] (97) at (-3, -4) {};
		\node [style=vert] (98) at (-3, 4) {};
		\node [style=vert] (100) at (-3, -5) {};
		\node [style=none] (101) at (-1, 4.5) {};
		\node [style=none] (102) at (-1, -5.5) {};
		\node [style=label-small] (103) at (5.25, 2.25) {$6$};
		\node [style=label-small] (104) at (6.25, 2.25) {$7$};
		\node [style=label-small] (106) at (-0.75, -3.75) {$6$};
		\node [style=label-small] (107) at (-0.75, -4.75) {$7$};
		\node [style=none] (111) at (-3.5, 4) {};
		\node [style=none] (112) at (-3.5, 3) {};
		\node [style=none] (113) at (-3.5, 1) {};
		\node [style=none] (114) at (-3.5, 0) {};
		\node [style=none] (115) at (-3.5, -1) {};
		\node [style=none] (116) at (-3.5, -2) {};
		\node [style=none] (117) at (-3.5, -3) {};
		\node [style=none] (118) at (-3.5, -4) {};
		\node [style=none] (119) at (-3.5, -5) {};
		\node [style=none] (120) at (6.5, 2) {};
		\node [style=none] (122) at (6.5, 4) {};
		\node [style=none] (123) at (6.5, 3) {};
		\node [style=none] (124) at (6.5, 1) {};
		\node [style=none] (125) at (6.5, 0) {};
		\node [style=none] (126) at (6.5, -1) {};
		\node [style=none] (127) at (6.5, -2) {};
		\node [style=none] (128) at (6.5, -3) {};
		\node [style=none] (129) at (6.5, -4) {};
		\node [style=none] (130) at (6.5, -5) {};
		\node [style=vert-white] (131) at (-2, 4) {};
		\node [style=vert-white] (132) at (-1, 4) {};
		\node [style=vert-white] (133) at (-1, 3) {};
		\node [style=vert-white] (134) at (0, 3) {};
		\node [style=vert-white] (135) at (0, 2) {};
		\node [style=vert-white] (136) at (0, 1) {};
		\node [style=vert-white] (137) at (0, 4) {};
		\node [style=vert-white] (138) at (1, 4) {};
		\node [style=vert-white] (139) at (1, 3) {};
		\node [style=vert-white] (140) at (1, 2) {};
		\node [style=vert-white] (141) at (2, 1) {};
		\node [style=vert-white] (142) at (2, 0) {};
		\node [style=vert-white] (143) at (2, -1) {};
		\node [style=vert-white] (144) at (4, -1) {};
		\node [style=vert-white] (145) at (4, -2) {};
		\node [style=vert-white] (146) at (2, -3) {};
		\node [style=vert-white] (147) at (2, -2) {};
		\node [style=vert-white] (148) at (1, -1) {};
		\node [style=vert-white] (149) at (1, 0) {};
		\node [style=vert-white] (150) at (4, 0) {};
		\node [style=vert-white] (151) at (0, 0) {};
		\node [style=vert-white] (152) at (1, 1) {};
		\node [style=vert-white] (153) at (4, 1) {};
		\node [style=vert-white] (154) at (5, 1) {};
		\node [style=vert-white] (155) at (2, 2) {};
		\node [style=vert-white] (156) at (3, 2) {};
		\node [style=vert-white] (157) at (4, 2) {};
		\node [style=vert-white] (158) at (5, 2) {};
		\node [style=vert-white] (159) at (6, -4) {};
		\node [style=vert-white] (160) at (6, -3) {};
		\node [style=vert-white] (161) at (6, -2) {};
		\node [style=vert-white] (162) at (6, -1) {};
		\node [style=vert-white] (163) at (6, 0) {};
		\node [style=vert-white] (164) at (6, 1) {};
		\node [style=vert-white] (165) at (6, 2) {};
		\node [style=vert-white] (166) at (6, 3) {};
		\node [style=vert-white] (167) at (6, 4) {};
		\node [style=vert-white] (168) at (5, 4) {};
		\node [style=vert-white] (169) at (5, 3) {};
		\node [style=vert-white] (170) at (4, 3) {};
		\node [style=vert-white] (171) at (4, 4) {};
		\node [style=vert-white] (172) at (3, 3) {};
		\node [style=vert-white] (173) at (3, 4) {};
		\node [style=vert-white] (174) at (2, 4) {};
		\node [style=vert-white] (175) at (2, 3) {};
	\end{pgfonlayer}
	\begin{pgfonlayer}{edgelayer}
		\draw [style=blue] (33.center) to (34.center);
		\draw [style=blue] (101.center) to (102.center);
		\draw [style=none] (111.center) to (122.center);
		\draw [style=none] (112.center) to (123.center);
		\draw [style=none] (113.center) to (124.center);
		\draw [style=none] (114.center) to (125.center);
		\draw [style=none] (115.center) to (126.center);
		\draw [style=none] (127.center) to (116.center);
		\draw [style=none] (117.center) to (128.center);
		\draw [style=none] (129.center) to (118.center);
		\draw [style=none] (119.center) to (130.center);
	\end{pgfonlayer}
\end{tikzpicture}
\caption{The analogue of Maya diagram for $w = [\cdots \bar 2\ \bar 1\ 0\ 2\ 4\ 1\ 5\ 3\ 6\ 7\ \cdots ]$, where $\bar i = -i$. The two blue lines represent the $x$-axis and $y$-axis of $\zz^2$.}
\label{fig:2d-maya-eg}
\end{figure}

From the definition, we know that if $f$ is a 2D Maya diagram, then the restriction to each row is a 1D Maya diagram, i.e. a partition. This can be interpreted as follow. Let $\pi_i$ denote this restriction and let ${G}_i(w)$ denote the permutation obtained from $w$ by sorting it into a $i$-Grassmanian permutation. For example $G_3(423561)=234156$. Geometrically, $G_i$ is the projection from the flag manifold to an $i$-Grassmanian. For any $w\in \sz$ as a 2D Maya diagram, we have that $\pi_i(w) = \lambda(G_i(w^{-1}))$. 

\begin{figure}[h]
\centering
	\begin{tikzpicture}[scale = 0.7]
	\begin{pgfonlayer}{nodelayer}
		\node [style=label-small] (0) at (-0.75, 2.25) {$0$};
		\node [style=label-small] (1) at (0.25, 2.25) {$1$};
		\node [style=label-small] (2) at (1.25, 2.25) {$2$};
		\node [style=label-small] (3) at (2.25, 2.25) {$3$};
		\node [style=label-small] (4) at (3.25, 2.25) {$4$};
		\node [style=label-small] (5) at (4.25, 2.25) {$5$};
		\node [style=label-small] (8) at (-1.75, 2.25) {$\bar 1$};
		\node [style=label-small] (9) at (-2.75, 2.25) {$\bar 2$};
		\node [style=label-small] (14) at (-0.75, 3.25) {$\bar 1$};
		\node [style=label-small] (15) at (-0.75, 4.25) {$\bar 2$};
		\node [style=label-small] (24) at (-0.75, 1.25) {$1$};
		\node [style=label-small] (25) at (-0.75, 0.25) {$2$};
		\node [style=label-small] (26) at (-0.75, -0.75) {$3$};
		\node [style=label-small] (27) at (-0.75, -1.75) {$4$};
		\node [style=label-small] (28) at (-0.75, -2.75) {$5$};
		\node [style=none] (33) at (-3.5, 2) {};
		\node [style=none] (34) at (6.5, 2) {};
		\node [style=vert] (35) at (3, 0) {};
		\node [style=vert] (36) at (3, -1) {};
		\node [style=vert] (37) at (1, -2) {};
		\node [style=vert] (38) at (4, -3) {};
		\node [style=vert] (42) at (5, -3) {};
		\node [style=vert] (45) at (3, -3) {};
		\node [style=vert] (46) at (5, -2) {};
		\node [style=vert] (47) at (5, -1) {};
		\node [style=vert] (48) at (3, -2) {};
		\node [style=vert] (49) at (1, -3) {};
		\node [style=vert] (51) at (5, 1) {};
		\node [style=vert] (52) at (5, 0) {};
		\node [style=vert] (53) at (1, -1) {};
		\node [style=vert] (54) at (0, -2) {};
		\node [style=vert] (55) at (0, -3) {};
		\node [style=vert] (56) at (1, 2) {};
		\node [style=vert] (57) at (1, 1) {};
		\node [style=vert] (58) at (1, 0) {};
		\node [style=vert] (59) at (0, -1) {};
		\node [style=vert] (60) at (-1, -2) {};
		\node [style=vert] (61) at (-1, -3) {};
		\node [style=vert] (62) at (-2, 2) {};
		\node [style=vert] (63) at (-2, 1) {};
		\node [style=vert] (64) at (-2, 0) {};
		\node [style=vert] (65) at (-2, -1) {};
		\node [style=vert] (66) at (-2, -2) {};
		\node [style=vert] (67) at (-2, -3) {};
		\node [style=vert] (68) at (3, -4) {};
		\node [style=vert] (70) at (5, -4) {};
		\node [style=vert] (71) at (4, -4) {};
		\node [style=vert] (72) at (2, -4) {};
		\node [style=vert] (73) at (1, -4) {};
		\node [style=vert] (74) at (0, -4) {};
		\node [style=vert] (75) at (-1, -4) {};
		\node [style=vert] (76) at (4, -5) {};
		\node [style=vert] (78) at (6, -5) {};
		\node [style=vert] (79) at (5, -5) {};
		\node [style=vert] (80) at (3, -5) {};
		\node [style=vert] (81) at (2, -5) {};
		\node [style=vert] (82) at (1, -5) {};
		\node [style=vert] (83) at (0, -5) {};
		\node [style=vert] (84) at (-2, -4) {};
		\node [style=vert] (85) at (-1, -5) {};
		\node [style=vert] (86) at (-2, -5) {};
		\node [style=vert] (88) at (1, 3) {};
		\node [style=vert] (90) at (-3, 3) {};
		\node [style=vert] (91) at (-3, 2) {};
		\node [style=vert] (92) at (-3, 1) {};
		\node [style=vert] (93) at (-3, 0) {};
		\node [style=vert] (94) at (-3, -1) {};
		\node [style=vert] (95) at (-3, -2) {};
		\node [style=vert] (96) at (-3, -3) {};
		\node [style=vert] (97) at (-3, -4) {};
		\node [style=vert] (98) at (-3, 4) {};
		\node [style=vert] (100) at (-3, -5) {};
		\node [style=none] (101) at (-1, 4.5) {};
		\node [style=none] (102) at (-1, -5.5) {};
		\node [style=label-small] (103) at (5.25, 2.25) {$6$};
		\node [style=label-small] (104) at (6.25, 2.25) {$7$};
		\node [style=label-small] (106) at (-0.75, -3.75) {$6$};
		\node [style=label-small] (107) at (-0.75, -4.75) {$7$};
		\node [style=none] (111) at (-3.5, 4) {};
		\node [style=none] (112) at (-3.5, 3) {};
		\node [style=none] (113) at (-3.5, 1) {};
		\node [style=none] (114) at (-3.5, 0) {};
		\node [style=none] (115) at (-3.5, -1) {};
		\node [style=none] (116) at (-3.5, -2) {};
		\node [style=none] (117) at (-3.5, -3) {};
		\node [style=none] (118) at (-3.5, -4) {};
		\node [style=none] (119) at (-3.5, -5) {};
		\node [style=none] (120) at (6.5, 2) {};
		\node [style=none] (122) at (6.5, 4) {};
		\node [style=none] (123) at (6.5, 3) {};
		\node [style=none] (124) at (6.5, 1) {};
		\node [style=none] (125) at (6.5, 0) {};
		\node [style=none] (126) at (6.5, -1) {};
		\node [style=none] (127) at (6.5, -2) {};
		\node [style=none] (128) at (6.5, -3) {};
		\node [style=none] (129) at (6.5, -4) {};
		\node [style=none] (130) at (6.5, -5) {};
		\node [style=vert-white] (131) at (-2, 4) {};
		\node [style=vert-white] (132) at (-1, 4) {};
		\node [style=vert-white] (133) at (-1, 3) {};
		\node [style=vert-white] (134) at (-2, 3) {};
		\node [style=vert-white] (135) at (-1, 2) {};
		\node [style=vert-white] (136) at (-1, 1) {};
		\node [style=vert-white] (137) at (0, 4) {};
		\node [style=vert-white] (138) at (1, 4) {};
		\node [style=vert-white] (139) at (0, 3) {};
		\node [style=vert-white] (140) at (0, 2) {};
		\node [style=vert-white] (141) at (2, 1) {};
		\node [style=vert-white] (142) at (2, 0) {};
		\node [style=vert-white] (143) at (2, -1) {};
		\node [style=vert-white] (144) at (4, -1) {};
		\node [style=vert-white] (145) at (4, -2) {};
		\node [style=vert-white] (146) at (2, -3) {};
		\node [style=vert-white] (147) at (2, -2) {};
		\node [style=vert-white] (148) at (-1, -1) {};
		\node [style=vert-white] (149) at (0, 0) {};
		\node [style=vert-white] (150) at (4, 0) {};
		\node [style=vert-white] (151) at (-1, 0) {};
		\node [style=vert-white] (152) at (0, 1) {};
		\node [style=vert-white] (153) at (3, 1) {};
		\node [style=vert-white] (154) at (4, 1) {};
		\node [style=vert-white] (155) at (2, 2) {};
		\node [style=vert-white] (156) at (3, 2) {};
		\node [style=vert-white] (157) at (4, 2) {};
		\node [style=vert-white] (158) at (5, 2) {};
		\node [style=vert-white] (159) at (6, -4) {};
		\node [style=vert-white] (160) at (6, -3) {};
		\node [style=vert-white] (161) at (6, -2) {};
		\node [style=vert-white] (162) at (6, -1) {};
		\node [style=vert-white] (163) at (6, 0) {};
		\node [style=vert-white] (164) at (6, 1) {};
		\node [style=vert-white] (165) at (6, 2) {};
		\node [style=vert-white] (166) at (6, 3) {};
		\node [style=vert-white] (167) at (6, 4) {};
		\node [style=vert-white] (168) at (5, 4) {};
		\node [style=vert-white] (169) at (5, 3) {};
		\node [style=vert-white] (170) at (4, 3) {};
		\node [style=vert-white] (171) at (4, 4) {};
		\node [style=vert-white] (172) at (3, 3) {};
		\node [style=vert-white] (173) at (3, 4) {};
		\node [style=vert-white] (174) at (2, 4) {};
		\node [style=vert-white] (175) at (2, 3) {};
		\node [style=label-small] (176) at (12.75, 2.25) {$0$};
		\node [style=label-small] (177) at (13.75, 2.25) {$1$};
		\node [style=label-small] (178) at (14.75, 2.25) {$2$};
		\node [style=label-small] (179) at (15.75, 2.25) {$3$};
		\node [style=label-small] (180) at (16.75, 2.25) {$4$};
		\node [style=label-small] (181) at (17.75, 2.25) {$5$};
		\node [style=label-small] (182) at (11.75, 2.25) {$\bar 1$};
		\node [style=label-small] (183) at (10.75, 2.25) {$\bar 2$};
		\node [style=label-small] (184) at (12.75, 3.25) {$\bar 1$};
		\node [style=label-small] (185) at (12.75, 4.25) {$\bar 2$};
		\node [style=label-small] (186) at (12.75, 1.25) {$1$};
		\node [style=label-small] (187) at (12.75, 0.25) {$2$};
		\node [style=label-small] (188) at (12.75, -0.75) {$3$};
		\node [style=label-small] (189) at (12.75, -1.75) {$4$};
		\node [style=label-small] (190) at (12.75, -2.75) {$5$};
		\node [style=none] (191) at (10, 2) {};
		\node [style=none] (192) at (20, 2) {};
		\node [style=vert] (193) at (16.5, 0) {};
		\node [style=vert] (194) at (16.5, -1) {};
		\node [style=vert] (195) at (14.5, -2) {};
		\node [style=vert] (196) at (17.5, -3) {};
		\node [style=vert] (197) at (18.5, -3) {};
		\node [style=vert] (198) at (16.5, -3) {};
		\node [style=vert] (199) at (17.5, -2) {};
		\node [style=vert] (200) at (17.5, -1) {};
		\node [style=vert] (201) at (16.5, -2) {};
		\node [style=vert] (202) at (14.5, -3) {};
		\node [style=vert] (203) at (16.5, 1) {};
		\node [style=vert] (204) at (17.5, 0) {};
		\node [style=vert] (205) at (14.5, -1) {};
		\node [style=vert] (206) at (13.5, -2) {};
		\node [style=vert] (207) at (13.5, -3) {};
		\node [style=vert] (208) at (13.5, 2) {};
		\node [style=vert] (209) at (13.5, 1) {};
		\node [style=vert] (210) at (13.5, 0) {};
		\node [style=vert] (211) at (13.5, -1) {};
		\node [style=vert] (212) at (12.5, -2) {};
		\node [style=vert] (213) at (12.5, -3) {};
		\node [style=vert] (214) at (11.5, 2) {};
		\node [style=vert] (215) at (11.5, 1) {};
		\node [style=vert] (216) at (11.5, 0) {};
		\node [style=vert] (217) at (11.5, -1) {};
		\node [style=vert] (218) at (11.5, -2) {};
		\node [style=vert] (219) at (11.5, -3) {};
		\node [style=vert] (220) at (16.5, -4) {};
		\node [style=vert] (221) at (18.5, -4) {};
		\node [style=vert] (222) at (17.5, -4) {};
		\node [style=vert] (223) at (15.5, -4) {};
		\node [style=vert] (224) at (14.5, -4) {};
		\node [style=vert] (225) at (13.5, -4) {};
		\node [style=vert] (226) at (12.5, -4) {};
		\node [style=vert] (227) at (17.5, -5) {};
		\node [style=vert] (228) at (19.5, -5) {};
		\node [style=vert] (229) at (18.5, -5) {};
		\node [style=vert] (230) at (16.5, -5) {};
		\node [style=vert] (231) at (15.5, -5) {};
		\node [style=vert] (232) at (14.5, -5) {};
		\node [style=vert] (233) at (13.5, -5) {};
		\node [style=vert] (234) at (11.5, -4) {};
		\node [style=vert] (235) at (12.5, -5) {};
		\node [style=vert] (236) at (11.5, -5) {};
		\node [style=vert] (237) at (13.5, 3) {};
		\node [style=vert] (238) at (10.5, 3) {};
		\node [style=vert] (239) at (10.5, 2) {};
		\node [style=vert] (240) at (10.5, 1) {};
		\node [style=vert] (241) at (10.5, 0) {};
		\node [style=vert] (242) at (10.5, -1) {};
		\node [style=vert] (243) at (10.5, -2) {};
		\node [style=vert] (244) at (10.5, -3) {};
		\node [style=vert] (245) at (10.5, -4) {};
		\node [style=vert] (246) at (10.5, 4) {};
		\node [style=vert] (247) at (10.5, -5) {};
		\node [style=none] (248) at (12.5, 4.5) {};
		\node [style=none] (249) at (12.5, -5.5) {};
		\node [style=label-small] (250) at (18.75, 2.25) {$6$};
		\node [style=label-small] (251) at (19.75, 2.25) {$7$};
		\node [style=label-small] (252) at (12.75, -3.75) {$6$};
		\node [style=label-small] (253) at (12.75, -4.75) {$7$};
		\node [style=none] (254) at (10, 4) {};
		\node [style=none] (255) at (10, 3) {};
		\node [style=none] (256) at (10, 1) {};
		\node [style=none] (257) at (10, 0) {};
		\node [style=none] (258) at (10, -1) {};
		\node [style=none] (259) at (10, -2) {};
		\node [style=none] (260) at (10, -3) {};
		\node [style=none] (261) at (10, -4) {};
		\node [style=none] (262) at (10, -5) {};
		\node [style=none] (263) at (20, 2) {};
		\node [style=none] (264) at (20, 4) {};
		\node [style=none] (265) at (20, 3) {};
		\node [style=none] (266) at (20, 1) {};
		\node [style=none] (267) at (20, 0) {};
		\node [style=none] (268) at (20, -1) {};
		\node [style=none] (269) at (20, -2) {};
		\node [style=none] (270) at (20, -3) {};
		\node [style=none] (271) at (20, -4) {};
		\node [style=none] (272) at (20, -5) {};
		\node [style=vert-white] (273) at (11.5, 4) {};
		\node [style=vert-white] (274) at (12.5, 4) {};
		\node [style=vert-white] (275) at (12.5, 3) {};
		\node [style=vert-white] (276) at (11.5, 3) {};
		\node [style=vert-white] (277) at (12.5, 2) {};
		\node [style=vert-white] (278) at (12.5, 1) {};
		\node [style=vert-white] (279) at (13.5, 4) {};
		\node [style=vert-white] (280) at (14.5, 4) {};
		\node [style=vert-white] (281) at (14.5, 3) {};
		\node [style=vert-white] (282) at (14.5, 2) {};
		\node [style=vert-white] (283) at (15.5, 1) {};
		\node [style=vert-white] (284) at (15.5, 0) {};
		\node [style=vert-white] (285) at (15.5, -1) {};
		\node [style=vert-white] (286) at (18.5, -1) {};
		\node [style=vert-white] (287) at (18.5, -2) {};
		\node [style=vert-white] (288) at (15.5, -3) {};
		\node [style=vert-white] (289) at (15.5, -2) {};
		\node [style=vert-white] (290) at (12.5, -1) {};
		\node [style=vert-white] (291) at (14.5, 0) {};
		\node [style=vert-white] (292) at (18.5, 0) {};
		\node [style=vert-white] (293) at (12.5, 0) {};
		\node [style=vert-white] (294) at (14.5, 1) {};
		\node [style=vert-white] (295) at (18.5, 1) {};
		\node [style=vert-white] (296) at (17.5, 1) {};
		\node [style=vert-white] (297) at (15.5, 2) {};
		\node [style=vert-white] (298) at (16.5, 2) {};
		\node [style=vert-white] (299) at (17.5, 2) {};
		\node [style=vert-white] (300) at (18.5, 2) {};
		\node [style=vert-white] (301) at (19.5, -4) {};
		\node [style=vert-white] (302) at (19.5, -3) {};
		\node [style=vert-white] (303) at (19.5, -2) {};
		\node [style=vert-white] (304) at (19.5, -1) {};
		\node [style=vert-white] (305) at (19.5, 0) {};
		\node [style=vert-white] (306) at (19.5, 1) {};
		\node [style=vert-white] (307) at (19.5, 2) {};
		\node [style=vert-white] (308) at (19.5, 3) {};
		\node [style=vert-white] (309) at (19.5, 4) {};
		\node [style=vert-white] (310) at (18.5, 4) {};
		\node [style=vert-white] (311) at (18.5, 3) {};
		\node [style=vert-white] (312) at (17.5, 3) {};
		\node [style=vert-white] (313) at (17.5, 4) {};
		\node [style=vert-white] (314) at (16.5, 3) {};
		\node [style=vert-white] (315) at (16.5, 4) {};
		\node [style=vert-white] (316) at (15.5, 4) {};
		\node [style=vert-white] (317) at (15.5, 3) {};
		\node [style=none] (318) at (7.25, -1) {};
		\node [style=none] (320) at (9.25, -1) {};
		\node [style=none] (321) at (1.5, -6.5) {$\cdots\ \bar 2\ 0\ 4\ 3\ \bar 1\ 6\ 2\ 5\ 1\ 7\ \cdots$};
		\node [style=none] (322) at (15, -6.5) {$\cdots\ \bar 2\ 0\ 4\ \bar 1\ 3\ 6\ 1\ 2\ 5\ 7\ \cdots$};
	\end{pgfonlayer}
	\begin{pgfonlayer}{edgelayer}
		\draw [style=blue] (33.center) to (34.center);
		\draw [style=blue] (101.center) to (102.center);
		\draw [style=fade] (111.center) to (122.center);
		\draw [style=fade] (112.center) to (123.center);
		\draw [style=fade] (113.center) to (124.center);
		\draw [style=fade] (114.center) to (125.center);
		\draw [style=fade] (115.center) to (126.center);
		\draw [style=fade] (127.center) to (116.center);
		\draw [style=fade] (117.center) to (128.center);
		\draw [style=fade] (129.center) to (118.center);
		\draw [style=fade] (119.center) to (130.center);
		\draw [style=blue] (191.center) to (192.center);
		\draw [style=blue] (248.center) to (249.center);
		\draw [style=fade] (254.center) to (264.center);
		\draw [style=fade] (255.center) to (265.center);
		\draw [style=fade] (256.center) to (266.center);
		\draw [style=fade] (257.center) to (267.center);
		\draw [style=fade] (258.center) to (268.center);
		\draw [style=fade] (269.center) to (259.center);
		\draw [style=fade] (260.center) to (270.center);
		\draw [style=fade] (271.center) to (261.center);
		\draw [style=fade] (262.center) to (272.center);
		\draw [style=red arrow thick, draw = black] (318.center) to (320.center);
		\draw [style=red arrow] (51) to (153);
		\draw [style=red arrow] (52) to (150);
		\draw [style=red arrow] (47) to (144);
		\draw [style=red arrow] (46) to (145);
		\draw [style=red arrow] (58) to (149);
		\draw [style=red arrow] (57) to (152);
		\draw [style=red arrow] (56) to (140);
		\draw [style=red arrow] (88) to (139);
		\draw [style=red arrow] (281) to (237);
		\draw [style=red arrow] (282) to (208);
		\draw [style=red arrow] (294) to (209);
		\draw [style=red arrow] (291) to (210);
		\draw [style=red arrow] (295) to (203);
		\draw [style=red arrow] (292) to (204);
		\draw [style=red arrow] (286) to (200);
		\draw [style=red arrow] (287) to (199);
	\end{pgfonlayer}
\end{tikzpicture}

\caption{Example of a one-step evolution of 2D fermions (from left to right). The red arrows represent the trajectories of the particle movements.}
\label{fig:2d-evolution-one-step}
\end{figure}
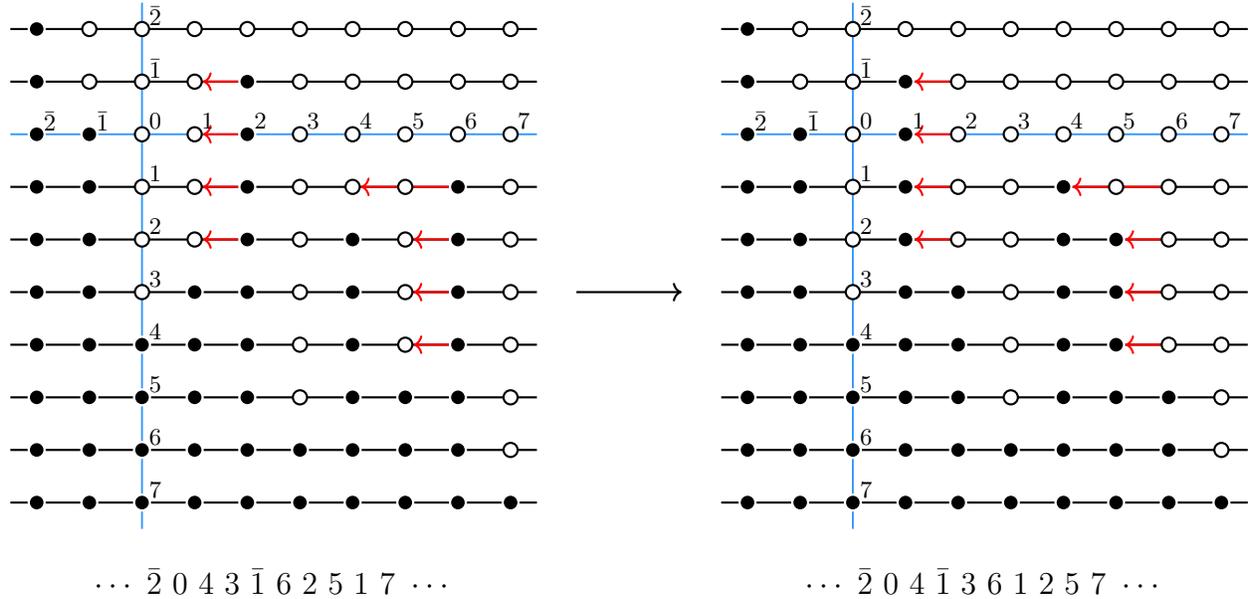

We now describe the evolutions of these 2D fermions. The particles can move to the left without exiting its row, and obey the same rules as the 1D fermions. Moreover, since an evolution has to preserve the axioms of 2D fermions (\Cref{def:2d-fermion}), it is automatic that if one particle is moving then all other particles in the same column have to move when possible.
\Cref{fig:2d-evolution-one-step} gives an example of such an evolution. This is in fact a three dimensional analogue of the lattice paths of \Cref{fig:vertex-one-step}. To make these 3D lattice paths more tractable, we give simplified model using certain tilings.

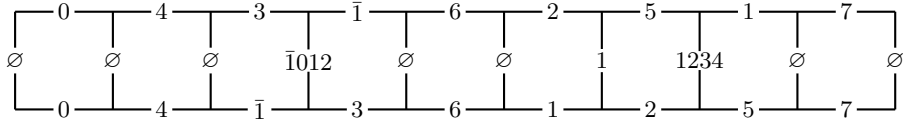
\begin{figure}[h]
\centering
\begin{tikzpicture}[scale = 0.65]
	\begin{pgfonlayer}{nodelayer}
		\node [style=none] (331) at (3, -13) {};
		\node [style=none] (332) at (5, -13) {};
		\node [style=none] (333) at (7, -13) {};
		\node [style=none] (334) at (9, -13) {};
		\node [style=none] (335) at (11, -13) {};
		\node [style=none] (336) at (13, -13) {};
		\node [style=none] (337) at (15, -13) {};
		\node [style=none] (338) at (17, -13) {};
		\node [style=none] (339) at (3, -15) {};
		\node [style=none] (340) at (5, -15) {};
		\node [style=none] (341) at (7, -15) {};
		\node [style=none] (342) at (9, -15) {};
		\node [style=none] (343) at (11, -15) {};
		\node [style=none] (344) at (13, -15) {};
		\node [style=none] (345) at (15, -15) {};
		\node [style=none] (346) at (17, -15) {};
		\node [style=none] (347) at (1, -13) {};
		\node [style=none] (348) at (1, -15) {};
		\node [style=none] (349) at (19, -13) {};
		\node [style=none] (350) at (19, -15) {};
		\node [style=none] (351) at (1, -13) {};
		\node [style=none] (352) at (1, -15) {};
		\node [style=none] (353) at (19, -15) {};
		\node [style=none] (354) at (19, -13) {};
		\node [style=label-white] (355) at (2, -13) {$0$};
		\node [style=label-white] (356) at (4, -13) {};
		\node [style=label-white] (357) at (4, -13) {$4$};
		\node [style=label-white] (358) at (6, -13) {};
		\node [style=label-white] (359) at (6, -13) {};
		\node [style=label-white] (360) at (6, -13) {$3$};
		\node [style=label-white] (361) at (8, -13) {$\bar 1$};
		\node [style=label-white] (362) at (10, -13) {$6$};
		\node [style=label-white] (363) at (12, -13) {$2$};
		\node [style=label-white] (364) at (14, -13) {$5$};
		\node [style=label-white] (365) at (16, -13) {$1$};
		\node [style=label-white] (366) at (18, -13) {$7$};
		\node [style=label-white] (367) at (2, -15) {$0$};
		\node [style=label-white] (368) at (4, -15) {$4$};
		\node [style=label-white] (369) at (6, -15) {$\bar 1$};
		\node [style=label-white] (370) at (8, -15) {$3$};
		\node [style=label-white] (371) at (10, -15) {$6$};
		\node [style=label-white] (372) at (12, -15) {$1$};
		\node [style=label-white] (373) at (14, -15) {$2$};
		\node [style=label-white] (374) at (16, -15) {$5$};
		\node [style=label-white] (375) at (18, -15) {$7$};
		\node [style=label-white] (376) at (1, -14) {$\emptyset$};
		\node [style=label-white] (377) at (3, -14) {$\emptyset$};
		\node [style=label-white] (378) at (5, -14) {$\emptyset$};
		\node [style=label-white] (379) at (19, -14) {$\emptyset$};
		\node [style=label-white] (380) at (17, -14) {$\emptyset$};
		\node [style=label-white] (381) at (11, -14) {$\emptyset$};
		\node [style=label-white] (382) at (9, -14) {$\emptyset$};
		\node [style=label-small] (384) at (7, -14) {$\bar 1 0 1 2$};
		\node [style=label-white] (386) at (13, -14) {$1$};
		\node [style=label-small] (387) at (15, -14) {$1234$};
	\end{pgfonlayer}
	\begin{pgfonlayer}{edgelayer}
		\draw [style=fade] (347.center) to (349.center);
		\draw [style=fade] (350.center) to (348.center);
		\draw [style=fade] (339.center) to (331.center);
		\draw [style=fade] (332.center) to (340.center);
		\draw [style=fade] (334.center) to (342.center);
		\draw [style=fade] (335.center) to (343.center);
		\draw [style=fade] (338.center) to (346.center);
		\draw [style=fade] (351.center) to (352.center);
		\draw [style=fade] (353.center) to (354.center);
		\draw [style=fade] (333.center) to (384);
		\draw [style=fade] (384) to (341.center);
		\draw [style=fade] (336.center) to (386);
		\draw [style=fade] (386) to (344.center);
		\draw [style=fade] (387) to (337.center);
		\draw [style=fade] (387) to (345.center);
	\end{pgfonlayer}
\end{tikzpicture}
\caption{Single-row tiling corresponding to the fermionic evolution in \Cref{fig:2d-evolution-one-step}. The non-empty intervals $[-1,2]$, $\{1\}$, and $[1,4]$ are exactly the row-indices of the trajectories of particle movements in \Cref{fig:2d-evolution-one-step}.}	
\label{fig:one-row-tiling}
\end{figure}

Generalizing the tiles in \Cref{fig:5-vertex-wt}, we define new tiles in which the particles/holes are replaced with numbers and lattice paths are replaced with intervals of $\zz$. The simplest way to define these tiles is the following
\[\begin{tikzpicture}[scale = 0.8]
	\begin{pgfonlayer}{nodelayer}
		\node [style=label-white, scale = 1.2, inner sep = 0.5mm] (323) at (7, -11) {$a$};
		\node [style=label-white, scale = 1.2, inner sep = 0.5mm] (324) at (7, -13) {$b$};
		\node [style=label-white, scale = 1.2, inner sep = 0.5mm] (325) at (6, -12) {$I$};
		\node [style=label-white, scale = 1.2, inner sep = 0.5mm] (326) at (8, -12) {$J$};
		\node [style=none] (327) at (6, -11) {};
		\node [style=none] (328) at (6, -13) {};
		\node [style=none] (329) at (8, -13) {};
		\node [style=none] (330) at (8, -11) {};
	\end{pgfonlayer}
	\begin{pgfonlayer}{edgelayer}
		\draw [style = fade] (330.center)
			 to (329.center)
			 to (328.center)
			 to (327.center)
			 to cycle;
	\end{pgfonlayer}
\end{tikzpicture}
\]
where $a,b\in\zz$ and $I,J\subset \zz$ are intervals satisfying $I\sqcup [b,\infty) =J\sqcup [a,\infty)=[\min(a,b),\infty)$. The weight is $x$ if $I\neq \emptyset$ and $1$ otherwise. Here the numbers $a,b$ represent the columns of the 2D Maya diagrams, and the intervals $I,J$ indicate particles1 moving along the rows. In particular, if $i\in I$ then there is a particle passing by the column along the $i$-th row. One can check that the rules of evolutions of 2D fermions exactly match the conditions of the vertex configuration. For example, the particle movement in \Cref{fig:2d-evolution-one-step} is equivalent to the tiling in \Cref{fig:one-row-tiling}. Abusing notation, we will define transfer matrices $T(x)$ in a similar fashion.

\begin{theorem}\label{thm:stanley_transfer_matrix}
	Let $T(x)$ be the transfer matrices for 2D fermions. We have that
	\[\langle u|T(x_n)\cdots T(x_2)T(x_1)|w\rangle =\mf_{w/u}(x_1,\cdots, x_n)\]
	the skew Stanley symmetric polynomials in $n$ variables.
\end{theorem}
\begin{proof}
	By definition of Stanley symmetric functions, it suffices to show that 
	\[\langle u|T(x)|w\rangle = \sum_{v} x^{\ell(v)}\]
	where the sum is over $v$ such that $w\doteq uv$ and $v$ has an increasing reduced word.
	Observe that if $w s_i < w$, then one can obtain $w s_i$ from $w$ via a fermionic evolution by moving all possible particles in the $(i+1)$-st column to the $i$-th column. For example,
	
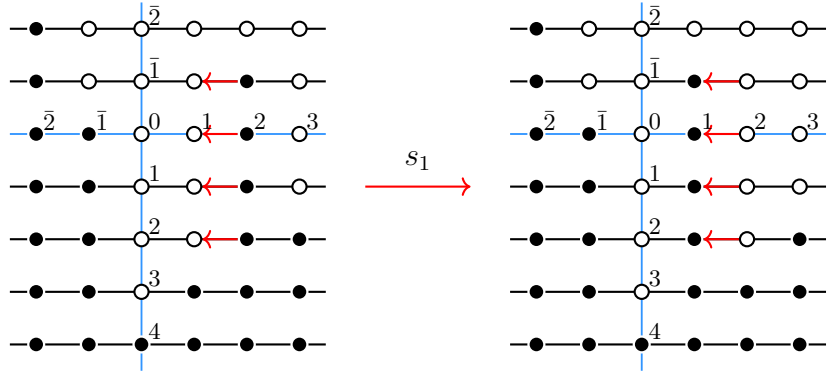
\[\begin{tikzpicture}[scale = 0.7]
	\begin{pgfonlayer}{nodelayer}
		\node [style=label-small] (0) at (1.5, 2.25) {$0$};
		\node [style=label-small] (1) at (2.5, 2.25) {$1$};
		\node [style=label-small] (2) at (3.5, 2.25) {$2$};
		\node [style=label-small] (3) at (4.5, 2.25) {$3$};
		\node [style=label-small] (8) at (0.5, 2.25) {$\bar 1$};
		\node [style=label-small] (9) at (-0.5, 2.25) {$\bar 2$};
		\node [style=label-small] (14) at (1.5, 3.25) {$\bar 1$};
		\node [style=label-small] (15) at (1.5, 4.25) {$\bar 2$};
		\node [style=label-small] (24) at (1.5, 1.25) {$1$};
		\node [style=label-small] (25) at (1.5, 0.25) {$2$};
		\node [style=label-small] (26) at (1.5, -0.75) {$3$};
		\node [style=label-small] (27) at (1.5, -1.75) {$4$};
		\node [style=none] (33) at (-1.25, 2) {};
		\node [style=none] (34) at (4.75, 2) {};
		\node [style=vert] (35) at (4.25, 0) {};
		\node [style=vert] (36) at (4.25, -1) {};
		\node [style=vert] (37) at (3.25, -2) {};
		\node [style=vert] (48) at (4.25, -2) {};
		\node [style=vert] (53) at (3.25, -1) {};
		\node [style=vert] (54) at (2.25, -2) {};
		\node [style=vert] (56) at (3.25, 2) {};
		\node [style=vert] (57) at (3.25, 1) {};
		\node [style=vert] (58) at (3.25, 0) {};
		\node [style=vert] (59) at (2.25, -1) {};
		\node [style=vert] (60) at (1.25, -2) {};
		\node [style=vert] (62) at (0.25, 2) {};
		\node [style=vert] (63) at (0.25, 1) {};
		\node [style=vert] (64) at (0.25, 0) {};
		\node [style=vert] (65) at (0.25, -1) {};
		\node [style=vert] (66) at (0.25, -2) {};
		\node [style=vert] (88) at (3.25, 3) {};
		\node [style=vert] (90) at (-0.75, 3) {};
		\node [style=vert] (91) at (-0.75, 2) {};
		\node [style=vert] (92) at (-0.75, 1) {};
		\node [style=vert] (93) at (-0.75, 0) {};
		\node [style=vert] (94) at (-0.75, -1) {};
		\node [style=vert] (95) at (-0.75, -2) {};
		\node [style=vert] (98) at (-0.75, 4) {};
		\node [style=none] (101) at (1.25, 4.5) {};
		\node [style=none] (102) at (1.25, -2.5) {};
		\node [style=none] (111) at (-1.25, 4) {};
		\node [style=none] (112) at (-1.25, 3) {};
		\node [style=none] (113) at (-1.25, 1) {};
		\node [style=none] (114) at (-1.25, 0) {};
		\node [style=none] (115) at (-1.25, -1) {};
		\node [style=none] (116) at (-1.25, -2) {};
		\node [style=none] (120) at (4.75, 2) {};
		\node [style=none] (122) at (4.75, 4) {};
		\node [style=none] (123) at (4.75, 3) {};
		\node [style=none] (124) at (4.75, 1) {};
		\node [style=none] (125) at (4.75, 0) {};
		\node [style=none] (126) at (4.75, -1) {};
		\node [style=none] (127) at (4.75, -2) {};
		\node [style=vert-white] (131) at (0.25, 4) {};
		\node [style=vert-white] (132) at (1.25, 4) {};
		\node [style=vert-white] (133) at (1.25, 3) {};
		\node [style=vert-white] (134) at (0.25, 3) {};
		\node [style=vert-white] (135) at (1.25, 2) {};
		\node [style=vert-white] (136) at (1.25, 1) {};
		\node [style=vert-white] (137) at (2.25, 4) {};
		\node [style=vert-white] (138) at (3.25, 4) {};
		\node [style=vert-white] (139) at (2.25, 3) {};
		\node [style=vert-white] (140) at (2.25, 2) {};
		\node [style=vert-white] (141) at (4.25, 1) {};
		\node [style=vert-white] (148) at (1.25, -1) {};
		\node [style=vert-white] (149) at (2.25, 0) {};
		\node [style=vert-white] (151) at (1.25, 0) {};
		\node [style=vert-white] (152) at (2.25, 1) {};
		\node [style=vert-white] (155) at (4.25, 2) {};
		\node [style=vert-white] (174) at (4.25, 4) {};
		\node [style=vert-white] (175) at (4.25, 3) {};
		\node [style=none] (318) at (5.5, 1) {};
		\node [style=none] (320) at (7.5, 1) {};
		\node [style=label-small] (321) at (11, 2.25) {$0$};
		\node [style=label-small] (322) at (12, 2.25) {$1$};
		\node [style=label-small] (323) at (13, 2.25) {$2$};
		\node [style=label-small] (324) at (14, 2.25) {$3$};
		\node [style=label-small] (325) at (10, 2.25) {$\bar 1$};
		\node [style=label-small] (326) at (9, 2.25) {$\bar 2$};
		\node [style=label-small] (327) at (11, 3.25) {$\bar 1$};
		\node [style=label-small] (328) at (11, 4.25) {$\bar 2$};
		\node [style=label-small] (329) at (11, 1.25) {$1$};
		\node [style=label-small] (330) at (11, 0.25) {$2$};
		\node [style=label-small] (331) at (11, -0.75) {$3$};
		\node [style=label-small] (332) at (11, -1.75) {$4$};
		\node [style=none] (333) at (8.25, 2) {};
		\node [style=none] (334) at (14.25, 2) {};
		\node [style=vert] (335) at (13.75, 0) {};
		\node [style=vert] (336) at (13.75, -1) {};
		\node [style=vert] (337) at (12.75, -2) {};
		\node [style=vert] (338) at (13.75, -2) {};
		\node [style=vert] (339) at (12.75, -1) {};
		\node [style=vert] (340) at (11.75, -2) {};
		\node [style=vert-white] (341) at (12.75, 2) {};
		\node [style=vert-white] (342) at (12.75, 1) {};
		\node [style=vert-white] (343) at (12.75, 0) {};
		\node [style=vert] (344) at (11.75, -1) {};
		\node [style=vert] (345) at (10.75, -2) {};
		\node [style=vert] (346) at (9.75, 2) {};
		\node [style=vert] (347) at (9.75, 1) {};
		\node [style=vert] (348) at (9.75, 0) {};
		\node [style=vert] (349) at (9.75, -1) {};
		\node [style=vert] (350) at (9.75, -2) {};
		\node [style=vert-white] (351) at (12.75, 3) {};
		\node [style=vert] (352) at (8.75, 3) {};
		\node [style=vert] (353) at (8.75, 2) {};
		\node [style=vert] (354) at (8.75, 1) {};
		\node [style=vert] (355) at (8.75, 0) {};
		\node [style=vert] (356) at (8.75, -1) {};
		\node [style=vert] (357) at (8.75, -2) {};
		\node [style=vert] (358) at (8.75, 4) {};
		\node [style=none] (359) at (10.75, 4.5) {};
		\node [style=none] (360) at (10.75, -2.5) {};
		\node [style=none] (361) at (8.25, 4) {};
		\node [style=none] (362) at (8.25, 3) {};
		\node [style=none] (363) at (8.25, 1) {};
		\node [style=none] (364) at (8.25, 0) {};
		\node [style=none] (365) at (8.25, -1) {};
		\node [style=none] (366) at (8.25, -2) {};
		\node [style=none] (367) at (14.25, 2) {};
		\node [style=none] (368) at (14.25, 4) {};
		\node [style=none] (369) at (14.25, 3) {};
		\node [style=none] (370) at (14.25, 1) {};
		\node [style=none] (371) at (14.25, 0) {};
		\node [style=none] (372) at (14.25, -1) {};
		\node [style=none] (373) at (14.25, -2) {};
		\node [style=vert-white] (374) at (9.75, 4) {};
		\node [style=vert-white] (375) at (10.75, 4) {};
		\node [style=vert-white] (376) at (10.75, 3) {};
		\node [style=vert-white] (377) at (9.75, 3) {};
		\node [style=vert-white] (378) at (10.75, 2) {};
		\node [style=vert-white] (379) at (10.75, 1) {};
		\node [style=vert-white] (380) at (11.75, 4) {};
		\node [style=vert-white] (381) at (12.75, 4) {};
		\node [style=vert] (382) at (11.75, 3) {};
		\node [style=vert] (383) at (11.75, 2) {};
		\node [style=vert-white] (384) at (13.75, 1) {};
		\node [style=vert-white] (385) at (10.75, -1) {};
		\node [style=vert] (386) at (11.75, 0) {};
		\node [style=vert-white] (387) at (10.75, 0) {};
		\node [style=vert] (388) at (11.75, 1) {};
		\node [style=vert-white] (389) at (13.75, 2) {};
		\node [style=vert-white] (390) at (13.75, 4) {};
		\node [style=vert-white] (391) at (13.75, 3) {};
		\node [style=none] (392) at (6.5, 1.5) {$s_1$};
	\end{pgfonlayer}
	\begin{pgfonlayer}{edgelayer}
		\draw [style=blue] (33.center) to (34.center);
		\draw [style=blue] (101.center) to (102.center);
		\draw [style=fade] (111.center) to (122.center);
		\draw [style=fade] (112.center) to (123.center);
		\draw [style=fade] (113.center) to (124.center);
		\draw [style=fade] (114.center) to (125.center);
		\draw [style=fade] (115.center) to (126.center);
		\draw [style=fade] (127.center) to (116.center);
		\draw [style=red arrow thick] (318.center) to (320.center);
		\draw [style=red arrow] (58) to (149);
		\draw [style=red arrow] (57) to (152);
		\draw [style=red arrow] (56) to (140);
		\draw [style=red arrow] (88) to (139);
		\draw [style=blue] (333.center) to (334.center);
		\draw [style=blue] (359.center) to (360.center);
		\draw [style=fade] (361.center) to (368.center);
		\draw [style=fade] (362.center) to (369.center);
		\draw [style=fade] (363.center) to (370.center);
		\draw [style=fade] (364.center) to (371.center);
		\draw [style=fade] (365.center) to (372.center);
		\draw [style=fade] (373.center) to (366.center);
		\draw [style=red arrow] (343) to (386);
		\draw [style=red arrow] (342) to (388);
		\draw [style=red arrow] (341) to (383);
		\draw [style=red arrow] (351) to (382);
	\end{pgfonlayer}
\end{tikzpicture}
\]
Now if particle movements take place at multiple columns, the non-intersecting property of the fermionic paths guarantees that the evolution is equivalent to right-multiplying $s_i$ from left to right for each column-$i$ with a particle movement, which is equivalent to $v$ being increasing. For example, in \Cref{fig:2d-evolution-one-step}, the corresponding increasing reduced word is $s_1s_4s_5$. The theorem then follows by induction by iteratively applying the transfer matrix $T$.
\end{proof}

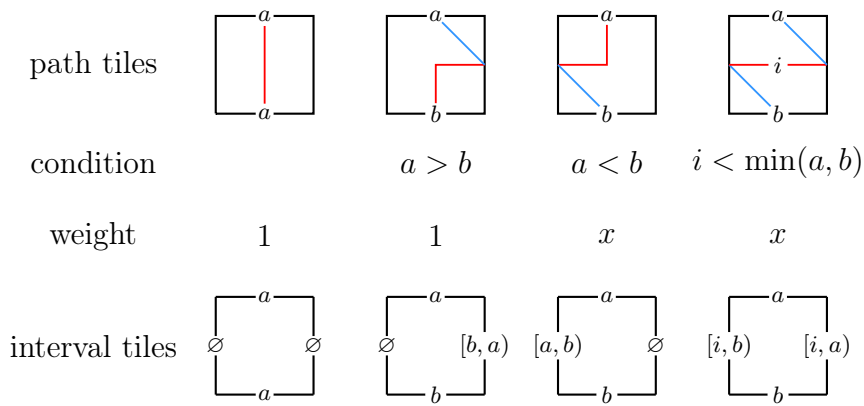
\begin{figure}[h]\centering
\begin{tikzpicture}[scale = 0.65]
	\begin{pgfonlayer}{nodelayer}
		\node [style=none] (0) at (-2, 1) {};
		\node [style=none] (1) at (0, 1) {};
		\node [style=none] (2) at (0, -1) {};
		\node [style=none] (3) at (-2, -1) {};
		\node [style=none] (4) at (1.5, 1) {};
		\node [style=none] (5) at (3.5, 1) {};
		\node [style=none] (6) at (3.5, -1) {};
		\node [style=none] (7) at (1.5, -1) {};
		\node [style=none] (8) at (5, 1) {};
		\node [style=none] (9) at (7, 1) {};
		\node [style=none] (10) at (7, -1) {};
		\node [style=none] (11) at (5, -1) {};
		\node [style=none] (12) at (-5.5, 1) {};
		\node [style=none] (13) at (-3.5, 1) {};
		\node [style=none] (14) at (-3.5, -1) {};
		\node [style=none] (15) at (-5.5, -1) {};
		\node [style=none] (16) at (-8, -2) {condition};
		\node [style=none] (17) at (-8, -3.5) {weight};
		\node [style=none] (18) at (-1, 1) {};
		\node [style=none] (19) at (0, 0) {};
		\node [style=none] (20) at (-1, -1) {};
		\node [style=none] (21) at (-1, 0) {};
		\node [style=none] (22) at (2.5, 0) {};
		\node [style=none] (23) at (2.5, 1) {};
		\node [style=none] (24) at (1.5, 0) {};
		\node [style=none] (25) at (2.5, -1) {};
		\node [style=none] (26) at (5, 0) {};
		\node [style=none] (27) at (7, 0) {};
		\node [style=none] (28) at (6, -1) {};
		\node [style=none] (29) at (6, 1) {};
		\node [style=label-white] (30) at (-4.5, 1) {$a$};
		\node [style=label-white] (31) at (-4.5, -1) {$a$};
		\node [style=label-white] (32) at (-1, -1) {$b$};
		\node [style=label-white] (33) at (-1, 1) {$a$};
		\node [style=label-white] (34) at (2.5, 1) {$a$};
		\node [style=label-white] (35) at (2.5, -1) {$b$};
		\node [style=label-white] (36) at (6, 1) {$a$};
		\node [style=label-white] (37) at (6, 0) {$i$};
		\node [style=label-white] (38) at (6, -1) {$b$};
		\node [style=none] (39) at (-1, -2) {$a>b$};
		\node [style=none] (40) at (2.5, -2) {$a<b$};
		\node [style=none] (41) at (6, -2) {$i<\min(a,b)$};
		\node [style=none] (42) at (-1, -3.5) {$1$};
		\node [style=none] (43) at (2.5, -3.5) {$x$};
		\node [style=none] (44) at (-4.5, -3.5) {$1$};
		\node [style=none] (45) at (6, -3.5) {$x$};
		\node [style=none] (47) at (0, -4.75) {};
		\node [style=none] (49) at (-2, -6.75) {};
		\node [style=none] (62) at (-1, -4.75) {};
		\node [style=none] (64) at (-1, -6.75) {};
		\node [style=label-white] (76) at (-1, -6.75) {$b$};
		\node [style=label-white] (77) at (-1, -4.75) {$a$};
		\node [style=none] (83) at (-8, -5.75) {interval tiles};
		\node [style=none] (84) at (-8, 0) {path tiles};
		\node [style=none] (85) at (-2, -4.75) {};
		\node [style=none] (86) at (0, -6.75) {};
		\node [style=label-small] (87) at (0, -5.75) {$[b,a)$};
		\node [style=label-small] (88) at (-2, -5.75) {$\emptyset$};
		\node [style=none] (89) at (3.5, -4.75) {};
		\node [style=none] (90) at (1.5, -6.75) {};
		\node [style=none] (91) at (2.5, -4.75) {};
		\node [style=none] (92) at (2.5, -6.75) {};
		\node [style=label-white] (93) at (2.5, -6.75) {$b$};
		\node [style=label-white] (94) at (2.5, -4.75) {$a$};
		\node [style=none] (95) at (1.5, -4.75) {};
		\node [style=none] (96) at (3.5, -6.75) {};
		\node [style=label-small] (97) at (3.5, -5.75) {$\emptyset$};
		\node [style=label-small] (98) at (1.5, -5.75) {$[a,b)$};
		\node [style=none] (99) at (-3.5, -4.75) {};
		\node [style=none] (100) at (-5.5, -6.75) {};
		\node [style=none] (101) at (-4.5, -4.75) {};
		\node [style=none] (102) at (-4.5, -6.75) {};
		\node [style=label-white] (103) at (-4.5, -6.75) {$a$};
		\node [style=label-white] (104) at (-4.5, -4.75) {$a$};
		\node [style=none] (105) at (-5.5, -4.75) {};
		\node [style=none] (106) at (-3.5, -6.75) {};
		\node [style=label-small] (107) at (-3.5, -5.75) {$\emptyset$};
		\node [style=label-small] (108) at (-5.5, -5.75) {$\emptyset$};
		\node [style=none] (109) at (7, -4.75) {};
		\node [style=none] (110) at (5, -6.75) {};
		\node [style=none] (111) at (6, -4.75) {};
		\node [style=none] (112) at (6, -6.75) {};
		\node [style=label-white] (113) at (6, -6.75) {$b$};
		\node [style=label-white] (114) at (6, -4.75) {$a$};
		\node [style=none] (115) at (5, -4.75) {};
		\node [style=none] (116) at (7, -6.75) {};
		\node [style=label-small] (117) at (7, -5.75) {$[i,a)$};
		\node [style=label-small] (118) at (5, -5.75) {$[i,b)$};
	\end{pgfonlayer}
	\begin{pgfonlayer}{edgelayer}
		\draw [style=fade] (1.center)
			 to (2.center)
			 to (3.center)
			 to (0.center)
			 to cycle;
		\draw [style=fade] (7.center)
			 to (4.center)
			 to (5.center)
			 to (6.center)
			 to cycle;
		\draw [style=fade] (10.center)
			 to (11.center)
			 to (8.center)
			 to (9.center)
			 to cycle;
		\draw [style=fade] (15.center)
			 to (12.center)
			 to (13.center)
			 to (14.center)
			 to cycle;
		\draw [style=red] (30) to (31);
		\draw [style=blue] (26.center) to (28.center);
		\draw [style=red] (24.center)
			 to (22.center)
			 to (23.center);
		\draw [style=red] (20.center)
			 to (21.center)
			 to (19.center);
		\draw [style=blue] (18.center) to (19.center);
		\draw [style=blue] (24.center) to (25.center);
		\draw [style=red] (26.center) to (27.center);
		\draw [style=blue] (29.center) to (27.center);
		\draw [style=fade] (85.center) to (47.center);
		\draw [style=fade] (49.center) to (86.center);
		\draw [style=fade] (86.center) to (87);
		\draw [style=fade] (47.center) to (87);
		\draw [style=fade] (85.center) to (88);
		\draw [style=fade] (88) to (49.center);
		\draw [style=fade] (95.center) to (89.center);
		\draw [style=fade] (90.center) to (96.center);
		\draw [style=fade] (96.center) to (97);
		\draw [style=fade] (89.center) to (97);
		\draw [style=fade] (95.center) to (98);
		\draw [style=fade] (98) to (90.center);
		\draw [style=fade] (105.center) to (99.center);
		\draw [style=fade] (100.center) to (106.center);
		\draw [style=fade] (106.center) to (107);
		\draw [style=fade] (99.center) to (107);
		\draw [style=fade] (105.center) to (108);
		\draw [style=fade] (108) to (100.center);
		\draw [style=fade] (115.center) to (109.center);
		\draw [style=fade] (110.center) to (116.center);
		\draw [style=fade] (116.center) to (117);
		\draw [style=fade] (109.center) to (117);
		\draw [style=fade] (115.center) to (118);
		\draw [style=fade] (118) to (110.center);
	\end{pgfonlayer}
\end{tikzpicture}
\caption{1st row: New tiles of labeled paths where path of different labels are colored differently. 2nd row: A tile is valid if the condition is satisfied. 3rd row: The weight of the tiles. 4th row: bijection with the interval tiles. }
\label{fig:new-tiles}
\end{figure}

We now give another alternative description of the vertex model, using tiles of labeled paths. There are four types of tiles given in \Cref{fig:new-tiles}. There is a simple bijection between the path version and the interval version of the tiles, also given in \Cref{fig:new-tiles}. When two paths meet, we use red (resp. blue) color to indicate the path with smaller (resp. larger) label. The connectivity depends on path labels rather than the colors. 

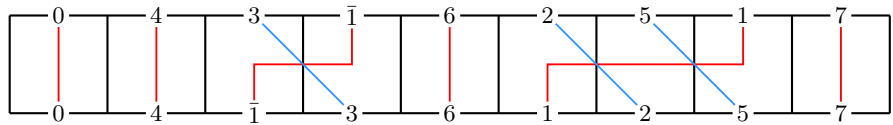
\begin{figure}[h]
	\begin{tikzpicture}[scale=0.65]
	\begin{pgfonlayer}{nodelayer}
		\node [style=none] (331) at (3, -13) {};
		\node [style=none] (332) at (5, -13) {};
		\node [style=none] (333) at (7, -13) {};
		\node [style=none] (334) at (9, -13) {};
		\node [style=none] (335) at (11, -13) {};
		\node [style=none] (336) at (13, -13) {};
		\node [style=none] (337) at (15, -13) {};
		\node [style=none] (338) at (17, -13) {};
		\node [style=none] (339) at (3, -15) {};
		\node [style=none] (340) at (5, -15) {};
		\node [style=none] (341) at (7, -15) {};
		\node [style=none] (342) at (9, -15) {};
		\node [style=none] (343) at (11, -15) {};
		\node [style=none] (344) at (13, -15) {};
		\node [style=none] (345) at (15, -15) {};
		\node [style=none] (346) at (17, -15) {};
		\node [style=none] (347) at (1, -13) {};
		\node [style=none] (348) at (1, -15) {};
		\node [style=none] (349) at (19, -13) {};
		\node [style=none] (350) at (19, -15) {};
		\node [style=none] (351) at (1, -13) {};
		\node [style=none] (352) at (1, -15) {};
		\node [style=none] (353) at (19, -15) {};
		\node [style=none] (354) at (19, -13) {};
		\node [style=label-white] (355) at (2, -13) {$0$};
		\node [style=label-white] (356) at (4, -13) {};
		\node [style=label-white] (357) at (4, -13) {$4$};
		\node [style=label-white] (358) at (6, -13) {};
		\node [style=label-white] (359) at (6, -13) {};
		\node [style=label-white] (360) at (6, -13) {$3$};
		\node [style=label-white] (361) at (8, -13) {$\bar 1$};
		\node [style=label-white] (362) at (10, -13) {$6$};
		\node [style=label-white] (363) at (12, -13) {$2$};
		\node [style=label-white] (364) at (14, -13) {$5$};
		\node [style=label-white] (365) at (16, -13) {$1$};
		\node [style=label-white] (366) at (18, -13) {$7$};
		\node [style=label-white] (367) at (2, -15) {$0$};
		\node [style=label-white] (368) at (4, -15) {$4$};
		\node [style=label-white] (369) at (6, -15) {$\bar 1$};
		\node [style=label-white] (370) at (8, -15) {$3$};
		\node [style=label-white] (371) at (10, -15) {$6$};
		\node [style=label-white] (372) at (12, -15) {$1$};
		\node [style=label-white] (373) at (14, -15) {$2$};
		\node [style=label-white] (374) at (16, -15) {$5$};
		\node [style=label-white] (375) at (18, -15) {$7$};
		\node [style=none] (376) at (6, -14) {};
		\node [style=none] (377) at (7, -14) {};
		\node [style=none] (378) at (8, -14) {};
		\node [style=none] (379) at (12, -14) {};
		\node [style=none] (380) at (13, -14) {};
		\node [style=none] (381) at (14, -14) {};
		\node [style=none] (382) at (15, -14) {};
		\node [style=none] (383) at (16, -14) {};
	\end{pgfonlayer}
	\begin{pgfonlayer}{edgelayer}
		\draw [style=fade] (347.center) to (349.center);
		\draw [style=fade] (350.center) to (348.center);
		\draw [style=fade] (339.center) to (331.center);
		\draw [style=fade] (332.center) to (340.center);
		\draw [style=fade] (334.center) to (342.center);
		\draw [style=fade] (335.center) to (343.center);
		\draw [style=fade] (338.center) to (346.center);
		\draw [style=fade] (351.center) to (352.center);
		\draw [style=fade] (353.center) to (354.center);
		\draw [style=fade] (333.center) to (341.center);
		\draw [style=fade] (336.center) to (344.center);
		\draw [style=fade] (337.center) to (345.center);
		\draw [style=red] (355) to (367);
		\draw [style=red] (357) to (368);
		\draw [style=red] (369)
			 to (376.center)
			 to (378.center)
			 to (361);
		\draw [style=red] (372)
			 to (379.center)
			 to (383.center)
			 to (365);
		\draw [style=red] (371) to (362);
		\draw [style=red] (375) to (366);
		\draw [style=blue] (360) to (370);
		\draw [style=blue] (363) to (373);
		\draw [style=blue] (364) to (374);
	\end{pgfonlayer}
\end{tikzpicture}
	\caption{Reinterpretation of \Cref{fig:one-row-tiling} using labeled paths.}
\label{fig:path-one-row-tiling}
\end{figure}

For example, with the new tiles, the one-row tiling in \Cref{fig:one-row-tiling} can be reinterpreted as follows in \Cref{fig:path-one-row-tiling}

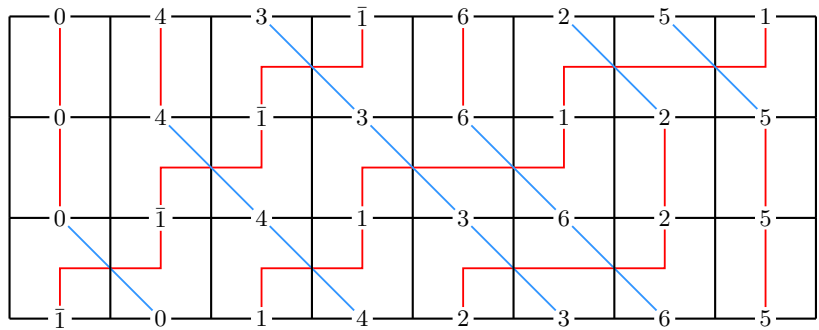
\begin{figure}[h]
\centering
\begin{tikzpicture}[scale = 0.67]
	\begin{pgfonlayer}{nodelayer}
		\node [style=none] (331) at (3, -13) {};
		\node [style=none] (332) at (5, -13) {};
		\node [style=none] (333) at (7, -13) {};
		\node [style=none] (334) at (9, -13) {};
		\node [style=none] (335) at (11, -13) {};
		\node [style=none] (336) at (13, -13) {};
		\node [style=none] (337) at (15, -13) {};
		\node [style=none] (338) at (17, -13) {};
		\node [style=none] (339) at (3, -15) {};
		\node [style=none] (340) at (5, -15) {};
		\node [style=none] (341) at (7, -15) {};
		\node [style=none] (342) at (9, -15) {};
		\node [style=none] (343) at (11, -15) {};
		\node [style=none] (344) at (13, -15) {};
		\node [style=none] (345) at (15, -15) {};
		\node [style=none] (346) at (17, -15) {};
		\node [style=none] (347) at (1, -13) {};
		\node [style=none] (348) at (1, -15) {};
		\node [style=none] (351) at (1, -13) {};
		\node [style=none] (352) at (1, -15) {};
		\node [style=label-white] (355) at (2, -13) {$0$};
		\node [style=label-white] (356) at (4, -13) {};
		\node [style=label-white] (357) at (4, -13) {$4$};
		\node [style=label-white] (358) at (6, -13) {};
		\node [style=label-white] (359) at (6, -13) {};
		\node [style=label-white] (360) at (6, -13) {$3$};
		\node [style=label-white] (361) at (8, -13) {$\bar 1$};
		\node [style=label-white] (362) at (10, -13) {$6$};
		\node [style=label-white] (363) at (12, -13) {$2$};
		\node [style=label-white] (364) at (14, -13) {$5$};
		\node [style=label-white] (365) at (16, -13) {$1$};
		\node [style=label-white] (367) at (2, -15) {$0$};
		\node [style=label-white] (368) at (4, -15) {$4$};
		\node [style=label-white] (369) at (6, -15) {$\bar 1$};
		\node [style=label-white] (370) at (8, -15) {$3$};
		\node [style=label-white] (371) at (10, -15) {$6$};
		\node [style=label-white] (372) at (12, -15) {$1$};
		\node [style=label-white] (373) at (14, -15) {$2$};
		\node [style=label-white] (374) at (16, -15) {$5$};
		\node [style=none] (376) at (6, -14) {};
		\node [style=none] (377) at (7, -14) {};
		\node [style=none] (378) at (8, -14) {};
		\node [style=none] (379) at (12, -14) {};
		\node [style=none] (380) at (13, -14) {};
		\node [style=none] (381) at (14, -14) {};
		\node [style=none] (382) at (15, -14) {};
		\node [style=none] (383) at (16, -14) {};
		\node [style=none] (384) at (3, -17) {};
		\node [style=none] (385) at (5, -17) {};
		\node [style=none] (386) at (7, -17) {};
		\node [style=none] (387) at (9, -17) {};
		\node [style=none] (388) at (11, -17) {};
		\node [style=none] (389) at (13, -17) {};
		\node [style=none] (390) at (15, -17) {};
		\node [style=none] (391) at (17, -17) {};
		\node [style=none] (392) at (1, -17) {};
		\node [style=none] (394) at (1, -17) {};
		\node [style=label-white] (396) at (2, -17) {$0$};
		\node [style=label-white] (397) at (4, -17) {$\bar 1$};
		\node [style=label-white] (398) at (6, -17) {$4$};
		\node [style=label-white] (399) at (8, -17) {$1$};
		\node [style=label-white] (400) at (10, -17) {$3$};
		\node [style=label-white] (401) at (12, -17) {$6$};
		\node [style=label-white] (402) at (14, -17) {$2$};
		\node [style=label-white] (403) at (16, -17) {$5$};
		\node [style=none] (405) at (12, -16) {};
		\node [style=none] (406) at (8, -16) {};
		\node [style=none] (407) at (6, -16) {};
		\node [style=none] (408) at (4, -16) {};
		\node [style=none] (409) at (3, -19) {};
		\node [style=none] (410) at (5, -19) {};
		\node [style=none] (411) at (7, -19) {};
		\node [style=none] (412) at (9, -19) {};
		\node [style=none] (413) at (11, -19) {};
		\node [style=none] (414) at (13, -19) {};
		\node [style=none] (415) at (15, -19) {};
		\node [style=none] (416) at (17, -19) {};
		\node [style=none] (417) at (1, -19) {};
		\node [style=none] (418) at (1, -19) {};
		\node [style=label-white] (419) at (2, -19) {$\bar 1$};
		\node [style=label-white] (420) at (4, -19) {$0$};
		\node [style=label-white] (421) at (6, -19) {$1$};
		\node [style=label-white] (422) at (8, -19) {$4$};
		\node [style=label-white] (423) at (10, -19) {$2$};
		\node [style=label-white] (424) at (12, -19) {$3$};
		\node [style=label-white] (425) at (14, -19) {$6$};
		\node [style=label-white] (426) at (16, -19) {$5$};
		\node [style=none] (427) at (14, -18) {};
		\node [style=none] (428) at (10, -18) {};
		\node [style=none] (429) at (8, -18) {};
		\node [style=none] (430) at (6, -18) {};
		\node [style=none] (431) at (4, -18) {};
		\node [style=none] (432) at (2, -18) {};
	\end{pgfonlayer}
	\begin{pgfonlayer}{edgelayer}
		\draw [style=fade] (339.center) to (331.center);
		\draw [style=fade] (332.center) to (340.center);
		\draw [style=fade] (334.center) to (342.center);
		\draw [style=fade] (335.center) to (343.center);
		\draw [style=fade] (338.center) to (346.center);
		\draw [style=fade] (351.center) to (352.center);
		\draw [style=fade] (333.center) to (341.center);
		\draw [style=fade] (336.center) to (344.center);
		\draw [style=fade] (337.center) to (345.center);
		\draw [style=red] (355) to (367);
		\draw [style=red] (357) to (368);
		\draw [style=red] (369)
			 to (376.center)
			 to (378.center)
			 to (361);
		\draw [style=red] (372)
			 to (379.center)
			 to (383.center)
			 to (365);
		\draw [style=red] (371) to (362);
		\draw [style=blue, in=135, out=-45] (360) to (370);
		\draw [style=blue] (363) to (373);
		\draw [style=blue] (364) to (374);
		\draw [style=fade] (352.center) to (394.center);
		\draw [style=fade] (339.center) to (384.center);
		\draw [style=fade] (340.center) to (385.center);
		\draw [style=fade] (341.center) to (386.center);
		\draw [style=fade] (342.center) to (387.center);
		\draw [style=fade] (343.center) to (388.center);
		\draw [style=fade] (344.center) to (389.center);
		\draw [style=fade] (345.center) to (390.center);
		\draw [style=fade] (346.center) to (391.center);
		\draw [style=red] (372)
			 to (405.center)
			 to (406.center)
			 to (399);
		\draw [style=blue] (370) to (400);
		\draw [style=blue] (371) to (401);
		\draw [style=red] (373) to (402);
		\draw [style=red] (374) to (403);
		\draw [style=fade] (351.center) to (338.center);
		\draw [style=fade] (394.center) to (391.center);
		\draw [style=fade] (352.center) to (346.center);
		\draw [style=red] (369)
			 to (407.center)
			 to (408.center)
			 to (397);
		\draw [style=blue] (368) to (398);
		\draw [style=red] (367) to (396);
		\draw [style=fade] (418.center) to (416.center);
		\draw [style=red] (402)
			 to (427.center)
			 to (428.center)
			 to (423);
		\draw [style=blue] (400) to (424);
		\draw [style=blue] (401) to (425);
		\draw [style=red] (403) to (426);
		\draw [style=red] (399)
			 to (429.center)
			 to (430.center)
			 to (421);
		\draw [style=red] (397)
			 to (431.center)
			 to (432.center)
			 to (419);
		\draw [style=blue] (396) to (420);
		\draw [style=blue] (422) to (398);
		\draw [style=fade] (385.center) to (410.center);
		\draw [style=fade] (409.center) to (384.center);
		\draw [style=fade] (418.center) to (394.center);
		\draw [style=fade] (411.center) to (386.center);
		\draw [style=fade] (387.center) to (412.center);
		\draw [style=fade] (388.center) to (413.center);
		\draw [style=fade] (389.center) to (414.center);
		\draw [style=fade] (390.center) to (415.center);
		\draw [style=fade] (391.center) to (416.center);
	\end{pgfonlayer}
\end{tikzpicture}
\caption{A 3-row tiling whose weight is $x_1^3x_2^3x_3^4$.}
\label{fig:3row-tiling}
\end{figure}

\begin{figure}[h]
\centering
\begin{tikzpicture}[scale=0.55]
	\begin{pgfonlayer}{nodelayer}
		\node [style=none] (433) at (27, -13) {};
		\node [style=none] (434) at (29, -13) {};
		\node [style=none] (435) at (31, -13) {};
		\node [style=none] (436) at (33, -13) {};
		\node [style=none] (437) at (35, -13) {};
		\node [style=none] (438) at (37, -13) {};
		\node [style=none] (439) at (39, -13) {};
		\node [style=none] (440) at (41, -13) {};
		\node [style=none] (441) at (25, -15) {};
		\node [style=none] (442) at (27, -15) {};
		\node [style=none] (443) at (29, -15) {};
		\node [style=none] (444) at (31, -15) {};
		\node [style=none] (445) at (33, -15) {};
		\node [style=none] (446) at (35, -15) {};
		\node [style=none] (447) at (37, -15) {};
		\node [style=none] (448) at (39, -15) {};
		\node [style=none] (449) at (25, -13) {};
		\node [style=none] (450) at (23, -15) {};
		\node [style=none] (451) at (25, -13) {};
		\node [style=none] (452) at (23, -15) {};
		\node [style=label-white] (453) at (26, -13) {$0$};
		\node [style=label-white] (454) at (28, -13) {};
		\node [style=label-white] (455) at (28, -13) {$4$};
		\node [style=label-white] (456) at (30, -13) {};
		\node [style=label-white] (457) at (30, -13) {};
		\node [style=label-white] (458) at (30, -13) {$3$};
		\node [style=label-white] (459) at (32, -13) {$\bar 1$};
		\node [style=label-white] (460) at (34, -13) {$6$};
		\node [style=label-white] (461) at (36, -13) {$2$};
		\node [style=label-white] (462) at (38, -13) {$5$};
		\node [style=label-white] (463) at (40, -13) {$1$};
		\node [style=label-white] (464) at (24, -15) {$0$};
		\node [style=label-white] (465) at (26, -15) {$4$};
		\node [style=label-white] (466) at (28, -15) {$\bar 1$};
		\node [style=label-white] (467) at (30, -15) {$3$};
		\node [style=label-white] (468) at (32, -15) {$6$};
		\node [style=label-white] (469) at (34, -15) {$1$};
		\node [style=label-white] (470) at (36, -15) {$2$};
		\node [style=label-white] (471) at (38, -15) {$5$};
		\node [style=none] (472) at (28, -14) {};
		\node [style=none] (473) at (30, -14) {};
		\node [style=none] (474) at (32, -14) {};
		\node [style=none] (475) at (34, -14) {};
		\node [style=none] (476) at (36, -14) {};
		\node [style=none] (477) at (37, -14) {};
		\node [style=none] (478) at (38, -14) {};
		\node [style=none] (479) at (40, -14) {};
		\node [style=none] (480) at (23, -17) {};
		\node [style=none] (481) at (25, -17) {};
		\node [style=none] (482) at (27, -17) {};
		\node [style=none] (483) at (29, -17) {};
		\node [style=none] (484) at (31, -17) {};
		\node [style=none] (485) at (33, -17) {};
		\node [style=none] (486) at (35, -17) {};
		\node [style=none] (487) at (37, -17) {};
		\node [style=none] (488) at (21, -17) {};
		\node [style=none] (489) at (21, -17) {};
		\node [style=label-white] (490) at (22, -17) {$0$};
		\node [style=label-white] (491) at (24, -17) {$\bar 1$};
		\node [style=label-white] (492) at (26, -17) {$4$};
		\node [style=label-white] (493) at (28, -17) {$1$};
		\node [style=label-white] (494) at (30, -17) {$3$};
		\node [style=label-white] (495) at (32, -17) {$6$};
		\node [style=label-white] (496) at (34, -17) {$2$};
		\node [style=label-white] (497) at (36, -17) {$5$};
		\node [style=none] (498) at (34, -16) {};
		\node [style=none] (499) at (28, -16) {};
		\node [style=none] (500) at (28, -16) {};
		\node [style=none] (501) at (24, -16) {};
		\node [style=none] (502) at (21, -19) {};
		\node [style=none] (503) at (23, -19) {};
		\node [style=none] (504) at (25, -19) {};
		\node [style=none] (505) at (27, -19) {};
		\node [style=none] (506) at (29, -19) {};
		\node [style=none] (507) at (31, -19) {};
		\node [style=none] (508) at (33, -19) {};
		\node [style=none] (509) at (35, -19) {};
		\node [style=none] (510) at (19, -19) {};
		\node [style=none] (511) at (19, -19) {};
		\node [style=label-white] (512) at (20, -19) {$\bar 1$};
		\node [style=label-white] (513) at (22, -19) {$0$};
		\node [style=label-white] (514) at (24, -19) {$1$};
		\node [style=label-white] (515) at (26, -19) {$4$};
		\node [style=label-white] (516) at (28, -19) {$2$};
		\node [style=label-white] (517) at (30, -19) {$3$};
		\node [style=label-white] (518) at (32, -19) {$6$};
		\node [style=label-white] (519) at (34, -19) {$5$};
		\node [style=none] (520) at (34, -18) {};
		\node [style=none] (521) at (28, -18) {};
		\node [style=none] (522) at (28, -18) {};
		\node [style=none] (523) at (24, -18) {};
		\node [style=none] (524) at (24, -18) {};
		\node [style=none] (525) at (20, -18) {};
		\node [style=none] (526) at (36, -17) {};
		\node [style=none] (527) at (34, -19) {};
		\node [style=none] (528) at (36, -18) {};
		\node [style=none] (529) at (34, -18) {};
		\node [style=none] (530) at (38, -15) {};
		\node [style=none] (531) at (36, -17) {};
		\node [style=none] (532) at (38, -16) {};
		\node [style=none] (533) at (36, -16) {};
		\node [style=none] (534) at (36, -15) {};
		\node [style=none] (535) at (34, -17) {};
		\node [style=none] (536) at (36, -16) {};
		\node [style=none] (537) at (34, -16) {};
		\node [style=none] (538) at (28, -13) {};
		\node [style=none] (539) at (26, -15) {};
		\node [style=none] (540) at (28, -14) {};
		\node [style=none] (541) at (26, -14) {};
		\node [style=none] (542) at (24, -15) {};
		\node [style=none] (543) at (22, -17) {};
		\node [style=none] (544) at (24, -16) {};
		\node [style=none] (545) at (22, -16) {};
		\node [style=none] (546) at (34, -13) {};
		\node [style=none] (547) at (32, -15) {};
		\node [style=none] (548) at (34, -14) {};
		\node [style=none] (549) at (32, -14) {};
		\node [style=none] (550) at (26, -13) {};
		\node [style=none] (551) at (24, -15) {};
		\node [style=none] (552) at (26, -14) {};
		\node [style=none] (553) at (24, -14) {};
	\end{pgfonlayer}
	\begin{pgfonlayer}{edgelayer}
		\draw [style=fade] (441.center) to (433.center);
		\draw [style=fade] (434.center) to (442.center);
		\draw [style=fade] (436.center) to (444.center);
		\draw [style=fade] (437.center) to (445.center);
		\draw [style=fade] (440.center) to (448.center);
		\draw [style=fade] (451.center) to (452.center);
		\draw [style=fade] (435.center) to (443.center);
		\draw [style=fade] (438.center) to (446.center);
		\draw [style=fade] (439.center) to (447.center);
		\draw [style=red, rounded corners=0.3cm] (466)
			 to (472.center)
			 to (474.center)
			 to (459);
		\draw [style=red, rounded corners=0.3cm] (469)
			 to (475.center)
			 to (479.center)
			 to (463);
		\draw [style=red] (458) to (467);
		\draw [style=red] (461) to (470);
		\draw [style=red] (462) to (471);
		\draw [style=fade] (452.center) to (489.center);
		\draw [style=fade] (441.center) to (480.center);
		\draw [style=fade] (442.center) to (481.center);
		\draw [style=fade] (443.center) to (482.center);
		\draw [style=fade] (444.center) to (483.center);
		\draw [style=fade] (445.center) to (484.center);
		\draw [style=fade] (446.center) to (485.center);
		\draw [style=fade] (447.center) to (486.center);
		\draw [style=fade] (448.center) to (487.center);
		\draw [style=red, rounded corners=0.3cm] (469)
			 to (498.center)
			 to (499.center)
			 to (493);
		\draw [style=red] (467) to (494);
		\draw [style=red] (468) to (495);
		\draw [style=fade] (451.center) to (440.center);
		\draw [style=fade] (489.center) to (487.center);
		\draw [style=fade] (452.center) to (448.center);
		\draw [style=red, rounded corners=0.3cm] (466)
			 to (500.center)
			 to (501.center)
			 to (491);
		\draw [style=red] (465) to (492);
		\draw [style=fade] (511.center) to (509.center);
		\draw [style=red, rounded corners=0.3cm] (496)
			 to (520.center)
			 to (521.center)
			 to (516);
		\draw [style=red] (494) to (517);
		\draw [style=red] (495) to (518);
		\draw [style=red, rounded corners=0.3cm] (493)
			 to (522.center)
			 to (523.center)
			 to (514);
		\draw [style=red, rounded corners=0.3cm] (491)
			 to (524.center)
			 to (525.center)
			 to (512);
		\draw [style=red] (490) to (513);
		\draw [style=red] (515) to (492);
		\draw [style=fade] (481.center) to (503.center);
		\draw [style=fade, in=225, out=45] (502.center) to (480.center);
		\draw [style=fade] (511.center) to (489.center);
		\draw [style=fade] (504.center) to (482.center);
		\draw [style=fade] (483.center) to (505.center);
		\draw [style=fade] (484.center) to (506.center);
		\draw [style=fade] (485.center) to (507.center);
		\draw [style=fade] (486.center) to (508.center);
		\draw [style=fade] (487.center) to (509.center);
		\draw [style=red, rounded corners=0.3cm] (526.center)
			 to (528.center)
			 to (529.center)
			 to (527.center);
		\draw [style=red, rounded corners=0.3cm] (530.center)
			 to (532.center)
			 to (533.center)
			 to (531.center);
		\draw [style=red, rounded corners=0.3cm] (534.center)
			 to (536.center)
			 to (537.center)
			 to (535.center);
		\draw [style=red, rounded corners=0.3cm] (538.center)
			 to (540.center)
			 to (541.center)
			 to (539.center);
		\draw [style=red, rounded corners=0.3cm] (542.center)
			 to (544.center)
			 to (545.center)
			 to (543.center);
		\draw [style=red, rounded corners=0.3cm] (546.center)
			 to (548.center)
			 to (549.center)
			 to (547.center);
		\draw [style=red, rounded corners=0.3cm] (550.center)
			 to (552.center)
			 to (553.center)
			 to (551.center);
	\end{pgfonlayer}
\end{tikzpicture}
\caption{The pipe-dream model corresponding to the tiling of \Cref{fig:3row-tiling}}
\label{fig:3row-pipedream}	
\end{figure}
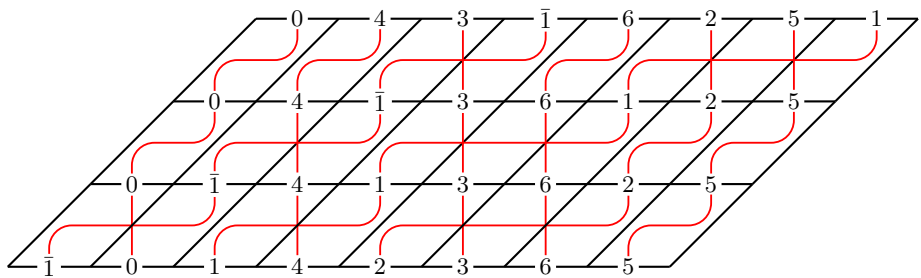

It turns out that the labeled path configurations of our vertex models are essentially a ``shifted'' version of the \emph{pipe-dreams} or RC-graphs of Bergeron-Billey \cite{bergeron1993rc} and Fomin-Kirrilov \cite{fomin1996yang}, where pipes are labeled by inverse permutations. For example, the tiling in \Cref{fig:3row-tiling} is in bijection with the piped-ream obtained by twisting the diagram then bending the straight edges into `bumps', see \Cref{fig:3row-pipedream}.

\section{Proof of \Cref{thm:heisenbergS,thm:stanley_hamiltonian}}\label{sec:proof}
We will now prove \Cref{thm:heisenbergS,thm:stanley_hamiltonian}. To streamline the proof, we will first introduce some diagrammatic presentation of the $H^{(k)}_i$ operators using a variation of the tiling model given in \Cref{fig:new-tiles}.

\subsection{Tiling Models for Bosonic Operators}
We first consider the $m>0$ case. As we have already seen in the proof of \Cref{thm:stanley_transfer_matrix}, the path configurations in \Cref{fig:new-tiles} account for the increasing words. Therefore we shall include the tiles in \Cref{fig:new-tiles} and introduce more.

\begin{definition}\label{def:Hplustiles}
	Define the $\nabla$-tiles to be the following path configurations.
	
	{\centering\begin{tikzpicture}[scale=0.7]
	\begin{pgfonlayer}{nodelayer}
		\node [style=none] (0) at (-4.5, 1.5) {};
		\node [style=none] (1) at (-2, 1.5) {};
		\node [style=none] (2) at (-2, -1) {};
		\node [style=none] (3) at (-4.5, -1) {};
		\node [style=none] (4) at (-1, 1.5) {};
		\node [style=none] (5) at (1.5, 1.5) {};
		\node [style=none] (6) at (1.5, -1) {};
		\node [style=none] (7) at (-1, -1) {};
		\node [style=none] (8) at (2.5, 1.5) {};
		\node [style=none] (9) at (5, 1.5) {};
		\node [style=none] (10) at (5, -1) {};
		\node [style=none] (11) at (2.5, -1) {};
		\node [style=none] (12) at (-8, 1.5) {};
		\node [style=none] (13) at (-5.5, 1.5) {};
		\node [style=none] (14) at (-5.5, -1) {};
		\node [style=none] (15) at (-8, -1) {};
		\node [style=label] (16) at (-9.5, -2.5) {Condition};
		\node [style=label] (17) at (-9.5, -4) {Weight};
		\node [style=none] (18) at (-3.25, 1.5) {};
		\node [style=none] (19) at (-2, 0.25) {};
		\node [style=none] (20) at (-3.25, -1) {};
		\node [style=none] (21) at (-3.25, 0.25) {};
		\node [style=none] (22) at (0.25, 0.25) {};
		\node [style=none] (23) at (0.25, 1.5) {};
		\node [style=none] (24) at (-1, 0.25) {};
		\node [style=none] (25) at (0.25, -1) {};
		\node [style=none] (26) at (2.5, 0.25) {};
		\node [style=none] (27) at (5, 0.25) {};
		\node [style=none] (28) at (3.75, -1) {};
		\node [style=none] (29) at (3.75, 1.5) {};
		\node [style=label-white] (30) at (-6.75, 1.5) {$a$};
		\node [style=label-white] (31) at (-6.75, -1) {$a$};
		\node [style=label-white] (32) at (-3.25, -1) {$b$};
		\node [style=label-white] (33) at (-3.25, 1.5) {$a$};
		\node [style=label-white] (34) at (0.25, 1.5) {$a$};
		\node [style=label-white] (35) at (0.25, -1) {$b$};
		\node [style=label-white] (36) at (3.75, 1.5) {$a$};
		\node [style=label-white] (37) at (3.75, 0.25) {$i$};
		\node [style=label-white] (38) at (3.75, -1) {$b$};
		\node [style=label] (39) at (-3.25, -2.5) {$a>b$};
		\node [style=label] (40) at (0.25, -2.5) {$a<b$};
		\node [style=label] (41) at (3.75, -2.5) {$i<\min(a,b)$};
		\node [style=label] (42) at (-3.25, -4) {$x$};
		\node [style=label] (43) at (0.25, -4) {$1$};
		\node [style=label] (44) at (-6.75, -4) {$1$};
		\node [style=label] (45) at (3.75, -4) {$x$};
		\node [style=label] (84) at (-9.5, 0) {$\nabla$-tiles};
		\node [style=none] (119) at (-3.25, 1.5) {};
		\node [style=none] (120) at (-2, 1.5) {};
		\node [style=none] (121) at (-2, -1) {};
		\node [style=none] (122) at (-3.25, -1) {};
		\node [style=none] (123) at (-1, 1.5) {};
		\node [style=none] (124) at (0.25, 1.5) {};
		\node [style=none] (125) at (0.25, -1) {};
		\node [style=none] (126) at (-1, -1) {};
		\node [style=none] (127) at (6, 1.5) {};
		\node [style=none] (128) at (8.5, 1.5) {};
		\node [style=none] (129) at (8.5, -1) {};
		\node [style=none] (130) at (6, -1) {};
		\node [style=none] (139) at (6, 1.5) {};
		\node [style=none] (140) at (8.5, 1.5) {};
		\node [style=none] (141) at (8.5, -1) {};
		\node [style=none] (142) at (6, -1) {};
		\node [style=none] (143) at (7.25, 0.25) {};
		\node [style=none] (144) at (7.25, 1.5) {};
		\node [style=none] (145) at (6, 0.25) {};
		\node [style=none] (146) at (7.25, -1) {};
		\node [style=label-white] (147) at (7.25, 1.5) {$a$};
		\node [style=label-white] (148) at (7.25, -1) {$b$};
		\node [style=none] (149) at (6, 1.5) {};
		\node [style=none] (151) at (7.25, -1) {};
		\node [style=none] (152) at (6, -1) {};
		\node [style=none] (153) at (7.25, 0.25) {};
		\node [style=none] (154) at (8.5, 0.25) {};
		\node [style=none] (155) at (6, 0.25) {};
		\node [style=none] (156) at (8.5, 0.25) {};
		\node [style=none] (157) at (9.5, 1.5) {};
		\node [style=none] (158) at (12, 1.5) {};
		\node [style=none] (159) at (12, -1) {};
		\node [style=none] (160) at (9.5, -1) {};
		\node [style=none] (162) at (10.75, 1.5) {};
		\node [style=none] (163) at (9.5, 0.25) {};
		\node [style=none] (164) at (10.75, -1) {};
		\node [style=label-white] (165) at (10.75, -1) {$a$};
		\node [style=label-white] (166) at (12, 0.25) {$b$};
		\node [style=none] (167) at (9.5, 1.5) {};
		\node [style=none] (170) at (9.5, -1) {};
		\node [style=none] (172) at (12, 0.25) {};
		\node [style=none] (173) at (9.5, 0.25) {};
		\node [style=none] (174) at (12, 0.25) {};
		\node [style=none] (175) at (2.5, 1.5) {};
		\node [style=none] (176) at (5, 1.5) {};
		\node [style=none] (177) at (5, -1) {};
		\node [style=none] (178) at (2.5, -1) {};
		\node [style=none] (187) at (9.5, 1.5) {};
		\node [style=none] (188) at (12, 1.5) {};
		\node [style=none] (189) at (12, -1) {};
		\node [style=none] (190) at (9.5, -1) {};
		\node [style=label-white] (191) at (7.25, 0.25) {$i$};
		\node [style=label-white] (192) at (8.5, 0.25) {$b$};
		\node [style=label-white] (193) at (6, 0.25) {$a$};
		\node [style=label-white] (194) at (10.75, 1.5) {$a$};
		\node [style=label-white] (195) at (9.5, 0.25) {$b$};
		\node [style=label-white] (196) at (11.25, 0.5) {$i$};
		\node [style=label] (197) at (7, -4) {$-x$};
		\node [style=label] (198) at (10.75, -4) {$-x^2$};
		\node [style=label] (199) at (7.25, -2.5) {$\max(a,b)<i$};
		\node [style=label] (200) at (10.75, -2.5) {$b<a<i$};
		\node [style=label] (201) at (-9.5, 2.5) {Type};
		\node [style=none] (202) at (-6.75, 2.5) {$0$};
		\node [style=none] (203) at (-3.25, 2.5) {$1$};
		\node [style=none] (204) at (0.25, 2.5) {$2$};
		\node [style=none] (205) at (3.75, 2.5) {$3$};
		\node [style=none] (206) at (7.25, 2.5) {$4$};
		\node [style=none] (207) at (10.75, 2.5) {$5$};
	\end{pgfonlayer}
	\begin{pgfonlayer}{edgelayer}
		\draw [style=fade] (1.center)
			 to (2.center)
			 to (3.center)
			 to (0.center)
			 to cycle;
		\draw [style=fade] (7.center)
			 to (4.center)
			 to (5.center)
			 to (6.center)
			 to cycle;
		\draw [style=fade] (10.center)
			 to (11.center)
			 to (8.center)
			 to (9.center)
			 to cycle;
		\draw [style=fade] (15.center)
			 to (12.center)
			 to (13.center)
			 to (14.center)
			 to cycle;
		\draw [style=red] (30) to (31);
		\draw [style=blue] (26.center) to (28.center);
		\draw [style=red, rounded corners=0.4cm] (24.center)
			 to (22.center)
			 to (23.center);
		\draw [style=red, rounded corners=0.4cm] (20.center)
			 to (21.center)
			 to (19.center);
		\draw [style=blue] (18.center) to (19.center);
		\draw [style=blue] (24.center) to (25.center);
		\draw [style=red] (26.center) to (27.center);
		\draw [style=blue] (29.center) to (27.center);
		\draw [style=shade] (119.center)
			 to (120.center)
			 to (121.center)
			 to (122.center)
			 to cycle;
		\draw [style=shade] (123.center)
			 to (124.center)
			 to (125.center)
			 to (126.center)
			 to cycle;
		\draw [style=shade] (127.center)
			 to (128.center)
			 to (129.center)
			 to (130.center)
			 to cycle;
		\draw [style=fade] (142.center)
			 to (139.center)
			 to (140.center)
			 to (141.center)
			 to cycle;
		\draw [style=red, rounded corners=0.4cm] (145.center)
			 to (143.center)
			 to (144.center);
		\draw [style=shade] (151.center) to (152.center);
		\draw [style=shade] (152.center) to (149.center);
		\draw [style=red, rounded corners=0.4cm] (151.center)
			 to (153.center)
			 to (154.center);
		\draw [style=blue, in=150, out=-30, looseness=1.75] (155.center) to (156.center);
		\draw [style=fade] (160.center)
			 to (157.center)
			 to (158.center)
			 to (159.center)
			 to cycle;
		\draw [style=shade] (170.center) to (167.center);
		\draw [style=blue, in=150, out=-30, looseness=1.75] (173.center) to (174.center);
		\draw [style=shade] (175.center)
			 to (176.center)
			 to (177.center)
			 to (178.center)
			 to cycle;
		\draw [style=shade] (187.center)
			 to (188.center)
			 to (189.center)
			 to (190.center)
			 to cycle;
		\draw [style=red] (163.center) to (172.center);
		\draw [style=red] (162.center) to (164.center);
	\end{pgfonlayer}
\end{tikzpicture}}
	
	Let $\mathcal{Z}^{\nabla}_{w/u}(x)$ be the weighted sum of all one-row tilings whose top labels are given by $w$ and bottom labels are given by $u$, such that there exists only one connected component of shaded tiles.
	We then define the transfer matrix $\tilde T_{+}(x)$ to be
	\[\tilde T_+(x)|w\rangle =\sum_{u}\mathcal{Z}^{\nabla}_{w/u}(x)|u\rangle\]
\end{definition}

\begin{lemma}
	\label{lem:transfer_equal_H_plus}
	The transfer matrix of $\nabla$-tiling model agrees with the $H_{n}$ operators, i.e.
	\[\tilde T_+(x)  = \sum_{n=1}^{\infty}x^n \alpha_n^\fff\]
\end{lemma}

\begin{proof}
	We shall construct a bijection between ribbons $v$ with $w=uv$ and one-row tilings of boundary $w,u$ whose weight equals to $(-1)^{\spin(v)}x^{\ell(v)}$. 
	Recall that $\tilde T_+$ only returns tilings with exactly one shaded region, therefore we observer immediately by definition of the $\nabla$-tiles that a shaded region consists of several red path segments that start with a $\rtile$ and end with a $\jtile$. For each such path-segment we define its reading word to be the word obtained by reading from left to right while each tile of type $1,3,4$ will contribute to a $s_i$ if it's in the $i$-th column and a type $5$ tile will contribute to a $s_{i+1}s_{i}$ if it's in the $i$-th column. We then obtain the reading word of an entire shaded component by going through all path-segments from right to left and concatenate the reading words of each of them. See \Cref{fig:H-plus-ribbon} for example.	
	\begin{figure}[h]
		\centerline{
	\begin{tikzpicture}[scale = 0.6]
	\begin{pgfonlayer}{nodelayer}
		\node [style=none] (208) at (-8, -6) {};
		\node [style=none] (209) at (-5.5, -6) {};
		\node [style=none] (210) at (-5.5, -8.5) {};
		\node [style=none] (211) at (-8, -8.5) {};
		\node [style=none] (212) at (-6.75, -6) {};
		\node [style=none] (213) at (-5.5, -7.25) {};
		\node [style=none] (214) at (-6.75, -8.5) {};
		\node [style=none] (215) at (-6.75, -7.25) {};
		\node [style=label-white] (217) at (-6.75, -6) {$3$};
		\node [style=none] (218) at (-6.75, -6) {};
		\node [style=none] (219) at (-5.5, -6) {};
		\node [style=none] (220) at (-5.5, -8.5) {};
		\node [style=none] (221) at (-6.75, -8.5) {};
		\node [style=none] (222) at (-5.5, -6) {};
		\node [style=none] (223) at (-3, -6) {};
		\node [style=none] (224) at (-3, -8.5) {};
		\node [style=none] (225) at (-5.5, -8.5) {};
		\node [style=none] (226) at (-5.5, -7.25) {};
		\node [style=none] (227) at (-3, -7.25) {};
		\node [style=none] (228) at (-4.25, -8.5) {};
		\node [style=none] (229) at (-4.25, -6) {};
		\node [style=label-white] (230) at (-4.25, -6) {$9$};
		\node [style=label-white] (232) at (-4.25, -8.5) {$3$};
		\node [style=none] (233) at (-5.5, -6) {};
		\node [style=none] (234) at (-3, -6) {};
		\node [style=none] (235) at (-3, -8.5) {};
		\node [style=none] (236) at (-5.5, -8.5) {};
		\node [style=none] (237) at (-3, -6) {};
		\node [style=none] (238) at (-0.5, -6) {};
		\node [style=none] (239) at (-0.5, -8.5) {};
		\node [style=none] (240) at (-3, -8.5) {};
		\node [style=none] (241) at (-1.75, -6) {};
		\node [style=none] (242) at (-3, -7.25) {};
		\node [style=none] (243) at (-1.75, -8.5) {};
		\node [style=label-white] (244) at (-1.75, -8.5) {$4$};
		\node [style=none] (246) at (-3, -6) {};
		\node [style=none] (247) at (-3, -8.5) {};
		\node [style=none] (248) at (-0.5, -7.25) {};
		\node [style=none] (249) at (-3, -7.25) {};
		\node [style=none] (250) at (-0.5, -7.25) {};
		\node [style=none] (251) at (-3, -6) {};
		\node [style=none] (252) at (-0.5, -6) {};
		\node [style=none] (253) at (-0.5, -8.5) {};
		\node [style=none] (254) at (-3, -8.5) {};
		\node [style=label-white] (255) at (-1.75, -6) {$4$};
		\node [style=none] (258) at (-0.5, -6) {};
		\node [style=none] (259) at (2, -6) {};
		\node [style=none] (260) at (2, -8.5) {};
		\node [style=none] (261) at (-0.5, -8.5) {};
		\node [style=none] (262) at (-0.5, -6) {};
		\node [style=none] (263) at (2, -6) {};
		\node [style=none] (264) at (2, -8.5) {};
		\node [style=none] (265) at (-0.5, -8.5) {};
		\node [style=none] (266) at (0.75, -7.25) {};
		\node [style=none] (267) at (0.75, -6) {};
		\node [style=none] (268) at (-0.5, -7.25) {};
		\node [style=none] (269) at (0.75, -8.5) {};
		\node [style=label-white] (270) at (0.75, -6) {$1$};
		\node [style=label-white] (271) at (0.75, -8.5) {$5$};
		\node [style=none] (272) at (-0.5, -6) {};
		\node [style=none] (273) at (0.75, -8.5) {};
		\node [style=none] (274) at (-0.5, -8.5) {};
		\node [style=none] (275) at (0.75, -7.25) {};
		\node [style=none] (276) at (2, -7.25) {};
		\node [style=none] (277) at (-0.5, -7.25) {};
		\node [style=none] (278) at (2, -7.25) {};
		\node [style=none] (282) at (2, -6) {};
		\node [style=none] (283) at (4.5, -6) {};
		\node [style=none] (284) at (4.5, -8.5) {};
		\node [style=none] (285) at (2, -8.5) {};
		\node [style=none] (286) at (3.25, -6) {};
		\node [style=none] (287) at (2, -7.25) {};
		\node [style=none] (288) at (3.25, -8.5) {};
		\node [style=label-white] (289) at (3.25, -8.5) {$7$};
		\node [style=none] (291) at (2, -6) {};
		\node [style=none] (292) at (2, -8.5) {};
		\node [style=none] (293) at (4.5, -7.25) {};
		\node [style=none] (294) at (2, -7.25) {};
		\node [style=none] (295) at (4.5, -7.25) {};
		\node [style=none] (296) at (2, -6) {};
		\node [style=none] (297) at (4.5, -6) {};
		\node [style=none] (298) at (4.5, -8.5) {};
		\node [style=none] (299) at (2, -8.5) {};
		\node [style=label-white] (300) at (3.25, -6) {$7$};
		\node [style=none] (303) at (4.5, -6) {};
		\node [style=none] (304) at (7, -6) {};
		\node [style=none] (305) at (7, -8.5) {};
		\node [style=none] (306) at (4.5, -8.5) {};
		\node [style=none] (307) at (4.5, -6) {};
		\node [style=none] (308) at (7, -6) {};
		\node [style=none] (309) at (7, -8.5) {};
		\node [style=none] (310) at (4.5, -8.5) {};
		\node [style=none] (311) at (5.75, -7.25) {};
		\node [style=none] (312) at (5.75, -6) {};
		\node [style=none] (313) at (4.5, -7.25) {};
		\node [style=none] (314) at (5.75, -8.5) {};
		\node [style=label-white] (315) at (5.75, -6) {$5$};
		\node [style=label-white] (316) at (5.75, -8.5) {$2$};
		\node [style=none] (317) at (4.5, -6) {};
		\node [style=none] (318) at (5.75, -8.5) {};
		\node [style=none] (319) at (4.5, -8.5) {};
		\node [style=none] (320) at (5.75, -7.25) {};
		\node [style=none] (321) at (7, -7.25) {};
		\node [style=none] (322) at (4.5, -7.25) {};
		\node [style=none] (323) at (7, -7.25) {};
		\node [style=none] (331) at (7, -7.25) {};
		\node [style=none] (332) at (7, -7.25) {};
		\node [style=none] (342) at (7, -6) {};
		\node [style=none] (343) at (9.5, -6) {};
		\node [style=none] (344) at (9.5, -8.5) {};
		\node [style=none] (345) at (7, -8.5) {};
		\node [style=none] (346) at (8.25, -7.25) {};
		\node [style=none] (347) at (8.25, -6) {};
		\node [style=none] (348) at (7, -7.25) {};
		\node [style=none] (349) at (8.25, -8.5) {};
		\node [style=label-white] (350) at (8.25, -6) {$2$};
		\node [style=label-white] (351) at (8.25, -8.5) {$9$};
		\node [style=none] (352) at (7, -6) {};
		\node [style=none] (353) at (8.25, -6) {};
		\node [style=none] (354) at (8.25, -8.5) {};
		\node [style=none] (355) at (7, -8.5) {};
		\node [style=none] (357) at (-8, -6) {};
		\node [style=none] (358) at (-8, -8.5) {};
		\node [style=none] (362) at (9.5, -6) {};
		\node [style=none] (365) at (9.5, -8.5) {};
		\node [style=label-white] (368) at (-6.75, -8.5) {$1$};
		\node [style=none] (369) at (-9.25, -7.25) {$\cdots$};
		\node [style=none] (370) at (10.75, -7.25) {$\cdots$};
	\end{pgfonlayer}
	\begin{pgfonlayer}{edgelayer}
		\draw [style=fade] (209.center)
			 to (210.center)
			 to (211.center)
			 to (208.center)
			 to cycle;
		\draw [style=red, rounded corners=0.4cm] (214.center)
			 to (215.center)
			 to (213.center);
		\draw [style=blue] (212.center) to (213.center);
		\draw [style=shade] (218.center)
			 to (219.center)
			 to (220.center)
			 to (221.center)
			 to cycle;
		\draw [style=fade] (224.center)
			 to (225.center)
			 to (222.center)
			 to (223.center)
			 to cycle;
		\draw [style=blue] (226.center) to (228.center);
		\draw [style=red] (226.center) to (227.center);
		\draw [style=blue] (229.center) to (227.center);
		\draw [style=shade] (233.center)
			 to (234.center)
			 to (235.center)
			 to (236.center)
			 to cycle;
		\draw [style=fade] (240.center)
			 to (237.center)
			 to (238.center)
			 to (239.center)
			 to cycle;
		\draw [style=shade] (247.center) to (246.center);
		\draw [style=blue, in=150, out=-30, looseness=1.75] (249.center) to (250.center);
		\draw [style=shade] (251.center)
			 to (252.center)
			 to (253.center)
			 to (254.center)
			 to cycle;
		\draw [style=red] (242.center) to (248.center);
		\draw [style=red, in=90, out=-90] (241.center) to (243.center);
		\draw [style=shade] (258.center)
			 to (259.center)
			 to (260.center)
			 to (261.center)
			 to cycle;
		\draw [style=fade] (265.center)
			 to (262.center)
			 to (263.center)
			 to (264.center)
			 to cycle;
		\draw [style=red, rounded corners=0.4cm] (268.center)
			 to (266.center)
			 to (267.center);
		\draw [style=shade] (273.center) to (274.center);
		\draw [style=shade] (274.center) to (272.center);
		\draw [style=red, rounded corners=0.4cm] (273.center)
			 to (275.center)
			 to (276.center);
		\draw [style=blue, in=150, out=-30, looseness=1.75] (277.center) to (278.center);
		\draw [style=fade] (285.center)
			 to (282.center)
			 to (283.center)
			 to (284.center)
			 to cycle;
		\draw [style=shade] (292.center) to (291.center);
		\draw [style=blue, in=150, out=-30, looseness=1.75] (294.center) to (295.center);
		\draw [style=shade] (296.center)
			 to (297.center)
			 to (298.center)
			 to (299.center)
			 to cycle;
		\draw [style=red] (287.center) to (293.center);
		\draw [style=red] (286.center) to (288.center);
		\draw [style=shade] (303.center)
			 to (304.center)
			 to (305.center)
			 to (306.center)
			 to cycle;
		\draw [style=fade] (310.center)
			 to (307.center)
			 to (308.center)
			 to (309.center)
			 to cycle;
		\draw [style=red, rounded corners=0.4cm] (313.center)
			 to (311.center)
			 to (312.center);
		\draw [style=shade] (318.center) to (319.center);
		\draw [style=shade] (319.center) to (317.center);
		\draw [style=red, rounded corners=0.4cm] (318.center)
			 to (320.center)
			 to (321.center);
		\draw [style=blue, in=150, out=-30, looseness=1.75] (322.center) to (323.center);
		\draw [style=fade] (345.center)
			 to (342.center)
			 to (343.center)
			 to (344.center)
			 to cycle;
		\draw [style=red, rounded corners=0.4cm] (348.center)
			 to (346.center)
			 to (347.center);
		\draw [style=blue] (348.center) to (349.center);
		\draw [style=shade] (352.center)
			 to (353.center)
			 to (354.center)
			 to (355.center)
			 to cycle;
		\draw [style=fade] (357.center) to (358.center);
		\draw [style=fade] (365.center) to (362.center);
	\end{pgfonlayer}
\end{tikzpicture}
	}
	\caption{Example of a single shaded component of an $\nabla$-tiling.  In this example, the boundaries are $u=[\cdots  1\ 3\ 4\ 5\ 7\ 2\ 9 \cdots]$, $w = [\cdots 3\ 9\ 4\ 1\ 7\ 5\ 2\cdots ]$ and the reading word is $v=s_6s_4s_5s_4s_1s_2s_3s_2$. We have that $w = uv$ with $\ell(w) = \ell(u)+\ell(v)$.}
	\label{fig:H-plus-ribbon}
	\end{figure}
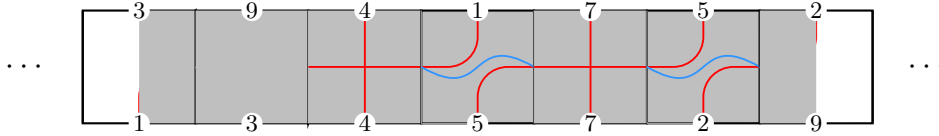

 We claim that the reading word $v$ obtained in this way is a ribbon and satisfy that $w=uv$ with $\ell(w)=\ell(u)+\ell(v)$. It is clear from the construction that the reading word obtained in this way necessarily avoids $2132$ and $2312$. To see that $\ell(w)=\ell(u)+\ell(v)$, we simply observe that by the conditions in \Cref{def:Hplustiles}, every time we apply a simply reflection according to the tiles we must increase the length. Furthermore, we need to show that the weight of such tiling is equal to $(-1)^{\spin(v)}x^{\ell(v)}$. Since every path-segment is increasing and hence contributes $0$ to the spin, we know that the $\spin(v)$ is exactly counted by the number of type $4_+$ and $5_+$ tiles. And it's clear that the exponent equals to the length of $v$ by construction.
	
We will next show that other direction of this bijection. Given $w=uv$ with $\ell(w)=\ell(u)+\ell(v)$ with $v$ being a ribbon. We will first find the reduced word of $v$ written maximally in the reverse lexicographical order. This is exactly the reading word obtained from a $H_+$-tiling with one shaded component, therefore we can construct such a tiling by reverse-engineering the exact same rules. That fact that this is weight preserving follows from the same argument as before.\end{proof} 

From here, we are ready to prove \Cref{thm:stanley_hamiltonian}, which we now restate an equivalent version.

\begin{theorem}\label{thm:stanley_hamiltonian_one_row}
	Let $T$ be the transfer matrix defined in \Cref{sec:2dfermion}. Then
	\begin{equation}\label{eq:one-row-exp}\exp\left(\sum_{i>o} {x^i\over i}\alpha_i^\fff\right ) =T(x)\end{equation}
\end{theorem}
\begin{proof}

	Writing the left hand side can be written as a triple product, and consider $$ \sum_{n\geq 0} {1\over n!}\left(\sum_{i>o} {x^i\over i}\alpha_i^\fff\right )|w\rangle$$ for some $w\in S_{\zz}$. Note that $|u\rangle$ can occur in the expansion (possibly with 0 coefficient), only if $u<w$ in the weak Bruhat order. In fact, every term in the expansion is a multi-row $\nabla$-tiling. We first note that the coefficient of $|u\rangle $ is zero unless $u^{-1}w$ has an increasing reduced word. When $u^{-1}w$ doesn't have an increasing reduced word, then we have $u^{-1}w=\cdots s_{m+1}s_m\cdots$ for some $m$. Then the weighted sum of $\nabla$-titles (multiplied by the coefficient from the infinite series) can be grouped into two different sets based on whether the $s_{m+1}s_m$ show up as a type 4 $\nabla$-tile or two tiles in different rows. These two sets must have opposite weight, therefore they sum to zero.
	
	Now consider the case when $u^{-1}w$ have an increasing reduced word. Suppose that the support of $u^{-1}w$ has $l$ connected component of length $n_1,\cdots,n_l$. In this case the appearances of $|u\rangle$ in \Cref{eq:one-row-exp} are in bijection with permuted partitions of each $n_i$. Therefore the coefficient of $|u\rangle$ is $$\sum_{
	\lambda^{(1)}\vdash n_1,\cdots,\lambda^{(l)}\vdash n_l} {1\over (n_1+\cdots+n_l)! \prod \lambda^{(1)}_i \cdots \prod \lambda^{(l)}_i},$$
	where each $\lambda^{(i)}$ is a partition of $n_i$, and $\prod \lambda_j^{(i)}$ is the product of all parts of $\lambda^{(i)}$.
	This alwas equals to $1$ for any numbers $n_1,\cdots,n_l$ (c.f. \cite[I.2.14']{macdonald1998symmetric}). This completes the proof.
\end{proof}

We will next introduce $\Delta_{k}$-tiles, which will characterize $\alpha_{m,k}^\fff$ operators for $m<0$.

\begin{definition}\label{def:Hminustiles}
Define $\Delta_k$-tiles as follows.
\vspace{1em}

\centerline{
\begin{tikzpicture}[scale=0.7]
	\begin{pgfonlayer}{nodelayer}
		\node [style=none] (0) at (-4.75, -1) {};
		\node [style=none] (1) at (-2.25, -1) {};
		\node [style=none] (2) at (-2.25, 1.5) {};
		\node [style=none] (3) at (-4.75, 1.5) {};
		\node [style=none] (4) at (-1.5, -1) {};
		\node [style=none] (5) at (1, -1) {};
		\node [style=none] (6) at (1, 1.5) {};
		\node [style=none] (7) at (-1.5, 1.5) {};
		\node [style=none] (8) at (1.75, -1) {};
		\node [style=none] (9) at (4.25, -1) {};
		\node [style=none] (10) at (4.25, 1.5) {};
		\node [style=none] (11) at (1.75, 1.5) {};
		\node [style=none] (12) at (-8, -1) {};
		\node [style=none] (13) at (-5.5, -1) {};
		\node [style=none] (14) at (-5.5, 1.5) {};
		\node [style=none] (15) at (-8, 1.5) {};
		\node [style=label] (16) at (-9.5, -2.5) {Condition};
		\node [style=label] (17) at (-9.5, -4) {Weight};
		\node [style=none] (18) at (-3.5, -1) {};
		\node [style=none] (19) at (-2.25, 0.25) {};
		\node [style=none] (20) at (-3.5, 1.5) {};
		\node [style=none] (21) at (-3.5, 0.25) {};
		\node [style=none] (22) at (-0.25, 0.25) {};
		\node [style=none] (23) at (-0.25, -1) {};
		\node [style=none] (24) at (-1.5, 0.25) {};
		\node [style=none] (25) at (-0.25, 1.5) {};
		\node [style=none] (26) at (1.75, 0.25) {};
		\node [style=none] (27) at (4.25, 0.25) {};
		\node [style=none] (28) at (3, 1.5) {};
		\node [style=none] (29) at (3, -1) {};
		\node [style=label-white] (30) at (-6.75, -1) {$a$};
		\node [style=label-white] (31) at (-6.75, 1.5) {$a$};
		\node [style=label-white] (32) at (-3.5, 1.5) {$a$};
		\node [style=label-white] (33) at (-3.5, -1) {$b$};
		\node [style=label-white] (34) at (-0.25, -1) {$a$};
		\node [style=label-white] (35) at (-0.25, 1.5) {$b$};
		\node [style=label-white] (36) at (3, -1) {$b$};
		\node [style=label-white] (37) at (3, 0.25) {$a$};
		\node [style=label-white] (38) at (3, 1.5) {$c$};
		\node [style=label, scale=0.8] (39) at (-3.5, -2.5) {$a\leq k<b$};
		\node [style=label, scale=0.8] (40) at (-0.25, -2.5) {$a\leq k<b$};
		\node [style=label, scale=0.8] (41) at (3, -2.5) {$a\leq k<b,c$};
		\node [style=label] (42) at (-3.5, -4) {$1$};
		\node [style=label] (43) at (-0.25, -4) {$x$};
		\node [style=label] (44) at (-6.75, -4) {$1$};
		\node [style=label] (45) at (3, -4) {$x$};
		\node [style=label] (84) at (-9.5, 0) {$\Delta_k$-tiles};
		\node [style=none] (119) at (-3.5, -1) {};
		\node [style=none] (120) at (-2.25, -1) {};
		\node [style=none] (121) at (-2.25, 1.5) {};
		\node [style=none] (122) at (-3.5, 1.5) {};
		\node [style=none] (123) at (-1.5, -1) {};
		\node [style=none] (124) at (-0.25, -1) {};
		\node [style=none] (125) at (-0.25, 1.5) {};
		\node [style=none] (126) at (-1.5, 1.5) {};
		\node [style=none] (127) at (5, -1) {};
		\node [style=none] (128) at (7.5, -1) {};
		\node [style=none] (129) at (7.5, 1.5) {};
		\node [style=none] (130) at (5, 1.5) {};
		\node [style=none] (139) at (5, -1) {};
		\node [style=none] (140) at (7.5, -1) {};
		\node [style=none] (141) at (7.5, 1.5) {};
		\node [style=none] (142) at (5, 1.5) {};
		\node [style=none] (143) at (6.25, 0.25) {};
		\node [style=none] (144) at (6.25, -1) {};
		\node [style=none] (145) at (5, 0.25) {};
		\node [style=none] (146) at (6.25, 1.5) {};
		\node [style=label-white] (147) at (6.25, -1) {$i$};
		\node [style=label-white] (148) at (6.25, 1.5) {$a$};
		\node [style=none] (149) at (5, -1) {};
		\node [style=none] (151) at (6.25, 1.5) {};
		\node [style=none] (152) at (5, 1.5) {};
		\node [style=none] (153) at (6.25, 0.25) {};
		\node [style=none] (154) at (7.5, 0.25) {};
		\node [style=none] (155) at (5, 0.25) {};
		\node [style=none] (156) at (7.5, 0.25) {};
		\node [style=none] (157) at (8.25, -1) {};
		\node [style=none] (158) at (10.75, -1) {};
		\node [style=none] (159) at (10.75, 1.5) {};
		\node [style=none] (160) at (8.25, 1.5) {};
		\node [style=none] (162) at (9.5, -1) {};
		\node [style=none] (163) at (8.25, 0.25) {};
		\node [style=none] (164) at (9.5, 1.5) {};
		\node [style=label-white] (165) at (9.5, 1.5) {$c$};
		\node [style=label-white] (166) at (10.75, 0.25) {$a$};
		\node [style=none] (167) at (8.25, -1) {};
		\node [style=none] (170) at (8.25, 1.5) {};
		\node [style=none] (172) at (10.75, 0.25) {};
		\node [style=none] (173) at (8.25, 0.25) {};
		\node [style=none] (174) at (10.75, 0.25) {};
		\node [style=none] (175) at (1.75, -1) {};
		\node [style=none] (176) at (4.25, -1) {};
		\node [style=none] (177) at (4.25, 1.5) {};
		\node [style=none] (178) at (1.75, 1.5) {};
		\node [style=none] (187) at (8.25, -1) {};
		\node [style=none] (188) at (10.75, -1) {};
		\node [style=none] (189) at (10.75, 1.5) {};
		\node [style=none] (190) at (8.25, 1.5) {};
		\node [style=label-white] (191) at (6.25, 0.25) {$b$};
		\node [style=label-white] (192) at (7.5, 0.25) {$a$};
		\node [style=label-white] (193) at (5, 0.25) {$i$};
		\node [style=label-white] (194) at (9.5, -1) {$c$};
		\node [style=label-white] (195) at (8.25, 0.25) {$a$};
		\node [style=label-white] (196) at (10, 0) {$b$};
		\node [style=label] (197) at (6, -4) {$-x$};
		\node [style=label] (198) at (9.5, -4) {$1$};
		\node [style=label, scale=0.8] (199) at (6.25, -2.5) {$a,i\leq k<b$};
		\node [style=label, scale=0.8] (200) at (9.5, -2.5) {$a\leq k<b<c$};
		\node [style=label] (201) at (-9.5, 2.5) {Type};
		\node [style=none] (202) at (-6.75, 2.5) {$0$};
		\node [style=none] (203) at (-3.5, 2.5) {$1$};
		\node [style=none] (204) at (-0.25, 2.5) {$2$};
		\node [style=none] (205) at (3, 2.5) {$3$};
		\node [style=none] (206) at (6.25, 2.5) {$4$};
		\node [style=none] (207) at (9.5, 2.5) {$5$};
		\node [style=label, scale=0.8] (208) at (-6.75, -2.5) {$a\leq k$};
		\node [style=none] (209) at (11.5, -1) {};
		\node [style=none] (210) at (14, -1) {};
		\node [style=none] (211) at (14, 1.5) {};
		\node [style=none] (212) at (11.5, 1.5) {};
		\node [style=none] (213) at (12.75, -1) {};
		\node [style=none] (214) at (11.5, 0.25) {};
		\node [style=none] (215) at (12.75, 1.5) {};
		\node [style=label-white] (216) at (12.75, 1.5) {$c$};
		\node [style=label-white] (217) at (14, 0.25) {$a$};
		\node [style=none] (218) at (11.5, -1) {};
		\node [style=none] (219) at (11.5, 1.5) {};
		\node [style=none] (220) at (14, 0.25) {};
		\node [style=none] (221) at (11.5, 0.25) {};
		\node [style=none] (222) at (14, 0.25) {};
		\node [style=none] (223) at (11.5, -1) {};
		\node [style=none] (224) at (14, -1) {};
		\node [style=none] (225) at (14, 1.5) {};
		\node [style=none] (226) at (11.5, 1.5) {};
		\node [style=label-white] (227) at (12.75, -1) {$c$};
		\node [style=label-white] (228) at (11.5, 0.25) {$a$};
		\node [style=label-white] (229) at (13.25, 0) {$b$};
		\node [style=label] (230) at (12.75, -4) {$1$};
		\node [style=label, scale=0.8] (231) at (12.75, -2.5) {$c<a\leq k<b$};
		\node [style=none] (232) at (12.75, 2.5) {$6$};
	\end{pgfonlayer}
	\begin{pgfonlayer}{edgelayer}
		\draw [style=fade] (1.center)
			 to (2.center)
			 to (3.center)
			 to (0.center)
			 to cycle;
		\draw [style=fade] (7.center)
			 to (4.center)
			 to (5.center)
			 to (6.center)
			 to cycle;
		\draw [style=fade] (10.center)
			 to (11.center)
			 to (8.center)
			 to (9.center)
			 to cycle;
		\draw [style=fade] (15.center)
			 to (12.center)
			 to (13.center)
			 to (14.center)
			 to cycle;
		\draw [style=red] (30) to (31);
		\draw [style=blue] (26.center) to (28.center);
		\draw [style=red, rounded corners=0.4cm] (24.center)
			 to (22.center)
			 to (23.center);
		\draw [style=red, rounded corners=0.4cm] (20.center)
			 to (21.center)
			 to (19.center);
		\draw [style=blue] (18.center) to (19.center);
		\draw [style=blue] (24.center) to (25.center);
		\draw [style=red] (26.center) to (27.center);
		\draw [style=blue] (29.center) to (27.center);
		\draw [style=shade] (119.center)
			 to (120.center)
			 to (121.center)
			 to (122.center)
			 to cycle;
		\draw [style=shade] (123.center)
			 to (124.center)
			 to (125.center)
			 to (126.center)
			 to cycle;
		\draw [style=shade] (127.center)
			 to (128.center)
			 to (129.center)
			 to (130.center)
			 to cycle;
		\draw [style=fade] (142.center)
			 to (139.center)
			 to (140.center)
			 to (141.center)
			 to cycle;
		\draw [style=red, rounded corners=0.4cm] (145.center)
			 to (143.center)
			 to (144.center);
		\draw [style=shade] (151.center) to (152.center);
		\draw [style=shade] (152.center) to (149.center);
		\draw [style=red, rounded corners=0.4cm] (151.center)
			 to (153.center)
			 to (154.center);
		\draw [style=blue, in=-150, out=30, looseness=1.75] (155.center) to (156.center);
		\draw [style=fade] (160.center)
			 to (157.center)
			 to (158.center)
			 to (159.center)
			 to cycle;
		\draw [style=shade] (170.center) to (167.center);
		\draw [style=blue, in=-150, out=30, looseness=1.75] (173.center) to (174.center);
		\draw [style=shade] (175.center)
			 to (176.center)
			 to (177.center)
			 to (178.center)
			 to cycle;
		\draw [style=shade] (187.center)
			 to (188.center)
			 to (189.center)
			 to (190.center)
			 to cycle;
		\draw [style=red] (163.center) to (172.center);
		\draw [style=blue] (162.center) to (164.center);
		\draw [style=fade] (212.center)
			 to (209.center)
			 to (210.center)
			 to (211.center)
			 to cycle;
		\draw [style=shade] (219.center) to (218.center);
		\draw [style=blue, in=-150, out=30, looseness=1.75] (221.center) to (222.center);
		\draw [style=shade] (223.center)
			 to (224.center)
			 to (225.center)
			 to (226.center)
			 to cycle;
		\draw [style=red] (214.center) to (220.center);
		\draw [style=red] (213.center) to (215.center);
	\end{pgfonlayer}
\end{tikzpicture}
}

Different from previous case, here we consider a tiling model with arbitrarily many rows (called \emph{layers} for convenience) of the same spectral parameter.  Let $\mathcal{Z}^{\Delta_k}_{w/u}(x)$ be the weighted sum of all tilings of arbitrarily many layers with top and bottom boundary given by $w$ and $u$, such that there exists only one connected component of shaded tiles and that no two weighted tiles (types 2,3,4) can appear in the same column\footnote{We can think of this as an one-row model consisting of infinitesimal rows.}.

	We then define the transfer matrix $\tilde T_{-}^{(k)}(x)$ to be
	\[\tilde T_-^{(k)}(x)|w\rangle =\sum_{u}\mathcal{Z}^{\Delta_k}_{w/u}(x)|u\rangle\]
\end{definition}
\begin{figure}[h]
		\begin{tikzpicture}[scale = 0.6]
	\begin{pgfonlayer}{nodelayer}
		\node [style=none] (208) at (-8, -8.5) {};
		\node [style=none] (209) at (-5.5, -8.5) {};
		\node [style=none] (210) at (-5.5, -6) {};
		\node [style=none] (211) at (-8, -6) {};
		\node [style=none] (212) at (-6.75, -8.5) {};
		\node [style=none] (213) at (-5.5, -7.25) {};
		\node [style=none] (214) at (-6.75, -6) {};
		\node [style=none] (215) at (-6.75, -7.25) {};
		\node [style=label-white] (217) at (-6.75, -8.5) {$5$};
		\node [style=none] (218) at (-6.75, -8.5) {};
		\node [style=none] (219) at (-5.5, -8.5) {};
		\node [style=none] (220) at (-5.5, -6) {};
		\node [style=none] (221) at (-6.75, -6) {};
		\node [style=none] (222) at (-5.5, -8.5) {};
		\node [style=none] (223) at (-3, -8.5) {};
		\node [style=none] (224) at (-3, -6) {};
		\node [style=none] (225) at (-5.5, -6) {};
		\node [style=none] (226) at (-5.5, -7.25) {};
		\node [style=none] (227) at (-3, -7.25) {};
		\node [style=none] (228) at (-4.25, -6) {};
		\node [style=none] (229) at (-4.25, -8.5) {};
		\node [style=label-white] (230) at (-4.25, -8.5) {$4$};
		\node [style=label-white] (232) at (-4.25, -6) {$5$};
		\node [style=none] (233) at (-5.5, -8.5) {};
		\node [style=none] (234) at (-3, -8.5) {};
		\node [style=none] (235) at (-3, -6) {};
		\node [style=none] (236) at (-5.5, -6) {};
		\node [style=none] (237) at (-3, -8.5) {};
		\node [style=none] (238) at (-0.5, -8.5) {};
		\node [style=none] (239) at (-0.5, -6) {};
		\node [style=none] (240) at (-3, -6) {};
		\node [style=none] (241) at (-1.75, -8.5) {};
		\node [style=none] (242) at (-3, -7.25) {};
		\node [style=none] (243) at (-1.75, -6) {};
		\node [style=label-white] (244) at (-1.75, -6) {$7$};
		\node [style=none] (246) at (-3, -8.5) {};
		\node [style=none] (247) at (-3, -6) {};
		\node [style=none] (248) at (-0.5, -7.25) {};
		\node [style=none] (249) at (-3, -7.25) {};
		\node [style=none] (250) at (-0.5, -7.25) {};
		\node [style=none] (251) at (-3, -8.5) {};
		\node [style=none] (252) at (-0.5, -8.5) {};
		\node [style=none] (253) at (-0.5, -6) {};
		\node [style=none] (254) at (-3, -6) {};
		\node [style=label-white] (255) at (-1.75, -8.5) {$7$};
		\node [style=none] (258) at (-0.5, -8.5) {};
		\node [style=none] (259) at (2, -8.5) {};
		\node [style=none] (260) at (2, -6) {};
		\node [style=none] (261) at (-0.5, -6) {};
		\node [style=none] (262) at (-0.5, -8.5) {};
		\node [style=none] (263) at (2, -8.5) {};
		\node [style=none] (264) at (2, -6) {};
		\node [style=none] (265) at (-0.5, -6) {};
		\node [style=none] (266) at (0.75, -7.25) {};
		\node [style=none] (267) at (0.75, -8.5) {};
		\node [style=none] (268) at (-0.5, -7.25) {};
		\node [style=none] (269) at (0.75, -6) {};
		\node [style=label-white] (270) at (0.75, -8.5) {$1$};
		\node [style=label-white] (271) at (0.75, -6) {$3$};
		\node [style=none] (272) at (-0.5, -8.5) {};
		\node [style=none] (273) at (0.75, -6) {};
		\node [style=none] (274) at (-0.5, -6) {};
		\node [style=none] (275) at (0.75, -7.25) {};
		\node [style=none] (276) at (2, -7.25) {};
		\node [style=none] (277) at (-0.5, -7.25) {};
		\node [style=none] (278) at (2, -7.25) {};
		\node [style=none] (282) at (2, -8.5) {};
		\node [style=none] (283) at (4.5, -8.5) {};
		\node [style=none] (284) at (4.5, -6) {};
		\node [style=none] (285) at (2, -6) {};
		\node [style=none] (286) at (3.25, -8.5) {};
		\node [style=none] (287) at (2, -7.25) {};
		\node [style=none] (288) at (3.25, -6) {};
		\node [style=label-white] (289) at (3.25, -6) {$6$};
		\node [style=none] (291) at (2, -8.5) {};
		\node [style=none] (292) at (2, -6) {};
		\node [style=none] (293) at (4.5, -7.25) {};
		\node [style=none] (294) at (2, -7.25) {};
		\node [style=none] (295) at (4.5, -7.25) {};
		\node [style=none] (296) at (2, -8.5) {};
		\node [style=none] (297) at (4.5, -8.5) {};
		\node [style=none] (298) at (4.5, -6) {};
		\node [style=none] (299) at (2, -6) {};
		\node [style=label-white] (300) at (3.25, -8.5) {$6$};
		\node [style=none] (303) at (4.5, -8.5) {};
		\node [style=none] (304) at (7, -8.5) {};
		\node [style=none] (305) at (7, -6) {};
		\node [style=none] (306) at (4.5, -6) {};
		\node [style=none] (307) at (4.5, -8.5) {};
		\node [style=none] (308) at (7, -8.5) {};
		\node [style=none] (309) at (7, -6) {};
		\node [style=none] (310) at (4.5, -6) {};
		\node [style=none] (311) at (5.75, -7.25) {};
		\node [style=none] (312) at (5.75, -8.5) {};
		\node [style=none] (313) at (4.5, -7.25) {};
		\node [style=none] (314) at (5.75, -6) {};
		\node [style=label-white] (315) at (5.75, -8.5) {$3$};
		\node [style=label-white] (316) at (5.75, -6) {$2$};
		\node [style=none] (317) at (4.5, -8.5) {};
		\node [style=none] (318) at (5.75, -6) {};
		\node [style=none] (319) at (4.5, -6) {};
		\node [style=none] (320) at (5.75, -7.25) {};
		\node [style=none] (321) at (7, -7.25) {};
		\node [style=none] (322) at (4.5, -7.25) {};
		\node [style=none] (323) at (7, -7.25) {};
		\node [style=none] (331) at (7, -7.25) {};
		\node [style=none] (332) at (7, -7.25) {};
		\node [style=none] (342) at (7, -8.5) {};
		\node [style=none] (343) at (9.5, -8.5) {};
		\node [style=none] (344) at (9.5, -6) {};
		\node [style=none] (345) at (7, -6) {};
		\node [style=none] (346) at (8.25, -7.25) {};
		\node [style=none] (347) at (8.25, -8.5) {};
		\node [style=none] (348) at (7, -7.25) {};
		\node [style=none] (349) at (8.25, -6) {};
		\node [style=label-white] (350) at (8.25, -8.5) {$2$};
		\node [style=label-white] (351) at (8.25, -6) {$4$};
		\node [style=none] (352) at (7, -8.5) {};
		\node [style=none] (353) at (8.25, -8.5) {};
		\node [style=none] (354) at (8.25, -6) {};
		\node [style=none] (355) at (7, -6) {};
		\node [style=none] (357) at (-8, -8.5) {};
		\node [style=none] (358) at (-8, -6) {};
		\node [style=none] (362) at (9.5, -8.5) {};
		\node [style=none] (365) at (9.5, -6) {};
		\node [style=label-white] (368) at (-6.75, -6) {$1$};
		\node [style=none] (369) at (-9.25, -7.25) {$\cdots$};
		\node [style=none] (370) at (10.75, -7.25) {$\cdots$};
	\end{pgfonlayer}
	\begin{pgfonlayer}{edgelayer}
		\draw [style=fade] (209.center)
			 to (210.center)
			 to (211.center)
			 to (208.center)
			 to cycle;
		\draw [style=red, rounded corners=0.4cm] (214.center)
			 to (215.center)
			 to (213.center);
		\draw [style=blue] (212.center) to (213.center);
		\draw [style=shade] (218.center)
			 to (219.center)
			 to (220.center)
			 to (221.center)
			 to cycle;
		\draw [style=fade] (224.center)
			 to (225.center)
			 to (222.center)
			 to (223.center)
			 to cycle;
		\draw [style=blue] (226.center) to (228.center);
		\draw [style=red] (226.center) to (227.center);
		\draw [style=blue] (229.center) to (227.center);
		\draw [style=shade] (233.center)
			 to (234.center)
			 to (235.center)
			 to (236.center)
			 to cycle;
		\draw [style=fade] (240.center)
			 to (237.center)
			 to (238.center)
			 to (239.center)
			 to cycle;
		\draw [style=shade] (247.center) to (246.center);
		\draw [style=blue, in=-150, out=30, looseness=1.75] (249.center) to (250.center);
		\draw [style=shade] (251.center)
			 to (252.center)
			 to (253.center)
			 to (254.center)
			 to cycle;
		\draw [style=red] (242.center) to (248.center);
		\draw [style=blue, in=-90, out=90] (241.center) to (243.center);
		\draw [style=shade] (258.center)
			 to (259.center)
			 to (260.center)
			 to (261.center)
			 to cycle;
		\draw [style=fade] (265.center)
			 to (262.center)
			 to (263.center)
			 to (264.center)
			 to cycle;
		\draw [style=red, rounded corners=0.4cm] (268.center)
			 to (266.center)
			 to (267.center);
		\draw [style=shade] (273.center) to (274.center);
		\draw [style=shade] (274.center) to (272.center);
		\draw [style=red, rounded corners=0.4cm] (273.center)
			 to (275.center)
			 to (276.center);
		\draw [style=blue, in=-150, out=30, looseness=1.75] (277.center) to (278.center);
		\draw [style=fade] (285.center)
			 to (282.center)
			 to (283.center)
			 to (284.center)
			 to cycle;
		\draw [style=shade] (292.center) to (291.center);
		\draw [style=blue, in=-150, out=30, looseness=1.75] (294.center) to (295.center);
		\draw [style=shade] (296.center)
			 to (297.center)
			 to (298.center)
			 to (299.center)
			 to cycle;
		\draw [style=red] (287.center) to (293.center);
		\draw [style=blue] (286.center) to (288.center);
		\draw [style=shade] (303.center)
			 to (304.center)
			 to (305.center)
			 to (306.center)
			 to cycle;
		\draw [style=fade] (310.center)
			 to (307.center)
			 to (308.center)
			 to (309.center)
			 to cycle;
		\draw [style=red, rounded corners=0.4cm] (313.center)
			 to (311.center)
			 to (312.center);
		\draw [style=shade] (318.center) to (319.center);
		\draw [style=shade] (319.center) to (317.center);
		\draw [style=red, rounded corners=0.4cm] (318.center)
			 to (320.center)
			 to (321.center);
		\draw [style=blue, in=-150, out=30, looseness=1.75] (322.center) to (323.center);
		\draw [style=fade] (345.center)
			 to (342.center)
			 to (343.center)
			 to (344.center)
			 to cycle;
		\draw [style=red, rounded corners=0.4cm] (348.center)
			 to (346.center)
			 to (347.center);
		\draw [style=blue] (348.center) to (349.center);
		\draw [style=shade] (352.center)
			 to (353.center)
			 to (354.center)
			 to (355.center)
			 to cycle;
		\draw [style=fade] (357.center) to (358.center);
		\draw [style=fade] (365.center) to (362.center);
	\end{pgfonlayer}
\end{tikzpicture}
		\caption{Example of a $\Delta_3$-tiling. Here $w = [\cdots 1\ 5\ 7\ 3\ 6\ 2\ 4\cdots]$ and $u=[\cdots 5\ 4\ 7\ 1\ 6\ 3\ 2\cdots]$. We have $u = \sigma w$ where $\sigma = t_{14}t_{15}t_{34}t_{24}=(1,4,6,7,2)$.}
		\label{fig:H-minus-ribbon}
	\end{figure}
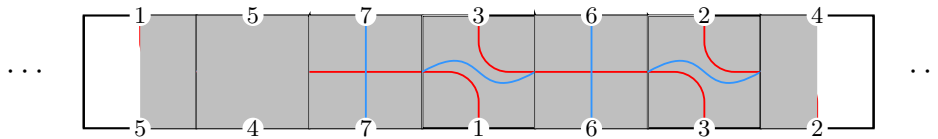

Analogs to \Cref{lem:transfer_equal_H_plus}, we show the following.
\begin{lemma}
	\label{lem:transfer_equal_H_minus}
	The transfer matrix of $\Delta_k$-tiling model agrees with the $\alpha_{-n,k}^\fff$ operators, i.e.
	\[\tilde T_-^{(k)}(x)  = \sum_{n=1}^{\infty}x^n \alpha_{-n,k}^\fff\]
\end{lemma}

\begin{proof}

	\begin{figure}[h]
		\input{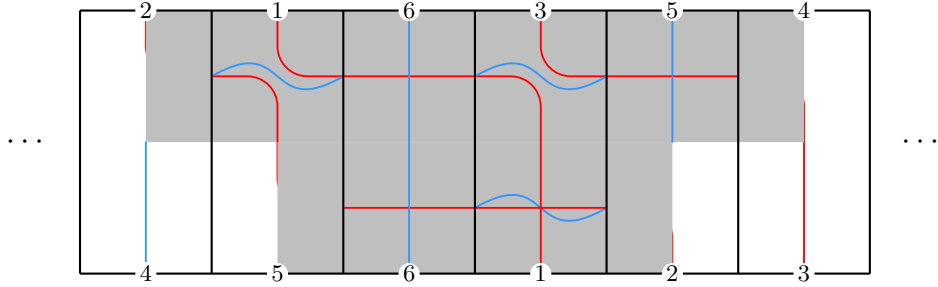}
		\caption{Example of a $\Delta_3$-tiling with multiple (two) layers. Here $w = [\cdots 2\ 1\ 6\ 3\ 5\ 4\cdots]$ and $u=[\cdots 4\ 5\ 6\ 1\ 2\ 3\cdots]$. We have $u = \sigma w$ where $\sigma = t_{34}t_{14}t_{24}t_{25}=(1,6,4,2,5)$.}
		\label{fig:H-minus-ribbon-2}
	\end{figure}

The conditions of $\Delta_k$ tiles guarantees that any tiling model obtained in this way will be an increasing $k$-Bruhat chain. And it is not hard to see that any $k$-Bruhat chain can be represented by a tiling. The requirement of having one shaded component ensures that the non-fixed paths pair-wise intersect, which means that the resulting permutation $w^{-1}u$ is a cycle. Furthermore, the number of weighted tiles equals to $\len(u)-\len(w)$, and since we do not allow multiple weighted tiles to appear in the same column, we must have $r=N-1$ where $N$ is the number of non-fixed paths. Therefore the resulting cycle must be a $(r+1)$-cycle. Thus the $\Delta_k$ tilings are in bijection with $k$ strong ribbons.\end{proof}
\begin{remark}
	Implicitly, the red (blue) coloring represents particle (hole) like movement. In terms of Gessel-Viennot non-intersecting lattice paths, our red paths corresponds to the NILP coming from Jacobi-Trudi determinant of $h_n$'s while the blue paths corresponds to those came from $e_n$'s. Note that for $\nabla$ tiles the red/blue coloring is local, but for $\Delta_k$ tiles the coloring is global (red when $\leq k$ and blue for $>k$).
\end{remark}
\subsection{Proof of \Cref{thm:heisenbergS,thm:stanley_hamiltonian}}
We will now use the $\nabla$ and $\Delta_k$ tiling diagrams to prove the commutation relations of the $\alpha_m^\fff$ operators.
\begin{prop}\label{prop: alpha_pos_commute}
	For any $m,n>0$, we have
	\([\alpha_m^\fff,\alpha_n^\fff] = 0\).
\end{prop}
\begin{proof}
  Take $w\in\sz$.	By \Cref{lem:transfer_equal_H_plus}, the expansion of $\alpha_m \alpha_n|w\rangle$ are given by all the two-row $\nabla$-tilings satisfying: (1) the top boundary is given by $w$ (2) each row has one shaded component (3) the weight is $\pm x_1^nx_2^m$. We shall now regard $\alpha_m\alpha_n|w\rangle$, and similarly $\alpha_n\alpha_m|w\rangle$, as sets of tilings. In particular, we will define an involution $\rho$, which maps an element of $\alpha_m\alpha_n|w\rangle$ to either (i) an element in $\alpha_n\alpha_m|w\rangle$ of the same sign or (ii) an element in $\alpha_m\alpha_n|w\rangle$ of opposite sign.
  
  Taking any tiling in $\alpha_m\alpha_n|w\rangle$. First off, if the shaded components of the two rows are disjoint, then $\rho$ will simply swap the two rows as follows.
 
  
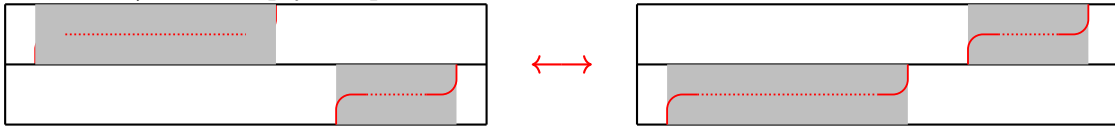
\centerline{\begin{tikzpicture}[scale = 0.4]
	\begin{pgfonlayer}{nodelayer}
		\node [style=none] (495) at (-7, -3) {};
		\node [style=none] (496) at (-8, -4) {};
		\node [style=none] (497) at (-8, -3) {};
		\node [style=none] (498) at (0, -3) {};
		\node [style=none] (499) at (0, -2) {};
		\node [style=none] (500) at (-1, -3) {};
		\node [style=none] (501) at (-9, -2) {};
		\node [style=none] (502) at (7, -2) {};
		\node [style=none] (503) at (-9, -4) {};
		\node [style=none] (504) at (7, -4) {};
		\node [style=none] (505) at (-9, -6) {};
		\node [style=none] (506) at (7, -6) {};
		\node [style=none] (507) at (-8, -2) {};
		\node [style=none] (508) at (0, -2) {};
		\node [style=none] (509) at (0, -4) {};
		\node [style=none] (510) at (-8, -4) {};
		\node [style=none] (515) at (2, -4) {};
		\node [style=none] (516) at (6, -4) {};
		\node [style=none] (517) at (6, -6) {};
		\node [style=none] (518) at (2, -6) {};
		\node [style=none] (525) at (3, -5) {};
		\node [style=none] (526) at (2, -6) {};
		\node [style=none] (527) at (2, -5) {};
		\node [style=none] (528) at (6, -5) {};
		\node [style=none] (529) at (6, -4) {};
		\node [style=none] (530) at (5, -5) {};
		\node [style=none] (531) at (12, -2) {};
		\node [style=none] (532) at (28, -2) {};
		\node [style=none] (533) at (12, -4) {};
		\node [style=none] (534) at (28, -4) {};
		\node [style=none] (535) at (12, -6) {};
		\node [style=none] (536) at (28, -6) {};
		\node [style=none] (537) at (13, -4) {};
		\node [style=none] (538) at (21, -4) {};
		\node [style=none] (539) at (21, -6) {};
		\node [style=none] (540) at (13, -6) {};
		\node [style=none] (541) at (23, -2) {};
		\node [style=none] (542) at (27, -2) {};
		\node [style=none] (543) at (27, -4) {};
		\node [style=none] (544) at (23, -4) {};
		\node [style=none] (551) at (14, -5) {};
		\node [style=none] (552) at (13, -6) {};
		\node [style=none] (553) at (13, -5) {};
		\node [style=none] (554) at (21, -5) {};
		\node [style=none] (555) at (21, -4) {};
		\node [style=none] (556) at (20, -5) {};
		\node [style=none] (557) at (24, -3) {};
		\node [style=none] (558) at (23, -4) {};
		\node [style=none] (559) at (23, -3) {};
		\node [style=none] (560) at (27, -3) {};
		\node [style=none] (561) at (27, -2) {};
		\node [style=none] (562) at (26, -3) {};
		\node [style=none] (563) at (8.5, -4) {};
		\node [style=none] (564) at (9.5, -4) {};
		\node [style=none] (565) at (10.5, -4) {};
	\end{pgfonlayer}
	\begin{pgfonlayer}{edgelayer}
		\draw [style=red, rounded corners=0.2cm] (496.center)
			 to (497.center)
			 to (495.center);
		\draw [style=red, rounded corners=0.2cm] (500.center)
			 to (498.center)
			 to (499.center);
		\draw [style=fade] (501.center) to (502.center);
		\draw [style=fade] (502.center) to (506.center);
		\draw [style=fade] (506.center) to (505.center);
		\draw [style=fade] (505.center) to (501.center);
		\draw [style=fade] (503.center) to (504.center);
		\draw [style=shade] (509.center)
			 to (510.center)
			 to (507.center)
			 to (508.center)
			 to cycle;
		\draw [style=shade] (517.center)
			 to (518.center)
			 to (515.center)
			 to (516.center)
			 to cycle;
		\draw [style=red, dashed] (495.center) to (500.center);
		\draw [style=red, rounded corners=0.2cm] (526.center)
			 to (527.center)
			 to (525.center);
		\draw [style=red, rounded corners=0.2cm] (530.center)
			 to (528.center)
			 to (529.center);
		\draw [style=red, dashed] (525.center) to (530.center);
		\draw [style=fade] (531.center) to (532.center);
		\draw [style=fade] (532.center) to (536.center);
		\draw [style=fade] (536.center) to (535.center);
		\draw [style=fade] (535.center) to (531.center);
		\draw [style=fade] (533.center) to (534.center);
		\draw [style=shade] (539.center)
			 to (540.center)
			 to (537.center)
			 to (538.center)
			 to cycle;
		\draw [style=shade] (543.center)
			 to (544.center)
			 to (541.center)
			 to (542.center)
			 to cycle;
		\draw [style=red, rounded corners=0.2cm] (552.center)
			 to (553.center)
			 to (551.center);
		\draw [style=red, rounded corners=0.2cm] (556.center)
			 to (554.center)
			 to (555.center);
		\draw [style=red, dashed] (551.center) to (556.center);
		\draw [style=red, rounded corners=0.2cm] (558.center)
			 to (559.center)
			 to (557.center);
		\draw [style=red, rounded corners=0.2cm] (562.center)
			 to (560.center)
			 to (561.center);
		\draw [style=red, dashed] (557.center) to (562.center);
		\draw [style=red arrow] (564.center) to (563.center);
		\draw [style=red arrow] (564.center) to (565.center);
	\end{pgfonlayer}
\end{tikzpicture}}
  
  We will complete the definition of $\rho$ in a case by case manner using diagrams similar to the above one. We will omit the blue path since they can be uniquely recovered. In what follows, we shall label the length of the shaded component, where a tile of type $5$ is considered to have length $2$.

  We next consider the cases when the two shaded component touches but don't overlap. 
  
  \centerline{
  \begin{tikzpicture}[scale = 0.4]
	\begin{pgfonlayer}{nodelayer}
		\node [style=none] (531) at (12, -2) {};
		\node [style=none] (532) at (28, -2) {};
		\node [style=none] (533) at (12, -4) {};
		\node [style=none] (534) at (28, -4) {};
		\node [style=none] (538) at (22, -4) {};
		\node [style=none] (541) at (22, -2) {};
		\node [style=none] (542) at (27, -2) {};
		\node [style=none] (543) at (27, -4) {};
		\node [style=none] (544) at (22, -4) {};
		\node [style=none] (555) at (22, -4) {};
		\node [style=none] (557) at (23, -3) {};
		\node [style=none] (558) at (22, -4) {};
		\node [style=none] (559) at (22, -3) {};
		\node [style=none] (560) at (27, -3) {};
		\node [style=none] (561) at (27, -2) {};
		\node [style=none] (562) at (26, -3) {};
		\node [style=none] (563) at (8.5, -4) {};
		\node [style=none] (564) at (9.5, -4) {};
		\node [style=none] (565) at (10.5, -4) {};
		\node [style=none] (566) at (22, -0.5) {};
		\node [style=none] (567) at (22, -1.5) {};
		\node [style=none] (568) at (27, -0.5) {};
		\node [style=none] (569) at (27, -1.5) {};
		\node [style=label-small] (570) at (24.5, -1) {$n$};
		\node [style=none] (571) at (22, -1) {};
		\node [style=none] (572) at (27, -1) {};
		\node [style=none] (575) at (12, -4) {};
		\node [style=none] (576) at (28, -4) {};
		\node [style=none] (577) at (12, -6) {};
		\node [style=none] (578) at (28, -6) {};
		\node [style=none] (579) at (13, -4) {};
		\node [style=none] (580) at (22, -4) {};
		\node [style=none] (581) at (22, -6) {};
		\node [style=none] (582) at (13, -6) {};
		\node [style=none] (586) at (22, -4) {};
		\node [style=none] (587) at (14, -5) {};
		\node [style=none] (588) at (13, -6) {};
		\node [style=none] (589) at (13, -5) {};
		\node [style=none] (590) at (22, -5) {};
		\node [style=none] (591) at (22, -4) {};
		\node [style=none] (592) at (21, -5) {};
		\node [style=none] (594) at (22, -4) {};
		\node [style=none] (595) at (-9, -2) {};
		\node [style=none] (596) at (7, -2) {};
		\node [style=none] (597) at (-9, -4) {};
		\node [style=none] (598) at (7, -4) {};
		\node [style=none] (599) at (-3, -2) {};
		\node [style=none] (600) at (6, -2) {};
		\node [style=none] (601) at (6, -4) {};
		\node [style=none] (602) at (-3, -4) {};
		\node [style=none] (603) at (6, -2) {};
		\node [style=none] (604) at (-2, -3) {};
		\node [style=none] (605) at (-3, -4) {};
		\node [style=none] (606) at (-3, -3) {};
		\node [style=none] (607) at (6, -3) {};
		\node [style=none] (608) at (6, -2) {};
		\node [style=none] (609) at (5, -3) {};
		\node [style=none] (610) at (6, -2) {};
		\node [style=none] (611) at (-9, -4) {};
		\node [style=none] (612) at (7, -4) {};
		\node [style=none] (613) at (-9, -6) {};
		\node [style=none] (614) at (7, -6) {};
		\node [style=none] (615) at (-8, -6) {};
		\node [style=none] (616) at (-8, -4) {};
		\node [style=none] (617) at (-3, -4) {};
		\node [style=none] (618) at (-3, -6) {};
		\node [style=none] (619) at (-8, -6) {};
		\node [style=none] (620) at (-8, -6) {};
		\node [style=none] (621) at (-7, -5) {};
		\node [style=none] (622) at (-8, -6) {};
		\node [style=none] (623) at (-8, -5) {};
		\node [style=none] (624) at (-3, -5) {};
		\node [style=none] (625) at (-3, -4) {};
		\node [style=none] (626) at (-4, -5) {};
		\node [style=none] (627) at (1, -0.5) {};
		\node [style=none] (628) at (1, -1.5) {};
		\node [style=none] (629) at (6, -0.5) {};
		\node [style=none] (630) at (6, -1.5) {};
		\node [style=label-small] (631) at (3.5, -1) {$n$};
		\node [style=none] (632) at (1, -1) {};
		\node [style=none] (633) at (6, -1) {};
		\node [style=none] (634) at (13, -6.5) {};
		\node [style=none] (635) at (13, -7.5) {};
		\node [style=none] (636) at (18, -6.5) {};
		\node [style=none] (637) at (18, -7.5) {};
		\node [style=label-small] (638) at (15.5, -7) {$n$};
		\node [style=none] (639) at (13, -7) {};
		\node [style=none] (640) at (18, -7) {};
		\node [style=none] (641) at (17, -5) {};
		\node [style=none] (642) at (19, -5) {};
		\node [style=none] (643) at (14, -5) {};
		\node [style=none] (645) at (19, -5) {};
		\node [style=none] (646) at (17, -5) {};
		\node [style=none] (647) at (19, -5) {};
		\node [style=none] (650) at (21, -5) {};
		\node [style=none] (651) at (-2, -3) {};
		\node [style=none] (652) at (2, -3) {};
		\node [style=none] (653) at (0, -3) {};
		\node [style=none] (654) at (2, -3) {};
		\node [style=none] (655) at (5, -3) {};
		\node [style=none] (675) at (27, -6.5) {};
		\node [style=none] (676) at (27, -7.5) {};
		\node [style=label-small] (677) at (22.5, -7) {$m$};
		\node [style=none] (678) at (18, -7) {};
		\node [style=none] (679) at (27, -7) {};
		\node [style=none] (680) at (-8, -6.5) {};
		\node [style=none] (681) at (-8, -7.5) {};
		\node [style=none] (682) at (-3, -6.5) {};
		\node [style=none] (683) at (-3, -7.5) {};
		\node [style=label-small] (684) at (-5.5, -7) {$n$};
		\node [style=none] (685) at (-8, -7) {};
		\node [style=none] (686) at (-3, -7) {};
		\node [style=none] (687) at (6, -6.5) {};
		\node [style=none] (688) at (6, -7.5) {};
		\node [style=label-small] (689) at (1.5, -7) {$m$};
		\node [style=none] (690) at (-3, -7) {};
		\node [style=none] (691) at (6, -7) {};
	\end{pgfonlayer}
	\begin{pgfonlayer}{edgelayer}
		\draw [style=fade] (531.center) to (532.center);
		\draw [style=fade] (533.center) to (534.center);
		\draw [style=shade] (543.center)
			 to (544.center)
			 to (541.center)
			 to (542.center)
			 to cycle;
		\draw [style=red, rounded corners=0.2cm] (558.center)
			 to (559.center)
			 to (557.center);
		\draw [style=red, rounded corners=0.2cm] (562.center)
			 to (560.center)
			 to (561.center);
		\draw [style=red, dashed] (557.center) to (562.center);
		\draw [style=red arrow] (564.center) to (563.center);
		\draw [style=red arrow] (564.center) to (565.center);
		\draw [style=fade] (566.center) to (567.center);
		\draw [style=fade] (568.center) to (569.center);
		\draw [style=red arrow] (570) to (571.center);
		\draw [style=red arrow] (570) to (572.center);
		\draw [style=fade] (578.center) to (577.center);
		\draw [style=fade] (575.center) to (576.center);
		\draw [style=shade] (581.center)
			 to (582.center)
			 to (579.center)
			 to (580.center)
			 to cycle;
		\draw [style=red, rounded corners=0.2cm] (588.center)
			 to (589.center)
			 to (587.center);
		\draw [style=red, rounded corners=0.2cm] (592.center)
			 to (590.center)
			 to (591.center);
		\draw [style=fade] (531.center) to (533.center);
		\draw [style=fade] (534.center) to (532.center);
		\draw [style=fade] (577.center) to (575.center);
		\draw [style=fade] (578.center) to (576.center);
		\draw [style=fade] (598.center) to (597.center);
		\draw [style=fade] (595.center) to (596.center);
		\draw [style=shade] (599.center)
			 to (600.center)
			 to (601.center)
			 to (602.center)
			 to cycle;
		\draw [style=red, rounded corners=0.2cm] (605.center)
			 to (606.center)
			 to (604.center);
		\draw [style=red, rounded corners=0.2cm] (609.center)
			 to (607.center)
			 to (608.center);
		\draw [style=fade] (597.center) to (595.center);
		\draw [style=fade] (598.center) to (596.center);
		\draw [style=fade] (611.center) to (612.center);
		\draw [style=fade] (613.center) to (614.center);
		\draw [style=shade] (618.center)
			 to (619.center)
			 to (616.center)
			 to (617.center)
			 to cycle;
		\draw [style=red, rounded corners=0.2cm] (622.center)
			 to (623.center)
			 to (621.center);
		\draw [style=red, rounded corners=0.2cm] (626.center)
			 to (624.center)
			 to (625.center);
		\draw [style=red, dashed] (621.center) to (626.center);
		\draw [style=fade] (611.center) to (613.center);
		\draw [style=fade] (614.center) to (612.center);
		\draw [style=fade] (627.center) to (628.center);
		\draw [style=fade] (629.center) to (630.center);
		\draw [style=red arrow] (631) to (632.center);
		\draw [style=red arrow] (631) to (633.center);
		\draw [style=fade] (634.center) to (635.center);
		\draw [style=fade] (636.center) to (637.center);
		\draw [style=red arrow] (638) to (639.center);
		\draw [style=red arrow] (638) to (640.center);
		\draw [style=red, rounded corners=0.2cm] (646.center) to (645.center);
		\draw [style=red, dashed] (643.center) to (646.center);
		\draw [style=red, dashed] (647.center) to (650.center);
		\draw [style=red, rounded corners=0.2cm] (653.center) to (652.center);
		\draw [style=red, dashed] (651.center) to (653.center);
		\draw [style=red, dashed] (654.center) to (655.center);
		\draw [style=fade] (675.center) to (676.center);
		\draw [style=fade, ->] (677) to (678.center);
		\draw [style=fade, ->] (677) to (679.center);
		\draw [style=fade] (680.center) to (681.center);
		\draw [style=fade] (682.center) to (683.center);
		\draw [style=red arrow] (684) to (685.center);
		\draw [style=red arrow] (684) to (686.center);
		\draw [style=fade] (687.center) to (688.center);
		\draw [style=red arrow] (689) to (690.center);
		\draw [style=red arrow] (689) to (691.center);
	\end{pgfonlayer}
\end{tikzpicture}
  }
  
  \centerline{
  \begin{tikzpicture}[scale = 0.4]
	\begin{pgfonlayer}{nodelayer}
		\node [style=none] (531) at (12, -2) {};
		\node [style=none] (532) at (28, -2) {};
		\node [style=none] (533) at (12, -4) {};
		\node [style=none] (534) at (28, -4) {};
		\node [style=none] (538) at (21, -4) {};
		\node [style=none] (541) at (21, -2) {};
		\node [style=none] (542) at (27, -2) {};
		\node [style=none] (543) at (27, -4) {};
		\node [style=none] (544) at (21, -4) {};
		\node [style=none] (555) at (21, -4) {};
		\node [style=none] (557) at (22, -3) {};
		\node [style=none] (558) at (21, -4) {};
		\node [style=none] (559) at (21, -3) {};
		\node [style=none] (560) at (27, -3) {};
		\node [style=none] (561) at (27, -2) {};
		\node [style=none] (562) at (26, -3) {};
		\node [style=none] (563) at (9, -4) {};
		\node [style=none] (564) at (10, -4) {};
		\node [style=none] (565) at (11, -4) {};
		\node [style=none] (566) at (21, -0.25) {};
		\node [style=none] (567) at (21, -1.25) {};
		\node [style=none] (568) at (27, -0.25) {};
		\node [style=none] (569) at (27, -1.25) {};
		\node [style=label-small] (570) at (24, -0.75) {$n$};
		\node [style=none] (571) at (21, -0.75) {};
		\node [style=none] (572) at (27, -0.75) {};
		\node [style=none] (575) at (12, -4) {};
		\node [style=none] (576) at (28, -4) {};
		\node [style=none] (577) at (12, -6) {};
		\node [style=none] (578) at (28, -6) {};
		\node [style=none] (579) at (13, -4) {};
		\node [style=none] (580) at (22, -4) {};
		\node [style=none] (581) at (22, -6) {};
		\node [style=none] (582) at (13, -6) {};
		\node [style=none] (587) at (14, -5) {};
		\node [style=none] (588) at (13, -6) {};
		\node [style=none] (589) at (13, -5) {};
		\node [style=none] (590) at (22, -5) {};
		\node [style=none] (591) at (22, -4) {};
		\node [style=none] (592) at (20, -5) {};
		\node [style=none] (595) at (-9, -2) {};
		\node [style=none] (596) at (8, -2) {};
		\node [style=none] (597) at (-9, -4) {};
		\node [style=none] (598) at (8, -4) {};
		\node [style=none] (599) at (-2, -2) {};
		\node [style=none] (600) at (7, -2) {};
		\node [style=none] (601) at (7, -4) {};
		\node [style=none] (602) at (-2, -4) {};
		\node [style=none] (603) at (7, -2) {};
		\node [style=none] (604) at (-1, -3) {};
		\node [style=none] (605) at (-2, -4) {};
		\node [style=none] (606) at (-2, -3) {};
		\node [style=none] (607) at (7, -3) {};
		\node [style=none] (608) at (7, -2) {};
		\node [style=none] (609) at (6, -3) {};
		\node [style=none] (610) at (7, -2) {};
		\node [style=none] (611) at (-9, -4) {};
		\node [style=none] (612) at (8, -4) {};
		\node [style=none] (613) at (-9, -6) {};
		\node [style=none] (614) at (8, -6) {};
		\node [style=none] (615) at (-8, -6) {};
		\node [style=none] (616) at (-8, -4) {};
		\node [style=none] (617) at (-2, -4) {};
		\node [style=none] (618) at (-2, -6) {};
		\node [style=none] (619) at (-8, -6) {};
		\node [style=none] (620) at (-8, -6) {};
		\node [style=none] (621) at (-7, -5) {};
		\node [style=none] (622) at (-8, -6) {};
		\node [style=none] (623) at (-8, -5) {};
		\node [style=none] (624) at (-2, -5) {};
		\node [style=none] (625) at (-2, -4) {};
		\node [style=none] (626) at (-3, -5) {};
		\node [style=none] (627) at (2, -0.5) {};
		\node [style=none] (628) at (2, -1.5) {};
		\node [style=none] (629) at (7, -0.5) {};
		\node [style=none] (630) at (7, -1.5) {};
		\node [style=label-small] (631) at (4.5, -1) {$n$};
		\node [style=none] (632) at (2, -1) {};
		\node [style=none] (633) at (7, -1) {};
		\node [style=none] (643) at (14, -5) {};
		\node [style=none] (646) at (20, -5) {};
		\node [style=none] (651) at (-1, -3) {};
		\node [style=none] (652) at (4, -3) {};
		\node [style=none] (653) at (1, -3) {};
		\node [style=none] (654) at (4, -3) {};
		\node [style=none] (655) at (6, -3) {};
		\node [style=none] (675) at (22, -6.5) {};
		\node [style=none] (676) at (22, -7.5) {};
		\node [style=label-small] (677) at (17.5, -7) {$m$};
		\node [style=none] (678) at (13, -7) {};
		\node [style=none] (679) at (22, -7) {};
		\node [style=none] (680) at (-8, -6.5) {};
		\node [style=none] (681) at (-8, -7.5) {};
		\node [style=none] (682) at (-2, -6.5) {};
		\node [style=none] (683) at (-2, -7.5) {};
		\node [style=label-small] (684) at (-5, -7) {$n$};
		\node [style=none] (685) at (-8, -7) {};
		\node [style=none] (686) at (-2, -7) {};
		\node [style=none] (687) at (7, -6.5) {};
		\node [style=none] (688) at (7, -7.5) {};
		\node [style=label-small] (689) at (2.5, -7) {$m$};
		\node [style=none] (690) at (-2, -7) {};
		\node [style=none] (691) at (7, -7) {};
		\node [style=none] (692) at (3, -2) {};
		\node [style=none] (693) at (3, -4) {};
		\node [style=none] (694) at (3, -4) {};
		\node [style=none] (695) at (13, -6.5) {};
		\node [style=none] (696) at (13, -7.5) {};
		\node [style=none] (697) at (23, -3) {};
		\node [style=none] (700) at (26, -3) {};
		\node [style=none] (702) at (22, -2) {};
		\node [style=none] (704) at (23, -3) {};
		\node [style=none] (705) at (22, -4) {};
		\node [style=none] (706) at (22, -3) {};
	\end{pgfonlayer}
	\begin{pgfonlayer}{edgelayer}
		\draw [style=fade] (531.center) to (532.center);
		\draw [style=fade] (533.center) to (534.center);
		\draw [style=shade] (543.center)
			 to (544.center)
			 to (541.center)
			 to (542.center)
			 to cycle;
		\draw [style=red, rounded corners=0.2cm] (558.center)
			 to (559.center)
			 to (557.center)
			 to (702.center);
		\draw [style=red, rounded corners=0.2cm] (562.center)
			 to (560.center)
			 to (561.center);
		\draw [style=red arrow] (564.center) to (563.center);
		\draw [style=red arrow] (564.center) to (565.center);
		\draw [style=fade] (566.center) to (567.center);
		\draw [style=fade] (568.center) to (569.center);
		\draw [style=red arrow] (570) to (571.center);
		\draw [style=red arrow] (570) to (572.center);
		\draw [style=fade] (578.center) to (577.center);
		\draw [style=shade] (581.center)
			 to (582.center)
			 to (579.center)
			 to (580.center)
			 to cycle;
		\draw [style=red, rounded corners=0.2cm] (588.center)
			 to (589.center)
			 to (587.center);
		\draw [style=red, rounded corners=0.2cm] (592.center)
			 to (590.center)
			 to (591.center);
		\draw [style=faded] (531.center) to (533.center);
		\draw [style=faded] (534.center) to (532.center);
		\draw [style=faded] (577.center) to (575.center);
		\draw [style=faded] (578.center) to (576.center);
		\draw [style=fade] (598.center) to (597.center);
		\draw [style=fade] (595.center) to (596.center);
		\draw [style=shade] (599.center)
			 to (600.center)
			 to (601.center)
			 to (602.center)
			 to cycle;
		\draw [style=red, rounded corners=0.2cm] (605.center)
			 to (606.center)
			 to (604.center);
		\draw [style=red, rounded corners=0.2cm] (609.center)
			 to (607.center)
			 to (608.center);
		\draw [style=faded] (597.center) to (595.center);
		\draw [style=faded] (598.center) to (596.center);
		\draw [style=fade] (613.center) to (614.center);
		\draw [style=shade] (618.center)
			 to (619.center)
			 to (616.center)
			 to (617.center)
			 to cycle;
		\draw [style=red, rounded corners=0.2cm] (622.center)
			 to (623.center)
			 to (621.center);
		\draw [style=red, rounded corners=0.2cm] (626.center)
			 to (624.center)
			 to (625.center);
		\draw [style=red, dashed] (621.center) to (626.center);
		\draw [style=faded] (611.center) to (613.center);
		\draw [style=faded] (614.center) to (612.center);
		\draw [style=fade] (627.center) to (628.center);
		\draw [style=fade] (629.center) to (630.center);
		\draw [style=red arrow] (631) to (632.center);
		\draw [style=red arrow] (631) to (633.center);
		\draw [style=red, dashed] (643.center) to (646.center);
		\draw [style=red, rounded corners=0.2cm] (653.center) to (652.center);
		\draw [style=red, dashed] (651.center) to (653.center);
		\draw [style=red, dashed] (654.center) to (655.center);
		\draw [style=fade] (675.center) to (676.center);
		\draw [style=red arrow] (677) to (678.center);
		\draw [style=red arrow] (677) to (679.center);
		\draw [style=fade] (680.center) to (681.center);
		\draw [style=fade] (682.center) to (683.center);
		\draw [style=red arrow] (684) to (685.center);
		\draw [style=red arrow] (684) to (686.center);
		\draw [style=fade] (687.center) to (688.center);
		\draw [style=red arrow] (689) to (690.center);
		\draw [style=red arrow] (689) to (691.center);
		\draw [style=red] (692.center) to (694.center);
		\draw [style=fade] (695.center) to (696.center);
		\draw [style=red, dashed] (697.center) to (700.center);
		\draw [style=red, rounded corners=0.2cm] (705.center)
			 to (706.center)
			 to (704.center);
	\end{pgfonlayer}
\end{tikzpicture}
  }
  
  \centerline{
  \begin{tikzpicture}[scale = 0.4]
	\begin{pgfonlayer}{nodelayer}
		\node [style=none] (531) at (12, -2) {};
		\node [style=none] (532) at (28, -2) {};
		\node [style=none] (533) at (12, -4) {};
		\node [style=none] (534) at (28, -4) {};
		\node [style=none] (538) at (22, -4) {};
		\node [style=none] (541) at (22, -2) {};
		\node [style=none] (542) at (27, -2) {};
		\node [style=none] (543) at (27, -4) {};
		\node [style=none] (544) at (22, -4) {};
		\node [style=none] (555) at (22, -4) {};
		\node [style=none] (557) at (23, -5) {};
		\node [style=none] (558) at (22, -6) {};
		\node [style=none] (559) at (22, -5) {};
		\node [style=none] (560) at (27, -3) {};
		\node [style=none] (561) at (27, -2) {};
		\node [style=none] (562) at (26, -3) {};
		\node [style=none] (563) at (9, -4) {};
		\node [style=none] (564) at (10, -4) {};
		\node [style=none] (565) at (11, -4) {};
		\node [style=none] (566) at (22, -0.25) {};
		\node [style=none] (567) at (22, -1.25) {};
		\node [style=none] (568) at (27, -0.25) {};
		\node [style=none] (569) at (27, -1.25) {};
		\node [style=label-small] (570) at (24.5, -0.75) {$n$};
		\node [style=none] (571) at (22, -0.75) {};
		\node [style=none] (572) at (27, -0.75) {};
		\node [style=none] (575) at (12, -4) {};
		\node [style=none] (576) at (28, -4) {};
		\node [style=none] (577) at (12, -6) {};
		\node [style=none] (578) at (28, -6) {};
		\node [style=none] (579) at (13, -4) {};
		\node [style=none] (580) at (23, -4) {};
		\node [style=none] (581) at (23, -6) {};
		\node [style=none] (582) at (13, -6) {};
		\node [style=none] (587) at (14, -5) {};
		\node [style=none] (588) at (13, -6) {};
		\node [style=none] (589) at (13, -5) {};
		\node [style=none] (590) at (22, -5) {};
		\node [style=none] (591) at (22, -4) {};
		\node [style=none] (592) at (21, -5) {};
		\node [style=none] (595) at (-9, -2) {};
		\node [style=none] (596) at (8, -2) {};
		\node [style=none] (597) at (-9, -4) {};
		\node [style=none] (598) at (8, -4) {};
		\node [style=none] (599) at (-3, -2) {};
		\node [style=none] (600) at (7, -2) {};
		\node [style=none] (601) at (7, -4) {};
		\node [style=none] (602) at (-3, -4) {};
		\node [style=none] (603) at (7, -2) {};
		\node [style=none] (604) at (-2, -3) {};
		\node [style=none] (605) at (-3, -4) {};
		\node [style=none] (606) at (-3, -3) {};
		\node [style=none] (607) at (7, -3) {};
		\node [style=none] (608) at (7, -2) {};
		\node [style=none] (609) at (6, -3) {};
		\node [style=none] (610) at (7, -2) {};
		\node [style=none] (611) at (-9, -4) {};
		\node [style=none] (612) at (8, -4) {};
		\node [style=none] (613) at (-9, -6) {};
		\node [style=none] (614) at (8, -6) {};
		\node [style=none] (615) at (-8, -6) {};
		\node [style=none] (616) at (-8, -4) {};
		\node [style=none] (617) at (-3, -4) {};
		\node [style=none] (618) at (-3, -6) {};
		\node [style=none] (619) at (-8, -6) {};
		\node [style=none] (620) at (-8, -6) {};
		\node [style=none] (621) at (-7, -5) {};
		\node [style=none] (622) at (-8, -6) {};
		\node [style=none] (623) at (-8, -5) {};
		\node [style=none] (624) at (-3, -5) {};
		\node [style=none] (625) at (-3, -4) {};
		\node [style=none] (626) at (-4, -5) {};
		\node [style=none] (627) at (2, -0.5) {};
		\node [style=none] (628) at (2, -1.5) {};
		\node [style=none] (629) at (7, -0.5) {};
		\node [style=none] (630) at (7, -1.5) {};
		\node [style=label-small] (631) at (4.5, -1) {$n$};
		\node [style=none] (632) at (2, -1) {};
		\node [style=none] (633) at (7, -1) {};
		\node [style=none] (643) at (14, -5) {};
		\node [style=none] (646) at (21, -5) {};
		\node [style=none] (651) at (-2, -3) {};
		\node [style=none] (652) at (3, -3) {};
		\node [style=none] (653) at (1, -3) {};
		\node [style=none] (654) at (3, -3) {};
		\node [style=none] (655) at (6, -3) {};
		\node [style=none] (675) at (23, -6.5) {};
		\node [style=none] (676) at (23, -7.5) {};
		\node [style=label-small] (677) at (18, -7) {$m$};
		\node [style=none] (678) at (13, -7) {};
		\node [style=none] (679) at (23, -7) {};
		\node [style=none] (680) at (-8, -6.5) {};
		\node [style=none] (681) at (-8, -7.5) {};
		\node [style=none] (682) at (-3, -6.5) {};
		\node [style=none] (683) at (-3, -7.5) {};
		\node [style=label-small] (684) at (-5.5, -7) {$n$};
		\node [style=none] (685) at (-8, -7) {};
		\node [style=none] (686) at (-3, -7) {};
		\node [style=none] (687) at (7, -6.5) {};
		\node [style=none] (688) at (7, -7.5) {};
		\node [style=label-small] (689) at (2.5, -7) {$m$};
		\node [style=none] (690) at (-3, -7) {};
		\node [style=none] (691) at (7, -7) {};
		\node [style=none] (692) at (2, -2) {};
		\node [style=none] (693) at (2, -4) {};
		\node [style=none] (694) at (2, -4) {};
		\node [style=none] (695) at (13, -6.5) {};
		\node [style=none] (696) at (13, -7.5) {};
		\node [style=none] (697) at (23, -3) {};
		\node [style=none] (700) at (26, -3) {};
		\node [style=none] (702) at (23, -4) {};
		\node [style=none] (704) at (23, -3) {};
		\node [style=none] (705) at (22, -4) {};
		\node [style=none] (706) at (22, -3) {};
	\end{pgfonlayer}
	\begin{pgfonlayer}{edgelayer}
		\draw [style=fade] (531.center) to (532.center);
		\draw [style=fade] (533.center) to (534.center);
		\draw [style=shade] (543.center)
			 to (544.center)
			 to (541.center)
			 to (542.center)
			 to cycle;
		\draw [style=red, rounded corners=0.2cm] (558.center)
			 to (559.center)
			 to (557.center)
			 to (702.center);
		\draw [style=red, rounded corners=0.2cm] (562.center)
			 to (560.center)
			 to (561.center);
		\draw [style=red arrow] (564.center) to (563.center);
		\draw [style=red arrow] (564.center) to (565.center);
		\draw [style=fade] (566.center) to (567.center);
		\draw [style=fade] (568.center) to (569.center);
		\draw [style=red arrow] (570) to (571.center);
		\draw [style=red arrow] (570) to (572.center);
		\draw [style=fade] (578.center) to (577.center);
		\draw [style=shade] (581.center)
			 to (582.center)
			 to (579.center)
			 to (580.center)
			 to cycle;
		\draw [style=red, rounded corners=0.2cm] (588.center)
			 to (589.center)
			 to (587.center);
		\draw [style=red, rounded corners=0.2cm] (592.center)
			 to (590.center)
			 to (591.center);
		\draw [style=faded] (531.center) to (533.center);
		\draw [style=faded] (534.center) to (532.center);
		\draw [style=faded] (577.center) to (575.center);
		\draw [style=faded] (578.center) to (576.center);
		\draw [style=fade] (598.center) to (597.center);
		\draw [style=fade] (595.center) to (596.center);
		\draw [style=shade] (599.center)
			 to (600.center)
			 to (601.center)
			 to (602.center)
			 to cycle;
		\draw [style=red, rounded corners=0.2cm] (605.center)
			 to (606.center)
			 to (604.center);
		\draw [style=red, rounded corners=0.2cm] (609.center)
			 to (607.center)
			 to (608.center);
		\draw [style=faded] (597.center) to (595.center);
		\draw [style=faded] (598.center) to (596.center);
		\draw [style=fade] (613.center) to (614.center);
		\draw [style=shade] (618.center)
			 to (619.center)
			 to (616.center)
			 to (617.center)
			 to cycle;
		\draw [style=red, rounded corners=0.2cm] (622.center)
			 to (623.center)
			 to (621.center);
		\draw [style=red, rounded corners=0.2cm] (626.center)
			 to (624.center)
			 to (625.center);
		\draw [style=red, dashed] (621.center) to (626.center);
		\draw [style=faded] (611.center) to (613.center);
		\draw [style=faded] (614.center) to (612.center);
		\draw [style=fade] (627.center) to (628.center);
		\draw [style=fade] (629.center) to (630.center);
		\draw [style=red arrow] (631) to (632.center);
		\draw [style=red arrow] (631) to (633.center);
		\draw [style=red, dashed] (643.center) to (646.center);
		\draw [style=red, rounded corners=0.2cm] (653.center) to (652.center);
		\draw [style=red, dashed] (651.center) to (653.center);
		\draw [style=red, dashed] (654.center) to (655.center);
		\draw [style=fade] (675.center) to (676.center);
		\draw [style=red arrow] (677) to (678.center);
		\draw [style=red arrow] (677) to (679.center);
		\draw [style=fade] (680.center) to (681.center);
		\draw [style=fade] (682.center) to (683.center);
		\draw [style=red arrow] (684) to (685.center);
		\draw [style=red arrow] (684) to (686.center);
		\draw [style=fade] (687.center) to (688.center);
		\draw [style=red arrow] (689) to (690.center);
		\draw [style=red arrow] (689) to (691.center);
		\draw [style=red] (692.center) to (694.center);
		\draw [style=fade] (695.center) to (696.center);
		\draw [style=red, dashed] (697.center) to (700.center);
		\draw [style=red, rounded corners=0.2cm] (705.center)
			 to (706.center)
			 to (704.center);
	\end{pgfonlayer}
\end{tikzpicture}
  }
  In the above three cases, $\rho$ sends $\alpha_n\alpha_m|w\rangle $ to $\alpha_m\alpha_n|w\rangle$. The next one will send $\alpha_m\alpha_n|w\rangle$ to itself while changing the spin by $1$ thereby ending up cancelling with each other.
  
  \centerline{
  \begin{tikzpicture}[scale = 0.4]
	\begin{pgfonlayer}{nodelayer}
		\node [style=none] (531) at (12, -2) {};
		\node [style=none] (532) at (29, -2) {};
		\node [style=none] (533) at (12, -4) {};
		\node [style=none] (534) at (29, -4) {};
		\node [style=none] (538) at (23, -6) {};
		\node [style=none] (541) at (23, -4) {};
		\node [style=none] (542) at (28, -4) {};
		\node [style=none] (543) at (28, -6) {};
		\node [style=none] (544) at (23, -6) {};
		\node [style=none] (555) at (23, -6) {};
		\node [style=none] (560) at (28, -5) {};
		\node [style=none] (561) at (28, -4) {};
		\node [style=none] (562) at (27, -5) {};
		\node [style=none] (563) at (9, -4) {};
		\node [style=none] (564) at (10, -4) {};
		\node [style=none] (565) at (11, -4) {};
		\node [style=none] (566) at (23, -6.5) {};
		\node [style=none] (567) at (23, -7.5) {};
		\node [style=none] (568) at (28, -6.5) {};
		\node [style=none] (569) at (28, -7.5) {};
		\node [style=label-small] (570) at (25.5, -7) {$n$};
		\node [style=none] (571) at (23, -7) {};
		\node [style=none] (572) at (28, -7) {};
		\node [style=none] (575) at (12, -4) {};
		\node [style=none] (576) at (29, -4) {};
		\node [style=none] (577) at (12, -6) {};
		\node [style=none] (578) at (29, -6) {};
		\node [style=none] (579) at (13, -2) {};
		\node [style=none] (580) at (23, -2) {};
		\node [style=none] (581) at (23, -4) {};
		\node [style=none] (582) at (13, -4) {};
		\node [style=none] (587) at (14, -3) {};
		\node [style=none] (588) at (13, -4) {};
		\node [style=none] (589) at (13, -3) {};
		\node [style=none] (590) at (23, -3) {};
		\node [style=none] (591) at (23, -2) {};
		\node [style=none] (592) at (21, -3) {};
		\node [style=none] (595) at (-9, -2) {};
		\node [style=none] (596) at (8, -2) {};
		\node [style=none] (597) at (-9, -4) {};
		\node [style=none] (598) at (8, -4) {};
		\node [style=none] (599) at (-3, -2) {};
		\node [style=none] (600) at (7, -2) {};
		\node [style=none] (601) at (7, -4) {};
		\node [style=none] (602) at (-3, -4) {};
		\node [style=none] (603) at (7, -2) {};
		\node [style=none] (604) at (-2, -3) {};
		\node [style=none] (605) at (-3, -4) {};
		\node [style=none] (606) at (-3, -3) {};
		\node [style=none] (607) at (7, -3) {};
		\node [style=none] (608) at (7, -2) {};
		\node [style=none] (609) at (6, -3) {};
		\node [style=none] (610) at (7, -2) {};
		\node [style=none] (611) at (-9, -4) {};
		\node [style=none] (612) at (8, -4) {};
		\node [style=none] (613) at (-9, -6) {};
		\node [style=none] (614) at (8, -6) {};
		\node [style=none] (615) at (-8, -6) {};
		\node [style=none] (616) at (-8, -4) {};
		\node [style=none] (617) at (-3, -4) {};
		\node [style=none] (618) at (-3, -6) {};
		\node [style=none] (619) at (-8, -6) {};
		\node [style=none] (620) at (-8, -6) {};
		\node [style=none] (621) at (-7, -5) {};
		\node [style=none] (622) at (-8, -6) {};
		\node [style=none] (623) at (-8, -5) {};
		\node [style=none] (624) at (-3, -5) {};
		\node [style=none] (625) at (-3, -4) {};
		\node [style=none] (626) at (-4, -5) {};
		\node [style=none] (627) at (2, -0.5) {};
		\node [style=none] (628) at (2, -1.5) {};
		\node [style=none] (629) at (7, -0.5) {};
		\node [style=none] (630) at (7, -1.5) {};
		\node [style=label-small] (631) at (4.5, -1) {$n$};
		\node [style=none] (632) at (2, -1) {};
		\node [style=none] (633) at (7, -1) {};
		\node [style=none] (643) at (14, -3) {};
		\node [style=none] (646) at (21, -3) {};
		\node [style=none] (651) at (-2, -3) {};
		\node [style=none] (652) at (3, -3) {};
		\node [style=none] (653) at (1, -3) {};
		\node [style=none] (654) at (3, -3) {};
		\node [style=none] (655) at (6, -3) {};
		\node [style=none] (675) at (23, -0.5) {};
		\node [style=none] (676) at (23, -1.5) {};
		\node [style=label-small] (677) at (18, -1) {$m$};
		\node [style=none] (678) at (13, -1) {};
		\node [style=none] (679) at (23, -1) {};
		\node [style=none] (680) at (-8, -6.5) {};
		\node [style=none] (681) at (-8, -7.5) {};
		\node [style=none] (682) at (-3, -6.5) {};
		\node [style=none] (683) at (-3, -7.5) {};
		\node [style=label-small] (684) at (-5.5, -7) {$n$};
		\node [style=none] (685) at (-8, -7) {};
		\node [style=none] (686) at (-3, -7) {};
		\node [style=none] (687) at (7, -6.5) {};
		\node [style=none] (688) at (7, -7.5) {};
		\node [style=label-small] (689) at (2, -7) {$m$};
		\node [style=none] (690) at (-3, -7) {};
		\node [style=none] (691) at (7, -7) {};
		\node [style=none] (692) at (2, -2) {};
		\node [style=none] (693) at (2, -4) {};
		\node [style=none] (694) at (2, -4) {};
		\node [style=none] (695) at (13, -0.5) {};
		\node [style=none] (696) at (13, -1.5) {};
		\node [style=none] (697) at (24, -5) {};
		\node [style=none] (700) at (27, -5) {};
		\node [style=none] (704) at (24, -5) {};
		\node [style=none] (705) at (23, -6) {};
		\node [style=none] (706) at (23, -5) {};
		\node [style=none] (707) at (2, -3) {};
		\node [style=none] (708) at (2, -2) {};
		\node [style=none] (709) at (1, -3) {};
		\node [style=none] (710) at (3, -3) {};
		\node [style=none] (711) at (2, -4) {};
		\node [style=none] (712) at (2, -3) {};
	\end{pgfonlayer}
	\begin{pgfonlayer}{edgelayer}
		\draw [style=fade] (531.center) to (532.center);
		\draw [style=fade] (533.center) to (534.center);
		\draw [style=shade] (543.center)
			 to (544.center)
			 to (541.center)
			 to (542.center)
			 to cycle;
		\draw [style=red, rounded corners=0.2cm] (562.center)
			 to (560.center)
			 to (561.center);
		\draw [style=red arrow] (564.center) to (563.center);
		\draw [style=red arrow] (564.center) to (565.center);
		\draw [style=fade] (566.center) to (567.center);
		\draw [style=fade] (568.center) to (569.center);
		\draw [style=red arrow] (570) to (571.center);
		\draw [style=red arrow] (570) to (572.center);
		\draw [style=fade] (578.center) to (577.center);
		\draw [style=shade] (581.center)
			 to (582.center)
			 to (579.center)
			 to (580.center)
			 to cycle;
		\draw [style=red, rounded corners=0.2cm] (588.center)
			 to (589.center)
			 to (587.center);
		\draw [style=red, rounded corners=0.2cm] (592.center)
			 to (590.center)
			 to (591.center);
		\draw [style=faded] (531.center) to (533.center);
		\draw [style=faded] (534.center) to (532.center);
		\draw [style=faded] (577.center) to (575.center);
		\draw [style=faded] (578.center) to (576.center);
		\draw [style=fade] (598.center) to (597.center);
		\draw [style=fade] (595.center) to (596.center);
		\draw [style=shade] (599.center)
			 to (600.center)
			 to (601.center)
			 to (602.center)
			 to cycle;
		\draw [style=red, rounded corners=0.2cm] (605.center)
			 to (606.center)
			 to (604.center);
		\draw [style=red, rounded corners=0.2cm] (609.center)
			 to (607.center)
			 to (608.center);
		\draw [style=faded] (597.center) to (595.center);
		\draw [style=faded] (598.center) to (596.center);
		\draw [style=fade] (613.center) to (614.center);
		\draw [style=shade] (618.center)
			 to (619.center)
			 to (616.center)
			 to (617.center)
			 to cycle;
		\draw [style=red, rounded corners=0.2cm] (622.center)
			 to (623.center)
			 to (621.center);
		\draw [style=red, rounded corners=0.2cm] (626.center)
			 to (624.center)
			 to (625.center);
		\draw [style=red, dashed] (621.center) to (626.center);
		\draw [style=faded] (611.center) to (613.center);
		\draw [style=faded] (614.center) to (612.center);
		\draw [style=fade] (627.center) to (628.center);
		\draw [style=fade] (629.center) to (630.center);
		\draw [style=red arrow] (631) to (632.center);
		\draw [style=red arrow] (631) to (633.center);
		\draw [style=red, dashed] (643.center) to (646.center);
		\draw [style=red, dashed] (651.center) to (653.center);
		\draw [style=red, dashed] (654.center) to (655.center);
		\draw [style=fade] (675.center) to (676.center);
		\draw [style=red arrow] (677) to (678.center);
		\draw [style=red arrow] (677) to (679.center);
		\draw [style=fade] (680.center) to (681.center);
		\draw [style=fade] (682.center) to (683.center);
		\draw [style=red arrow] (684) to (685.center);
		\draw [style=red arrow] (684) to (686.center);
		\draw [style=fade] (687.center) to (688.center);
		\draw [style=red arrow] (689) to (690.center);
		\draw [style=red arrow] (689) to (691.center);
		\draw [style=fade] (695.center) to (696.center);
		\draw [style=red, dashed] (697.center) to (700.center);
		\draw [style=red, rounded corners=0.2cm] (705.center)
			 to (706.center)
			 to (704.center);
		\draw [style=red, rounded corners=0.2cm] (709.center)
			 to (707.center)
			 to (708.center);
		\draw [style=red, rounded corners=0.2cm] (711.center)
			 to (712.center)
			 to (710.center);
	\end{pgfonlayer}
\end{tikzpicture}
  }
  
  We next look at the case when the two shaded components overlap. In the next case $\rho$ exchanges $\alpha_n\alpha_m$ with $\alpha_m\alpha_n$ while changing the spin by $2$.
  
  \centerline{
  \begin{tikzpicture}[scale = 0.4]
	\begin{pgfonlayer}{nodelayer}
		\node [style=none] (531) at (12, -2) {};
		\node [style=none] (532) at (29, -2) {};
		\node [style=none] (533) at (12, -4) {};
		\node [style=none] (534) at (29, -4) {};
		\node [style=none] (538) at (24, -4) {};
		\node [style=none] (541) at (20, -2) {};
		\node [style=none] (542) at (28, -2) {};
		\node [style=none] (543) at (28, -4) {};
		\node [style=none] (544) at (20, -4) {};
		\node [style=none] (560) at (28, -3) {};
		\node [style=none] (561) at (28, -2) {};
		\node [style=none] (562) at (27, -3) {};
		\node [style=none] (563) at (9, -4) {};
		\node [style=none] (564) at (10, -4) {};
		\node [style=none] (565) at (11, -4) {};
		\node [style=none] (566) at (20, -0.25) {};
		\node [style=none] (567) at (20, -1.25) {};
		\node [style=none] (568) at (28, -0.25) {};
		\node [style=none] (569) at (28, -1.25) {};
		\node [style=label-small] (570) at (24, -0.75) {$n$};
		\node [style=none] (571) at (20, -0.75) {};
		\node [style=none] (572) at (28, -0.75) {};
		\node [style=none] (575) at (12, -4) {};
		\node [style=none] (576) at (29, -4) {};
		\node [style=none] (577) at (12, -6) {};
		\node [style=none] (578) at (29, -6) {};
		\node [style=none] (579) at (13, -4) {};
		\node [style=none] (580) at (24, -4) {};
		\node [style=none] (581) at (24, -6) {};
		\node [style=none] (582) at (13, -6) {};
		\node [style=none] (587) at (14, -5) {};
		\node [style=none] (588) at (13, -6) {};
		\node [style=none] (589) at (13, -5) {};
		\node [style=none] (590) at (24, -5) {};
		\node [style=none] (591) at (24, -4) {};
		\node [style=none] (592) at (23, -5) {};
		\node [style=none] (595) at (-9, -2) {};
		\node [style=none] (596) at (8, -2) {};
		\node [style=none] (597) at (-9, -4) {};
		\node [style=none] (598) at (8, -4) {};
		\node [style=none] (599) at (-8, -2) {};
		\node [style=none] (600) at (3, -2) {};
		\node [style=none] (601) at (3, -4) {};
		\node [style=none] (602) at (-8, -4) {};
		\node [style=none] (603) at (3, -2) {};
		\node [style=none] (604) at (-7, -3) {};
		\node [style=none] (605) at (-8, -4) {};
		\node [style=none] (606) at (-8, -3) {};
		\node [style=none] (607) at (3, -3) {};
		\node [style=none] (608) at (3, -2) {};
		\node [style=none] (609) at (2, -3) {};
		\node [style=none] (610) at (3, -2) {};
		\node [style=none] (611) at (-9, -4) {};
		\node [style=none] (612) at (8, -4) {};
		\node [style=none] (613) at (-9, -6) {};
		\node [style=none] (614) at (8, -6) {};
		\node [style=none] (615) at (-1, -6) {};
		\node [style=none] (616) at (-1, -4) {};
		\node [style=none] (617) at (7, -4) {};
		\node [style=none] (618) at (7, -6) {};
		\node [style=none] (619) at (-1, -6) {};
		\node [style=none] (620) at (-1, -6) {};
		\node [style=none] (621) at (0, -5) {};
		\node [style=none] (622) at (-1, -6) {};
		\node [style=none] (623) at (-1, -5) {};
		\node [style=none] (627) at (-1, -6.5) {};
		\node [style=none] (628) at (-1, -7.5) {};
		\node [style=none] (629) at (7, -6.5) {};
		\node [style=none] (630) at (7, -7.5) {};
		\node [style=label-small] (631) at (3, -7) {$n$};
		\node [style=none] (632) at (-1, -7) {};
		\node [style=none] (633) at (7, -7) {};
		\node [style=none] (643) at (14, -5) {};
		\node [style=none] (646) at (19, -5) {};
		\node [style=none] (653) at (2, -3) {};
		\node [style=none] (655) at (2, -3) {};
		\node [style=none] (675) at (24, -6.5) {};
		\node [style=none] (676) at (24, -7.5) {};
		\node [style=label-small] (677) at (18.5, -7) {$m$};
		\node [style=none] (678) at (13, -7) {};
		\node [style=none] (679) at (24, -7) {};
		\node [style=none] (682) at (-8, -0.5) {};
		\node [style=none] (683) at (-8, -1.5) {};
		\node [style=none] (686) at (-8, -1) {};
		\node [style=none] (687) at (3, -0.5) {};
		\node [style=none] (688) at (3, -1.5) {};
		\node [style=label-small] (689) at (-2.5, -1) {$m$};
		\node [style=none] (690) at (-8, -1) {};
		\node [style=none] (691) at (3, -1) {};
		\node [style=none] (695) at (13, -6.5) {};
		\node [style=none] (696) at (13, -7.5) {};
		\node [style=none] (697) at (25, -3) {};
		\node [style=none] (700) at (27, -3) {};
		\node [style=none] (704) at (25, -3) {};
		\node [style=none] (705) at (24, -4) {};
		\node [style=none] (706) at (24, -3) {};
		\node [style=none] (716) at (7, -5) {};
		\node [style=none] (717) at (7, -4) {};
		\node [style=none] (718) at (6, -5) {};
		\node [style=none] (720) at (20, -4) {};
		\node [style=none] (721) at (21, -3) {};
		\node [style=none] (722) at (20, -4) {};
		\node [style=none] (723) at (20, -3) {};
		\node [style=none] (725) at (0, -5) {};
		\node [style=none] (726) at (6, -5) {};
		\node [style=none] (737) at (-7, -3) {};
		\node [style=none] (738) at (2, -3) {};
		\node [style=none] (743) at (20, -4) {};
		\node [style=none] (744) at (20, -6) {};
		\node [style=none] (745) at (23, -5) {};
		\node [style=none] (746) at (21, -5) {};
		\node [style=none] (747) at (23, -5) {};
		\node [style=none] (748) at (24, -3) {};
		\node [style=none] (749) at (24, -2) {};
		\node [style=none] (750) at (23, -3) {};
		\node [style=none] (752) at (21, -3) {};
		\node [style=none] (753) at (23, -3) {};
		\node [style=none] (754) at (20, -6) {};
		\node [style=none] (755) at (21, -5) {};
		\node [style=none] (756) at (21, -5) {};
		\node [style=none] (757) at (20, -6) {};
		\node [style=none] (758) at (20, -5) {};
		\node [style=none] (759) at (20, -5) {};
		\node [style=none] (760) at (20, -4) {};
		\node [style=none] (761) at (19, -5) {};
		\node [style=none] (762) at (19, -5) {};
	\end{pgfonlayer}
	\begin{pgfonlayer}{edgelayer}
		\draw [style=fade] (531.center) to (532.center);
		\draw [style=fade] (533.center) to (534.center);
		\draw [style=shade] (543.center)
			 to (544.center)
			 to (541.center)
			 to (542.center)
			 to cycle;
		\draw [style=red, rounded corners=0.2cm] (562.center)
			 to (560.center)
			 to (561.center);
		\draw [style=red arrow] (564.center) to (563.center);
		\draw [style=red arrow] (564.center) to (565.center);
		\draw [style=fade] (566.center) to (567.center);
		\draw [style=fade] (568.center) to (569.center);
		\draw [style=red arrow] (570) to (571.center);
		\draw [style=red arrow] (570) to (572.center);
		\draw [style=fade] (578.center) to (577.center);
		\draw [style=shade] (581.center)
			 to (582.center)
			 to (579.center)
			 to (580.center)
			 to cycle;
		\draw [style=red, rounded corners=0.2cm] (588.center)
			 to (589.center)
			 to (587.center);
		\draw [style=red, rounded corners=0.2cm] (592.center)
			 to (590.center)
			 to (591.center);
		\draw [style=faded] (531.center) to (533.center);
		\draw [style=faded] (534.center) to (532.center);
		\draw [style=faded] (577.center) to (575.center);
		\draw [style=faded] (578.center) to (576.center);
		\draw [style=fade] (598.center) to (597.center);
		\draw [style=fade] (595.center) to (596.center);
		\draw [style=shade] (599.center)
			 to (600.center)
			 to (601.center)
			 to (602.center)
			 to cycle;
		\draw [style=red, rounded corners=0.2cm] (605.center)
			 to (606.center)
			 to (604.center);
		\draw [style=red, rounded corners=0.2cm] (609.center)
			 to (607.center)
			 to (608.center);
		\draw [style=faded] (597.center) to (595.center);
		\draw [style=faded] (598.center) to (596.center);
		\draw [style=fade] (613.center) to (614.center);
		\draw [style=shade] (618.center)
			 to (619.center)
			 to (616.center)
			 to (617.center)
			 to cycle;
		\draw [style=red, rounded corners=0.2cm] (622.center)
			 to (623.center)
			 to (621.center);
		\draw [style=faded] (611.center) to (613.center);
		\draw [style=faded] (614.center) to (612.center);
		\draw [style=fade] (627.center) to (628.center);
		\draw [style=fade] (629.center) to (630.center);
		\draw [style=red arrow] (631) to (632.center);
		\draw [style=red arrow] (631) to (633.center);
		\draw [style=red, dashed] (643.center) to (646.center);
		\draw [style=fade] (675.center) to (676.center);
		\draw [style=red arrow] (677) to (678.center);
		\draw [style=red arrow] (677) to (679.center);
		\draw [style=fade] (682.center) to (683.center);
		\draw [style=fade] (687.center) to (688.center);
		\draw [style=red arrow] (689) to (690.center);
		\draw [style=red arrow] (689) to (691.center);
		\draw [style=fade] (695.center) to (696.center);
		\draw [style=red, dashed] (697.center) to (700.center);
		\draw [style=red, rounded corners=0.2cm] (705.center)
			 to (706.center)
			 to (704.center);
		\draw [style=red, rounded corners=0.2cm] (718.center)
			 to (716.center)
			 to (717.center);
		\draw [style=red, rounded corners=0.2cm] (722.center)
			 to (723.center)
			 to (721.center);
		\draw [style=red, dashed] (725.center) to (726.center);
		\draw [style=red, dashed] (737.center) to (738.center);
		\draw [style=red, dashed] (746.center) to (747.center);
		\draw [style=red, rounded corners=0.2cm] (750.center)
			 to (748.center)
			 to (749.center);
		\draw [style=red, dashed] (752.center) to (753.center);
		\draw [style=red, rounded corners=0.2cm] (757.center)
			 to (758.center)
			 to (756.center);
		\draw [style=red, rounded corners=0.2cm] (761.center)
			 to (759.center)
			 to (760.center);
	\end{pgfonlayer}
\end{tikzpicture}
  }
  
 In the next case $\rho$ exchanges $\alpha_m\alpha_n$ with itself while changing the spin by $1$.
  
  \centerline{
  \begin{tikzpicture}[scale = 0.4]
	\begin{pgfonlayer}{nodelayer}
		\node [style=none] (531) at (12, -2) {};
		\node [style=none] (532) at (29, -2) {};
		\node [style=none] (533) at (12, -4) {};
		\node [style=none] (534) at (29, -4) {};
		\node [style=none] (538) at (22, -4) {};
		\node [style=none] (541) at (17, -2) {};
		\node [style=none] (542) at (28, -2) {};
		\node [style=none] (543) at (28, -4) {};
		\node [style=none] (544) at (17, -4) {};
		\node [style=none] (560) at (28, -3) {};
		\node [style=none] (561) at (28, -2) {};
		\node [style=none] (562) at (27, -3) {};
		\node [style=none] (563) at (9, -4) {};
		\node [style=none] (564) at (10, -4) {};
		\node [style=none] (565) at (11, -4) {};
		\node [style=none] (566) at (17, -0.25) {};
		\node [style=none] (567) at (17, -1.25) {};
		\node [style=none] (568) at (28, -0.25) {};
		\node [style=none] (569) at (28, -1.25) {};
		\node [style=label-small] (570) at (22.5, -0.75) {$m$};
		\node [style=none] (571) at (17, -0.75) {};
		\node [style=none] (572) at (28, -0.75) {};
		\node [style=none] (575) at (12, -4) {};
		\node [style=none] (576) at (29, -4) {};
		\node [style=none] (577) at (12, -6) {};
		\node [style=none] (578) at (29, -6) {};
		\node [style=none] (579) at (13, -4) {};
		\node [style=none] (580) at (22, -4) {};
		\node [style=none] (581) at (22, -6) {};
		\node [style=none] (582) at (13, -6) {};
		\node [style=none] (587) at (14, -5) {};
		\node [style=none] (588) at (13, -6) {};
		\node [style=none] (589) at (13, -5) {};
		\node [style=none] (590) at (22, -5) {};
		\node [style=none] (591) at (22, -4) {};
		\node [style=none] (592) at (21, -5) {};
		\node [style=none] (595) at (-9, -2) {};
		\node [style=none] (596) at (8, -2) {};
		\node [style=none] (597) at (-9, -4) {};
		\node [style=none] (598) at (8, -4) {};
		\node [style=none] (599) at (-4, -2) {};
		\node [style=none] (600) at (7, -2) {};
		\node [style=none] (601) at (7, -4) {};
		\node [style=none] (602) at (-4, -4) {};
		\node [style=none] (603) at (7, -2) {};
		\node [style=none] (604) at (-3, -3) {};
		\node [style=none] (605) at (-4, -4) {};
		\node [style=none] (606) at (-4, -3) {};
		\node [style=none] (607) at (7, -3) {};
		\node [style=none] (608) at (7, -2) {};
		\node [style=none] (609) at (6, -3) {};
		\node [style=none] (610) at (7, -2) {};
		\node [style=none] (611) at (-9, -4) {};
		\node [style=none] (612) at (8, -4) {};
		\node [style=none] (613) at (-9, -6) {};
		\node [style=none] (614) at (8, -6) {};
		\node [style=none] (615) at (-8, -6) {};
		\node [style=none] (616) at (-8, -4) {};
		\node [style=none] (617) at (2, -4) {};
		\node [style=none] (618) at (2, -6) {};
		\node [style=none] (619) at (-8, -6) {};
		\node [style=none] (620) at (-8, -6) {};
		\node [style=none] (621) at (-7, -5) {};
		\node [style=none] (622) at (-8, -6) {};
		\node [style=none] (623) at (-8, -5) {};
		\node [style=none] (624) at (-4, -5) {};
		\node [style=none] (625) at (-4, -4) {};
		\node [style=none] (626) at (-5, -5) {};
		\node [style=none] (627) at (-8, -6.5) {};
		\node [style=none] (628) at (-8, -7.5) {};
		\node [style=none] (629) at (2, -6.5) {};
		\node [style=none] (630) at (2, -7.5) {};
		\node [style=label-small] (631) at (-3, -7) {$n$};
		\node [style=none] (632) at (-8, -7) {};
		\node [style=none] (633) at (2, -7) {};
		\node [style=none] (643) at (14, -5) {};
		\node [style=none] (646) at (16, -5) {};
		\node [style=none] (651) at (2, -3) {};
		\node [style=none] (653) at (6, -3) {};
		\node [style=none] (655) at (6, -3) {};
		\node [style=none] (675) at (22, -6.5) {};
		\node [style=none] (676) at (22, -7.5) {};
		\node [style=label-small] (677) at (17.5, -7) {$n$};
		\node [style=none] (678) at (13, -7) {};
		\node [style=none] (679) at (22, -7) {};
		\node [style=none] (682) at (-4, -0.5) {};
		\node [style=none] (683) at (-4, -1.5) {};
		\node [style=none] (686) at (-4, -1) {};
		\node [style=none] (687) at (7, -0.5) {};
		\node [style=none] (688) at (7, -1.5) {};
		\node [style=label-small] (689) at (1.5, -1) {$m$};
		\node [style=none] (690) at (-4, -1) {};
		\node [style=none] (691) at (7, -1) {};
		\node [style=none] (695) at (13, -6.5) {};
		\node [style=none] (696) at (13, -7.5) {};
		\node [style=none] (697) at (23, -3) {};
		\node [style=none] (700) at (27, -3) {};
		\node [style=none] (704) at (23, -3) {};
		\node [style=none] (705) at (22, -4) {};
		\node [style=none] (706) at (22, -3) {};
		\node [style=none] (713) at (-3, -5) {};
		\node [style=none] (714) at (-4, -6) {};
		\node [style=none] (715) at (-4, -5) {};
		\node [style=none] (716) at (2, -5) {};
		\node [style=none] (717) at (2, -4) {};
		\node [style=none] (718) at (1, -5) {};
		\node [style=none] (720) at (17, -4) {};
		\node [style=none] (721) at (18, -3) {};
		\node [style=none] (722) at (17, -4) {};
		\node [style=none] (723) at (17, -3) {};
		\node [style=none] (725) at (-3, -5) {};
		\node [style=none] (726) at (1, -5) {};
		\node [style=none] (737) at (-3, -3) {};
		\node [style=none] (738) at (0, -3) {};
		\node [style=none] (743) at (17, -4) {};
		\node [style=none] (744) at (17, -6) {};
		\node [style=none] (745) at (21, -5) {};
		\node [style=none] (746) at (18, -5) {};
		\node [style=none] (747) at (21, -5) {};
		\node [style=none] (748) at (22, -3) {};
		\node [style=none] (749) at (22, -2) {};
		\node [style=none] (750) at (21, -3) {};
		\node [style=none] (752) at (18, -3) {};
		\node [style=none] (753) at (21, -3) {};
	\end{pgfonlayer}
	\begin{pgfonlayer}{edgelayer}
		\draw [style=fade] (531.center) to (532.center);
		\draw [style=fade] (533.center) to (534.center);
		\draw [style=shade] (543.center)
			 to (544.center)
			 to (541.center)
			 to (542.center)
			 to cycle;
		\draw [style=red, rounded corners=0.2cm] (562.center)
			 to (560.center)
			 to (561.center);
		\draw [style=red arrow] (564.center) to (563.center);
		\draw [style=red arrow] (564.center) to (565.center);
		\draw [style=fade] (566.center) to (567.center);
		\draw [style=fade] (568.center) to (569.center);
		\draw [style=red arrow] (570) to (571.center);
		\draw [style=red arrow] (570) to (572.center);
		\draw [style=fade] (578.center) to (577.center);
		\draw [style=shade] (581.center)
			 to (582.center)
			 to (579.center)
			 to (580.center)
			 to cycle;
		\draw [style=red, rounded corners=0.2cm] (588.center)
			 to (589.center)
			 to (587.center);
		\draw [style=red, rounded corners=0.2cm] (592.center)
			 to (590.center)
			 to (591.center);
		\draw [style=faded] (531.center) to (533.center);
		\draw [style=faded] (534.center) to (532.center);
		\draw [style=faded] (577.center) to (575.center);
		\draw [style=faded] (578.center) to (576.center);
		\draw [style=fade] (598.center) to (597.center);
		\draw [style=fade] (595.center) to (596.center);
		\draw [style=shade] (599.center)
			 to (600.center)
			 to (601.center)
			 to (602.center)
			 to cycle;
		\draw [style=red, rounded corners=0.2cm] (605.center)
			 to (606.center)
			 to (604.center);
		\draw [style=red, rounded corners=0.2cm] (609.center)
			 to (607.center)
			 to (608.center);
		\draw [style=faded] (597.center) to (595.center);
		\draw [style=faded] (598.center) to (596.center);
		\draw [style=fade] (613.center) to (614.center);
		\draw [style=shade] (618.center)
			 to (619.center)
			 to (616.center)
			 to (617.center)
			 to cycle;
		\draw [style=red, rounded corners=0.2cm] (622.center)
			 to (623.center)
			 to (621.center);
		\draw [style=red, rounded corners=0.2cm] (626.center)
			 to (624.center)
			 to (625.center);
		\draw [style=red, dashed] (621.center) to (626.center);
		\draw [style=faded] (611.center) to (613.center);
		\draw [style=faded] (614.center) to (612.center);
		\draw [style=fade] (627.center) to (628.center);
		\draw [style=fade] (629.center) to (630.center);
		\draw [style=red arrow] (631) to (632.center);
		\draw [style=red arrow] (631) to (633.center);
		\draw [style=red, dashed] (643.center) to (646.center);
		\draw [style=red, dashed] (651.center) to (653.center);
		\draw [style=fade] (675.center) to (676.center);
		\draw [style=red arrow] (677) to (678.center);
		\draw [style=red arrow] (677) to (679.center);
		\draw [style=fade] (682.center) to (683.center);
		\draw [style=fade] (687.center) to (688.center);
		\draw [style=red arrow] (689) to (690.center);
		\draw [style=red arrow] (689) to (691.center);
		\draw [style=fade] (695.center) to (696.center);
		\draw [style=red, dashed] (697.center) to (700.center);
		\draw [style=red, rounded corners=0.2cm] (705.center)
			 to (706.center)
			 to (704.center);
		\draw [style=red, rounded corners=0.2cm] (714.center)
			 to (715.center)
			 to (713.center);
		\draw [style=red, rounded corners=0.2cm] (718.center)
			 to (716.center)
			 to (717.center);
		\draw [style=red, rounded corners=0.2cm] (722.center)
			 to (723.center)
			 to (721.center);
		\draw [style=red, dashed] (725.center) to (726.center);
		\draw [style=red, dashed] (737.center) to (738.center);
		\draw [style=red] (738.center) to (651.center);
		\draw [style=red] (743.center) to (744.center);
		\draw [style=red, dashed] (746.center) to (747.center);
		\draw [style=red] (646.center) to (746.center);
		\draw [style=red, rounded corners=0.2cm] (750.center)
			 to (748.center)
			 to (749.center);
		\draw [style=red, dashed] (752.center) to (753.center);
	\end{pgfonlayer}
\end{tikzpicture}
  }
 
And finally in the rest of the cases, $\rho$ exchanges $\alpha_m\alpha_n$ with $H_nH_m$ with the same spin.

  
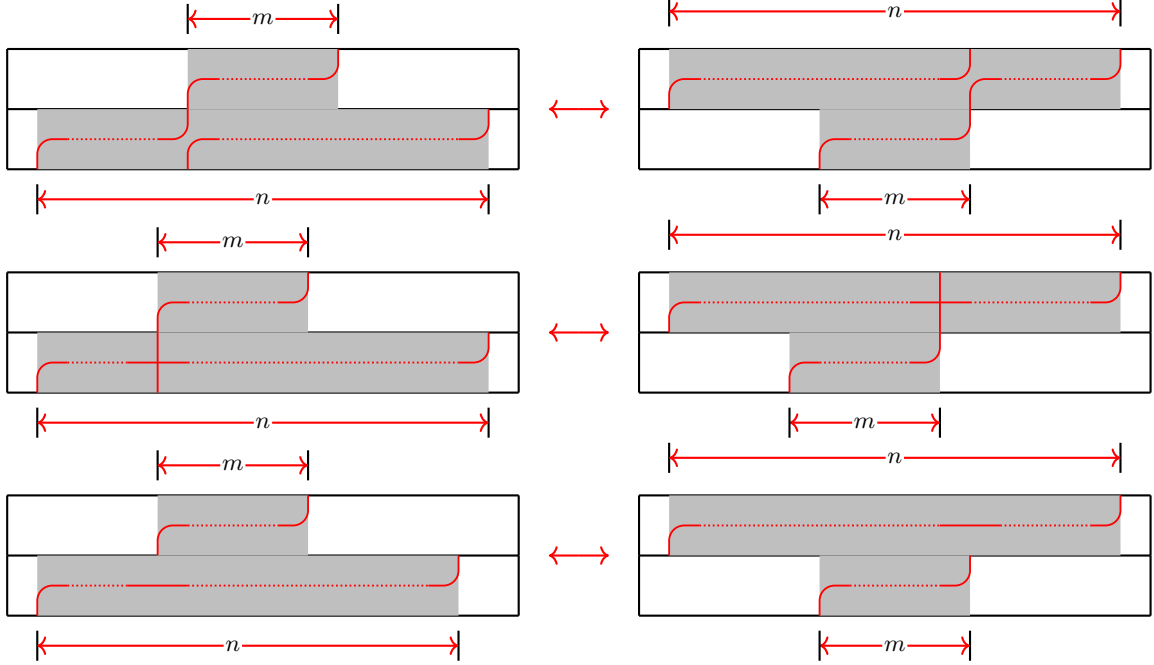
\centerline{
  \begin{tikzpicture}[scale = 0.4]
	\begin{pgfonlayer}{nodelayer}
		\node [style=none] (531) at (12, -2) {};
		\node [style=none] (532) at (29, -2) {};
		\node [style=none] (533) at (12, -4) {};
		\node [style=none] (534) at (29, -4) {};
		\node [style=none] (538) at (23, -4) {};
		\node [style=none] (541) at (13, -2) {};
		\node [style=none] (542) at (28, -2) {};
		\node [style=none] (543) at (28, -4) {};
		\node [style=none] (544) at (13, -4) {};
		\node [style=none] (560) at (28, -3) {};
		\node [style=none] (561) at (28, -2) {};
		\node [style=none] (562) at (27, -3) {};
		\node [style=none] (563) at (9, -4) {};
		\node [style=none] (564) at (10, -4) {};
		\node [style=none] (565) at (11, -4) {};
		\node [style=none] (566) at (13, -0.25) {};
		\node [style=none] (567) at (13, -1.25) {};
		\node [style=none] (568) at (28, -0.25) {};
		\node [style=none] (569) at (28, -1.25) {};
		\node [style=label-small] (570) at (20.5, -0.75) {$n$};
		\node [style=none] (571) at (13, -0.75) {};
		\node [style=none] (572) at (28, -0.75) {};
		\node [style=none] (575) at (12, -4) {};
		\node [style=none] (576) at (29, -4) {};
		\node [style=none] (577) at (12, -6) {};
		\node [style=none] (578) at (29, -6) {};
		\node [style=none] (579) at (18, -4) {};
		\node [style=none] (580) at (23, -4) {};
		\node [style=none] (581) at (23, -6) {};
		\node [style=none] (582) at (18, -6) {};
		\node [style=none] (590) at (23, -5) {};
		\node [style=none] (591) at (23, -4) {};
		\node [style=none] (592) at (22, -5) {};
		\node [style=none] (595) at (-9, -2) {};
		\node [style=none] (596) at (8, -2) {};
		\node [style=none] (597) at (-9, -4) {};
		\node [style=none] (598) at (8, -4) {};
		\node [style=none] (599) at (-3, -2) {};
		\node [style=none] (600) at (2, -2) {};
		\node [style=none] (601) at (2, -4) {};
		\node [style=none] (602) at (-3, -4) {};
		\node [style=none] (603) at (2, -2) {};
		\node [style=none] (604) at (-2, -3) {};
		\node [style=none] (605) at (-3, -4) {};
		\node [style=none] (606) at (-3, -3) {};
		\node [style=none] (607) at (2, -3) {};
		\node [style=none] (608) at (2, -2) {};
		\node [style=none] (609) at (1, -3) {};
		\node [style=none] (610) at (2, -2) {};
		\node [style=none] (611) at (-9, -4) {};
		\node [style=none] (612) at (8, -4) {};
		\node [style=none] (613) at (-9, -6) {};
		\node [style=none] (614) at (8, -6) {};
		\node [style=none] (616) at (-8, -4) {};
		\node [style=none] (617) at (7, -4) {};
		\node [style=none] (618) at (7, -6) {};
		\node [style=none] (619) at (-8, -6) {};
		\node [style=none] (621) at (-2, -5) {};
		\node [style=none] (622) at (-3, -6) {};
		\node [style=none] (623) at (-3, -5) {};
		\node [style=none] (627) at (-8, -6.5) {};
		\node [style=none] (628) at (-8, -7.5) {};
		\node [style=none] (629) at (7, -6.5) {};
		\node [style=none] (630) at (7, -7.5) {};
		\node [style=label-small] (631) at (-0.5, -7) {$n$};
		\node [style=none] (632) at (-8, -7) {};
		\node [style=none] (633) at (7, -7) {};
		\node [style=none] (653) at (1, -3) {};
		\node [style=none] (655) at (1, -3) {};
		\node [style=none] (675) at (23, -6.5) {};
		\node [style=none] (676) at (23, -7.5) {};
		\node [style=label-small] (677) at (20.5, -7) {$m$};
		\node [style=none] (678) at (18, -7) {};
		\node [style=none] (679) at (23, -7) {};
		\node [style=none] (682) at (-3, -0.5) {};
		\node [style=none] (683) at (-3, -1.5) {};
		\node [style=none] (686) at (-3, -1) {};
		\node [style=none] (687) at (2, -0.5) {};
		\node [style=none] (688) at (2, -1.5) {};
		\node [style=label-small] (689) at (-0.5, -1) {$m$};
		\node [style=none] (690) at (-3, -1) {};
		\node [style=none] (691) at (2, -1) {};
		\node [style=none] (695) at (18, -6.5) {};
		\node [style=none] (696) at (18, -7.5) {};
		\node [style=none] (697) at (24, -3) {};
		\node [style=none] (700) at (27, -3) {};
		\node [style=none] (704) at (24, -3) {};
		\node [style=none] (705) at (23, -4) {};
		\node [style=none] (706) at (23, -3) {};
		\node [style=none] (716) at (7, -5) {};
		\node [style=none] (717) at (7, -4) {};
		\node [style=none] (718) at (6, -5) {};
		\node [style=none] (720) at (13, -4) {};
		\node [style=none] (721) at (14, -3) {};
		\node [style=none] (722) at (13, -4) {};
		\node [style=none] (723) at (13, -3) {};
		\node [style=none] (725) at (-2, -5) {};
		\node [style=none] (726) at (6, -5) {};
		\node [style=none] (737) at (-2, -3) {};
		\node [style=none] (738) at (1, -3) {};
		\node [style=none] (744) at (18, -6) {};
		\node [style=none] (745) at (22, -5) {};
		\node [style=none] (746) at (19, -5) {};
		\node [style=none] (747) at (22, -5) {};
		\node [style=none] (748) at (23, -3) {};
		\node [style=none] (749) at (23, -2) {};
		\node [style=none] (750) at (22, -3) {};
		\node [style=none] (752) at (14, -3) {};
		\node [style=none] (753) at (22, -3) {};
		\node [style=none] (754) at (18, -6) {};
		\node [style=none] (755) at (19, -5) {};
		\node [style=none] (756) at (19, -5) {};
		\node [style=none] (757) at (18, -6) {};
		\node [style=none] (758) at (18, -5) {};
		\node [style=none] (764) at (-8, -6) {};
		\node [style=none] (765) at (-7, -5) {};
		\node [style=none] (766) at (-8, -6) {};
		\node [style=none] (767) at (-8, -5) {};
		\node [style=none] (768) at (-7, -5) {};
		\node [style=none] (769) at (-4, -5) {};
		\node [style=none] (770) at (-3, -4) {};
		\node [style=none] (774) at (-3, -5) {};
		\node [style=none] (775) at (-3, -5) {};
		\node [style=none] (776) at (-3, -4) {};
		\node [style=none] (777) at (-4, -5) {};
		\node [style=none] (778) at (-4, -5) {};
	\end{pgfonlayer}
	\begin{pgfonlayer}{edgelayer}
		\draw [style=fade] (531.center) to (532.center);
		\draw [style=fade] (533.center) to (534.center);
		\draw [style=shade] (543.center)
			 to (544.center)
			 to (541.center)
			 to (542.center)
			 to cycle;
		\draw [style=red, rounded corners=0.2cm] (562.center)
			 to (560.center)
			 to (561.center);
		\draw [style=red arrow] (564.center) to (563.center);
		\draw [style=red arrow] (564.center) to (565.center);
		\draw [style=fade] (566.center) to (567.center);
		\draw [style=fade] (568.center) to (569.center);
		\draw [style=red arrow] (570) to (571.center);
		\draw [style=red arrow] (570) to (572.center);
		\draw [style=fade] (578.center) to (577.center);
		\draw [style=shade] (581.center)
			 to (582.center)
			 to (579.center)
			 to (580.center)
			 to cycle;
		\draw [style=red, rounded corners=0.2cm] (592.center)
			 to (590.center)
			 to (591.center);
		\draw [style=faded] (531.center) to (533.center);
		\draw [style=faded] (534.center) to (532.center);
		\draw [style=faded] (577.center) to (575.center);
		\draw [style=faded] (578.center) to (576.center);
		\draw [style=fade] (598.center) to (597.center);
		\draw [style=fade] (595.center) to (596.center);
		\draw [style=shade] (599.center)
			 to (600.center)
			 to (601.center)
			 to (602.center)
			 to cycle;
		\draw [style=red, rounded corners=0.2cm] (605.center)
			 to (606.center)
			 to (604.center);
		\draw [style=red, rounded corners=0.2cm] (609.center)
			 to (607.center)
			 to (608.center);
		\draw [style=faded] (597.center) to (595.center);
		\draw [style=faded] (598.center) to (596.center);
		\draw [style=fade] (613.center) to (614.center);
		\draw [style=shade] (618.center)
			 to (619.center)
			 to (616.center)
			 to (617.center)
			 to cycle;
		\draw [style=red, rounded corners=0.2cm] (622.center)
			 to (623.center)
			 to (621.center);
		\draw [style=faded] (611.center) to (613.center);
		\draw [style=faded] (614.center) to (612.center);
		\draw [style=fade] (627.center) to (628.center);
		\draw [style=fade] (629.center) to (630.center);
		\draw [style=red arrow] (631) to (632.center);
		\draw [style=red arrow] (631) to (633.center);
		\draw [style=fade] (675.center) to (676.center);
		\draw [style=red arrow] (677) to (678.center);
		\draw [style=red arrow] (677) to (679.center);
		\draw [style=fade] (682.center) to (683.center);
		\draw [style=fade] (687.center) to (688.center);
		\draw [style=red arrow] (689) to (690.center);
		\draw [style=red arrow] (689) to (691.center);
		\draw [style=fade] (695.center) to (696.center);
		\draw [style=red, dashed] (697.center) to (700.center);
		\draw [style=red, rounded corners=0.2cm] (705.center)
			 to (706.center)
			 to (704.center);
		\draw [style=red, rounded corners=0.2cm] (718.center)
			 to (716.center)
			 to (717.center);
		\draw [style=red, rounded corners=0.2cm] (722.center)
			 to (723.center)
			 to (721.center);
		\draw [style=red, dashed] (725.center) to (726.center);
		\draw [style=red, dashed] (737.center) to (738.center);
		\draw [style=red, dashed] (746.center) to (747.center);
		\draw [style=red, rounded corners=0.2cm] (750.center)
			 to (748.center)
			 to (749.center);
		\draw [style=red, dashed] (752.center) to (753.center);
		\draw [style=red, rounded corners=0.2cm] (757.center)
			 to (758.center)
			 to (756.center);
		\draw [style=red, rounded corners=0.2cm] (766.center)
			 to (767.center)
			 to (765.center);
		\draw [style=red, dashed] (768.center) to (769.center);
		\draw [style=red, rounded corners=0.2cm] (777.center)
			 to (775.center)
			 to (776.center);
	\end{pgfonlayer}
\end{tikzpicture}
  }
  
  
  \centerline{
  \begin{tikzpicture}[scale=0.4]
	\begin{pgfonlayer}{nodelayer}
		\node [style=none] (531) at (12, -2) {};
		\node [style=none] (532) at (29, -2) {};
		\node [style=none] (533) at (12, -4) {};
		\node [style=none] (534) at (29, -4) {};
		\node [style=none] (538) at (22, -4) {};
		\node [style=none] (541) at (13, -2) {};
		\node [style=none] (542) at (28, -2) {};
		\node [style=none] (543) at (28, -4) {};
		\node [style=none] (544) at (13, -4) {};
		\node [style=none] (560) at (28, -3) {};
		\node [style=none] (561) at (28, -2) {};
		\node [style=none] (562) at (27, -3) {};
		\node [style=none] (563) at (9, -4) {};
		\node [style=none] (564) at (10, -4) {};
		\node [style=none] (565) at (11, -4) {};
		\node [style=none] (566) at (13, -0.25) {};
		\node [style=none] (567) at (13, -1.25) {};
		\node [style=none] (568) at (28, -0.25) {};
		\node [style=none] (569) at (28, -1.25) {};
		\node [style=label-small] (570) at (20.5, -0.75) {$n$};
		\node [style=none] (571) at (13, -0.75) {};
		\node [style=none] (572) at (28, -0.75) {};
		\node [style=none] (575) at (12, -4) {};
		\node [style=none] (576) at (29, -4) {};
		\node [style=none] (577) at (12, -6) {};
		\node [style=none] (578) at (29, -6) {};
		\node [style=none] (579) at (17, -4) {};
		\node [style=none] (580) at (22, -4) {};
		\node [style=none] (581) at (22, -6) {};
		\node [style=none] (582) at (17, -6) {};
		\node [style=none] (590) at (22, -5) {};
		\node [style=none] (591) at (22, -4) {};
		\node [style=none] (592) at (21, -5) {};
		\node [style=none] (595) at (-9, -2) {};
		\node [style=none] (596) at (8, -2) {};
		\node [style=none] (597) at (-9, -4) {};
		\node [style=none] (598) at (8, -4) {};
		\node [style=none] (599) at (-4, -2) {};
		\node [style=none] (600) at (1, -2) {};
		\node [style=none] (601) at (1, -4) {};
		\node [style=none] (602) at (-4, -4) {};
		\node [style=none] (603) at (1, -2) {};
		\node [style=none] (604) at (-3, -3) {};
		\node [style=none] (605) at (-4, -4) {};
		\node [style=none] (606) at (-4, -3) {};
		\node [style=none] (607) at (1, -3) {};
		\node [style=none] (608) at (1, -2) {};
		\node [style=none] (609) at (0, -3) {};
		\node [style=none] (610) at (1, -2) {};
		\node [style=none] (611) at (-9, -4) {};
		\node [style=none] (612) at (8, -4) {};
		\node [style=none] (613) at (-9, -6) {};
		\node [style=none] (614) at (8, -6) {};
		\node [style=none] (616) at (-8, -4) {};
		\node [style=none] (617) at (7, -4) {};
		\node [style=none] (618) at (7, -6) {};
		\node [style=none] (619) at (-8, -6) {};
		\node [style=none] (621) at (-3, -5) {};
		\node [style=none] (622) at (-4, -6) {};
		\node [style=none] (627) at (-8, -6.5) {};
		\node [style=none] (628) at (-8, -7.5) {};
		\node [style=none] (629) at (7, -6.5) {};
		\node [style=none] (630) at (7, -7.5) {};
		\node [style=label-small] (631) at (-0.5, -7) {$n$};
		\node [style=none] (632) at (-8, -7) {};
		\node [style=none] (633) at (7, -7) {};
		\node [style=none] (653) at (0, -3) {};
		\node [style=none] (655) at (0, -3) {};
		\node [style=none] (675) at (22, -6.5) {};
		\node [style=none] (676) at (22, -7.5) {};
		\node [style=label-small] (677) at (19.5, -7) {$m$};
		\node [style=none] (678) at (17, -7) {};
		\node [style=none] (679) at (22, -7) {};
		\node [style=none] (682) at (-4, -0.5) {};
		\node [style=none] (683) at (-4, -1.5) {};
		\node [style=none] (686) at (-4, -1) {};
		\node [style=none] (687) at (1, -0.5) {};
		\node [style=none] (688) at (1, -1.5) {};
		\node [style=label-small] (689) at (-1.5, -1) {$m$};
		\node [style=none] (690) at (-4, -1) {};
		\node [style=none] (691) at (1, -1) {};
		\node [style=none] (695) at (17, -6.5) {};
		\node [style=none] (696) at (17, -7.5) {};
		\node [style=none] (700) at (27, -3) {};
		\node [style=none] (705) at (22, -4) {};
		\node [style=none] (716) at (7, -5) {};
		\node [style=none] (717) at (7, -4) {};
		\node [style=none] (718) at (6, -5) {};
		\node [style=none] (720) at (13, -4) {};
		\node [style=none] (721) at (14, -3) {};
		\node [style=none] (722) at (13, -4) {};
		\node [style=none] (723) at (13, -3) {};
		\node [style=none] (725) at (-3, -5) {};
		\node [style=none] (726) at (6, -5) {};
		\node [style=none] (737) at (-3, -3) {};
		\node [style=none] (738) at (0, -3) {};
		\node [style=none] (744) at (17, -6) {};
		\node [style=none] (745) at (21, -5) {};
		\node [style=none] (746) at (18, -5) {};
		\node [style=none] (747) at (21, -5) {};
		\node [style=none] (752) at (14, -3) {};
		\node [style=none] (753) at (21, -3) {};
		\node [style=none] (754) at (17, -6) {};
		\node [style=none] (755) at (18, -5) {};
		\node [style=none] (756) at (18, -5) {};
		\node [style=none] (757) at (17, -6) {};
		\node [style=none] (758) at (17, -5) {};
		\node [style=none] (764) at (-8, -6) {};
		\node [style=none] (765) at (-7, -5) {};
		\node [style=none] (766) at (-8, -6) {};
		\node [style=none] (767) at (-8, -5) {};
		\node [style=none] (768) at (-7, -5) {};
		\node [style=none] (769) at (-5, -5) {};
		\node [style=none] (770) at (-4, -4) {};
		\node [style=none] (776) at (-4, -4) {};
		\node [style=none] (778) at (-5, -5) {};
		\node [style=none] (780) at (23, -3) {};
		\node [style=none] (781) at (27, -3) {};
		\node [style=none] (782) at (22, -2) {};
	\end{pgfonlayer}
	\begin{pgfonlayer}{edgelayer}
		\draw [style=fade] (531.center) to (532.center);
		\draw [style=fade] (533.center) to (534.center);
		\draw [style=shade] (543.center)
			 to (544.center)
			 to (541.center)
			 to (542.center)
			 to cycle;
		\draw [style=red, rounded corners=0.2cm] (562.center)
			 to (560.center)
			 to (561.center);
		\draw [style=red arrow] (564.center) to (563.center);
		\draw [style=red arrow] (564.center) to (565.center);
		\draw [style=fade] (566.center) to (567.center);
		\draw [style=fade] (568.center) to (569.center);
		\draw [style=red arrow] (570) to (571.center);
		\draw [style=red arrow] (570) to (572.center);
		\draw [style=fade] (578.center) to (577.center);
		\draw [style=shade] (581.center)
			 to (582.center)
			 to (579.center)
			 to (580.center)
			 to cycle;
		\draw [style=red, rounded corners=0.2cm] (592.center)
			 to (590.center)
			 to (591.center);
		\draw [style=faded] (531.center) to (533.center);
		\draw [style=faded] (534.center) to (532.center);
		\draw [style=faded] (577.center) to (575.center);
		\draw [style=faded] (578.center) to (576.center);
		\draw [style=fade] (598.center) to (597.center);
		\draw [style=fade] (595.center) to (596.center);
		\draw [style=shade] (599.center)
			 to (600.center)
			 to (601.center)
			 to (602.center)
			 to cycle;
		\draw [style=red, rounded corners=0.2cm] (605.center)
			 to (606.center)
			 to (604.center);
		\draw [style=red, rounded corners=0.2cm] (609.center)
			 to (607.center)
			 to (608.center);
		\draw [style=faded] (597.center) to (595.center);
		\draw [style=faded] (598.center) to (596.center);
		\draw [style=fade] (613.center) to (614.center);
		\draw [style=shade] (618.center)
			 to (619.center)
			 to (616.center)
			 to (617.center)
			 to cycle;
		\draw [style=faded] (611.center) to (613.center);
		\draw [style=faded] (614.center) to (612.center);
		\draw [style=fade] (627.center) to (628.center);
		\draw [style=fade] (629.center) to (630.center);
		\draw [style=red arrow] (631) to (632.center);
		\draw [style=red arrow] (631) to (633.center);
		\draw [style=fade] (675.center) to (676.center);
		\draw [style=red arrow] (677) to (678.center);
		\draw [style=red arrow] (677) to (679.center);
		\draw [style=fade] (682.center) to (683.center);
		\draw [style=fade] (687.center) to (688.center);
		\draw [style=red arrow] (689) to (690.center);
		\draw [style=red arrow] (689) to (691.center);
		\draw [style=fade] (695.center) to (696.center);
		\draw [style=red, rounded corners=0.2cm] (718.center)
			 to (716.center)
			 to (717.center);
		\draw [style=red, rounded corners=0.2cm] (722.center)
			 to (723.center)
			 to (721.center);
		\draw [style=red, dashed] (725.center) to (726.center);
		\draw [style=red, dashed] (737.center) to (738.center);
		\draw [style=red, dashed] (746.center) to (747.center);
		\draw [style=red, dashed] (752.center) to (753.center);
		\draw [style=red, rounded corners=0.2cm] (757.center)
			 to (758.center)
			 to (756.center);
		\draw [style=red, rounded corners=0.2cm] (766.center)
			 to (767.center)
			 to (765.center);
		\draw [style=red, dashed] (768.center) to (769.center);
		\draw [style=red] (778.center) to (725.center);
		\draw [style=red] (622.center) to (776.center);
		\draw [style=red, dashed] (780.center) to (781.center);
		\draw [style=red] (753.center) to (780.center);
		\draw [style=red] (782.center) to (705.center);
	\end{pgfonlayer}
\end{tikzpicture}
  }
  
  \centerline{
  \begin{tikzpicture}[scale = 0.4]
	\begin{pgfonlayer}{nodelayer}
		\node [style=none] (531) at (12, -2) {};
		\node [style=none] (532) at (29, -2) {};
		\node [style=none] (533) at (12, -4) {};
		\node [style=none] (534) at (29, -4) {};
		\node [style=none] (538) at (23, -4) {};
		\node [style=none] (541) at (13, -2) {};
		\node [style=none] (542) at (28, -2) {};
		\node [style=none] (543) at (28, -4) {};
		\node [style=none] (544) at (13, -4) {};
		\node [style=none] (560) at (28, -3) {};
		\node [style=none] (561) at (28, -2) {};
		\node [style=none] (562) at (27, -3) {};
		\node [style=none] (563) at (9, -4) {};
		\node [style=none] (564) at (10, -4) {};
		\node [style=none] (565) at (11, -4) {};
		\node [style=none] (566) at (13, -0.25) {};
		\node [style=none] (567) at (13, -1.25) {};
		\node [style=none] (568) at (28, -0.25) {};
		\node [style=none] (569) at (28, -1.25) {};
		\node [style=label-small] (570) at (20.5, -0.75) {$n$};
		\node [style=none] (571) at (13, -0.75) {};
		\node [style=none] (572) at (28, -0.75) {};
		\node [style=none] (575) at (12, -4) {};
		\node [style=none] (576) at (29, -4) {};
		\node [style=none] (577) at (12, -6) {};
		\node [style=none] (578) at (29, -6) {};
		\node [style=none] (579) at (18, -4) {};
		\node [style=none] (580) at (23, -4) {};
		\node [style=none] (581) at (23, -6) {};
		\node [style=none] (582) at (18, -6) {};
		\node [style=none] (590) at (23, -5) {};
		\node [style=none] (591) at (23, -4) {};
		\node [style=none] (592) at (22, -5) {};
		\node [style=none] (595) at (-9, -2) {};
		\node [style=none] (596) at (8, -2) {};
		\node [style=none] (597) at (-9, -4) {};
		\node [style=none] (598) at (8, -4) {};
		\node [style=none] (599) at (-4, -2) {};
		\node [style=none] (600) at (1, -2) {};
		\node [style=none] (601) at (1, -4) {};
		\node [style=none] (602) at (-4, -4) {};
		\node [style=none] (603) at (1, -2) {};
		\node [style=none] (604) at (-3, -3) {};
		\node [style=none] (605) at (-4, -4) {};
		\node [style=none] (606) at (-4, -3) {};
		\node [style=none] (607) at (1, -3) {};
		\node [style=none] (608) at (1, -2) {};
		\node [style=none] (609) at (0, -3) {};
		\node [style=none] (610) at (1, -2) {};
		\node [style=none] (611) at (-9, -4) {};
		\node [style=none] (612) at (8, -4) {};
		\node [style=none] (613) at (-9, -6) {};
		\node [style=none] (614) at (8, -6) {};
		\node [style=none] (616) at (-8, -4) {};
		\node [style=none] (617) at (6, -4) {};
		\node [style=none] (618) at (6, -6) {};
		\node [style=none] (619) at (-8, -6) {};
		\node [style=none] (621) at (-3, -5) {};
		\node [style=none] (622) at (-4, -6) {};
		\node [style=none] (627) at (-8, -6.5) {};
		\node [style=none] (628) at (-8, -7.5) {};
		\node [style=none] (629) at (6, -6.5) {};
		\node [style=none] (630) at (6, -7.5) {};
		\node [style=label-small] (631) at (-1.5, -7) {$n$};
		\node [style=none] (632) at (-8, -7) {};
		\node [style=none] (633) at (6, -7) {};
		\node [style=none] (653) at (0, -3) {};
		\node [style=none] (655) at (0, -3) {};
		\node [style=none] (675) at (23, -6.5) {};
		\node [style=none] (676) at (23, -7.5) {};
		\node [style=label-small] (677) at (20.5, -7) {$m$};
		\node [style=none] (678) at (18, -7) {};
		\node [style=none] (679) at (23, -7) {};
		\node [style=none] (682) at (-4, -0.5) {};
		\node [style=none] (683) at (-4, -1.5) {};
		\node [style=none] (686) at (-4, -1) {};
		\node [style=none] (687) at (1, -0.5) {};
		\node [style=none] (688) at (1, -1.5) {};
		\node [style=label-small] (689) at (-1.5, -1) {$m$};
		\node [style=none] (690) at (-4, -1) {};
		\node [style=none] (691) at (1, -1) {};
		\node [style=none] (695) at (18, -6.5) {};
		\node [style=none] (696) at (18, -7.5) {};
		\node [style=none] (700) at (24, -3) {};
		\node [style=none] (705) at (23, -4) {};
		\node [style=none] (716) at (6, -5) {};
		\node [style=none] (717) at (6, -4) {};
		\node [style=none] (718) at (5, -5) {};
		\node [style=none] (720) at (13, -4) {};
		\node [style=none] (721) at (14, -3) {};
		\node [style=none] (722) at (13, -4) {};
		\node [style=none] (723) at (13, -3) {};
		\node [style=none] (725) at (-3, -5) {};
		\node [style=none] (726) at (5, -5) {};
		\node [style=none] (737) at (-3, -3) {};
		\node [style=none] (738) at (0, -3) {};
		\node [style=none] (744) at (18, -6) {};
		\node [style=none] (745) at (22, -5) {};
		\node [style=none] (746) at (19, -5) {};
		\node [style=none] (747) at (22, -5) {};
		\node [style=none] (752) at (14, -3) {};
		\node [style=none] (753) at (22, -3) {};
		\node [style=none] (754) at (18, -6) {};
		\node [style=none] (755) at (19, -5) {};
		\node [style=none] (756) at (19, -5) {};
		\node [style=none] (757) at (18, -6) {};
		\node [style=none] (758) at (18, -5) {};
		\node [style=none] (764) at (-8, -6) {};
		\node [style=none] (765) at (-7, -5) {};
		\node [style=none] (766) at (-8, -6) {};
		\node [style=none] (767) at (-8, -5) {};
		\node [style=none] (768) at (-7, -5) {};
		\node [style=none] (769) at (-5, -5) {};
		\node [style=none] (770) at (-4, -4) {};
		\node [style=none] (776) at (-4, -4) {};
		\node [style=none] (778) at (-5, -5) {};
		\node [style=none] (779) at (23, -8.5) {};
		\node [style=none] (780) at (24, -3) {};
		\node [style=none] (781) at (27, -3) {};
	\end{pgfonlayer}
	\begin{pgfonlayer}{edgelayer}
		\draw [style=fade] (531.center) to (532.center);
		\draw [style=fade] (533.center) to (534.center);
		\draw [style=shade] (543.center)
			 to (544.center)
			 to (541.center)
			 to (542.center)
			 to cycle;
		\draw [style=red, rounded corners=0.2cm] (562.center)
			 to (560.center)
			 to (561.center);
		\draw [style=red arrow] (564.center) to (563.center);
		\draw [style=red arrow] (564.center) to (565.center);
		\draw [style=fade] (566.center) to (567.center);
		\draw [style=fade] (568.center) to (569.center);
		\draw [style=red arrow] (570) to (571.center);
		\draw [style=red arrow] (570) to (572.center);
		\draw [style=fade] (578.center) to (577.center);
		\draw [style=shade] (581.center)
			 to (582.center)
			 to (579.center)
			 to (580.center)
			 to cycle;
		\draw [style=red, rounded corners=0.2cm] (592.center)
			 to (590.center)
			 to (591.center);
		\draw [style=faded] (531.center) to (533.center);
		\draw [style=faded] (534.center) to (532.center);
		\draw [style=faded] (577.center) to (575.center);
		\draw [style=faded] (578.center) to (576.center);
		\draw [style=fade] (598.center) to (597.center);
		\draw [style=fade] (595.center) to (596.center);
		\draw [style=shade] (599.center)
			 to (600.center)
			 to (601.center)
			 to (602.center)
			 to cycle;
		\draw [style=red, rounded corners=0.2cm] (605.center)
			 to (606.center)
			 to (604.center);
		\draw [style=red, rounded corners=0.2cm] (609.center)
			 to (607.center)
			 to (608.center);
		\draw [style=faded] (597.center) to (595.center);
		\draw [style=faded] (598.center) to (596.center);
		\draw [style=fade] (613.center) to (614.center);
		\draw [style=shade] (618.center)
			 to (619.center)
			 to (616.center)
			 to (617.center)
			 to cycle;
		\draw [style=faded] (611.center) to (613.center);
		\draw [style=faded] (614.center) to (612.center);
		\draw [style=fade] (627.center) to (628.center);
		\draw [style=fade] (629.center) to (630.center);
		\draw [style=red arrow] (631) to (632.center);
		\draw [style=red arrow] (631) to (633.center);
		\draw [style=fade] (675.center) to (676.center);
		\draw [style=red arrow] (677) to (678.center);
		\draw [style=red arrow] (677) to (679.center);
		\draw [style=fade] (682.center) to (683.center);
		\draw [style=fade] (687.center) to (688.center);
		\draw [style=red arrow] (689) to (690.center);
		\draw [style=red arrow] (689) to (691.center);
		\draw [style=fade] (695.center) to (696.center);
		\draw [style=red, rounded corners=0.2cm] (718.center)
			 to (716.center)
			 to (717.center);
		\draw [style=red, rounded corners=0.2cm] (722.center)
			 to (723.center)
			 to (721.center);
		\draw [style=red, dashed] (725.center) to (726.center);
		\draw [style=red, dashed] (737.center) to (738.center);
		\draw [style=red, dashed] (746.center) to (747.center);
		\draw [style=red, dashed] (752.center) to (753.center);
		\draw [style=red, rounded corners=0.2cm] (757.center)
			 to (758.center)
			 to (756.center);
		\draw [style=red, rounded corners=0.2cm] (766.center)
			 to (767.center)
			 to (765.center);
		\draw [style=red, dashed] (768.center) to (769.center);
		\draw [style=red] (778.center) to (725.center);
		\draw [style=red, dashed] (780.center) to (781.center);
		\draw [style=red] (753.center) to (780.center);
	\end{pgfonlayer}
\end{tikzpicture}
  }
  
  This completes the proof.
  \end{proof}
 
The proof of the other commutation relations will be done in the same spirit.

\begin{prop}\label{prop: alpha_neg_commute}
	For any $m,n>0$ and $k$, we have $[\alpha_{-m,k}^\fff,\alpha_{-n,k}^\fff]=0$.
\end{prop}
\begin{proof}
The proof is similar to the previous case. The expansion of $\alpha_{-m,k} \alpha_{-n,k}|w\rangle$ is given by set of $\Delta_k$-tilings of two spectral parameters.
	 Regarding $\alpha_{-m,k}\alpha_{-n,k}|w\rangle$, and similarly $\alpha_{-n,k}\alpha_{-m,k}|w\rangle$, as sets of tilings, we will define an involution $\rho$, which maps an element of $\alpha_{-m,k}\alpha_{-n,k}|w\rangle$ to either (i) an element in $\alpha_{-n,k}\alpha_{-m,k}|w\rangle$ of the same sign or (ii) an element in $\alpha_{-m,k}\alpha_{-n,k}|w\rangle$ of opposite sign.  When the tilings have one layer, then the construction of $\rho$ can be describe exactly the same was as in the proof of \Cref{prop: alpha_pos_commute}, flipping all the diagrams upside down. We list these diagrams in the appendix.
	 
	 In the presence of multiple layers, we shall focus on the left-most (or right-most) path and apply the usual operation as if they were one layer. A few new operations are needed, which we listed as follows.
	 
	 
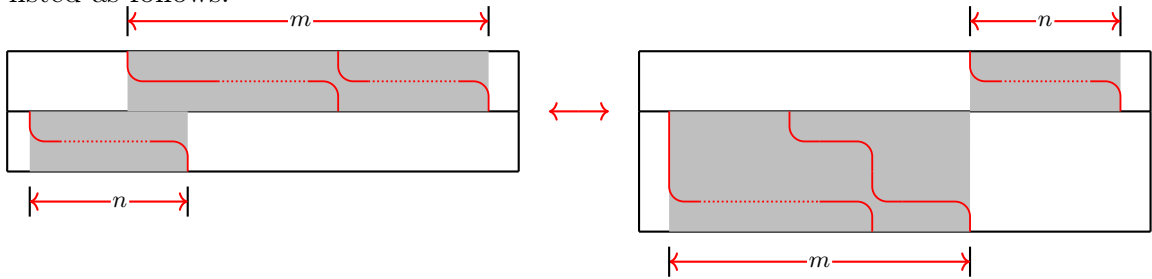
\centerline{\begin{tikzpicture}[xscale=0.4,yscale = - 0.4]
	\begin{pgfonlayer}{nodelayer}
		\node [style=none] (531) at (12, 0) {};
		\node [style=none] (532) at (29, 0) {};
		\node [style=none] (533) at (12, -4) {};
		\node [style=none] (534) at (29, -4) {};
		\node [style=none] (538) at (23, -6) {};
		\node [style=none] (541) at (23, -4) {};
		\node [style=none] (542) at (28, -4) {};
		\node [style=none] (543) at (28, -6) {};
		\node [style=none] (544) at (23, -6) {};
		\node [style=none] (555) at (23, -6) {};
		\node [style=none] (560) at (28, -5) {};
		\node [style=none] (561) at (28, -4) {};
		\node [style=none] (562) at (27, -5) {};
		\node [style=none] (563) at (9, -4) {};
		\node [style=none] (564) at (10, -4) {};
		\node [style=none] (565) at (11, -4) {};
		\node [style=none] (566) at (23, -6.5) {};
		\node [style=none] (567) at (23, -7.5) {};
		\node [style=none] (568) at (28, -6.5) {};
		\node [style=none] (569) at (28, -7.5) {};
		\node [style=label-small] (570) at (25.5, -7) {$n$};
		\node [style=none] (571) at (23, -7) {};
		\node [style=none] (572) at (28, -7) {};
		\node [style=none] (575) at (12, -4) {};
		\node [style=none] (576) at (29, -4) {};
		\node [style=none] (577) at (12, -6) {};
		\node [style=none] (578) at (29, -6) {};
		\node [style=none] (579) at (13, 0) {};
		\node [style=none] (580) at (23, 0) {};
		\node [style=none] (581) at (23, -4) {};
		\node [style=none] (582) at (13, -4) {};
		\node [style=none] (587) at (14, -1) {};
		\node [style=none] (588) at (13, -4) {};
		\node [style=none] (589) at (13, -1) {};
		\node [style=none] (590) at (19.75, -1) {};
		\node [style=none] (591) at (19.75, 0) {};
		\node [style=none] (592) at (18, -1) {};
		\node [style=none] (595) at (-9, -2) {};
		\node [style=none] (596) at (8, -2) {};
		\node [style=none] (597) at (-9, -4) {};
		\node [style=none] (598) at (8, -4) {};
		\node [style=none] (611) at (-9, -4) {};
		\node [style=none] (612) at (8, -4) {};
		\node [style=none] (613) at (-9, -6) {};
		\node [style=none] (614) at (8, -6) {};
		\node [style=none] (627) at (-8.25, -0.5) {};
		\node [style=none] (628) at (-8.25, -1.5) {};
		\node [style=none] (629) at (-3, -0.5) {};
		\node [style=none] (630) at (-3, -1.5) {};
		\node [style=label-small] (631) at (-5.25, -1) {$n$};
		\node [style=none] (632) at (-8.25, -1) {};
		\node [style=none] (633) at (-3, -1) {};
		\node [style=none] (643) at (14, -1) {};
		\node [style=none] (646) at (18, -1) {};
		\node [style=none] (675) at (23, 1.5) {};
		\node [style=none] (676) at (23, 0.5) {};
		\node [style=label-small] (677) at (18, 1) {$m$};
		\node [style=none] (678) at (13, 1) {};
		\node [style=none] (679) at (23, 1) {};
		\node [style=none] (682) at (-5, -6.5) {};
		\node [style=none] (683) at (-5, -7.5) {};
		\node [style=none] (686) at (-5, -7) {};
		\node [style=none] (687) at (7, -6.5) {};
		\node [style=none] (688) at (7, -7.5) {};
		\node [style=label-small] (689) at (0.75, -7) {$m$};
		\node [style=none] (690) at (-5, -7) {};
		\node [style=none] (691) at (7, -7) {};
		\node [style=none] (695) at (13, 1.5) {};
		\node [style=none] (696) at (13, 0.5) {};
		\node [style=none] (697) at (24, -5) {};
		\node [style=none] (700) at (27, -5) {};
		\node [style=none] (704) at (24, -5) {};
		\node [style=none] (705) at (23, -6) {};
		\node [style=none] (706) at (23, -5) {};
		\node [style=none] (713) at (19.75, -2) {};
		\node [style=none] (714) at (19.75, -3) {};
		\node [style=none] (715) at (19.75, -2) {};
		\node [style=none] (716) at (18, -3) {};
		\node [style=none] (719) at (20.75, -1) {};
		\node [style=none] (720) at (19.75, -2) {};
		\node [style=none] (721) at (19.75, -1) {};
		\node [style=none] (722) at (20.75, -1) {};
		\node [style=none] (728) at (23, 0) {};
		\node [style=none] (729) at (23, -1) {};
		\node [style=none] (730) at (23, 0) {};
		\node [style=none] (731) at (20.75, -1) {};
		\node [style=none] (738) at (18, -3) {};
		\node [style=none] (739) at (17, -4) {};
		\node [style=none] (740) at (17, -3) {};
		\node [style=none] (741) at (18, -3) {};
		\node [style=none] (742) at (-3, -2) {};
		\node [style=none] (743) at (-3, -2) {};
		\node [style=none] (744) at (-8.25, -4) {};
		\node [style=none] (745) at (-8.25, -2) {};
		\node [style=none] (746) at (-3, -2) {};
		\node [style=none] (747) at (-3, -4) {};
		\node [style=none] (748) at (-8.25, -4) {};
		\node [style=none] (749) at (-8.25, -4) {};
		\node [style=none] (750) at (-7.25, -3) {};
		\node [style=none] (751) at (-8.25, -4) {};
		\node [style=none] (752) at (-8.25, -3) {};
		\node [style=none] (753) at (-3, -3) {};
		\node [style=none] (754) at (-3, -2) {};
		\node [style=none] (755) at (-4.25, -3) {};
		\node [style=none] (756) at (-5, -4) {};
		\node [style=none] (757) at (7, -4) {};
		\node [style=none] (758) at (7, -6) {};
		\node [style=none] (759) at (-5, -6) {};
		\node [style=none] (760) at (7, -4) {};
		\node [style=none] (761) at (-2, -5) {};
		\node [style=none] (762) at (-5, -6) {};
		\node [style=none] (763) at (-5, -5) {};
		\node [style=none] (764) at (7, -5) {};
		\node [style=none] (765) at (7, -4) {};
		\node [style=none] (766) at (6, -5) {};
		\node [style=none] (767) at (7, -4) {};
		\node [style=none] (770) at (-2, -5) {};
		\node [style=none] (771) at (3, -5) {};
		\node [style=none] (772) at (1, -5) {};
		\node [style=none] (773) at (3, -5) {};
		\node [style=none] (774) at (6, -5) {};
		\node [style=none] (775) at (2, -4) {};
		\node [style=none] (776) at (2, -6) {};
		\node [style=none] (777) at (2, -6) {};
		\node [style=none] (778) at (2, -5) {};
		\node [style=none] (779) at (2, -4) {};
		\node [style=none] (780) at (1, -5) {};
		\node [style=none] (781) at (3, -5) {};
		\node [style=none] (782) at (2, -6) {};
		\node [style=none] (783) at (2, -5) {};
	\end{pgfonlayer}
	\begin{pgfonlayer}{edgelayer}
		\draw [style=fade] (531.center) to (532.center);
		\draw [style=fade] (533.center) to (534.center);
		\draw [style=shade] (543.center)
			 to (544.center)
			 to (541.center)
			 to (542.center)
			 to cycle;
		\draw [style=red, rounded corners=0.2cm] (562.center)
			 to (560.center)
			 to (561.center);
		\draw [style=red arrow] (564.center) to (563.center);
		\draw [style=red arrow] (564.center) to (565.center);
		\draw [style=fade] (566.center) to (567.center);
		\draw [style=fade] (568.center) to (569.center);
		\draw [style=red arrow] (570) to (571.center);
		\draw [style=red arrow] (570) to (572.center);
		\draw [style=fade] (578.center) to (577.center);
		\draw [style=shade] (581.center)
			 to (582.center)
			 to (579.center)
			 to (580.center)
			 to cycle;
		\draw [style=red, rounded corners=0.2cm] (588.center)
			 to (589.center)
			 to (587.center);
		\draw [style=red, rounded corners=0.2cm] (592.center)
			 to (590.center)
			 to (591.center);
		\draw [style=faded] (531.center) to (533.center);
		\draw [style=faded] (534.center) to (532.center);
		\draw [style=faded] (577.center) to (575.center);
		\draw [style=faded] (578.center) to (576.center);
		\draw [style=fade] (598.center) to (597.center);
		\draw [style=fade] (595.center) to (596.center);
		\draw [style=faded] (597.center) to (595.center);
		\draw [style=faded] (598.center) to (596.center);
		\draw [style=fade] (613.center) to (614.center);
		\draw [style=faded] (611.center) to (613.center);
		\draw [style=faded] (614.center) to (612.center);
		\draw [style=fade] (627.center) to (628.center);
		\draw [style=fade] (629.center) to (630.center);
		\draw [style=red arrow] (631) to (632.center);
		\draw [style=red arrow] (631) to (633.center);
		\draw [style=red, dashed] (643.center) to (646.center);
		\draw [style=fade] (675.center) to (676.center);
		\draw [style=red arrow] (677) to (678.center);
		\draw [style=red arrow] (677) to (679.center);
		\draw [style=fade] (682.center) to (683.center);
		\draw [style=fade] (687.center) to (688.center);
		\draw [style=red arrow] (689) to (690.center);
		\draw [style=red arrow] (689) to (691.center);
		\draw [style=fade] (695.center) to (696.center);
		\draw [style=red, dashed] (697.center) to (700.center);
		\draw [style=red, rounded corners=0.2cm] (705.center)
			 to (706.center)
			 to (704.center);
		\draw [style=red, rounded corners=0.2cm] (716.center)
			 to (714.center)
			 to (715.center);
		\draw [style=red, rounded corners=0.2cm] (720.center)
			 to (721.center)
			 to (719.center);
		\draw [style=red, rounded corners=0.2cm] (731.center)
			 to (729.center)
			 to (730.center);
		\draw [style=red, rounded corners=0.2cm] (739.center)
			 to (740.center)
			 to (738.center);
		\draw [style=shade] (747.center)
			 to (748.center)
			 to (745.center)
			 to (746.center)
			 to cycle;
		\draw [style=red, rounded corners=0.2cm] (751.center)
			 to (752.center)
			 to (750.center);
		\draw [style=red, rounded corners=0.2cm] (755.center)
			 to (753.center)
			 to (754.center);
		\draw [style=red, dashed] (750.center) to (755.center);
		\draw [style=shade] (756.center)
			 to (757.center)
			 to (758.center)
			 to (759.center)
			 to cycle;
		\draw [style=red, rounded corners=0.2cm] (762.center)
			 to (763.center)
			 to (761.center);
		\draw [style=red, rounded corners=0.2cm] (766.center)
			 to (764.center)
			 to (765.center);
		\draw [style=red, dashed] (770.center) to (772.center);
		\draw [style=red, dashed] (773.center) to (774.center);
		\draw [style=red, rounded corners=0.2cm] (780.center)
			 to (778.center)
			 to (779.center);
		\draw [style=red, rounded corners=0.2cm] (782.center)
			 to (783.center)
			 to (781.center);
	\end{pgfonlayer}
\end{tikzpicture}}
	 
	 
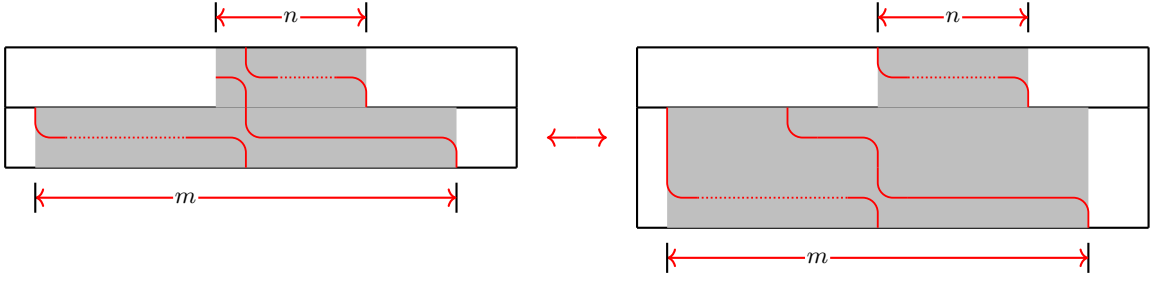
\centerline{\begin{tikzpicture}[xscale=0.4,yscale = -0.4]
	\begin{pgfonlayer}{nodelayer}
		\node [style=none] (531) at (12, 0) {};
		\node [style=none] (532) at (29, 0) {};
		\node [style=none] (533) at (12, -4) {};
		\node [style=none] (534) at (29, -4) {};
		\node [style=none] (538) at (20, -6) {};
		\node [style=none] (541) at (20, -4) {};
		\node [style=none] (542) at (25, -4) {};
		\node [style=none] (543) at (25, -6) {};
		\node [style=none] (544) at (20, -6) {};
		\node [style=none] (555) at (20, -6) {};
		\node [style=none] (560) at (25, -5) {};
		\node [style=none] (561) at (25, -4) {};
		\node [style=none] (562) at (24, -5) {};
		\node [style=none] (563) at (9, -3) {};
		\node [style=none] (564) at (10, -3) {};
		\node [style=none] (565) at (11, -3) {};
		\node [style=none] (566) at (20, -6.5) {};
		\node [style=none] (567) at (20, -7.5) {};
		\node [style=none] (568) at (25, -6.5) {};
		\node [style=none] (569) at (25, -7.5) {};
		\node [style=label-small] (570) at (22.5, -7) {$n$};
		\node [style=none] (571) at (20, -7) {};
		\node [style=none] (572) at (25, -7) {};
		\node [style=none] (575) at (12, -4) {};
		\node [style=none] (576) at (29, -4) {};
		\node [style=none] (577) at (12, -6) {};
		\node [style=none] (578) at (29, -6) {};
		\node [style=none] (579) at (13, 0) {};
		\node [style=none] (580) at (27, 0) {};
		\node [style=none] (581) at (27, -4) {};
		\node [style=none] (582) at (13, -4) {};
		\node [style=none] (587) at (14, -1) {};
		\node [style=none] (588) at (13, -4) {};
		\node [style=none] (589) at (13, -1) {};
		\node [style=none] (590) at (20, -1) {};
		\node [style=none] (591) at (20, 0) {};
		\node [style=none] (592) at (19, -1) {};
		\node [style=none] (643) at (14, -1) {};
		\node [style=none] (646) at (19, -1) {};
		\node [style=none] (675) at (27, 1.5) {};
		\node [style=none] (676) at (27, 0.5) {};
		\node [style=label-small] (677) at (18, 1) {$m$};
		\node [style=none] (678) at (13, 1) {};
		\node [style=none] (679) at (27, 1) {};
		\node [style=none] (695) at (13, 1.5) {};
		\node [style=none] (696) at (13, 0.5) {};
		\node [style=none] (697) at (21, -5) {};
		\node [style=none] (700) at (24, -5) {};
		\node [style=none] (704) at (21, -5) {};
		\node [style=none] (705) at (20, -6) {};
		\node [style=none] (706) at (20, -5) {};
		\node [style=none] (713) at (20, -2) {};
		\node [style=none] (714) at (20, -3) {};
		\node [style=none] (715) at (20, -2) {};
		\node [style=none] (716) at (18, -3) {};
		\node [style=none] (719) at (21, -1) {};
		\node [style=none] (720) at (20, -2) {};
		\node [style=none] (721) at (20, -1) {};
		\node [style=none] (722) at (21, -1) {};
		\node [style=none] (729) at (27, -1) {};
		\node [style=none] (730) at (27, 0) {};
		\node [style=none] (731) at (21, -1) {};
		\node [style=none] (738) at (18, -3) {};
		\node [style=none] (739) at (17, -4) {};
		\node [style=none] (740) at (17, -3) {};
		\node [style=none] (741) at (18, -3) {};
		\node [style=none] (742) at (-9, -2) {};
		\node [style=none] (743) at (8, -2) {};
		\node [style=none] (744) at (-9, -4) {};
		\node [style=none] (745) at (8, -4) {};
		\node [style=none] (746) at (-2, -6) {};
		\node [style=none] (747) at (-2, -4) {};
		\node [style=none] (748) at (3, -4) {};
		\node [style=none] (749) at (3, -6) {};
		\node [style=none] (750) at (-2, -6) {};
		\node [style=none] (751) at (-2, -6) {};
		\node [style=none] (752) at (3, -5) {};
		\node [style=none] (753) at (3, -4) {};
		\node [style=none] (754) at (2, -5) {};
		\node [style=none] (755) at (-2, -6.5) {};
		\node [style=none] (756) at (-2, -7.5) {};
		\node [style=none] (757) at (3, -6.5) {};
		\node [style=none] (758) at (3, -7.5) {};
		\node [style=label-small] (759) at (0.5, -7) {$n$};
		\node [style=none] (760) at (-2, -7) {};
		\node [style=none] (761) at (3, -7) {};
		\node [style=none] (762) at (-9, -4) {};
		\node [style=none] (763) at (8, -4) {};
		\node [style=none] (764) at (-9, -6) {};
		\node [style=none] (765) at (8, -6) {};
		\node [style=none] (766) at (-8, -2) {};
		\node [style=none] (767) at (6, -2) {};
		\node [style=none] (768) at (6, -4) {};
		\node [style=none] (769) at (-8, -4) {};
		\node [style=none] (770) at (-7, -3) {};
		\node [style=none] (771) at (-8, -4) {};
		\node [style=none] (772) at (-8, -3) {};
		\node [style=none] (773) at (-1, -3) {};
		\node [style=none] (774) at (-1, -2) {};
		\node [style=none] (775) at (-3, -3) {};
		\node [style=none] (776) at (-7, -3) {};
		\node [style=none] (777) at (-3, -3) {};
		\node [style=none] (778) at (6, -0.5) {};
		\node [style=none] (779) at (6, -1.5) {};
		\node [style=label-small] (780) at (-3, -1) {$m$};
		\node [style=none] (781) at (-8, -1) {};
		\node [style=none] (782) at (6, -1) {};
		\node [style=none] (783) at (-8, -0.5) {};
		\node [style=none] (784) at (-8, -1.5) {};
		\node [style=none] (785) at (0, -5) {};
		\node [style=none] (786) at (2, -5) {};
		\node [style=none] (787) at (0, -5) {};
		\node [style=none] (788) at (-1, -6) {};
		\node [style=none] (789) at (-1, -5) {};
		\node [style=none] (790) at (-2, -4) {};
		\node [style=none] (794) at (0, -3) {};
		\node [style=none] (795) at (-1, -4) {};
		\node [style=none] (796) at (-1, -3) {};
		\node [style=none] (797) at (0, -3) {};
		\node [style=none] (798) at (6, -3) {};
		\node [style=none] (799) at (6, -2) {};
		\node [style=none] (800) at (0, -3) {};
		\node [style=none] (804) at (-2, -3) {};
		\node [style=none] (809) at (-1, -4) {};
		\node [style=none] (810) at (-2, -5) {};
		\node [style=none] (811) at (-1, -5) {};
	\end{pgfonlayer}
	\begin{pgfonlayer}{edgelayer}
		\draw [style=fade] (531.center) to (532.center);
		\draw [style=fade] (533.center) to (534.center);
		\draw [style=shade] (543.center)
			 to (544.center)
			 to (541.center)
			 to (542.center)
			 to cycle;
		\draw [style=red, rounded corners=0.2cm] (562.center)
			 to (560.center)
			 to (561.center);
		\draw [style=red arrow] (564.center) to (563.center);
		\draw [style=red arrow] (564.center) to (565.center);
		\draw [style=fade] (566.center) to (567.center);
		\draw [style=fade] (568.center) to (569.center);
		\draw [style=red arrow] (570) to (571.center);
		\draw [style=red arrow] (570) to (572.center);
		\draw [style=fade] (578.center) to (577.center);
		\draw [style=shade] (581.center)
			 to (582.center)
			 to (579.center)
			 to (580.center)
			 to cycle;
		\draw [style=red, rounded corners=0.2cm] (588.center)
			 to (589.center)
			 to (587.center);
		\draw [style=red, rounded corners=0.2cm] (592.center)
			 to (590.center)
			 to (591.center);
		\draw [style=faded] (531.center) to (533.center);
		\draw [style=faded] (534.center) to (532.center);
		\draw [style=faded] (577.center) to (575.center);
		\draw [style=faded] (578.center) to (576.center);
		\draw [style=red, dashed] (643.center) to (646.center);
		\draw [style=fade] (675.center) to (676.center);
		\draw [style=red arrow] (677) to (678.center);
		\draw [style=red arrow] (677) to (679.center);
		\draw [style=fade] (695.center) to (696.center);
		\draw [style=red, dashed] (697.center) to (700.center);
		\draw [style=red, rounded corners=0.2cm] (705.center)
			 to (706.center)
			 to (704.center);
		\draw [style=red, rounded corners=0.2cm] (716.center)
			 to (714.center)
			 to (715.center);
		\draw [style=red, rounded corners=0.2cm] (720.center)
			 to (721.center)
			 to (719.center);
		\draw [style=red, rounded corners=0.2cm] (731.center)
			 to (729.center)
			 to (730.center);
		\draw [style=red, rounded corners=0.2cm] (739.center)
			 to (740.center)
			 to (738.center);
		\draw [style=fade] (742.center) to (743.center);
		\draw [style=fade] (744.center) to (745.center);
		\draw [style=shade] (749.center)
			 to (750.center)
			 to (747.center)
			 to (748.center)
			 to cycle;
		\draw [style=red, rounded corners=0.2cm] (754.center)
			 to (752.center)
			 to (753.center);
		\draw [style=fade] (755.center) to (756.center);
		\draw [style=fade] (757.center) to (758.center);
		\draw [style=red arrow] (759) to (760.center);
		\draw [style=red arrow] (759) to (761.center);
		\draw [style=fade] (765.center) to (764.center);
		\draw [style=shade] (768.center)
			 to (769.center)
			 to (766.center)
			 to (767.center)
			 to cycle;
		\draw [style=red, rounded corners=0.2cm] (771.center)
			 to (772.center)
			 to (770.center);
		\draw [style=red, rounded corners=0.2cm] (775.center)
			 to (773.center)
			 to (774.center);
		\draw [style=faded] (742.center) to (744.center);
		\draw [style=faded] (745.center) to (743.center);
		\draw [style=faded] (764.center) to (762.center);
		\draw [style=faded] (765.center) to (763.center);
		\draw [style=red, dashed] (776.center) to (777.center);
		\draw [style=fade] (778.center) to (779.center);
		\draw [style=red arrow] (780) to (781.center);
		\draw [style=red arrow] (780) to (782.center);
		\draw [style=fade] (783.center) to (784.center);
		\draw [style=red, dashed] (785.center) to (786.center);
		\draw [style=red, rounded corners=0.2cm] (788.center)
			 to (789.center)
			 to (787.center);
		\draw [style=red, rounded corners=0.2cm] (795.center)
			 to (796.center)
			 to (794.center);
		\draw [style=red, rounded corners=0.2cm] (800.center)
			 to (798.center)
			 to (799.center);
		\draw [style=red, rounded corners=0.2cm] (810.center)
			 to (811.center)
			 to (809.center);
	\end{pgfonlayer}
\end{tikzpicture}}
	 
	 The cases where both tows have multiple layers can be handled analogously by symmetry.
	 \end{proof}
	 
Before proving the next commutation relations, we first state a useful lemma. For $w,u\in S_{\zz}$, we say that $w/u$ is a \emph{primitive} $k$-strong-ribbon of size $r$ if $w/u$ is a $k$-strong ribbon and $u\cdot w^{-1}$ is a weak ribbon.
\begin{lemma}\label{lem:UDDUfixedpoints}
	For any $k\in\zz$ and $w\in S_{\zz}$, we have that
	\[\#\{u :w/u \text{ is primitive, } \ell(u)-\ell(w)=r \} - \#\{u :u/w \text{ is primitive, } \ell(w)-\ell(u)=r \} = r\]
	\end{lemma}
\begin{proof}
Let $w/u$ be a primitive $k$-strong ribbon. Then there exists a one row $\Delta_k$-tiling with boundaries labeled by $(w,u)$. But at the same time, there must also exist a one row $\nabla$-tiling with boundaries labeled by $(u,w)$, which is identical to the said $\Delta_k$-tiling up to a reflection.
The possible $\Delta_k$-tilings with such properties must not use tiles of types $5,6$ and thus contain only one single layer. 

Now identify numbers $\leq k$ with particles and identify numbers $>k$ with holes, we may realize a permutation $w$ as the partition $\pi_k(w)$\footnote{Recall that $\pi_k(w)$ is the partition obtained by projecting $w$ onto the $k$-Grassmanian}. Disallowing tiles of types $5$ and $6$, the $\Delta_k$-tilings exactly matches the operation of adding a ribbon on $\pi_k(w)$. This reduces the question to the case of partitions, which is a well known result (see e.g. \cite{lam2010quantized}).
\end{proof}

\begin{prop}
	For any $m,n>0$ and $k$, we have $[\alpha_{n}^\fff,\alpha_{-m,k}^\fff]=n \delta_{m,n}$. 
\end{prop}
\begin{proof}
	Similar to the previous proofs, we consider $\alpha_{n}\alpha_{-m,k}|w\rangle $ and $ \alpha_{-m,k}\alpha_{n}|w\rangle$ as two-row tilings (using $\nabla$ and $\Delta_k$ for each row), we will construct an involution $\rho$ on $\alpha_{n}\alpha_{-m,k}|w\rangle \cup \alpha_{-m,k}\alpha_{n}|w\rangle$. The operation $\rho$ either sends an element of $\alpha_{n}\alpha_{-m,k}|w\rangle$ to an element of $\alpha_{-m,k}\alpha_{n}|w\rangle$ of the same weight, or to an element of the same set with opposite weight. For readability, we describe the operation in the \Cref{sec:appendix}.
	
	When $m\neq n$, the involution $\rho$ has no fixed point, therefore $[\alpha_{n}^\fff,\alpha_{-m,k}^\fff]=0$. On the other hand, when $m=n$, the fixed points of $\rho$ are exactly the primitive $k$ strong ribbons, which are enumerated by \Cref{lem:UDDUfixedpoints}. Thus $[\alpha_{n}^\fff,\alpha_{-m,k}^\fff]=n$.
\end{proof}

\begin{corollary}
	For any choice of $k\in \mathbb{N}$, the operators $\alpha_{-m,k}^\fff$ and $\alpha_{n}^\fff$ defines an $\mathfrak{h}$ action on $\fff$.
\end{corollary}

We are now ready to prove \Cref{thm:stanley_hamiltonian}.
\begin{proof}[Proof of \Cref{thm:stanley_hamiltonian}]
	\Cref{thm:stanley_hamiltonian_one_row} shows that $e^{\phi(x)}|w\rangle = T(x)|w\rangle $, thus $$\langle u|\cdots e^{\phi(x_2)}e^{\phi(x_1)}|w\rangle = \langle u| \cdots T(x_2)T(x_1)|w\rangle,$$ which equals to $\mf_{w/u}(x_\infty)$ by \Cref{thm:stanley_transfer_matrix}.
	
	The second part of the theorem says that the map sending $\tilde\Phi:|w\rangle \mapsto \mf_{w} = \langle \id |e^{H(t)}|w\rangle$ preserves the $\mathfrak{h}$-actions. For this we follow the exact same proof of \cite[Theorem 5.1]{jimbo1983solitons}.\footnote{The proof of \cite[Theorem 5.1]{jimbo1983solitons} doesn't depend on the actual $\mathfrak h$-action, therefore can be generalized for any $\mathfrak h$-modules. This kind of generalization was also done in \cite{lam2006combinatorial} through different method.}
	
	Fix any $k\in\mathbb{N}$. By commutativity of $\alpha_m^\fff$ among themselves, we have
   $$\alpha_m^{\bb}\langle \id |e^{H(t)}|w\rangle ={\partial\over \partial t_i}\langle \id |e^{H(t)}|w\rangle = \langle \id|{\partial\over \partial t_i}e^{H(t)}|w\rangle =\langle \id|e^{H(t)}\alpha_m^\fff|w\rangle.$$
   
   On the other hand, we have
   \[\langle \id| e^{H(t)}\alpha_{-n,k}^\fff |u\rangle =\langle \id| e^{H(t)}\alpha_{-n,k}^\fff e^{-H(t)}e^{H(t)}|u\rangle \]
   and 
   \[ e^{H(t)}\alpha_{-n,k}^\fff e^{-H(t)}=\alpha_{-n,k}^\fff+[H(t),\alpha_{-n,t}^\fff]+ {1\over 2!}[H(t),[H(t),\alpha_{-n,t}^\fff]]+\cdots=\alpha_{-n,t}^\fff+nt_n\]
   
   Combining we get,
    \[\langle \id| e^{H(t)}\alpha_{-n,k}^\fff |u\rangle = \langle \id|(\alpha_{-n,t}^\fff+nt_n)e^{H(t)}|w\rangle = nt_n\langle \id|e^{H(t)}|w\rangle = \alpha^{\bb}_{-n,k}\qedhere\]
\end{proof}

\section{Back Symmetric Functions and Generalized Boson-Fermion Correspondence}
\label{sec:back-stable}
Throughout the rest of the paper, we use $\xx_+$ to abbreviate the variable set $\{x_1,x_2,\cdots\}$, $\xx_{\leq k}$ for $\{\cdots x_{k-1},x_k\}$ and $\xx_{\infty}$ for $\{x_i:i\in\zz\}$.

Following \cite{lam2021back}, let $R$ be the $\qq$-algebra of formal power series $f$ in $\xx_\infty$ with bounded total degree (there is an $M$ such that all
monomials in $f$ have total degree at most $M$) and above bounded support (there is an
$N$ such that the variables $x_i$ do not appear  for $i > N$). Say that $f\in R$ is \emph{back symmetric} if there is a $k \in \zz$ such that $s_i(f) = f$ for all $i \leq k$, where $s_i$ is the action of $\sz$ on $R$ by swapping $x_i$ and $x_{i+1}$.
Denote $\bR$ the elements of $R$ that are back-symmetric.

The algebra structure on $\bR$ is given by
\[\bR = \Lambda^{\leq k} \otimes \qq [ \xx_\infty ],\]
for any choice of $k\in\zz$,
which is viewed as the polynomial algebra generated by algebraically independent generators $\{p^{\leq k}_i,i\in\zz_{\geq 0}\}$ and $\xx_\infty$. 

The ring $\bR$ has basis given by the back-stable Schubert polynomials, which we define via a co-product formula of \cite{lam2021back}.
\begin{definition}
	For $w\in\sz$, define
	\[\bS_w=\sum_{\substack{w\doteq vu\\v\in S_{\neq 0}}}\mf_{u}^{\leq 0}\mathfrak{S}_v\]
	where $S_{\neq 0}$ is the subgroup of $\sz$ that do not use the simple reflection $s_0$, and $\mathfrak{S}_v$ is the usual Schubert polynomial.
\end{definition}

\begin{remark}
	This particular definition of back-stable Schubert polynomial is naturally an element of $\Lambda^{\leq 0}\otimes \qq[\xx_\infty]$, where $\mf_u^{\leq 0}\in\Lambda^{\leq 0}$ is the symmetric part and $\mathfrak{S}_v\in \qq[\xx_\infty]$ is the finite polynomial part.
\end{remark}
We will now give $\bR$ an action of $\fH$, and name it the higher bosonic space $\bbb=\bR$.
\begin{definition}
	For $n\in\zz_{\geq 0}$ and $k\in\zz$, define operators $\alpha_n^\bbb$ and $\alpha_{-n,k}^\bbb$ on $\bbb=\bR$ as follows.
	\begin{align*}
	\alpha_n^\bbb f=n{\partial\over\partial  p_n^{\leq 0}}f	\\
	\alpha_{-n,k}^\bbb f=p_n^{\leq k} f
	\end{align*}
\end{definition}
\begin{remark}
	The partial derivative ${\partial / \partial p_n^{\leq 0}}$ is taken with respect to the generators $\{p_n^{\leq 0},n\in\zz_{\geq 0}\}\cup\{x_i:i\in\zz\}$ under the identification $\bbb = \Lambda^{\leq 0}\otimes \qq[\xx_\infty]$. Note that if we choose a different isomorphism $\bbb =\Lambda^{\leq k}\otimes \qq[\xx_\infty]$, the partial derivative $ {\partial / \partial p_n^{\leq 0}}$ will always be the same, i.e. independent of the choice of $k$.
\end{remark}
\begin{prop}
	The operators $\{\alpha_n^\bbb,\alpha_{-n,k}^\bbb\}$ satisfy the relation of $\fH$.
\end{prop}
\begin{proof}
	Straightforward.
\end{proof}
\begin{theorem}
	The map $ \Psi:\fff\to\bbb,|w\rangle \mapsto \bS_w$ is an isomorphism of $\fH$-modules.
\end{theorem}
\begin{proof}
	This map is clearly a bijection since the back-stable Schubert polynomials is a basis of $\bbb$. We now prove that it preserves $\fH$-actions.
	First consider the action of $\alpha_n^\bbb$,
	\begin{align*}
		\alpha_n^\bbb \bS_w&=\sum_{\substack{w\doteq uv\\v\in S_{\neq 0}}}\alpha_n^\bbb\mf_{u}^{\leq 0}\mathfrak{S}_v\\
		&=\sum_{\substack{w\doteq vu\\v\in S_{\neq 0}}}
		\sum_{\substack{ u\doteq \theta \rho\\ \ell(\rho)=n\\ \rho \text{ weak-ribbon} } } (-1)^{\spin(\rho)}\mf_\theta^{\leq 0} \mathfrak{S}_v \\
		&=\sum_{\substack{ w \doteq v \theta \rho\\\ell(\rho)=n\\ v\in S_{\neq 0}}}(-1)^{\spin(\rho)}\mf_\theta^{\leq 0} \mathfrak{S}_v\\
		&=\sum_{\substack{ w \doteq v \theta \rho\\ \ell(\rho)=n\\ v\in S_{\neq 0}}}\sum (-1)^{\spin(\rho)} \mf_\theta^{\leq 0} \mathfrak{S}_v\\
		&=\sum_{\substack{ w \doteq u \rho\\ \ell(\rho)=n}} (-1)^{\spin(\rho)} \sum_{\substack{u\doteq v\theta\\v\notin S_{\neq 0}}} \mf_\theta^{\leq 0} \mathfrak{S}_v\\
		&=\sum_{\substack{ w \doteq u \rho\\ \ell(\rho)=n}} (-1)^{\spin(\rho)} \bS_u=\Psi(\alpha_n^\fff|w\rangle)
	\end{align*}
	
For the other part, we fix any $k\in\mathbb{N}$. We shall use the co-product formula for the $k$-Grassmanian instead.
\begin{align*}
	\alpha_{-n,k}^\bbb \bS_w &= \sum_{\substack{w\doteq uv\\v\in S_{\neq k}}}\alpha_{-n,k}^\bbb\mf_{u}^{\leq k}\mathfrak{S}_v\\&=
\sum_{\substack{w\doteq uv\\v\in S_{\neq k}}} \sum_{u/\theta\ k\text{-strong ribbon}}(-1)^{\text{spin}(u/\theta)}\mf_{\theta}^{\leq k}  \mathfrak{S}_v\\
&=\sum_{ \substack{w \doteq \rho\theta v \\ \rho \text{\ is a cycle}\\v \notin S_{\neq k} } }(-1)^{\text{spin}(w/\theta v)} \mf_{\theta}^{\leq k}  \mathfrak{S}_v\\
&=\sum_{w/\theta v\ k\text{-strong ribbon}}(-1)^{\text{spin}(w/\theta v)}\sum_{\substack{\sigma\doteq \theta v\\v\in S_{\neq k}}} \mf_{\theta}^{\leq k}  \mathfrak{S}_v\\
&= \sum_{w/\sigma\ k\text{-strong ribbon}}(-1)^{\text{spin}(w/\sigma)}\mf^{\leq k}_\sigma \qedhere
\end{align*}
\end{proof}
We have now completed the commutative diagram in \Cref{eq:comm-diagram}.

\section*{Acknowledgements}
I thank Daoji Huang for introducing me to Schubert calculus, Ben Brubaker for inspiring my interest in the boson-fermion correspondence, and Yibo Gao for his insights on the Bruhat order. I thank Phillipe Di Frencesco, Sergey Fomin, Shiliang Gao, Jiyang Gao, Alex Yong, Tianyi Yu, Pasha Pylyavskyy and Paul Zinn-Justin for interesting discussions. An earlier version of this work was presented as part of my doctoral preliminary oral examination at the University of Minnesota, for which I thank my committee members Pasha Pylyavskyy, Roy Cook, Gregg Musiker, and Anna Weigandt for their comments and feedback.

\bibliographystyle{alpha}
\bibliography{main.bib}
\appendix
\section{Commutation relations}\label{sec:appendix}
\subsection{Diagrams for \Cref{prop: alpha_neg_commute}}
The following are list of description for the $\rho$-operator in the proof of \Cref{prop: alpha_neg_commute}. 

\centerline{\begin{tikzpicture}[xscale = 0.4,yscale = -0.4]
	\begin{pgfonlayer}{nodelayer}
		\node [style=none] (531) at (12, -2) {};
		\node [style=none] (532) at (28, -2) {};
		\node [style=none] (533) at (12, -4) {};
		\node [style=none] (534) at (28, -4) {};
		\node [style=none] (538) at (22, -4) {};
		\node [style=none] (541) at (22, -2) {};
		\node [style=none] (542) at (27, -2) {};
		\node [style=none] (543) at (27, -4) {};
		\node [style=none] (544) at (22, -4) {};
		\node [style=none] (555) at (22, -4) {};
		\node [style=none] (557) at (23, -3) {};
		\node [style=none] (558) at (22, -4) {};
		\node [style=none] (559) at (22, -3) {};
		\node [style=none] (560) at (27, -3) {};
		\node [style=none] (561) at (27, -2) {};
		\node [style=none] (562) at (26, -3) {};
		\node [style=none] (563) at (8.5, -4) {};
		\node [style=none] (564) at (9.5, -4) {};
		\node [style=none] (565) at (10.5, -4) {};
		\node [style=none] (566) at (22, -0.5) {};
		\node [style=none] (567) at (22, -1.5) {};
		\node [style=none] (568) at (27, -0.5) {};
		\node [style=none] (569) at (27, -1.5) {};
		\node [style=label-small] (570) at (24.5, -1) {$n$};
		\node [style=none] (571) at (22, -1) {};
		\node [style=none] (572) at (27, -1) {};
		\node [style=none] (575) at (12, -4) {};
		\node [style=none] (576) at (28, -4) {};
		\node [style=none] (577) at (12, -6) {};
		\node [style=none] (578) at (28, -6) {};
		\node [style=none] (579) at (13, -4) {};
		\node [style=none] (580) at (22, -4) {};
		\node [style=none] (581) at (22, -6) {};
		\node [style=none] (582) at (13, -6) {};
		\node [style=none] (586) at (22, -4) {};
		\node [style=none] (587) at (14, -5) {};
		\node [style=none] (588) at (13, -6) {};
		\node [style=none] (589) at (13, -5) {};
		\node [style=none] (590) at (22, -5) {};
		\node [style=none] (591) at (22, -4) {};
		\node [style=none] (592) at (21, -5) {};
		\node [style=none] (594) at (22, -4) {};
		\node [style=none] (595) at (-9, -2) {};
		\node [style=none] (596) at (7, -2) {};
		\node [style=none] (597) at (-9, -4) {};
		\node [style=none] (598) at (7, -4) {};
		\node [style=none] (599) at (-3, -2) {};
		\node [style=none] (600) at (6, -2) {};
		\node [style=none] (601) at (6, -4) {};
		\node [style=none] (602) at (-3, -4) {};
		\node [style=none] (603) at (6, -2) {};
		\node [style=none] (604) at (-2, -3) {};
		\node [style=none] (605) at (-3, -4) {};
		\node [style=none] (606) at (-3, -3) {};
		\node [style=none] (607) at (6, -3) {};
		\node [style=none] (608) at (6, -2) {};
		\node [style=none] (609) at (5, -3) {};
		\node [style=none] (610) at (6, -2) {};
		\node [style=none] (611) at (-9, -4) {};
		\node [style=none] (612) at (7, -4) {};
		\node [style=none] (613) at (-9, -6) {};
		\node [style=none] (614) at (7, -6) {};
		\node [style=none] (615) at (-8, -6) {};
		\node [style=none] (616) at (-8, -4) {};
		\node [style=none] (617) at (-3, -4) {};
		\node [style=none] (618) at (-3, -6) {};
		\node [style=none] (619) at (-8, -6) {};
		\node [style=none] (620) at (-8, -6) {};
		\node [style=none] (621) at (-7, -5) {};
		\node [style=none] (622) at (-8, -6) {};
		\node [style=none] (623) at (-8, -5) {};
		\node [style=none] (624) at (-3, -5) {};
		\node [style=none] (625) at (-3, -4) {};
		\node [style=none] (626) at (-4, -5) {};
		\node [style=none] (627) at (1, -0.5) {};
		\node [style=none] (628) at (1, -1.5) {};
		\node [style=none] (629) at (6, -0.5) {};
		\node [style=none] (630) at (6, -1.5) {};
		\node [style=label-small] (631) at (3.5, -1) {$n$};
		\node [style=none] (632) at (1, -1) {};
		\node [style=none] (633) at (6, -1) {};
		\node [style=none] (634) at (13, -6.5) {};
		\node [style=none] (635) at (13, -7.5) {};
		\node [style=none] (636) at (18, -6.5) {};
		\node [style=none] (637) at (18, -7.5) {};
		\node [style=label-small] (638) at (15.5, -7) {$n$};
		\node [style=none] (639) at (13, -7) {};
		\node [style=none] (640) at (18, -7) {};
		\node [style=none] (641) at (17, -5) {};
		\node [style=none] (642) at (19, -5) {};
		\node [style=none] (643) at (14, -5) {};
		\node [style=none] (645) at (19, -5) {};
		\node [style=none] (646) at (17, -5) {};
		\node [style=none] (647) at (19, -5) {};
		\node [style=none] (650) at (21, -5) {};
		\node [style=none] (651) at (-2, -3) {};
		\node [style=none] (652) at (2, -3) {};
		\node [style=none] (653) at (0, -3) {};
		\node [style=none] (654) at (2, -3) {};
		\node [style=none] (655) at (5, -3) {};
		\node [style=none] (675) at (27, -6.5) {};
		\node [style=none] (676) at (27, -7.5) {};
		\node [style=label-small] (677) at (22.5, -7) {$m$};
		\node [style=none] (678) at (18, -7) {};
		\node [style=none] (679) at (27, -7) {};
		\node [style=none] (680) at (-8, -6.5) {};
		\node [style=none] (681) at (-8, -7.5) {};
		\node [style=none] (682) at (-3, -6.5) {};
		\node [style=none] (683) at (-3, -7.5) {};
		\node [style=label-small] (684) at (-5.5, -7) {$n$};
		\node [style=none] (685) at (-8, -7) {};
		\node [style=none] (686) at (-3, -7) {};
		\node [style=none] (687) at (6, -6.5) {};
		\node [style=none] (688) at (6, -7.5) {};
		\node [style=label-small] (689) at (1.5, -7) {$m$};
		\node [style=none] (690) at (-3, -7) {};
		\node [style=none] (691) at (6, -7) {};
	\end{pgfonlayer}
	\begin{pgfonlayer}{edgelayer}
		\draw [style=fade] (531.center) to (532.center);
		\draw [style=fade] (533.center) to (534.center);
		\draw [style=shade] (543.center)
			 to (544.center)
			 to (541.center)
			 to (542.center)
			 to cycle;
		\draw [style=red, rounded corners=0.2cm] (558.center)
			 to (559.center)
			 to (557.center);
		\draw [style=red, rounded corners=0.2cm] (562.center)
			 to (560.center)
			 to (561.center);
		\draw [style=red, dashed] (557.center) to (562.center);
		\draw [style=red arrow] (564.center) to (563.center);
		\draw [style=red arrow] (564.center) to (565.center);
		\draw [style=fade] (566.center) to (567.center);
		\draw [style=fade] (568.center) to (569.center);
		\draw [style=red arrow] (570) to (571.center);
		\draw [style=red arrow] (570) to (572.center);
		\draw [style=fade] (578.center) to (577.center);
		\draw [style=fade] (575.center) to (576.center);
		\draw [style=shade] (581.center)
			 to (582.center)
			 to (579.center)
			 to (580.center)
			 to cycle;
		\draw [style=red, rounded corners=0.2cm] (588.center)
			 to (589.center)
			 to (587.center);
		\draw [style=red, rounded corners=0.2cm] (592.center)
			 to (590.center)
			 to (591.center);
		\draw [style=fade] (531.center) to (533.center);
		\draw [style=fade] (534.center) to (532.center);
		\draw [style=fade] (577.center) to (575.center);
		\draw [style=fade] (578.center) to (576.center);
		\draw [style=fade] (598.center) to (597.center);
		\draw [style=fade] (595.center) to (596.center);
		\draw [style=shade] (599.center)
			 to (600.center)
			 to (601.center)
			 to (602.center)
			 to cycle;
		\draw [style=red, rounded corners=0.2cm] (605.center)
			 to (606.center)
			 to (604.center);
		\draw [style=red, rounded corners=0.2cm] (609.center)
			 to (607.center)
			 to (608.center);
		\draw [style=fade] (597.center) to (595.center);
		\draw [style=fade] (598.center) to (596.center);
		\draw [style=fade] (611.center) to (612.center);
		\draw [style=fade] (613.center) to (614.center);
		\draw [style=shade] (618.center)
			 to (619.center)
			 to (616.center)
			 to (617.center)
			 to cycle;
		\draw [style=red, rounded corners=0.2cm] (622.center)
			 to (623.center)
			 to (621.center);
		\draw [style=red, rounded corners=0.2cm] (626.center)
			 to (624.center)
			 to (625.center);
		\draw [style=red, dashed] (621.center) to (626.center);
		\draw [style=fade] (611.center) to (613.center);
		\draw [style=fade] (614.center) to (612.center);
		\draw [style=fade] (627.center) to (628.center);
		\draw [style=fade] (629.center) to (630.center);
		\draw [style=red arrow] (631) to (632.center);
		\draw [style=red arrow] (631) to (633.center);
		\draw [style=fade] (634.center) to (635.center);
		\draw [style=fade] (636.center) to (637.center);
		\draw [style=red arrow] (638) to (639.center);
		\draw [style=red arrow] (638) to (640.center);
		\draw [style=red, rounded corners=0.2cm] (646.center) to (645.center);
		\draw [style=red, dashed] (643.center) to (646.center);
		\draw [style=red, dashed] (647.center) to (650.center);
		\draw [style=red, rounded corners=0.2cm] (653.center) to (652.center);
		\draw [style=red, dashed] (651.center) to (653.center);
		\draw [style=red, dashed] (654.center) to (655.center);
		\draw [style=fade] (675.center) to (676.center);
		\draw [style=fade, ->] (677) to (678.center);
		\draw [style=fade, ->] (677) to (679.center);
		\draw [style=fade] (680.center) to (681.center);
		\draw [style=fade] (682.center) to (683.center);
		\draw [style=red arrow] (684) to (685.center);
		\draw [style=red arrow] (684) to (686.center);
		\draw [style=fade] (687.center) to (688.center);
		\draw [style=red arrow] (689) to (690.center);
		\draw [style=red arrow] (689) to (691.center);
	\end{pgfonlayer}
\end{tikzpicture}}
\centerline{\begin{tikzpicture}[xscale = 0.4,yscale = -0.4]
	\begin{pgfonlayer}{nodelayer}
		\node [style=none] (531) at (12, -2) {};
		\node [style=none] (532) at (29, -2) {};
		\node [style=none] (533) at (12, -4) {};
		\node [style=none] (534) at (29, -4) {};
		\node [style=none] (538) at (23, -6) {};
		\node [style=none] (541) at (23, -4) {};
		\node [style=none] (542) at (28, -4) {};
		\node [style=none] (543) at (28, -6) {};
		\node [style=none] (544) at (23, -6) {};
		\node [style=none] (555) at (23, -6) {};
		\node [style=none] (560) at (28, -5) {};
		\node [style=none] (561) at (28, -4) {};
		\node [style=none] (562) at (27, -5) {};
		\node [style=none] (563) at (9, -4) {};
		\node [style=none] (564) at (10, -4) {};
		\node [style=none] (565) at (11, -4) {};
		\node [style=none] (566) at (23, -6.5) {};
		\node [style=none] (567) at (23, -7.5) {};
		\node [style=none] (568) at (28, -6.5) {};
		\node [style=none] (569) at (28, -7.5) {};
		\node [style=label-small] (570) at (25.5, -7) {$n$};
		\node [style=none] (571) at (23, -7) {};
		\node [style=none] (572) at (28, -7) {};
		\node [style=none] (575) at (12, -4) {};
		\node [style=none] (576) at (29, -4) {};
		\node [style=none] (577) at (12, -6) {};
		\node [style=none] (578) at (29, -6) {};
		\node [style=none] (579) at (13, -2) {};
		\node [style=none] (580) at (23, -2) {};
		\node [style=none] (581) at (23, -4) {};
		\node [style=none] (582) at (13, -4) {};
		\node [style=none] (587) at (14, -3) {};
		\node [style=none] (588) at (13, -4) {};
		\node [style=none] (589) at (13, -3) {};
		\node [style=none] (590) at (23, -3) {};
		\node [style=none] (591) at (23, -2) {};
		\node [style=none] (592) at (21, -3) {};
		\node [style=none] (595) at (-9, -2) {};
		\node [style=none] (596) at (8, -2) {};
		\node [style=none] (597) at (-9, -4) {};
		\node [style=none] (598) at (8, -4) {};
		\node [style=none] (599) at (-3, -2) {};
		\node [style=none] (600) at (7, -2) {};
		\node [style=none] (601) at (7, -4) {};
		\node [style=none] (602) at (-3, -4) {};
		\node [style=none] (603) at (7, -2) {};
		\node [style=none] (604) at (-2, -3) {};
		\node [style=none] (605) at (-3, -4) {};
		\node [style=none] (606) at (-3, -3) {};
		\node [style=none] (607) at (7, -3) {};
		\node [style=none] (608) at (7, -2) {};
		\node [style=none] (609) at (6, -3) {};
		\node [style=none] (610) at (7, -2) {};
		\node [style=none] (611) at (-9, -4) {};
		\node [style=none] (612) at (8, -4) {};
		\node [style=none] (613) at (-9, -6) {};
		\node [style=none] (614) at (8, -6) {};
		\node [style=none] (615) at (-8, -6) {};
		\node [style=none] (616) at (-8, -4) {};
		\node [style=none] (617) at (-3, -4) {};
		\node [style=none] (618) at (-3, -6) {};
		\node [style=none] (619) at (-8, -6) {};
		\node [style=none] (620) at (-8, -6) {};
		\node [style=none] (621) at (-7, -5) {};
		\node [style=none] (622) at (-8, -6) {};
		\node [style=none] (623) at (-8, -5) {};
		\node [style=none] (624) at (-3, -5) {};
		\node [style=none] (625) at (-3, -4) {};
		\node [style=none] (626) at (-4, -5) {};
		\node [style=none] (627) at (2, -0.5) {};
		\node [style=none] (628) at (2, -1.5) {};
		\node [style=none] (629) at (7, -0.5) {};
		\node [style=none] (630) at (7, -1.5) {};
		\node [style=label-small] (631) at (4.5, -1) {$n$};
		\node [style=none] (632) at (2, -1) {};
		\node [style=none] (633) at (7, -1) {};
		\node [style=none] (643) at (14, -3) {};
		\node [style=none] (646) at (21, -3) {};
		\node [style=none] (651) at (-2, -3) {};
		\node [style=none] (652) at (3, -3) {};
		\node [style=none] (653) at (1, -3) {};
		\node [style=none] (654) at (3, -3) {};
		\node [style=none] (655) at (6, -3) {};
		\node [style=none] (675) at (23, -0.5) {};
		\node [style=none] (676) at (23, -1.5) {};
		\node [style=label-small] (677) at (18, -1) {$m$};
		\node [style=none] (678) at (13, -1) {};
		\node [style=none] (679) at (23, -1) {};
		\node [style=none] (680) at (-8, -6.5) {};
		\node [style=none] (681) at (-8, -7.5) {};
		\node [style=none] (682) at (-3, -6.5) {};
		\node [style=none] (683) at (-3, -7.5) {};
		\node [style=label-small] (684) at (-5.5, -7) {$n$};
		\node [style=none] (685) at (-8, -7) {};
		\node [style=none] (686) at (-3, -7) {};
		\node [style=none] (687) at (7, -6.5) {};
		\node [style=none] (688) at (7, -7.5) {};
		\node [style=label-small] (689) at (2, -7) {$m$};
		\node [style=none] (690) at (-3, -7) {};
		\node [style=none] (691) at (7, -7) {};
		\node [style=none] (692) at (2, -2) {};
		\node [style=none] (693) at (2, -4) {};
		\node [style=none] (694) at (2, -4) {};
		\node [style=none] (695) at (13, -0.5) {};
		\node [style=none] (696) at (13, -1.5) {};
		\node [style=none] (697) at (24, -5) {};
		\node [style=none] (700) at (27, -5) {};
		\node [style=none] (704) at (24, -5) {};
		\node [style=none] (705) at (23, -6) {};
		\node [style=none] (706) at (23, -5) {};
		\node [style=none] (707) at (2, -3) {};
		\node [style=none] (708) at (2, -2) {};
		\node [style=none] (709) at (1, -3) {};
		\node [style=none] (710) at (3, -3) {};
		\node [style=none] (711) at (2, -4) {};
		\node [style=none] (712) at (2, -3) {};
	\end{pgfonlayer}
	\begin{pgfonlayer}{edgelayer}
		\draw [style=fade] (531.center) to (532.center);
		\draw [style=fade] (533.center) to (534.center);
		\draw [style=shade] (543.center)
			 to (544.center)
			 to (541.center)
			 to (542.center)
			 to cycle;
		\draw [style=red, rounded corners=0.2cm] (562.center)
			 to (560.center)
			 to (561.center);
		\draw [style=red arrow] (564.center) to (563.center);
		\draw [style=red arrow] (564.center) to (565.center);
		\draw [style=fade] (566.center) to (567.center);
		\draw [style=fade] (568.center) to (569.center);
		\draw [style=red arrow] (570) to (571.center);
		\draw [style=red arrow] (570) to (572.center);
		\draw [style=fade] (578.center) to (577.center);
		\draw [style=shade] (581.center)
			 to (582.center)
			 to (579.center)
			 to (580.center)
			 to cycle;
		\draw [style=red, rounded corners=0.2cm] (588.center)
			 to (589.center)
			 to (587.center);
		\draw [style=red, rounded corners=0.2cm] (592.center)
			 to (590.center)
			 to (591.center);
		\draw [style=faded] (531.center) to (533.center);
		\draw [style=faded] (534.center) to (532.center);
		\draw [style=faded] (577.center) to (575.center);
		\draw [style=faded] (578.center) to (576.center);
		\draw [style=fade] (598.center) to (597.center);
		\draw [style=fade] (595.center) to (596.center);
		\draw [style=shade] (599.center)
			 to (600.center)
			 to (601.center)
			 to (602.center)
			 to cycle;
		\draw [style=red, rounded corners=0.2cm] (605.center)
			 to (606.center)
			 to (604.center);
		\draw [style=red, rounded corners=0.2cm] (609.center)
			 to (607.center)
			 to (608.center);
		\draw [style=faded] (597.center) to (595.center);
		\draw [style=faded] (598.center) to (596.center);
		\draw [style=fade] (613.center) to (614.center);
		\draw [style=shade] (618.center)
			 to (619.center)
			 to (616.center)
			 to (617.center)
			 to cycle;
		\draw [style=red, rounded corners=0.2cm] (622.center)
			 to (623.center)
			 to (621.center);
		\draw [style=red, rounded corners=0.2cm] (626.center)
			 to (624.center)
			 to (625.center);
		\draw [style=red, dashed] (621.center) to (626.center);
		\draw [style=faded] (611.center) to (613.center);
		\draw [style=faded] (614.center) to (612.center);
		\draw [style=fade] (627.center) to (628.center);
		\draw [style=fade] (629.center) to (630.center);
		\draw [style=red arrow] (631) to (632.center);
		\draw [style=red arrow] (631) to (633.center);
		\draw [style=red, dashed] (643.center) to (646.center);
		\draw [style=red, dashed] (651.center) to (653.center);
		\draw [style=red, dashed] (654.center) to (655.center);
		\draw [style=fade] (675.center) to (676.center);
		\draw [style=red arrow] (677) to (678.center);
		\draw [style=red arrow] (677) to (679.center);
		\draw [style=fade] (680.center) to (681.center);
		\draw [style=fade] (682.center) to (683.center);
		\draw [style=red arrow] (684) to (685.center);
		\draw [style=red arrow] (684) to (686.center);
		\draw [style=fade] (687.center) to (688.center);
		\draw [style=red arrow] (689) to (690.center);
		\draw [style=red arrow] (689) to (691.center);
		\draw [style=fade] (695.center) to (696.center);
		\draw [style=red, dashed] (697.center) to (700.center);
		\draw [style=red, rounded corners=0.2cm] (705.center)
			 to (706.center)
			 to (704.center);
		\draw [style=red, rounded corners=0.2cm] (709.center)
			 to (707.center)
			 to (708.center);
		\draw [style=red, rounded corners=0.2cm] (711.center)
			 to (712.center)
			 to (710.center);
	\end{pgfonlayer}
\end{tikzpicture}}
\centerline{\begin{tikzpicture}[xscale = 0.4,yscale = -0.4]
	\begin{pgfonlayer}{nodelayer}
		\node [style=none] (531) at (12, -2) {};
		\node [style=none] (532) at (29, -2) {};
		\node [style=none] (533) at (12, -4) {};
		\node [style=none] (534) at (29, -4) {};
		\node [style=none] (538) at (24, -4) {};
		\node [style=none] (541) at (20, -2) {};
		\node [style=none] (542) at (28, -2) {};
		\node [style=none] (543) at (28, -4) {};
		\node [style=none] (544) at (20, -4) {};
		\node [style=none] (560) at (28, -3) {};
		\node [style=none] (561) at (28, -2) {};
		\node [style=none] (562) at (27, -3) {};
		\node [style=none] (563) at (9, -4) {};
		\node [style=none] (564) at (10, -4) {};
		\node [style=none] (565) at (11, -4) {};
		\node [style=none] (566) at (20, -0.25) {};
		\node [style=none] (567) at (20, -1.25) {};
		\node [style=none] (568) at (28, -0.25) {};
		\node [style=none] (569) at (28, -1.25) {};
		\node [style=label-small] (570) at (24, -0.75) {$n$};
		\node [style=none] (571) at (20, -0.75) {};
		\node [style=none] (572) at (28, -0.75) {};
		\node [style=none] (575) at (12, -4) {};
		\node [style=none] (576) at (29, -4) {};
		\node [style=none] (577) at (12, -6) {};
		\node [style=none] (578) at (29, -6) {};
		\node [style=none] (579) at (13, -4) {};
		\node [style=none] (580) at (24, -4) {};
		\node [style=none] (581) at (24, -6) {};
		\node [style=none] (582) at (13, -6) {};
		\node [style=none] (587) at (14, -5) {};
		\node [style=none] (588) at (13, -6) {};
		\node [style=none] (589) at (13, -5) {};
		\node [style=none] (590) at (24, -5) {};
		\node [style=none] (591) at (24, -4) {};
		\node [style=none] (592) at (23, -5) {};
		\node [style=none] (595) at (-9, -2) {};
		\node [style=none] (596) at (8, -2) {};
		\node [style=none] (597) at (-9, -4) {};
		\node [style=none] (598) at (8, -4) {};
		\node [style=none] (599) at (-8, -2) {};
		\node [style=none] (600) at (3, -2) {};
		\node [style=none] (601) at (3, -4) {};
		\node [style=none] (602) at (-8, -4) {};
		\node [style=none] (603) at (3, -2) {};
		\node [style=none] (604) at (-7, -3) {};
		\node [style=none] (605) at (-8, -4) {};
		\node [style=none] (606) at (-8, -3) {};
		\node [style=none] (607) at (3, -3) {};
		\node [style=none] (608) at (3, -2) {};
		\node [style=none] (609) at (2, -3) {};
		\node [style=none] (610) at (3, -2) {};
		\node [style=none] (611) at (-9, -4) {};
		\node [style=none] (612) at (8, -4) {};
		\node [style=none] (613) at (-9, -6) {};
		\node [style=none] (614) at (8, -6) {};
		\node [style=none] (615) at (-1, -6) {};
		\node [style=none] (616) at (-1, -4) {};
		\node [style=none] (617) at (7, -4) {};
		\node [style=none] (618) at (7, -6) {};
		\node [style=none] (619) at (-1, -6) {};
		\node [style=none] (620) at (-1, -6) {};
		\node [style=none] (621) at (0, -5) {};
		\node [style=none] (622) at (-1, -6) {};
		\node [style=none] (623) at (-1, -5) {};
		\node [style=none] (627) at (-1, -6.5) {};
		\node [style=none] (628) at (-1, -7.5) {};
		\node [style=none] (629) at (7, -6.5) {};
		\node [style=none] (630) at (7, -7.5) {};
		\node [style=label-small] (631) at (3, -7) {$n$};
		\node [style=none] (632) at (-1, -7) {};
		\node [style=none] (633) at (7, -7) {};
		\node [style=none] (643) at (14, -5) {};
		\node [style=none] (646) at (19, -5) {};
		\node [style=none] (653) at (2, -3) {};
		\node [style=none] (655) at (2, -3) {};
		\node [style=none] (675) at (24, -6.5) {};
		\node [style=none] (676) at (24, -7.5) {};
		\node [style=label-small] (677) at (18.5, -7) {$m$};
		\node [style=none] (678) at (13, -7) {};
		\node [style=none] (679) at (24, -7) {};
		\node [style=none] (682) at (-8, -0.5) {};
		\node [style=none] (683) at (-8, -1.5) {};
		\node [style=none] (686) at (-8, -1) {};
		\node [style=none] (687) at (3, -0.5) {};
		\node [style=none] (688) at (3, -1.5) {};
		\node [style=label-small] (689) at (-2.5, -1) {$m$};
		\node [style=none] (690) at (-8, -1) {};
		\node [style=none] (691) at (3, -1) {};
		\node [style=none] (695) at (13, -6.5) {};
		\node [style=none] (696) at (13, -7.5) {};
		\node [style=none] (697) at (25, -3) {};
		\node [style=none] (700) at (27, -3) {};
		\node [style=none] (704) at (25, -3) {};
		\node [style=none] (705) at (24, -4) {};
		\node [style=none] (706) at (24, -3) {};
		\node [style=none] (716) at (7, -5) {};
		\node [style=none] (717) at (7, -4) {};
		\node [style=none] (718) at (6, -5) {};
		\node [style=none] (720) at (20, -4) {};
		\node [style=none] (721) at (21, -3) {};
		\node [style=none] (722) at (20, -4) {};
		\node [style=none] (723) at (20, -3) {};
		\node [style=none] (725) at (0, -5) {};
		\node [style=none] (726) at (6, -5) {};
		\node [style=none] (737) at (-7, -3) {};
		\node [style=none] (738) at (2, -3) {};
		\node [style=none] (743) at (20, -4) {};
		\node [style=none] (744) at (20, -6) {};
		\node [style=none] (745) at (23, -5) {};
		\node [style=none] (746) at (21, -5) {};
		\node [style=none] (747) at (23, -5) {};
		\node [style=none] (748) at (24, -3) {};
		\node [style=none] (749) at (24, -2) {};
		\node [style=none] (750) at (23, -3) {};
		\node [style=none] (752) at (21, -3) {};
		\node [style=none] (753) at (23, -3) {};
		\node [style=none] (754) at (20, -6) {};
		\node [style=none] (755) at (21, -5) {};
		\node [style=none] (756) at (21, -5) {};
		\node [style=none] (757) at (20, -6) {};
		\node [style=none] (758) at (20, -5) {};
		\node [style=none] (759) at (20, -5) {};
		\node [style=none] (760) at (20, -4) {};
		\node [style=none] (761) at (19, -5) {};
		\node [style=none] (762) at (19, -5) {};
	\end{pgfonlayer}
	\begin{pgfonlayer}{edgelayer}
		\draw [style=fade] (531.center) to (532.center);
		\draw [style=fade] (533.center) to (534.center);
		\draw [style=shade] (543.center)
			 to (544.center)
			 to (541.center)
			 to (542.center)
			 to cycle;
		\draw [style=red, rounded corners=0.2cm] (562.center)
			 to (560.center)
			 to (561.center);
		\draw [style=red arrow] (564.center) to (563.center);
		\draw [style=red arrow] (564.center) to (565.center);
		\draw [style=fade] (566.center) to (567.center);
		\draw [style=fade] (568.center) to (569.center);
		\draw [style=red arrow] (570) to (571.center);
		\draw [style=red arrow] (570) to (572.center);
		\draw [style=fade] (578.center) to (577.center);
		\draw [style=shade] (581.center)
			 to (582.center)
			 to (579.center)
			 to (580.center)
			 to cycle;
		\draw [style=red, rounded corners=0.2cm] (588.center)
			 to (589.center)
			 to (587.center);
		\draw [style=red, rounded corners=0.2cm] (592.center)
			 to (590.center)
			 to (591.center);
		\draw [style=faded] (531.center) to (533.center);
		\draw [style=faded] (534.center) to (532.center);
		\draw [style=faded] (577.center) to (575.center);
		\draw [style=faded] (578.center) to (576.center);
		\draw [style=fade] (598.center) to (597.center);
		\draw [style=fade] (595.center) to (596.center);
		\draw [style=shade] (599.center)
			 to (600.center)
			 to (601.center)
			 to (602.center)
			 to cycle;
		\draw [style=red, rounded corners=0.2cm] (605.center)
			 to (606.center)
			 to (604.center);
		\draw [style=red, rounded corners=0.2cm] (609.center)
			 to (607.center)
			 to (608.center);
		\draw [style=faded] (597.center) to (595.center);
		\draw [style=faded] (598.center) to (596.center);
		\draw [style=fade] (613.center) to (614.center);
		\draw [style=shade] (618.center)
			 to (619.center)
			 to (616.center)
			 to (617.center)
			 to cycle;
		\draw [style=red, rounded corners=0.2cm] (622.center)
			 to (623.center)
			 to (621.center);
		\draw [style=faded] (611.center) to (613.center);
		\draw [style=faded] (614.center) to (612.center);
		\draw [style=fade] (627.center) to (628.center);
		\draw [style=fade] (629.center) to (630.center);
		\draw [style=red arrow] (631) to (632.center);
		\draw [style=red arrow] (631) to (633.center);
		\draw [style=red, dashed] (643.center) to (646.center);
		\draw [style=fade] (675.center) to (676.center);
		\draw [style=red arrow] (677) to (678.center);
		\draw [style=red arrow] (677) to (679.center);
		\draw [style=fade] (682.center) to (683.center);
		\draw [style=fade] (687.center) to (688.center);
		\draw [style=red arrow] (689) to (690.center);
		\draw [style=red arrow] (689) to (691.center);
		\draw [style=fade] (695.center) to (696.center);
		\draw [style=red, dashed] (697.center) to (700.center);
		\draw [style=red, rounded corners=0.2cm] (705.center)
			 to (706.center)
			 to (704.center);
		\draw [style=red, rounded corners=0.2cm] (718.center)
			 to (716.center)
			 to (717.center);
		\draw [style=red, rounded corners=0.2cm] (722.center)
			 to (723.center)
			 to (721.center);
		\draw [style=red, dashed] (725.center) to (726.center);
		\draw [style=red, dashed] (737.center) to (738.center);
		\draw [style=red, dashed] (746.center) to (747.center);
		\draw [style=red, rounded corners=0.2cm] (750.center)
			 to (748.center)
			 to (749.center);
		\draw [style=red, dashed] (752.center) to (753.center);
		\draw [style=red, rounded corners=0.2cm] (757.center)
			 to (758.center)
			 to (756.center);
		\draw [style=red, rounded corners=0.2cm] (761.center)
			 to (759.center)
			 to (760.center);
	\end{pgfonlayer}
\end{tikzpicture}}
\centerline{\begin{tikzpicture}[xscale = 0.4,yscale = -0.4]
	\begin{pgfonlayer}{nodelayer}
		\node [style=none] (531) at (12, -2) {};
		\node [style=none] (532) at (29, -2) {};
		\node [style=none] (533) at (12, -4) {};
		\node [style=none] (534) at (29, -4) {};
		\node [style=none] (538) at (22, -4) {};
		\node [style=none] (541) at (17, -2) {};
		\node [style=none] (542) at (28, -2) {};
		\node [style=none] (543) at (28, -4) {};
		\node [style=none] (544) at (17, -4) {};
		\node [style=none] (560) at (28, -3) {};
		\node [style=none] (561) at (28, -2) {};
		\node [style=none] (562) at (27, -3) {};
		\node [style=none] (563) at (9, -4) {};
		\node [style=none] (564) at (10, -4) {};
		\node [style=none] (565) at (11, -4) {};
		\node [style=none] (566) at (17, -0.25) {};
		\node [style=none] (567) at (17, -1.25) {};
		\node [style=none] (568) at (28, -0.25) {};
		\node [style=none] (569) at (28, -1.25) {};
		\node [style=label-small] (570) at (22.5, -0.75) {$m$};
		\node [style=none] (571) at (17, -0.75) {};
		\node [style=none] (572) at (28, -0.75) {};
		\node [style=none] (575) at (12, -4) {};
		\node [style=none] (576) at (29, -4) {};
		\node [style=none] (577) at (12, -6) {};
		\node [style=none] (578) at (29, -6) {};
		\node [style=none] (579) at (13, -4) {};
		\node [style=none] (580) at (22, -4) {};
		\node [style=none] (581) at (22, -6) {};
		\node [style=none] (582) at (13, -6) {};
		\node [style=none] (587) at (14, -5) {};
		\node [style=none] (588) at (13, -6) {};
		\node [style=none] (589) at (13, -5) {};
		\node [style=none] (590) at (22, -5) {};
		\node [style=none] (591) at (22, -4) {};
		\node [style=none] (592) at (21, -5) {};
		\node [style=none] (595) at (-9, -2) {};
		\node [style=none] (596) at (8, -2) {};
		\node [style=none] (597) at (-9, -4) {};
		\node [style=none] (598) at (8, -4) {};
		\node [style=none] (599) at (-4, -2) {};
		\node [style=none] (600) at (7, -2) {};
		\node [style=none] (601) at (7, -4) {};
		\node [style=none] (602) at (-4, -4) {};
		\node [style=none] (603) at (7, -2) {};
		\node [style=none] (604) at (-3, -3) {};
		\node [style=none] (605) at (-4, -4) {};
		\node [style=none] (606) at (-4, -3) {};
		\node [style=none] (607) at (7, -3) {};
		\node [style=none] (608) at (7, -2) {};
		\node [style=none] (609) at (6, -3) {};
		\node [style=none] (610) at (7, -2) {};
		\node [style=none] (611) at (-9, -4) {};
		\node [style=none] (612) at (8, -4) {};
		\node [style=none] (613) at (-9, -6) {};
		\node [style=none] (614) at (8, -6) {};
		\node [style=none] (615) at (-8, -6) {};
		\node [style=none] (616) at (-8, -4) {};
		\node [style=none] (617) at (2, -4) {};
		\node [style=none] (618) at (2, -6) {};
		\node [style=none] (619) at (-8, -6) {};
		\node [style=none] (620) at (-8, -6) {};
		\node [style=none] (621) at (-7, -5) {};
		\node [style=none] (622) at (-8, -6) {};
		\node [style=none] (623) at (-8, -5) {};
		\node [style=none] (624) at (-4, -5) {};
		\node [style=none] (625) at (-4, -4) {};
		\node [style=none] (626) at (-5, -5) {};
		\node [style=none] (627) at (-8, -6.5) {};
		\node [style=none] (628) at (-8, -7.5) {};
		\node [style=none] (629) at (2, -6.5) {};
		\node [style=none] (630) at (2, -7.5) {};
		\node [style=label-small] (631) at (-3, -7) {$n$};
		\node [style=none] (632) at (-8, -7) {};
		\node [style=none] (633) at (2, -7) {};
		\node [style=none] (643) at (14, -5) {};
		\node [style=none] (646) at (16, -5) {};
		\node [style=none] (651) at (2, -3) {};
		\node [style=none] (653) at (6, -3) {};
		\node [style=none] (655) at (6, -3) {};
		\node [style=none] (675) at (22, -6.5) {};
		\node [style=none] (676) at (22, -7.5) {};
		\node [style=label-small] (677) at (17.5, -7) {$n$};
		\node [style=none] (678) at (13, -7) {};
		\node [style=none] (679) at (22, -7) {};
		\node [style=none] (682) at (-4, -0.5) {};
		\node [style=none] (683) at (-4, -1.5) {};
		\node [style=none] (686) at (-4, -1) {};
		\node [style=none] (687) at (7, -0.5) {};
		\node [style=none] (688) at (7, -1.5) {};
		\node [style=label-small] (689) at (1.5, -1) {$m$};
		\node [style=none] (690) at (-4, -1) {};
		\node [style=none] (691) at (7, -1) {};
		\node [style=none] (695) at (13, -6.5) {};
		\node [style=none] (696) at (13, -7.5) {};
		\node [style=none] (697) at (23, -3) {};
		\node [style=none] (700) at (27, -3) {};
		\node [style=none] (704) at (23, -3) {};
		\node [style=none] (705) at (22, -4) {};
		\node [style=none] (706) at (22, -3) {};
		\node [style=none] (713) at (-3, -5) {};
		\node [style=none] (714) at (-4, -6) {};
		\node [style=none] (715) at (-4, -5) {};
		\node [style=none] (716) at (2, -5) {};
		\node [style=none] (717) at (2, -4) {};
		\node [style=none] (718) at (1, -5) {};
		\node [style=none] (720) at (17, -4) {};
		\node [style=none] (721) at (18, -3) {};
		\node [style=none] (722) at (17, -4) {};
		\node [style=none] (723) at (17, -3) {};
		\node [style=none] (725) at (-3, -5) {};
		\node [style=none] (726) at (1, -5) {};
		\node [style=none] (737) at (-3, -3) {};
		\node [style=none] (738) at (0, -3) {};
		\node [style=none] (743) at (17, -4) {};
		\node [style=none] (744) at (17, -6) {};
		\node [style=none] (745) at (21, -5) {};
		\node [style=none] (746) at (18, -5) {};
		\node [style=none] (747) at (21, -5) {};
		\node [style=none] (748) at (22, -3) {};
		\node [style=none] (749) at (22, -2) {};
		\node [style=none] (750) at (21, -3) {};
		\node [style=none] (752) at (18, -3) {};
		\node [style=none] (753) at (21, -3) {};
	\end{pgfonlayer}
	\begin{pgfonlayer}{edgelayer}
		\draw [style=fade] (531.center) to (532.center);
		\draw [style=fade] (533.center) to (534.center);
		\draw [style=shade] (543.center)
			 to (544.center)
			 to (541.center)
			 to (542.center)
			 to cycle;
		\draw [style=red, rounded corners=0.2cm] (562.center)
			 to (560.center)
			 to (561.center);
		\draw [style=red arrow] (564.center) to (563.center);
		\draw [style=red arrow] (564.center) to (565.center);
		\draw [style=fade] (566.center) to (567.center);
		\draw [style=fade] (568.center) to (569.center);
		\draw [style=red arrow] (570) to (571.center);
		\draw [style=red arrow] (570) to (572.center);
		\draw [style=fade] (578.center) to (577.center);
		\draw [style=shade] (581.center)
			 to (582.center)
			 to (579.center)
			 to (580.center)
			 to cycle;
		\draw [style=red, rounded corners=0.2cm] (588.center)
			 to (589.center)
			 to (587.center);
		\draw [style=red, rounded corners=0.2cm] (592.center)
			 to (590.center)
			 to (591.center);
		\draw [style=faded] (531.center) to (533.center);
		\draw [style=faded] (534.center) to (532.center);
		\draw [style=faded] (577.center) to (575.center);
		\draw [style=faded] (578.center) to (576.center);
		\draw [style=fade] (598.center) to (597.center);
		\draw [style=fade] (595.center) to (596.center);
		\draw [style=shade] (599.center)
			 to (600.center)
			 to (601.center)
			 to (602.center)
			 to cycle;
		\draw [style=red, rounded corners=0.2cm] (605.center)
			 to (606.center)
			 to (604.center);
		\draw [style=red, rounded corners=0.2cm] (609.center)
			 to (607.center)
			 to (608.center);
		\draw [style=faded] (597.center) to (595.center);
		\draw [style=faded] (598.center) to (596.center);
		\draw [style=fade] (613.center) to (614.center);
		\draw [style=shade] (618.center)
			 to (619.center)
			 to (616.center)
			 to (617.center)
			 to cycle;
		\draw [style=red, rounded corners=0.2cm] (622.center)
			 to (623.center)
			 to (621.center);
		\draw [style=red, rounded corners=0.2cm] (626.center)
			 to (624.center)
			 to (625.center);
		\draw [style=red, dashed] (621.center) to (626.center);
		\draw [style=faded] (611.center) to (613.center);
		\draw [style=faded] (614.center) to (612.center);
		\draw [style=fade] (627.center) to (628.center);
		\draw [style=fade] (629.center) to (630.center);
		\draw [style=red arrow] (631) to (632.center);
		\draw [style=red arrow] (631) to (633.center);
		\draw [style=red, dashed] (643.center) to (646.center);
		\draw [style=red, dashed] (651.center) to (653.center);
		\draw [style=fade] (675.center) to (676.center);
		\draw [style=red arrow] (677) to (678.center);
		\draw [style=red arrow] (677) to (679.center);
		\draw [style=fade] (682.center) to (683.center);
		\draw [style=fade] (687.center) to (688.center);
		\draw [style=red arrow] (689) to (690.center);
		\draw [style=red arrow] (689) to (691.center);
		\draw [style=fade] (695.center) to (696.center);
		\draw [style=red, dashed] (697.center) to (700.center);
		\draw [style=red, rounded corners=0.2cm] (705.center)
			 to (706.center)
			 to (704.center);
		\draw [style=red, rounded corners=0.2cm] (714.center)
			 to (715.center)
			 to (713.center);
		\draw [style=red, rounded corners=0.2cm] (718.center)
			 to (716.center)
			 to (717.center);
		\draw [style=red, rounded corners=0.2cm] (722.center)
			 to (723.center)
			 to (721.center);
		\draw [style=red, dashed] (725.center) to (726.center);
		\draw [style=red, dashed] (737.center) to (738.center);
		\draw [style=red] (738.center) to (651.center);
		\draw [style=red] (743.center) to (744.center);
		\draw [style=red, dashed] (746.center) to (747.center);
		\draw [style=red] (646.center) to (746.center);
		\draw [style=red, rounded corners=0.2cm] (750.center)
			 to (748.center)
			 to (749.center);
		\draw [style=red, dashed] (752.center) to (753.center);
	\end{pgfonlayer}
\end{tikzpicture}}
\centerline{\begin{tikzpicture}[xscale = 0.4,yscale = -0.4]
	\begin{pgfonlayer}{nodelayer}
		\node [style=none] (531) at (12, -2) {};
		\node [style=none] (532) at (29, -2) {};
		\node [style=none] (533) at (12, -4) {};
		\node [style=none] (534) at (29, -4) {};
		\node [style=none] (538) at (23, -4) {};
		\node [style=none] (541) at (13, -2) {};
		\node [style=none] (542) at (28, -2) {};
		\node [style=none] (543) at (28, -4) {};
		\node [style=none] (544) at (13, -4) {};
		\node [style=none] (560) at (28, -3) {};
		\node [style=none] (561) at (28, -2) {};
		\node [style=none] (562) at (27, -3) {};
		\node [style=none] (563) at (9, -4) {};
		\node [style=none] (564) at (10, -4) {};
		\node [style=none] (565) at (11, -4) {};
		\node [style=none] (566) at (13, -0.25) {};
		\node [style=none] (567) at (13, -1.25) {};
		\node [style=none] (568) at (28, -0.25) {};
		\node [style=none] (569) at (28, -1.25) {};
		\node [style=label-small] (570) at (20.5, -0.75) {$n$};
		\node [style=none] (571) at (13, -0.75) {};
		\node [style=none] (572) at (28, -0.75) {};
		\node [style=none] (575) at (12, -4) {};
		\node [style=none] (576) at (29, -4) {};
		\node [style=none] (577) at (12, -6) {};
		\node [style=none] (578) at (29, -6) {};
		\node [style=none] (579) at (18, -4) {};
		\node [style=none] (580) at (23, -4) {};
		\node [style=none] (581) at (23, -6) {};
		\node [style=none] (582) at (18, -6) {};
		\node [style=none] (590) at (23, -5) {};
		\node [style=none] (591) at (23, -4) {};
		\node [style=none] (592) at (22, -5) {};
		\node [style=none] (595) at (-9, -2) {};
		\node [style=none] (596) at (8, -2) {};
		\node [style=none] (597) at (-9, -4) {};
		\node [style=none] (598) at (8, -4) {};
		\node [style=none] (599) at (-3, -2) {};
		\node [style=none] (600) at (2, -2) {};
		\node [style=none] (601) at (2, -4) {};
		\node [style=none] (602) at (-3, -4) {};
		\node [style=none] (603) at (2, -2) {};
		\node [style=none] (604) at (-2, -3) {};
		\node [style=none] (605) at (-3, -4) {};
		\node [style=none] (606) at (-3, -3) {};
		\node [style=none] (607) at (2, -3) {};
		\node [style=none] (608) at (2, -2) {};
		\node [style=none] (609) at (1, -3) {};
		\node [style=none] (610) at (2, -2) {};
		\node [style=none] (611) at (-9, -4) {};
		\node [style=none] (612) at (8, -4) {};
		\node [style=none] (613) at (-9, -6) {};
		\node [style=none] (614) at (8, -6) {};
		\node [style=none] (616) at (-8, -4) {};
		\node [style=none] (617) at (7, -4) {};
		\node [style=none] (618) at (7, -6) {};
		\node [style=none] (619) at (-8, -6) {};
		\node [style=none] (621) at (-2, -5) {};
		\node [style=none] (622) at (-3, -6) {};
		\node [style=none] (623) at (-3, -5) {};
		\node [style=none] (627) at (-8, -6.5) {};
		\node [style=none] (628) at (-8, -7.5) {};
		\node [style=none] (629) at (7, -6.5) {};
		\node [style=none] (630) at (7, -7.5) {};
		\node [style=label-small] (631) at (-0.5, -7) {$n$};
		\node [style=none] (632) at (-8, -7) {};
		\node [style=none] (633) at (7, -7) {};
		\node [style=none] (653) at (1, -3) {};
		\node [style=none] (655) at (1, -3) {};
		\node [style=none] (675) at (23, -6.5) {};
		\node [style=none] (676) at (23, -7.5) {};
		\node [style=label-small] (677) at (20.5, -7) {$m$};
		\node [style=none] (678) at (18, -7) {};
		\node [style=none] (679) at (23, -7) {};
		\node [style=none] (682) at (-3, -0.5) {};
		\node [style=none] (683) at (-3, -1.5) {};
		\node [style=none] (686) at (-3, -1) {};
		\node [style=none] (687) at (2, -0.5) {};
		\node [style=none] (688) at (2, -1.5) {};
		\node [style=label-small] (689) at (-0.5, -1) {$m$};
		\node [style=none] (690) at (-3, -1) {};
		\node [style=none] (691) at (2, -1) {};
		\node [style=none] (695) at (18, -6.5) {};
		\node [style=none] (696) at (18, -7.5) {};
		\node [style=none] (697) at (24, -3) {};
		\node [style=none] (700) at (27, -3) {};
		\node [style=none] (704) at (24, -3) {};
		\node [style=none] (705) at (23, -4) {};
		\node [style=none] (706) at (23, -3) {};
		\node [style=none] (716) at (7, -5) {};
		\node [style=none] (717) at (7, -4) {};
		\node [style=none] (718) at (6, -5) {};
		\node [style=none] (720) at (13, -4) {};
		\node [style=none] (721) at (14, -3) {};
		\node [style=none] (722) at (13, -4) {};
		\node [style=none] (723) at (13, -3) {};
		\node [style=none] (725) at (-2, -5) {};
		\node [style=none] (726) at (6, -5) {};
		\node [style=none] (737) at (-2, -3) {};
		\node [style=none] (738) at (1, -3) {};
		\node [style=none] (744) at (18, -6) {};
		\node [style=none] (745) at (22, -5) {};
		\node [style=none] (746) at (19, -5) {};
		\node [style=none] (747) at (22, -5) {};
		\node [style=none] (748) at (23, -3) {};
		\node [style=none] (749) at (23, -2) {};
		\node [style=none] (750) at (22, -3) {};
		\node [style=none] (752) at (14, -3) {};
		\node [style=none] (753) at (22, -3) {};
		\node [style=none] (754) at (18, -6) {};
		\node [style=none] (755) at (19, -5) {};
		\node [style=none] (756) at (19, -5) {};
		\node [style=none] (757) at (18, -6) {};
		\node [style=none] (758) at (18, -5) {};
		\node [style=none] (764) at (-8, -6) {};
		\node [style=none] (765) at (-7, -5) {};
		\node [style=none] (766) at (-8, -6) {};
		\node [style=none] (767) at (-8, -5) {};
		\node [style=none] (768) at (-7, -5) {};
		\node [style=none] (769) at (-4, -5) {};
		\node [style=none] (770) at (-3, -4) {};
		\node [style=none] (774) at (-3, -5) {};
		\node [style=none] (775) at (-3, -5) {};
		\node [style=none] (776) at (-3, -4) {};
		\node [style=none] (777) at (-4, -5) {};
		\node [style=none] (778) at (-4, -5) {};
	\end{pgfonlayer}
	\begin{pgfonlayer}{edgelayer}
		\draw [style=fade] (531.center) to (532.center);
		\draw [style=fade] (533.center) to (534.center);
		\draw [style=shade] (543.center)
			 to (544.center)
			 to (541.center)
			 to (542.center)
			 to cycle;
		\draw [style=red, rounded corners=0.2cm] (562.center)
			 to (560.center)
			 to (561.center);
		\draw [style=red arrow] (564.center) to (563.center);
		\draw [style=red arrow] (564.center) to (565.center);
		\draw [style=fade] (566.center) to (567.center);
		\draw [style=fade] (568.center) to (569.center);
		\draw [style=red arrow] (570) to (571.center);
		\draw [style=red arrow] (570) to (572.center);
		\draw [style=fade] (578.center) to (577.center);
		\draw [style=shade] (581.center)
			 to (582.center)
			 to (579.center)
			 to (580.center)
			 to cycle;
		\draw [style=red, rounded corners=0.2cm] (592.center)
			 to (590.center)
			 to (591.center);
		\draw [style=faded] (531.center) to (533.center);
		\draw [style=faded] (534.center) to (532.center);
		\draw [style=faded] (577.center) to (575.center);
		\draw [style=faded] (578.center) to (576.center);
		\draw [style=fade] (598.center) to (597.center);
		\draw [style=fade] (595.center) to (596.center);
		\draw [style=shade] (599.center)
			 to (600.center)
			 to (601.center)
			 to (602.center)
			 to cycle;
		\draw [style=red, rounded corners=0.2cm] (605.center)
			 to (606.center)
			 to (604.center);
		\draw [style=red, rounded corners=0.2cm] (609.center)
			 to (607.center)
			 to (608.center);
		\draw [style=faded] (597.center) to (595.center);
		\draw [style=faded] (598.center) to (596.center);
		\draw [style=fade] (613.center) to (614.center);
		\draw [style=shade] (618.center)
			 to (619.center)
			 to (616.center)
			 to (617.center)
			 to cycle;
		\draw [style=red, rounded corners=0.2cm] (622.center)
			 to (623.center)
			 to (621.center);
		\draw [style=faded] (611.center) to (613.center);
		\draw [style=faded] (614.center) to (612.center);
		\draw [style=fade] (627.center) to (628.center);
		\draw [style=fade] (629.center) to (630.center);
		\draw [style=red arrow] (631) to (632.center);
		\draw [style=red arrow] (631) to (633.center);
		\draw [style=fade] (675.center) to (676.center);
		\draw [style=red arrow] (677) to (678.center);
		\draw [style=red arrow] (677) to (679.center);
		\draw [style=fade] (682.center) to (683.center);
		\draw [style=fade] (687.center) to (688.center);
		\draw [style=red arrow] (689) to (690.center);
		\draw [style=red arrow] (689) to (691.center);
		\draw [style=fade] (695.center) to (696.center);
		\draw [style=red, dashed] (697.center) to (700.center);
		\draw [style=red, rounded corners=0.2cm] (705.center)
			 to (706.center)
			 to (704.center);
		\draw [style=red, rounded corners=0.2cm] (718.center)
			 to (716.center)
			 to (717.center);
		\draw [style=red, rounded corners=0.2cm] (722.center)
			 to (723.center)
			 to (721.center);
		\draw [style=red, dashed] (725.center) to (726.center);
		\draw [style=red, dashed] (737.center) to (738.center);
		\draw [style=red, dashed] (746.center) to (747.center);
		\draw [style=red, rounded corners=0.2cm] (750.center)
			 to (748.center)
			 to (749.center);
		\draw [style=red, dashed] (752.center) to (753.center);
		\draw [style=red, rounded corners=0.2cm] (757.center)
			 to (758.center)
			 to (756.center);
		\draw [style=red, rounded corners=0.2cm] (766.center)
			 to (767.center)
			 to (765.center);
		\draw [style=red, dashed] (768.center) to (769.center);
		\draw [style=red, rounded corners=0.2cm] (777.center)
			 to (775.center)
			 to (776.center);
	\end{pgfonlayer}
\end{tikzpicture}}
\centerline{\begin{tikzpicture}[xscale=0.4,yscale = -0.4]
	\begin{pgfonlayer}{nodelayer}
		\node [style=none] (531) at (12, 0) {};
		\node [style=none] (532) at (29, 0) {};
		\node [style=none] (533) at (12, -4) {};
		\node [style=none] (534) at (29, -4) {};
		\node [style=none] (538) at (20, -6) {};
		\node [style=none] (541) at (20, -4) {};
		\node [style=none] (542) at (25, -4) {};
		\node [style=none] (543) at (25, -6) {};
		\node [style=none] (544) at (20, -6) {};
		\node [style=none] (555) at (20, -6) {};
		\node [style=none] (560) at (25, -5) {};
		\node [style=none] (561) at (25, -4) {};
		\node [style=none] (562) at (24, -5) {};
		\node [style=none] (563) at (9, -3) {};
		\node [style=none] (564) at (10, -3) {};
		\node [style=none] (565) at (11, -3) {};
		\node [style=none] (566) at (20, -6.5) {};
		\node [style=none] (567) at (20, -7.5) {};
		\node [style=none] (568) at (25, -6.5) {};
		\node [style=none] (569) at (25, -7.5) {};
		\node [style=label-small] (570) at (22.5, -7) {$n$};
		\node [style=none] (571) at (20, -7) {};
		\node [style=none] (572) at (25, -7) {};
		\node [style=none] (575) at (12, -4) {};
		\node [style=none] (576) at (29, -4) {};
		\node [style=none] (577) at (12, -6) {};
		\node [style=none] (578) at (29, -6) {};
		\node [style=none] (579) at (13, 0) {};
		\node [style=none] (580) at (27, 0) {};
		\node [style=none] (581) at (27, -4) {};
		\node [style=none] (582) at (13, -4) {};
		\node [style=none] (587) at (14, -1) {};
		\node [style=none] (588) at (13, -4) {};
		\node [style=none] (589) at (13, -1) {};
		\node [style=none] (590) at (20, -1) {};
		\node [style=none] (591) at (20, 0) {};
		\node [style=none] (592) at (19, -1) {};
		\node [style=none] (643) at (14, -1) {};
		\node [style=none] (646) at (19, -1) {};
		\node [style=none] (675) at (27, 1.5) {};
		\node [style=none] (676) at (27, 0.5) {};
		\node [style=label-small] (677) at (18, 1) {$m$};
		\node [style=none] (678) at (13, 1) {};
		\node [style=none] (679) at (27, 1) {};
		\node [style=none] (695) at (13, 1.5) {};
		\node [style=none] (696) at (13, 0.5) {};
		\node [style=none] (697) at (21, -5) {};
		\node [style=none] (700) at (24, -5) {};
		\node [style=none] (704) at (21, -5) {};
		\node [style=none] (705) at (20, -6) {};
		\node [style=none] (706) at (20, -5) {};
		\node [style=none] (713) at (20, -2) {};
		\node [style=none] (714) at (20, -3) {};
		\node [style=none] (715) at (20, -2) {};
		\node [style=none] (716) at (18, -3) {};
		\node [style=none] (719) at (21, -1) {};
		\node [style=none] (720) at (20, -2) {};
		\node [style=none] (721) at (20, -1) {};
		\node [style=none] (722) at (21, -1) {};
		\node [style=none] (729) at (27, -1) {};
		\node [style=none] (730) at (27, 0) {};
		\node [style=none] (731) at (21, -1) {};
		\node [style=none] (738) at (18, -3) {};
		\node [style=none] (739) at (17, -4) {};
		\node [style=none] (740) at (17, -3) {};
		\node [style=none] (741) at (18, -3) {};
		\node [style=none] (742) at (-9, -2) {};
		\node [style=none] (743) at (8, -2) {};
		\node [style=none] (744) at (-9, -4) {};
		\node [style=none] (745) at (8, -4) {};
		\node [style=none] (746) at (-2, -6) {};
		\node [style=none] (747) at (-2, -4) {};
		\node [style=none] (748) at (3, -4) {};
		\node [style=none] (749) at (3, -6) {};
		\node [style=none] (750) at (-2, -6) {};
		\node [style=none] (751) at (-2, -6) {};
		\node [style=none] (752) at (3, -5) {};
		\node [style=none] (753) at (3, -4) {};
		\node [style=none] (754) at (2, -5) {};
		\node [style=none] (755) at (-2, -6.5) {};
		\node [style=none] (756) at (-2, -7.5) {};
		\node [style=none] (757) at (3, -6.5) {};
		\node [style=none] (758) at (3, -7.5) {};
		\node [style=label-small] (759) at (0.5, -7) {$n$};
		\node [style=none] (760) at (-2, -7) {};
		\node [style=none] (761) at (3, -7) {};
		\node [style=none] (762) at (-9, -4) {};
		\node [style=none] (763) at (8, -4) {};
		\node [style=none] (764) at (-9, -6) {};
		\node [style=none] (765) at (8, -6) {};
		\node [style=none] (766) at (-8, -2) {};
		\node [style=none] (767) at (6, -2) {};
		\node [style=none] (768) at (6, -4) {};
		\node [style=none] (769) at (-8, -4) {};
		\node [style=none] (770) at (-7, -3) {};
		\node [style=none] (771) at (-8, -4) {};
		\node [style=none] (772) at (-8, -3) {};
		\node [style=none] (773) at (-1, -3) {};
		\node [style=none] (774) at (-1, -2) {};
		\node [style=none] (775) at (-3, -3) {};
		\node [style=none] (776) at (-7, -3) {};
		\node [style=none] (777) at (-3, -3) {};
		\node [style=none] (778) at (6, -0.5) {};
		\node [style=none] (779) at (6, -1.5) {};
		\node [style=label-small] (780) at (-3, -1) {$m$};
		\node [style=none] (781) at (-8, -1) {};
		\node [style=none] (782) at (6, -1) {};
		\node [style=none] (783) at (-8, -0.5) {};
		\node [style=none] (784) at (-8, -1.5) {};
		\node [style=none] (785) at (0, -5) {};
		\node [style=none] (786) at (2, -5) {};
		\node [style=none] (787) at (0, -5) {};
		\node [style=none] (788) at (-1, -6) {};
		\node [style=none] (789) at (-1, -5) {};
		\node [style=none] (790) at (-2, -4) {};
		\node [style=none] (794) at (0, -3) {};
		\node [style=none] (795) at (-1, -4) {};
		\node [style=none] (796) at (-1, -3) {};
		\node [style=none] (797) at (0, -3) {};
		\node [style=none] (798) at (6, -3) {};
		\node [style=none] (799) at (6, -2) {};
		\node [style=none] (800) at (0, -3) {};
		\node [style=none] (804) at (-2, -3) {};
		\node [style=none] (809) at (-1, -4) {};
		\node [style=none] (810) at (-2, -5) {};
		\node [style=none] (811) at (-1, -5) {};
	\end{pgfonlayer}
	\begin{pgfonlayer}{edgelayer}
		\draw [style=fade] (531.center) to (532.center);
		\draw [style=fade] (533.center) to (534.center);
		\draw [style=shade] (543.center)
			 to (544.center)
			 to (541.center)
			 to (542.center)
			 to cycle;
		\draw [style=red, rounded corners=0.2cm] (562.center)
			 to (560.center)
			 to (561.center);
		\draw [style=red arrow] (564.center) to (563.center);
		\draw [style=red arrow] (564.center) to (565.center);
		\draw [style=fade] (566.center) to (567.center);
		\draw [style=fade] (568.center) to (569.center);
		\draw [style=red arrow] (570) to (571.center);
		\draw [style=red arrow] (570) to (572.center);
		\draw [style=fade] (578.center) to (577.center);
		\draw [style=shade] (581.center)
			 to (582.center)
			 to (579.center)
			 to (580.center)
			 to cycle;
		\draw [style=red, rounded corners=0.2cm] (588.center)
			 to (589.center)
			 to (587.center);
		\draw [style=red, rounded corners=0.2cm] (592.center)
			 to (590.center)
			 to (591.center);
		\draw [style=faded] (531.center) to (533.center);
		\draw [style=faded] (534.center) to (532.center);
		\draw [style=faded] (577.center) to (575.center);
		\draw [style=faded] (578.center) to (576.center);
		\draw [style=red, dashed] (643.center) to (646.center);
		\draw [style=fade] (675.center) to (676.center);
		\draw [style=red arrow] (677) to (678.center);
		\draw [style=red arrow] (677) to (679.center);
		\draw [style=fade] (695.center) to (696.center);
		\draw [style=red, dashed] (697.center) to (700.center);
		\draw [style=red, rounded corners=0.2cm] (705.center)
			 to (706.center)
			 to (704.center);
		\draw [style=red, rounded corners=0.2cm] (716.center)
			 to (714.center)
			 to (715.center);
		\draw [style=red, rounded corners=0.2cm] (720.center)
			 to (721.center)
			 to (719.center);
		\draw [style=red, rounded corners=0.2cm] (731.center)
			 to (729.center)
			 to (730.center);
		\draw [style=red, rounded corners=0.2cm] (739.center)
			 to (740.center)
			 to (738.center);
		\draw [style=fade] (742.center) to (743.center);
		\draw [style=fade] (744.center) to (745.center);
		\draw [style=shade] (749.center)
			 to (750.center)
			 to (747.center)
			 to (748.center)
			 to cycle;
		\draw [style=red, rounded corners=0.2cm] (754.center)
			 to (752.center)
			 to (753.center);
		\draw [style=fade] (755.center) to (756.center);
		\draw [style=fade] (757.center) to (758.center);
		\draw [style=red arrow] (759) to (760.center);
		\draw [style=red arrow] (759) to (761.center);
		\draw [style=fade] (765.center) to (764.center);
		\draw [style=shade] (768.center)
			 to (769.center)
			 to (766.center)
			 to (767.center)
			 to cycle;
		\draw [style=red, rounded corners=0.2cm] (771.center)
			 to (772.center)
			 to (770.center);
		\draw [style=red, rounded corners=0.2cm] (775.center)
			 to (773.center)
			 to (774.center);
		\draw [style=faded] (742.center) to (744.center);
		\draw [style=faded] (745.center) to (743.center);
		\draw [style=faded] (764.center) to (762.center);
		\draw [style=faded] (765.center) to (763.center);
		\draw [style=red, dashed] (776.center) to (777.center);
		\draw [style=fade] (778.center) to (779.center);
		\draw [style=red arrow] (780) to (781.center);
		\draw [style=red arrow] (780) to (782.center);
		\draw [style=fade] (783.center) to (784.center);
		\draw [style=red, dashed] (785.center) to (786.center);
		\draw [style=red, rounded corners=0.2cm] (788.center)
			 to (789.center)
			 to (787.center);
		\draw [style=red, rounded corners=0.2cm] (795.center)
			 to (796.center)
			 to (794.center);
		\draw [style=red, rounded corners=0.2cm] (800.center)
			 to (798.center)
			 to (799.center);
		\draw [style=red, rounded corners=0.2cm] (810.center)
			 to (811.center)
			 to (809.center);
	\end{pgfonlayer}
\end{tikzpicture}}
\centerline{\begin{tikzpicture}[xscale=0.4,yscale = - 0.4]
	\begin{pgfonlayer}{nodelayer}
		\node [style=none] (531) at (12, 0) {};
		\node [style=none] (532) at (29, 0) {};
		\node [style=none] (533) at (12, -4) {};
		\node [style=none] (534) at (29, -4) {};
		\node [style=none] (538) at (23, -6) {};
		\node [style=none] (541) at (23, -4) {};
		\node [style=none] (542) at (28, -4) {};
		\node [style=none] (543) at (28, -6) {};
		\node [style=none] (544) at (23, -6) {};
		\node [style=none] (555) at (23, -6) {};
		\node [style=none] (560) at (28, -5) {};
		\node [style=none] (561) at (28, -4) {};
		\node [style=none] (562) at (27, -5) {};
		\node [style=none] (563) at (9, -4) {};
		\node [style=none] (564) at (10, -4) {};
		\node [style=none] (565) at (11, -4) {};
		\node [style=none] (566) at (23, -6.5) {};
		\node [style=none] (567) at (23, -7.5) {};
		\node [style=none] (568) at (28, -6.5) {};
		\node [style=none] (569) at (28, -7.5) {};
		\node [style=label-small] (570) at (25.5, -7) {$n$};
		\node [style=none] (571) at (23, -7) {};
		\node [style=none] (572) at (28, -7) {};
		\node [style=none] (575) at (12, -4) {};
		\node [style=none] (576) at (29, -4) {};
		\node [style=none] (577) at (12, -6) {};
		\node [style=none] (578) at (29, -6) {};
		\node [style=none] (579) at (13, 0) {};
		\node [style=none] (580) at (23, 0) {};
		\node [style=none] (581) at (23, -4) {};
		\node [style=none] (582) at (13, -4) {};
		\node [style=none] (587) at (14, -1) {};
		\node [style=none] (588) at (13, -4) {};
		\node [style=none] (589) at (13, -1) {};
		\node [style=none] (590) at (19.75, -1) {};
		\node [style=none] (591) at (19.75, 0) {};
		\node [style=none] (592) at (18, -1) {};
		\node [style=none] (595) at (-9, -2) {};
		\node [style=none] (596) at (8, -2) {};
		\node [style=none] (597) at (-9, -4) {};
		\node [style=none] (598) at (8, -4) {};
		\node [style=none] (611) at (-9, -4) {};
		\node [style=none] (612) at (8, -4) {};
		\node [style=none] (613) at (-9, -6) {};
		\node [style=none] (614) at (8, -6) {};
		\node [style=none] (627) at (-8.25, -0.5) {};
		\node [style=none] (628) at (-8.25, -1.5) {};
		\node [style=none] (629) at (-3, -0.5) {};
		\node [style=none] (630) at (-3, -1.5) {};
		\node [style=label-small] (631) at (-5.25, -1) {$n$};
		\node [style=none] (632) at (-8.25, -1) {};
		\node [style=none] (633) at (-3, -1) {};
		\node [style=none] (643) at (14, -1) {};
		\node [style=none] (646) at (18, -1) {};
		\node [style=none] (675) at (23, 1.5) {};
		\node [style=none] (676) at (23, 0.5) {};
		\node [style=label-small] (677) at (18, 1) {$m$};
		\node [style=none] (678) at (13, 1) {};
		\node [style=none] (679) at (23, 1) {};
		\node [style=none] (682) at (-5, -6.5) {};
		\node [style=none] (683) at (-5, -7.5) {};
		\node [style=none] (686) at (-5, -7) {};
		\node [style=none] (687) at (7, -6.5) {};
		\node [style=none] (688) at (7, -7.5) {};
		\node [style=label-small] (689) at (0.75, -7) {$m$};
		\node [style=none] (690) at (-5, -7) {};
		\node [style=none] (691) at (7, -7) {};
		\node [style=none] (695) at (13, 1.5) {};
		\node [style=none] (696) at (13, 0.5) {};
		\node [style=none] (697) at (24, -5) {};
		\node [style=none] (700) at (27, -5) {};
		\node [style=none] (704) at (24, -5) {};
		\node [style=none] (705) at (23, -6) {};
		\node [style=none] (706) at (23, -5) {};
		\node [style=none] (713) at (19.75, -2) {};
		\node [style=none] (714) at (19.75, -3) {};
		\node [style=none] (715) at (19.75, -2) {};
		\node [style=none] (716) at (18, -3) {};
		\node [style=none] (719) at (20.75, -1) {};
		\node [style=none] (720) at (19.75, -2) {};
		\node [style=none] (721) at (19.75, -1) {};
		\node [style=none] (722) at (20.75, -1) {};
		\node [style=none] (728) at (23, 0) {};
		\node [style=none] (729) at (23, -1) {};
		\node [style=none] (730) at (23, 0) {};
		\node [style=none] (731) at (20.75, -1) {};
		\node [style=none] (738) at (18, -3) {};
		\node [style=none] (739) at (17, -4) {};
		\node [style=none] (740) at (17, -3) {};
		\node [style=none] (741) at (18, -3) {};
		\node [style=none] (742) at (-3, -2) {};
		\node [style=none] (743) at (-3, -2) {};
		\node [style=none] (744) at (-8.25, -4) {};
		\node [style=none] (745) at (-8.25, -2) {};
		\node [style=none] (746) at (-3, -2) {};
		\node [style=none] (747) at (-3, -4) {};
		\node [style=none] (748) at (-8.25, -4) {};
		\node [style=none] (749) at (-8.25, -4) {};
		\node [style=none] (750) at (-7.25, -3) {};
		\node [style=none] (751) at (-8.25, -4) {};
		\node [style=none] (752) at (-8.25, -3) {};
		\node [style=none] (753) at (-3, -3) {};
		\node [style=none] (754) at (-3, -2) {};
		\node [style=none] (755) at (-4.25, -3) {};
		\node [style=none] (756) at (-5, -4) {};
		\node [style=none] (757) at (7, -4) {};
		\node [style=none] (758) at (7, -6) {};
		\node [style=none] (759) at (-5, -6) {};
		\node [style=none] (760) at (7, -4) {};
		\node [style=none] (761) at (-2, -5) {};
		\node [style=none] (762) at (-5, -6) {};
		\node [style=none] (763) at (-5, -5) {};
		\node [style=none] (764) at (7, -5) {};
		\node [style=none] (765) at (7, -4) {};
		\node [style=none] (766) at (6, -5) {};
		\node [style=none] (767) at (7, -4) {};
		\node [style=none] (770) at (-2, -5) {};
		\node [style=none] (771) at (3, -5) {};
		\node [style=none] (772) at (1, -5) {};
		\node [style=none] (773) at (3, -5) {};
		\node [style=none] (774) at (6, -5) {};
		\node [style=none] (775) at (2, -4) {};
		\node [style=none] (776) at (2, -6) {};
		\node [style=none] (777) at (2, -6) {};
		\node [style=none] (778) at (2, -5) {};
		\node [style=none] (779) at (2, -4) {};
		\node [style=none] (780) at (1, -5) {};
		\node [style=none] (781) at (3, -5) {};
		\node [style=none] (782) at (2, -6) {};
		\node [style=none] (783) at (2, -5) {};
	\end{pgfonlayer}
	\begin{pgfonlayer}{edgelayer}
		\draw [style=fade] (531.center) to (532.center);
		\draw [style=fade] (533.center) to (534.center);
		\draw [style=shade] (543.center)
			 to (544.center)
			 to (541.center)
			 to (542.center)
			 to cycle;
		\draw [style=red, rounded corners=0.2cm] (562.center)
			 to (560.center)
			 to (561.center);
		\draw [style=red arrow] (564.center) to (563.center);
		\draw [style=red arrow] (564.center) to (565.center);
		\draw [style=fade] (566.center) to (567.center);
		\draw [style=fade] (568.center) to (569.center);
		\draw [style=red arrow] (570) to (571.center);
		\draw [style=red arrow] (570) to (572.center);
		\draw [style=fade] (578.center) to (577.center);
		\draw [style=shade] (581.center)
			 to (582.center)
			 to (579.center)
			 to (580.center)
			 to cycle;
		\draw [style=red, rounded corners=0.2cm] (588.center)
			 to (589.center)
			 to (587.center);
		\draw [style=red, rounded corners=0.2cm] (592.center)
			 to (590.center)
			 to (591.center);
		\draw [style=faded] (531.center) to (533.center);
		\draw [style=faded] (534.center) to (532.center);
		\draw [style=faded] (577.center) to (575.center);
		\draw [style=faded] (578.center) to (576.center);
		\draw [style=fade] (598.center) to (597.center);
		\draw [style=fade] (595.center) to (596.center);
		\draw [style=faded] (597.center) to (595.center);
		\draw [style=faded] (598.center) to (596.center);
		\draw [style=fade] (613.center) to (614.center);
		\draw [style=faded] (611.center) to (613.center);
		\draw [style=faded] (614.center) to (612.center);
		\draw [style=fade] (627.center) to (628.center);
		\draw [style=fade] (629.center) to (630.center);
		\draw [style=red arrow] (631) to (632.center);
		\draw [style=red arrow] (631) to (633.center);
		\draw [style=red, dashed] (643.center) to (646.center);
		\draw [style=fade] (675.center) to (676.center);
		\draw [style=red arrow] (677) to (678.center);
		\draw [style=red arrow] (677) to (679.center);
		\draw [style=fade] (682.center) to (683.center);
		\draw [style=fade] (687.center) to (688.center);
		\draw [style=red arrow] (689) to (690.center);
		\draw [style=red arrow] (689) to (691.center);
		\draw [style=fade] (695.center) to (696.center);
		\draw [style=red, dashed] (697.center) to (700.center);
		\draw [style=red, rounded corners=0.2cm] (705.center)
			 to (706.center)
			 to (704.center);
		\draw [style=red, rounded corners=0.2cm] (716.center)
			 to (714.center)
			 to (715.center);
		\draw [style=red, rounded corners=0.2cm] (720.center)
			 to (721.center)
			 to (719.center);
		\draw [style=red, rounded corners=0.2cm] (731.center)
			 to (729.center)
			 to (730.center);
		\draw [style=red, rounded corners=0.2cm] (739.center)
			 to (740.center)
			 to (738.center);
		\draw [style=shade] (747.center)
			 to (748.center)
			 to (745.center)
			 to (746.center)
			 to cycle;
		\draw [style=red, rounded corners=0.2cm] (751.center)
			 to (752.center)
			 to (750.center);
		\draw [style=red, rounded corners=0.2cm] (755.center)
			 to (753.center)
			 to (754.center);
		\draw [style=red, dashed] (750.center) to (755.center);
		\draw [style=shade] (756.center)
			 to (757.center)
			 to (758.center)
			 to (759.center)
			 to cycle;
		\draw [style=red, rounded corners=0.2cm] (762.center)
			 to (763.center)
			 to (761.center);
		\draw [style=red, rounded corners=0.2cm] (766.center)
			 to (764.center)
			 to (765.center);
		\draw [style=red, dashed] (770.center) to (772.center);
		\draw [style=red, dashed] (773.center) to (774.center);
		\draw [style=red, rounded corners=0.2cm] (780.center)
			 to (778.center)
			 to (779.center);
		\draw [style=red, rounded corners=0.2cm] (782.center)
			 to (783.center)
			 to (781.center);
	\end{pgfonlayer}
\end{tikzpicture}}

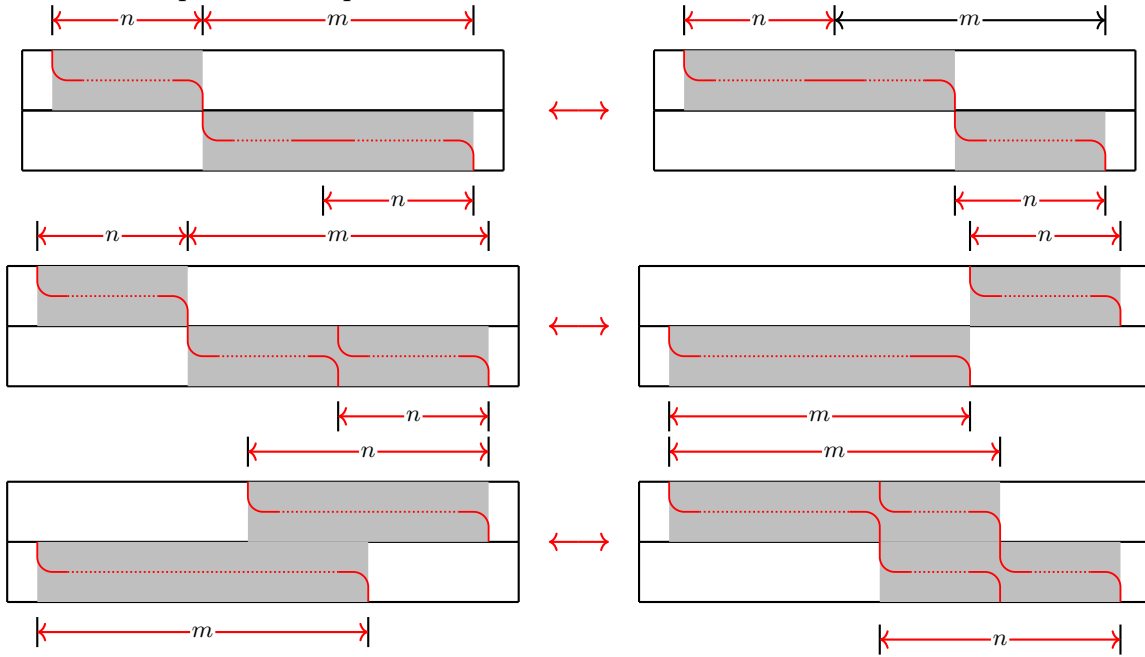
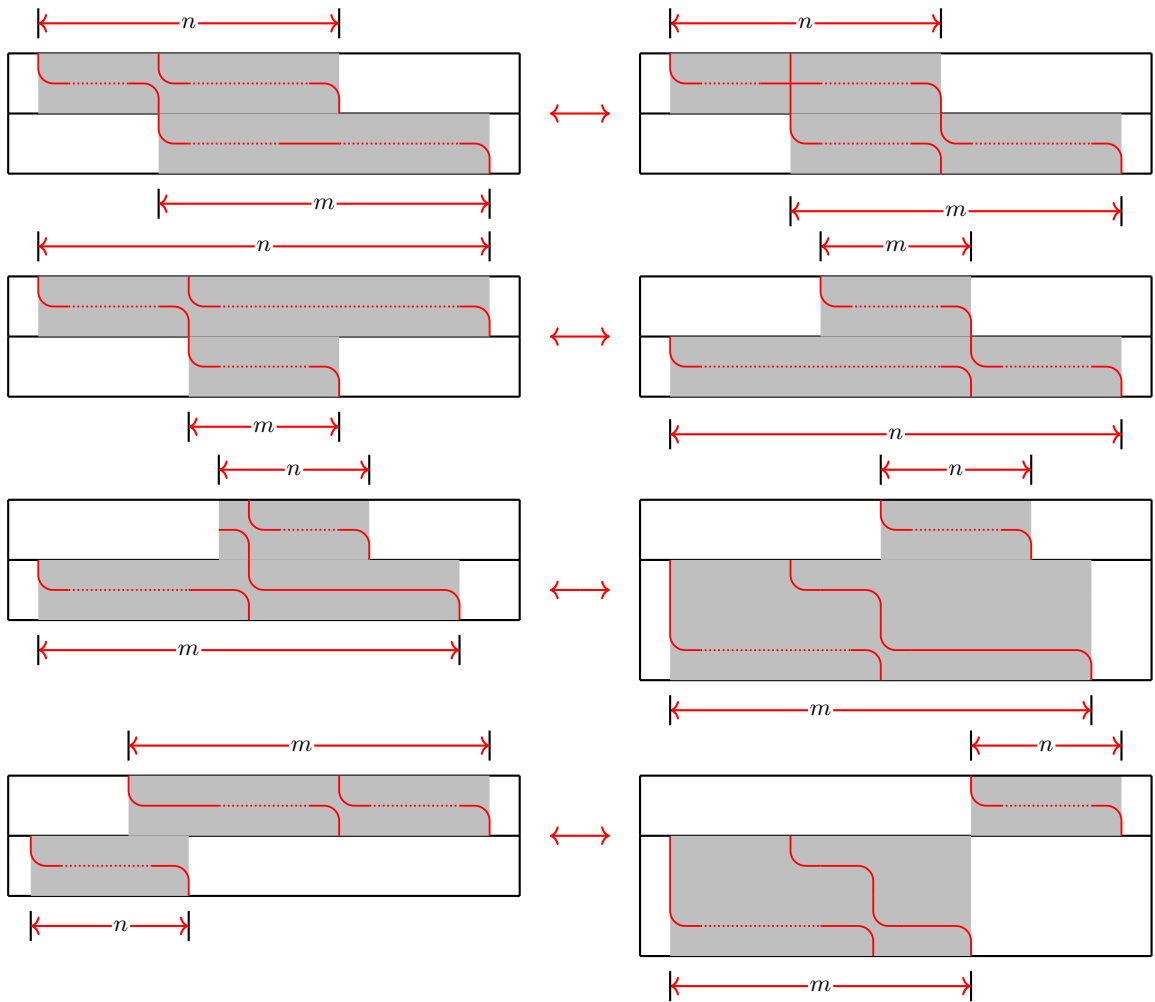
\end{document}